\documentclass[12pt]{amsart}
\usepackage{amssymb}
\usepackage{amsbsy}
\usepackage{amscd}
\usepackage{verbatim}
\makeatletter
\def\Bbb{\mathbb}
\def\frak{\mathfrak}

\newenvironment{pf*}[1]{\proof[#1]}{\endproof}
\newcommand{\rom}{\textup}
\makeatother
\newenvironment{aenume}{%
  \begin{enumerate}%
  }{\end{enumerate}}
\newtheorem{Theorem}[equation]{Theorem}
\newtheorem{Corollary}[equation]{Corollary}
\newtheorem{Lemma}[equation]{Lemma}
\newtheorem{Proposition}[equation]{Proposition}

\theoremstyle{definition}
\newtheorem{Definition}[equation]{Definition}
\newtheorem{Example}[equation]{Example}

\newtheorem{Notation}[equation]{Notation}
\newtheorem{Convention}[equation]{Convention}

\theoremstyle{remark}
\newtheorem{Remark}[equation]{Remark}
\newtheorem*{Claim}{Claim}

\numberwithin{equation}{subsection}

\newcommand{\thmref}[1]{Theorem~\ref{#1}}
\newcommand{\secref}[1]{\S\ref{#1}}
\newcommand{\lemref}[1]{Lemma~\ref{#1}}
\newcommand{\propref}[1]{Proposition~\ref{#1}}
\newcommand{\corref}[1]{Corollary~\ref{#1}}
\newcommand{\subsecref}[1]{\S\ref{#1}}

\newcommand{\defref}[1]{Definition~\ref{#1}}
\newcommand{\remref}[1]{Remark~\ref{#1}}
\newcommand{\lsp}[2]{\,{}^{#1}{#2}}
\newcommand{\defeq}{\overset{\operatorname{\scriptstyle def.}}{=}}
\newcommand{\C}{{\Bbb C}}
\newcommand{\Z}{{\Bbb Z}}
\newcommand{\Q}{{\Bbb Q}}
\newcommand{\R}{{\Bbb R}}

\newcommand{\SL}{\operatorname{\rm SL}}

\newcommand{\GL}{\operatorname{GL}}

\newcommand{\U}{\operatorname{\rm U}}

\newcommand{\algsl}{\operatorname{\frak{sl}}} 

\newcommand{\gl}{\operatorname{\frak{gl}}}
\newcommand{\algu}{\operatorname{\frak{u}}}

\newcommand{\Spec}{\operatorname{Spec}\nolimits}
\newcommand{\Proj}{\operatorname{Proj}\nolimits}
\newcommand{\End}{\operatorname{End}}
\newcommand{\Hom}{\operatorname{Hom}}
\newcommand{\Ext}{\operatorname{Ext}}
\newcommand{\Ker}{\operatorname{Ker}}
\newcommand{\Coker}{\operatorname{Coker}}
\newcommand{\Ima}{\operatorname{Im}}

\newcommand{\rank}{\operatorname{rank}}

\newcommand{\tr}{\operatorname{tr}}

\newcommand{\id}{\operatorname{id}}
\newcommand{\ve}{\varepsilon}
\newcommand{\vin}{\operatorname{in}} 
\newcommand{\vout}{\operatorname{out}} 
\newcommand{\bM}{{\mathbf M}} 
\newcommand{\M}{{\frak M}} 
\newcommand{\Mreg}{{\frak M_0^{\operatorname{reg}}}}
\newcommand{\La}{{\frak L}} 
\newcommand{\Fi}{{\frak F}} 
\newcommand{\dslash}{/\!\!/} 
\newcommand{\bv}{{\mathbf v}} 
\newcommand{\bw}{{\mathbf w}} 
\newcommand{\bC}{{\mathbf C}} 
\newcommand{\bA}{{\mathbf A}} 
\newcommand{\bI}{{\mathbf I}} 

\newcommand{\codim}{\operatorname{codim}} 
\newcommand{\gr}{\operatorname{gr}} 
\newcommand{\topdeg}{\operatorname{top}} 
\newcommand{\Uq}{{\mathbf U}_q(\mathfrak g)} 
\newcommand{\Un}{{\mathbf U}_q}
\newcommand{\Ua}{{\mathbf U}_q(\widehat{\mathfrak g})} 
\newcommand{\Ul}{{\mathbf U}_q({\mathbf L}{\mathfrak g})} 

\newcommand{\Ui}{{\mathbf U}^{\Z}_q(\mathfrak g)}
\newcommand{\Ue}{{\mathbf U}_{\varepsilon}(\mathfrak g)}
\newcommand{\Ueone}{{\mathbf U}_{\varepsilon = 1}(\mathfrak g)}
\newcommand{\Uli}{{\mathbf U}^{\Z}_q({\mathbf L}\mathfrak g)}
\newcommand{\Ule}{{\mathbf U}_{\varepsilon}({\mathbf L}\mathfrak g)}
\newcommand{\Pa}{{\frak P}} 
\newcommand{\HomE}{\operatorname{E}}
\newcommand{\HomL}{\operatorname{L}}

\newcommand{\hiw}{{L}} 

\newcommand{\shfO}{\mathcal O}
\newcommand{\Mp}{\widetilde\M}
\newcommand{\Mpp}{\widehat\M}
\newcommand{\Wedge}{{\textstyle \bigwedge}}
\newcommand{\Zw}{Z(\bw)}
\newcommand{\Mw}{\M(\bw)}
\newcommand{\Law}{\La(\bw)}
\newcommand{\ngpiso}{\overset{.}{=}}
\newcommand{\idl}{\mathfrak I}
\newcommand{\Ktop}{K_{\operatorname{top}}}
\newcommand{\Kztop}{K_{0,\operatorname{top}}}
\newcommand{\Kotop}{K_{1,\operatorname{top}}}

\begin{document}
\author{Hiraku Nakajima}
\title[quiver varieties and quantum affine algebras]
{Quiver varieties and finite dimensional representations of quantum
affine algebras
}
\address{Department of Mathematics, Kyoto University, Kyoto 606-8502,
Japan}
\email{nakajima@kusm.kyoto-u.ac.jp}
\thanks{Supported by the Grant-in-aid
for Scientific Research (No.11740011), the Ministry of Education,
Japan and the National Science Foundation Grants \#DMS 97-29992.}
\begin{abstract}
We study finite dimensional representations of the quantum affine
algebra $\Ua$ using geometry of quiver varieties introduced by the
author~\cite{KN,Na-quiver,Na-alg}. As an application, we obtain
character formulas expressed in terms of intersection cohomologies of
quiver varieties.
\end{abstract}

\maketitle

\section*{Introduction}\label{sec:intro}

Let $\mathfrak g$ be a simple finite dimensional Lie algebra of type
$ADE$,
$\widehat{\mathfrak g}$ be the corresponding (untwisted)
affine Lie algebra,
and $\Ua$ be its quantum enveloping algebra of Drinfel'd-Jimbo, or the
quantum affine algebra for short.
In this paper we study finite dimensional representations of $\Ua$,
using geometry of quiver varieties which were introduced in
\cite{KN,Na-quiver,Na-alg}.

There are immense literature on finite dimensional representations of
$\Ua$, see for example \cite{AK,CP,FR,Hata,Kl} and the reference
therein.
A basic result relevant to us is due to Chari-Pressley~\cite{CP2}:
irreducible finite dimensional representations are classified by
$n$-tuple of polynomials, where $n$ is the rank of $\mathfrak g$. This
result was announced for Yangian earlier by Drinfel'd~\cite{Drinfeld}.
Hence the polynomials are called Drinfel'd polynomials. However, not
much are known for properties of irreducible finite dimensional
representations, say their dimensions, tensor product decomposition,
etc.

Quiver varieties are generalization of moduli spaces of instantons
(anti-self-dual connections) on certain classes of real
$4$-dimensional hyper-K\"ahler manifolds, called ALE spaces~\cite{KN}.
They can be defined for any finite graph, but we are concerned with
the Dynkin graph of type $ADE$ corresponding to $\mathfrak g$ for a
moment.
Motivated by results of Ringel \cite{Ri} and Lusztig \cite{Lu},
the author has been studying their properties \cite{Na-quiver,Na-alg}.
In particular, it was shown that there is a homomorphism
\[
   {\mathbf U}(\mathfrak g) \to H_{\topdeg}(Z(\bw),\C),
\]
where ${\mathbf U}(\mathfrak g)$ is the universal enveloping algebra
of $\mathfrak g$, $Z(\bw)$ is a certain lagrangian subvariety of the
product of quiver varieties (the quiver variety depends on a choice of
a dominant weight $\bw$), and $H_{\topdeg}(\ ,\C)$ denotes the top
degree homology group with complex coefficients. The multiplication on
the right hand side is defined by the convolution product.

During the study, it became clear that the quiver varieties are
analogous to the cotangent bundle $T^*\mathcal B$ of the flag variety
$\mathcal B$.
The lagrangian subvariety $Z(\bw)$ is an analogue of Steinberg
variety $Z = T^*\mathcal B\times_{\mathcal N} T^*\mathcal B$, where
$\mathcal N$ is the nilpotent cone and $T^*\mathcal B\to \mathcal N$
is the Springer resolution.
The above mentioned result is an analogue of Ginzburg's lagrangian
construction of the Weyl group $W$ \cite{Gi-La}.
If we replace homology group by equivariant $K$-homology group in the
case of $T^*\mathcal B$, we get the affine Hecke algebra $H_q$ instead
of $W$ as was shown by Kazhdan-Lusztig~\cite{KL},
Ginzburg~\cite{Gi-book}.
Thus it became natural to conjecture that equivariant $K$-homology
group of the quiver variety gave us the quantum affine algebra $\Ua$.
After the author wrote \cite{Na-quiver}, many people suggested him
this conjecture, for example Kashiwara, Ginzburg, Lusztig, Vasserot.

A geometric approach to finite dimension representations of $\Ua$
(when $\mathfrak g = \algsl_n$) was given by
Ginzburg-Vasserot~\cite{GV,Vasserot}.
They used the cotangent bundle of the $n$-step partial flag variety,
which is an example of a quiver variety of type $A$. Thus their result
can be considered as a partial solution to the conjecture.

In \cite{Gr_aff} Grojnowski constructed the lower half part $\Ua^-$ of
$\Ua$ on equivariant $K$-homology of a certain lagrangian
subvariety of the cotangent bundle of a variety ${\mathbf E_{\mathbf
d}}$.
This ${\mathbf E_{\mathbf d}}$ was used earlier by Lusztig for the
construction of canonical bases on lower half part $\Uq^-$ of the
quantized enveloping algebra $\Uq$.
Grojnowski's construction was motivated in part by Tanisaki's
result~\cite{Tani}: a homomorphism from the finite Hecke algebra to
equivariant $K$-homology of Steinberg variety is defined by assigning
to perverse sheaves (or more precisely Hodge modules) on $\mathcal B$
their characteristic cycles.
In the same way, he considered characteristic cycles of perverse
sheaves on ${\mathbf E_{\mathbf d}}$. Thus he obtained a homomorphism
from $\Uq^-$ to $K$-homology of the lagrangian subvariety.
This lagrangian subvariety contains a lagrangian subvariety of the
quiver variety as an open subvariety.
Thus his construction was a solution to the `half' of the conjecture.

Later Grojnowski wrote an `advertisement' of his book on the full
conjecture~\cite{Gr}. Unfortunately, details were not explained, and
his book is not published yet.

The purpose of this paper is to solve the conjecture affirmatively,
and to derive results whose analogue are known for $H_q$.
Recall that Kazhdan-Lusztig \cite{KL} gave a classification of simple
modules of $H_q$, using above mentioned $K$-theoretic construction.
Our analogue is Drinfel'd-Chari-Pressley classification.
Also Ginzburg gave a character formula, called a $p$-adic analogue of
the Kazhdan-Lusztig multiplicity formula~\cite{Gi-book}.
(See introduction of \cite{Gi-book} for more detailed account and
historical comments.)
We prove a similar formula for $\Ua$ in this paper.

Let us describe the contents of this paper in more detail.
In \secref{sec:alg} we recall a new realization of $\Ua$, called
Drinfel'd realization~\cite{Drinfeld}.
It is more suitable than the original one for our purpose,
or rather, we can consider it as a definition of $\Ua$.
We also introduce the quantum loop algebra $\Ul$, which is a quotient
of $\Ua$, i.e., the quantum affine algebra without central
extension. Since the centeral extension acts trivially on finite
dimensional representations, we study $\Ul$ rather than $\Ua$.
Introducing a certain $\Z[q,q^{-1}]$-subalgebra $\Uli$ of $\Ul$, we
define a specialization $\Ule$ of $\Ul$ at $q = \varepsilon$.
This $\Uli$ was originally introduced by Chari-Pressley~\cite{CP3} for
the study of finite dimensional representations of $\Ule$ when
$\varepsilon$ is a root of unity.
Then we recall basic results on finite dimensional representations of
$\Ule$. We introduce several concepts, such as {\it l\/}-weights, {\it
l\/}-dominant, {\it l\/}-highest weight modules, {\it l\/}-fundamental
representation, etc. These are analogue of the same concepts without
{\it l\/} for $\Ue$-modules. `{\it l\/}' stands for the loop. In the
literature, some of these concepts were refered without `{\it l\/}'.

In \secref{sec:quiver_variety} we introduce two types of quiver
varieties $\Mw$, $\M_0(\infty,\bw)$ (both depend on a choice of a
dominant weight $\bw = \sum w_k \Lambda_k$). They are analogue of
$T^*\mathcal B$ and the nilpotent cone $\mathcal N$ respectively, and
have the following properties:
\begin{enumerate}
\item $\Mw$ is a nonsingular quasi-projective variety, having many 
components of various dimensions.
\item $\M_0(\infty,\bw)$ is an affine algebraic variety, not necessarily
irreducible.
\item Both $\Mw$ and $\M_0(\infty,\bw)$ have a
$G_\bw\times\C^*$-action, where $G_\bw = \prod \GL_{w_k}(\C)$.
\item There is a $G_\bw\times\C^*$-equivariant projective
morphism $\pi\colon\Mw\to \M_0(\infty,\bw)$.
\end{enumerate}

In \secref{sec:strat}$\sim$\secref{sec:convolution} we prepare some
results on quiver varieties and $K$-theory which we use in later
sections.

In \secref{sec:hom} we consider an analogue of the Steinberg variety
$\Zw = \M(\bw)\times_{\M_0(\infty,\bw)}\M(\bw)$ and its equivariant
$K$-homology $K^{G_\bw\times\C^*}(\Zw)$.
We construct an algebra homomorphism 
\[ 
   \Ul\to K^{G_\bw\times\C^*}(\Zw)\otimes_{\Z[q,q^{-1}]}\Q(q).
\]

In \secref{sec:hom} we define images of generators, and check the
defining relations in \secref{sec:rel1},\secref{sec:rel2}. Unlike the
case of the affine Hecke algebra, where $H_q$ is isomorphic to
$K^{G\times\C^*}(Z)$ ($Z =$ the Steinberg variety), this homomorphism
is {\it not\/} an isomorphism, neither injective nor surjective.

In \secref{sec:int} we show that the above homomorphism induces a
homomorphism between $\Uli$ and
$K^{G_\bw\times\C^*}(\Zw)/\!\operatorname{torsion}$.
(It is natural to expect that $\Uli$ is an integral form of $\Ul$ and
that $K^{G_\bw\times\C^*}(\Zw)$ is torsion-free, but we do not have
the proofs.)

In \secref{sec:std} we introduce a {\it standard module\/}
$M_{x,a}$. It depends on a choice of a point $x\in\M_0(\infty,\bw)$
and a semisimple element $a = (s,\varepsilon)\in G_\bw\times\C^*$ such
that $x$ is fixed by $a$.
The parameter $\varepsilon$ corresponds to the specialization $q =
\varepsilon$, while $s$ corresponds to Drinfel'd polynomials.
In this paper, we assume $\varepsilon$ is {\it not\/} a root of unity, 
though most of our results holds even in that case (see
\remref{rem:root of unity}).
Let $A$ be the Zariski closure of $a^{\Z}$. We define $M_{x,a}$ as the
specialized equivariant $K$-homology $K^A(\Mw_x)\otimes_{R(A)}\C_a$,
where $\Mw_x$ is a fiber of $\Mw\to\M_0(\infty,\bw)$ at $x$, and
$\C_a$ is an $R(A)$-algebra structure on $\C$ determined by $a$.
By the convolution product, $M_{x,a}$ has a
$K^A(\Zw)\otimes_{R(A)}\C_a$-module structure. Thus it has a
$\Ule$-module structure by the above homomorphism.
By the localization theorem of equivariant $K$-homology due to
Thomason~\cite{T-loc}, $M_{x,a}$ is isomorphic to the
complexified (non-equivariant) $K$-homology $K(\Mw_x^A)\otimes\C$ of
the fixed point set $\Mw_x^A$.
Moreover, it is isomorphic to $H_*(\Mw_x^A,\C)$ via the Chern
character homomorphism thanks to a result in \secref{sec:freeness}.
We also show that $M_{x,a}$ is a finite dimensional
{\it l\/}-highest weight module.
As a usual argument for Verma modules, $M_{x,a}$ has the unique
(nonzero) simple quotient.
The author conjectures that $M_{x,a}$ is a tensor product of {\it
l\/}-fundamental representations in some order. This is proved when
the parameter is generic in \subsecref{subsec:generic}.

In \secref{sec:simple} we show that the standard modules $M_{x,a}$ and
$M_{y,a}$ are isomorphic if and only if $x$ and $y$ are contained in
the same stratum.
Here the fixed point set $\M_0(\infty,\bw)^A$ has a stratification
$\M_0(\infty,\bw)^A = \bigsqcup_\rho \Mreg(\rho)$ defined in
\secref{sec:Fixed}.
Furthermore, we show that the index set $\{ \rho \}$ of stratum
coincides with the set $\mathcal P = \{ P \}$ of {\it l\/}-dominant
{\it l\/}-weights of $M_{0,a}$, the standard module corresponding to
the central fiber $\pi^{-1}(0)$. Let us denote by $\rho_P$ the index
corresponding to $P$.
Thus we may denote $M_{x,a}$ and its unique simple quotient by $M(P)$
and $L(P)$ respectively if $x$ is contained in the stratum
$\Mreg(\rho_P)$ corresponding to an {\it l\/}-dominant {\it
l\/}-weight $P$. We prove the multiplicity formula 
\[
   \left[ M(P) : L(Q) \right]
   = \dim H^*(i_x^! IC(\Mreg(\rho_Q))),
\]
where $x$ is a point in $\Mreg(\rho_P)$,
$i_x\colon \{x\}\to \M_0(\infty,\bw)^A$ is the inclusion,
and $IC(\Mreg(\rho_Q))$ is the intersection cohomology complex
attached to $\Mreg(\rho_Q)$ and the constant local system
$\C_{\Mreg(\rho_Q)}$.

Our result is simpler than the case of the affine Hecke algebra:
nonconstant local systems never appear. An algebraic reason for this
is that all modules are {\it l\/}-highest weight. It compensate for
the difference of $\Uli$ and $K^{G_\bw\times\C^*}(\Zw)$ during the
proof of the multiplicity formula.

If $\mathfrak g$ is of type $A$, then $\M_0(\infty,\bw)^A$ coincides
with a product of varieties ${\mathbf E_{\mathbf d}}$ studied by
Lusztig \cite{Lu}, where the underlying graph is of type $A$. In
particular, the Poincar\'e polynomial of 
\(
   H^*(i_x^! IC(\Mreg(\rho_Q)))
\)
is a Kazhdan-Lusztig polynomial for a Weyl group of type $A$.
We should have a combinatorial algorithm to compute Poincar\'e
polynomials of $H^*(i_x^! IC(\Mreg(\rho_Q)))$ for general $\mathfrak
g$.

Once we know $\dim H^*(i_x^! IC(\Mreg(\rho_Q)))$, information about
$L(P)$ can be deduced from information about $M(P)$, which is easier
to study.
For example, consider the following problems:
\begin{enumerate}
\item Compute Frenkel-Reshetikhin's $q$-characters~\cite{FR}.
\item Decompose restrictions of finite dimensional $\Ule$-modules
to $\Ue$-modules (see \cite{Kl}).
\end{enumerate}

These problems for $M(P)$ are easier than those for $L(P)$, and we
have the following answers.

Frenkel-Reshetikhin's $q$-characters are generating functions of
dimensions of {\it l\/}-weight spaces (see \subsecref{subsec:FR}).
In \subsecref{subsec:FR} we show that these dimensions are Euler
numbers of connected components of $\Mw^A$ for standard modules
$M_{0,a}$.
As an application, we prove a conjecture in \cite{FR} for $\mathfrak
g$ of type $ADE$.
These Euler numbers should be computable.

Let $\operatorname{Res}M(P)$ be the restriction of $M(P)$ to a
$\Ue$-module. In \secref{sec:Ue} we show the multiplicity formula
\[
   \left[ \operatorname{Res}M(P) : L(\bw-\bv) \right]
   = \dim H^*(i_x^! IC(\Mreg(\bv,\bw))),
\]
where $\bv$ is a weight such that $\bw - \bv$ is dominant,
$L(\bw-\bv)$ is the corresponding irreducible finite dimensional module
(these are concept for usual $\mathfrak g$ without `{\it l\/}'),
$x$ is a point in $\Mreg(\rho_P)$,
$i_x\colon \{x\}\to \M_0(\infty,\bw)$ is the inclusion,
$\Mreg(\bv,\bw)$ is a stratum of $\M(\infty,\bw)$,
and $IC(\Mreg(\bv,\bw))$ is the intersection cohomology complex
attached to $\Mreg(\bv,\bw)$ and the constant local system
$\C_{\Mreg(\bv,\bw)}$.

If $\mathfrak g$ is of type $A$, then $\Mreg(\bv,\bw)$ coincides with
a nilpotent orbit cut out by Slodowy's transversal slice by
\cite[8.4]{Na-quiver}. The Poincar\'e polynomials of $H^*(i_x^!
IC(\Mreg(\bv,\bw)))$ were calculated by Lusztig~\cite{Lu-Green} and
coincide with Kostka polynomials. This result is compatible with the
conjecture that $M(P)$ is a tensor product of {\it l\/}-fundamental
representations, for the restriction of an {\it l\/}-fundamental
representation is simple for type $A$, and Kostaka polynomials give
tensor product decompositions.
We should have a combinatorial algorithm to compute Poincar\'e
polynomials of $H^*(i_x^!IC(\Mreg(\bv,\bw)))$ for general $\mathfrak
g$.

We give two examples where $\Mreg(\bv,\bw)$ can be described
explicitly.

Consider the case that $\bw$ is a fundamental weight of type $A$, or
more generally a fundamental weight such that the label of the
corresponding vertex of the Dynkin diagram is $1$.
Then it is easy to see that the corresponding quiver variety
$\M_0(\infty,\bw)$ consists of a single point $0$. Thus
$\operatorname{Res}M(P)$ remains irreducible in this case.

If $\bw$ is the highest weight of the adjoint representation, the
corresponding $\M_0(\infty,\bw)$ is a simple singularity
$\C^2/\Gamma$, where $\Gamma$ is a finite subgroup of $\SL_2(\C)$
of the type corresponding to $\mathfrak g$. Then
$\M_0(\infty,\bw)$ has two strata $\{0\}$ and
$(\C^2\setminus\{0\})/\Gamma$.
The intersection cohomology complexes are constant sheaves. Hence we
have
\[
  \operatorname{Res}M(P) = L(\bw)\oplus L(0).
\]
These two results were shown by Chari-Pressley~\cite{CP-ex} by a totally
different method.

As we mentioned, the quantum affine algebra $\Ua$ has another
realization, called Drinfel'd new realization.
This Drinfel'd construction can be applied to any symmetrizable
Kac-Moody algebra $\mathfrak g$, not necessarily a finite dimensional
one.
This generalization also fit our result, for quiver varieties can be
defined for arbitrary finite graphs.
If we replace finite dimensional representations by {\it
l\/}-integrable representations, parts of our result can be
generalized to a Kac-Moody algebra $\mathfrak g$, at least when it is
symmetric.
For example, we generalize the Drinfel'd-Chari-Pressley parametrization.
A generalization of the multiplicity formula requires further study,
since we have extra strata which are not parametrized by $\mathcal P$.

If $\mathfrak g$ is an affine Lie algebra, then $\Ua$ is the quantum
affinization of the affine Lie algebra. It is called a {\it double
loop algebra}, or {\it toroidal algebra}, and has been studied by
various people, see for example \cite{GKV,Saito,STU,VV1} and the
reference therein. A first step to a geometric approach to the
toroidal algebra using quiver varieties for the affine Dynkin graph of
type $\widetilde A$ was given by M.~Varagnolo and
E.~Vasserot~\cite{VV2}.
In fact, quiver varieties for affine Dynkin graphs are moduli spaces
of instantons (or torsion free sheave) on ALE spaces. Thus these cases
are relevant to the original motivation, i.e., a study of relation
between the $4$-dimensional gauge theory and the representation
theory. In some cases, these quiver varieties coincide with Hilbert
schemes of points on ALE spaces, for which many results have been
obtained (see \cite{Lecture}). We would like to return back in future.

If we replace equivariant $K$-homology by equivariant homology, we
should get the Yangian $Y(\mathfrak g)$ instead of $\Ua$. This
conjecture is motivated again by the analogy of quiver varieties with
$T^*\mathcal B$. The equivariant homology of $T^*\mathcal B$ gives the
{\it graded\/} Hecke algebra~\cite{Lu:GH}, which is an analogue of
$Y(\mathfrak g)$ for $H_q$.
As an application, the conjecture implies that the representation
theory of $\Ua$ and that of the Yangian is the same. This has been
believed by many people, but there is no written proof.

While the author was preparing this paper, he was informed that
Frenkel-Mukhin \cite{FM} proved the conjecture in \cite{FR} for
general $\mathfrak g$.

\subsection*{Acknowledgement}
Part of this work was done while the author enjoyed the hospitality of
the Institute for Advanced Study. The author is grateful to G.~Lusztig 
for his interests and encouragements.

\tableofcontents

\section{Quantum affine algebra}\label{sec:alg}
In this section, we give a quick review for the definitions of the
quantized universal enveloping algebra $\Uq$ of the Kac-Moody algebra
$\mathfrak g$ associated with a symmetrizable generalized Cartan
matrix, its affinization $\Ua$, and the associated loop algebra $\Ul$.
Although the algebras defined via quiver varieties are automatically
symmetric, we treat the nonsymmetric case also for the
completeness.

\subsection{quantized universal enveloping algebra}\label{subsec:QUE}

Let $q$ be an indeterminate. For nonnegative integers $n\ge r$, define 
\begin{equation}\label{eq:q-binom}
  [n]_q \defeq \frac{q^n - q^{-n}}{q - q^{-1}}, \quad
  [n]_q ! \defeq 
  \begin{cases}
   [n]_q [n-1]_q \cdots [2]_q [1]_q &(n > 0),\\
   1 &(n=0),
  \end{cases}
\quad  
  \begin{bmatrix}
    n \\ r
  \end{bmatrix}_q \defeq \frac{[n]_q !}{[r]_q! [n-r]_q!}.
\end{equation}

Suppose that the following data are given:
\begin{enumerate}
  \item $P$ : free $\Z$-module (weight lattice),
  \item $P^* = \Hom_{\Z}(P, \Z)$ with a natural pairing 
    $\langle\ , \ \rangle\colon P\otimes P^*\to \Z$,
    \item an index set $I$ of simple roots
  \item $\alpha_k\in P$ ($k \in I$) (simple root),
  \item $h_k \in P^*$ ($k \in I$) (simple coroot),
  \item a symmetric bilinear form $(\ ,\ )$ on $P$.
\end{enumerate}
Those are required to satisfy the followings:
\begin{aenume}
\item $\langle h_k, \lambda\rangle 
= 2(\alpha_k, \lambda)/(\alpha_k,\alpha_k)$ for  
  $k \in I$ and $\lambda\in P$,
\item $\bC \defeq (\langle h_k, \alpha_l \rangle)_{k,l}$ is a symmetrizable
  generalized Cartan matrix, i.e., $\langle h_k, \alpha_k\rangle = 2$,
  and 
  $\langle h_k, \alpha_l\rangle\in\Z_{\le 0}$ and
  $\langle h_k, \alpha_l\rangle = 0 \Longleftrightarrow
   \langle h_l, \alpha_k\rangle = 0$ for $k\ne l$,
\item $(\alpha_k,\alpha_k)\in 2\Z_{> 0}$,
\item $\{\alpha_k\}_{k\in I}$ are linearly independent,
\item there exists $\Lambda_k\in P$ ($k \in I$) such that 
  $\langle h_l, \Lambda_k\rangle = \delta_{kl}$
  (fundamental weight).
\end{aenume}

The quantized universal enveloping algebra $\Uq$ of the Kac-Moody
algebra is the $\Q(q)$-algebra
generated by $e_k$, $f_k$ ($k \in I$), $q^h$ ($h\in P^*$)
with relations
{\allowdisplaybreaks
\begin{gather}
  q^0 = 1, \quad q^h q^{h'} = q^{h+h'},\\
  q^h e_k q^{-h} = q^{\langle h, \alpha_k\rangle} e_k,\quad
  q^h f_k q^{-h} = q^{-\langle h, \alpha_k\rangle} f_k,\\
  e_k f_l - f_l e_k = \delta_{kl}
   \frac{q^{(\alpha_k,\alpha_k)h_k/2}-q^{-(\alpha_k,\alpha_k)h_k/2}}
         {q_k - q_k^{-1}},\\
  \sum_{p=0}^{b}(-1)^p \begin{bmatrix} b \\ p\end{bmatrix}_{q_k}
                     e_k^p e_l e_k^{b-p} =
  \sum_{p=0}^{b}(-1)^p \begin{bmatrix} b \\ p\end{bmatrix}_{q_k}
                     f_k^p f_l f_k^{b-p} = 0 \quad
  \text{for $k\ne l$,}
\end{gather}
where $q_k = q^{(\alpha_k,\alpha_k)/2}$,
$b = 1 - \langle h_k,\alpha_l\rangle$.}

Let $\Uq^+$ (resp.\ $\Uq^-$) be the $\Q(q)$-subalgebra of $\Uq$
generated by elements $e_k$'s (resp.\ $f_k$'s).  Let $\Uq^0$
be the $\Q(q)$-subalgebra generated by elements $q^h$ ($h\in P^*$).
Then we have the triangle decomposition \cite[3.2.5]{Lu-Book}:
\begin{equation}
  \Uq \cong \Uq^+ \otimes \Uq^0 \otimes \Uq^- .
  \label{triangle}
\end{equation}

Let $e_k^{(n)}\defeq e_k^n/[n]_{q_k}!$, 
$f_k^{(n)}\defeq f_k^n/[n]_{q_k}!$.
Let $\Ui$ be the $\Z[q,q^{-1}]$-subalgebra of $\Uq$ generated by
elements
$e_k^{(n)}$, $f_k^{(n)}$, $q^h$ for $k\in I$, $n\in\Z_{> 0}$,
$h\in P^*$.
It is known that $\Ui$ is an {\it integral form\/} of $\Uq$, i.e.,
the natural map $\Ui\otimes_{\Z[q,q^{-1}]}\Q(q) \to \Uq$ is an
isomorphism. (See \cite[9.3.1]{CP}.)
For $\varepsilon\in\C^*$, Let us define $\Ue$ as
$\Ui\otimes_{\Z[q,q^{-1}]}\C$ via the
algebra homomorphism $\Z[q,q^{-1}]\to \C$ that takes $q$ to
$\varepsilon$. It will be called the {\it specialized\/} quantized
enveloping algebra.
We say a $\Uq$-module $M$ (defined over $\Q(q)$)
is a {\it highest weight module\/} with
{\it highest weight\/} $\Lambda\in P$ if
there exists a vector $m_0\in M$ such that
\begin{align}
  & e_k\ast m_0 = 0, \qquad \Uq^- \ast m_0 = M, \\
  & q^h\ast m_0 = q^{\langle h, \Lambda\rangle} m_0
\qquad\text{for any $h\in P^*$.}
\end{align}
Then there exists a direct sum decomposition $M =
\bigoplus_{\lambda\in P} M_\lambda$ (weight space decomposition)
where
$M_\lambda \defeq \{ m\mid 
q^h\cdot v = q^{\langle h, \lambda\rangle} m$ for any $h\in P^*$\}.
By using the triangular decomposition \eqref{triangle}, one can show
that the simple highest weight $\Uq$-module is determined
uniquely by $\Lambda$.

We say a $\Uq$-module $M$ (defined over $\Q(q)$) is {\it integrable\/}
if $M$ has a weight space decomposition $M = \bigoplus_{\lambda\in P}
M_\lambda$ with $\dim M_\lambda < \infty$,
and for any $m\in M$, there exists $n_0\ge 1$ such that $e_{k}^{n}
\ast m = f_{k}^{n}\ast m = 0$ for all $k\in I$ and $n\ge n_0$.

The (unique) simple highest weight $\Uq$-module with highest
weight $\Lambda$ is integrable if and only if
$\Lambda$ is a dominant integral weight $\Lambda$, i.e.,
$\langle\Lambda, h_k\rangle\in \Z_{\ge 0}$ for any $k\in I$
(\cite[3.5.6, 3.5.8]{Lu-Book}).
In this case, the integrable highest
weight $\Un$-module with highest weight $\Lambda$ 
is denoted by $\hiw(\Lambda)$.

For a $\Ue$-module $M$ (defined over $\C$), we define highest weight
modules, integrable modules, etc in a similar way.

Suppose $\Lambda$ is dominant. 
Let $\hiw(\Lambda)^\Z \defeq \Ui\ast m_0$, where $m_0$ is the highest
weight vector.
It is known that the natural map
$\hiw(\Lambda)^\Z\otimes_{\Z[q,q^{-1}]}\Q(q)\to \hiw(\Lambda)$
is an isomorphism and 
$\hiw(\Lambda)^\Z\otimes_{\Z[q,q^{-1}]}\C$ is the simple
integrable highest weight module of the corresponding Kac-Moody
algebra $\mathfrak g$ with highest weight $\Lambda$, where
$\Z[q,q^{-1}]\to \C$ is the homomorphism that sends $q$ to $1$
(\cite[Chapter~14 and 33.1.3]{Lu-Book}).
Unless $\varepsilon$ is a root of unity, the simple integrable highest
weight $\Ue$-module is the specialization of $\hiw(\Lambda)^\Z$
(\cite[10.1.14, 10.1.15]{CP}).

\subsection{quantum affine algebra}\label{subsec:quantum_affine}


The {\it quantum affinization\/} $\Ua$ of $\Uq$ (or simply {\it
quantum affine algebra\/}) is an associative algebra over $\Q(q)$
generated by $e_{k,r}, f_{k,r}$ ($k\in I$, $r\in\Z$), $q^h$ ($h \in
P^*$), $q^{\pm c/2}$, $q^{\pm d}$, $h_{k,m}$ ($k\in I$, $m\in
\Z\setminus\{0\}$) with the following defining relations
{\allowdisplaybreaks[4]
\begin{gather}
  \text{$q^{\pm c/2}$ is central,} \label{eq:relCcent}\\
  q^0 = 1, \quad q^h q^{h'} = q^{h+h'}, \quad
  [q^h, h_{k,m}] = 0, \quad q^{d} q^{-d} = 1, \quad
  q^{c/2} q^{-c/2} = 1, \label{eq:relHH}\\
%
%
  \psi_k^\pm(z) \psi_l^\pm(w) = \psi_l^\pm(w) \psi_k^\pm(z),
  \label{eq:relHH2}\\
  \psi_k^-(z) \psi_l^+(w) =
   \frac{(z - q^{-(\alpha_k,\alpha_l)}q^c w)
         (z - q^{(\alpha_k,\alpha_l)}q^{-c} w)}
        {(z - q^{(\alpha_k,\alpha_l)}q^c w)
         (z - q^{-(\alpha_k,\alpha_l)}q^{-c} w)}
   \psi_l^+(w) \psi_k^-(z)
  ,\label{eq:relHH3}\\
  [ q^{d}, q^h ] = 0, \quad q^d h_{k,m} q^{-d} = q^m h_{k,m},
  \quad
  q^d e_{k,r} q^{-d} = q^r e_{k,r}, \quad
  q^d f_{k,r} q^{-d} = q^r f_{k,r}, \label{eq:relD}
\\
  ( q^{\pm sc/2} z - q^{\pm \langle h_k,\alpha_l\rangle} w)
    \psi_l^s(z) x_k^\pm(w) =
  ( q^{\pm\langle h_k,\alpha_l\rangle} q^{\pm sc/2} z -  w)
    x^\pm_k(w) \psi_l^s(z), \quad
    \label{eq:relHE}
\\
%
%
%
  \left[x_{k}^+(z), x_{l}^-(w)\right] =
  \frac{\delta_{kl}}{q_k  - q_k^{-1}}
  \left\{\delta\left(q^c\frac{w}{z}\right)\psi^+_k(q^{c/2}w) -
        \delta\left(q^c\frac{z}{w}\right)\psi^-_k(q^{c/2}z)\right\},
    \label{eq:relEF}
\\
%
   (z - q^{\pm 2} w) x_{k}^\pm(z) x_{k}^\pm(w)
   = (q^{\pm 2} z - w) x_{k}^\pm(w) x_{k}^\pm(z),
   \label{eq:relexE2}
\\
%
  \prod_{p=1}^{-\langle \alpha_k,h_l\rangle}
   (z - q^{\pm(b'-2p)} w) x_{k}^\pm(z) x_{l}^\pm(w)
 = \prod_{p=1}^{-\langle \alpha_k, h_l\rangle}
   (q^{\pm(b'-2p)} z - w) x_{l}^\pm(w) x_{k}^\pm(z),
   \quad\text{if $k\neq l$},
   \label{eq:relexE1}
\\
%
%
%
%
  \sum_{\sigma\in S_b}
   \sum_{p=0}^{b}(-1)^p 
   \begin{bmatrix} b \\ p\end{bmatrix}_{q_k}
   x_{k}^\pm(z_{\sigma(1)})\cdots x_{k}^\pm(z_{\sigma(p)})
   x_{l}^\pm(w)
   x_{k}^\pm(z_{\sigma(p+1)})\cdots x_{k}^\pm(z_{\sigma(b)}) = 0,
   \quad \text{if $k\neq l$,}
 \label{eq:relDS}
%
\end{gather}
where}
$q_k = q^{(\alpha_k,\alpha_k)/2}$,
$s = \pm$,
$b = 1-\langle h_k, \alpha_l\rangle$, 
$b' = - (\alpha_k,\alpha_l)$,
and 
$S_b$ is the symmetric group of $b$ letters.
Here $\delta(z)$, $x_k^+(z)$, $x_k^-(z)$, $\psi^{\pm}_{k}(z)$ are
generating functions defined by
{\allowdisplaybreaks[4]
\begin{gather*}
   \delta(z) \defeq \sum_{r=-\infty}^\infty z^{r}, \qquad
   x_k^+(z) \defeq \sum_{r=-\infty}^\infty e_{k,r} z^{-r}, \qquad
   x_k^-(z) \defeq \sum_{r=-\infty}^\infty f_{k,r} z^{-r}, \\
   \psi^{\pm}_k(z)
  \defeq q^{\pm (\alpha_k,\alpha_k)h_k/2}
   \exp\left(\pm (q_k-q_k^{-1})\sum_{m=1}^\infty h_{k,\pm m} z^{\mp m}\right).
\end{gather*}
We} also need the following generating function later:
\begin{equation*}
   p_k^\pm(z) \defeq 
   \exp\left(
     - \sum_{m=1}^\infty \frac{h_{k,\pm m}}{[m]_{q_k}} z^{\mp m}
   \right).
\end{equation*}
We have
\(
  \psi^{\pm}_k(z) = q^{\pm (\alpha_k,\alpha_k)h_k/2}
  p_k^\pm(q_k z)/p_k^\pm(q_k^{-1} z).
\)

\begin{Remark}
When $\mathfrak g$ is finite dimensional, then $\min(\langle \alpha_k,
h_l\rangle, \langle \alpha_l, h_k\rangle) = \text{$0$ or $1$}$. Then
the relation~\eqref{eq:relexE1} reduces to the one in literature. Our
generalization seems natural since we will check it later, at least
for symmetric $\mathfrak g$.
\end{Remark}

Let $\Ua^+$ (resp.\ $\Ua^-$) be the $\Q(q)$-subalgebra of $\Ua$
generated by elements $e_{k,r}$'s (resp.\ $f_{k,r}$'s).  Let $\Ua^0$
be the $\Q(q)$-subalgebra generated by elements $q^h$, $h_{k,m}$.

The {\it quantum loop algebra\/} $\Ul$ is the subalgebra
of $\Ua/(q^{\pm c/2} - 1)$ generated by 
$e_{k,r}, f_{k,r}$ ($k\in I$, $r\in\Z$), $q^h$ ($h \in P^*$),
$h_{k,m}$ ($k\in I$, $m\in \Z\setminus\{0\}$),
i.e., generators other than $q^{\pm c/2}$, $q^{\pm d}$.
We will be concerned only with the quantum loop algebra, and not with
the quantum affine algebra in the sequel.

There is a homomorphism $\Uq\to \Ul$ defined by
\begin{equation*}
   q^h \mapsto q^h, \quad
   e_k \mapsto e_{k,0}, \quad f_k \mapsto f_{k,0}.
\end{equation*}

Let $e_{k,r}^{(n)} \defeq e_{k,r}^n / [n]_{q_k}!$,
$f_{k,r}^{(n)} \defeq f_{k,r}^n / [n]_{q_k}!$.
Let $\Uli$ be the $\Z[q,q^{-1}]$-subalgebra generated by
$e_{k,r}^{(n)}$, $f_{k,r}^{(n)}$, $q^h$
and the coefficients of $p_k^\pm(z)$
for $k\in I$, $r\in \Z$, $n\in\Z_{> 0}$, $h\in P^*$.
(It should be true that
$\Uli$ is free over $\Z[q,q^{-1}]$ and
that the natural map $\Uli\otimes_{\Z[q,q^{-1}]}\Q(q) \to \Ul$ is an
isomorphism. But the author does not know how to prove them.)
This subalgebra was introduced by Chari-Pressley \cite{CP3}.
Let $\Uli^+$ (resp.\ $\Uli^-$) be $\Z[q,q^{-1}]$-subalgebra generated
by $e_{k,r}^{(n)}$ (resp.\ $f_{k,r}^{(n)}$) for $k\in I$, $r\in \Z$,
$n\in Z_{> 0}$. We have $\Uli^\pm\subset \Uli$.
Let $\Uli^0$ be the $\Z[q,q^{-1}]$-subalgebra generated by $q^h$, the
coefficients of $p_k^\pm(z)$ and
\begin{equation*}
   \begin{bmatrix}
      q^{h_k}; n \\ r
   \end{bmatrix}
   \defeq
   \prod_{s=1}^r \frac{q^{(\alpha_k,\alpha_k)h_k/2} q_k^{n-s+1}
   - q^{-(\alpha_k,\alpha_k)h_k/2} q_k^{-n+s-1}}{q_k^s - q_k^{-s}}
\end{equation*}
for all $h\in P$, $k\in I$, $n\in \Z$, $r\in \Z_{> 0}$. One can
easily shown that $\Uli^0\subset \Uli$ (see e.g.,
\cite[3.1.9]{Lu-Book}).

For $\varepsilon\in\C^*$, let $\Ule$ be the {\it specialized quantum
loop algebra\/} defined by $\Uli\otimes_{\Z[q,q^{-1}]}\C$ via the
algebra homomorphism $\Z[q,q^{-1}]\to \C$ that takes $q$ to
$\varepsilon$.
We assume $\varepsilon$ is {\it not\/} a root of unity in this paper.
Let $\Ule^\pm$, $\Ule^0$ be the specialization of $\Uli^\pm$, $\Uli^0$ 
respectively. We have a weak form of the triangular decomposition
\begin{equation}
\label{eq:weaktri}
   \Ule = \Ule^-\cdot \Ule^0\cdot \Ule^+,
\end{equation}
which follows from the definition (cf.~\cite[6.1]{CP3}).

We say a $\Ule$-module $M$ is a {\it l-highest weight module\/} (`{\it
l\/}' stands for the loop) with {\it l-highest weight\/} 
$\Psi^\pm(z) = (\Psi^{\pm}_{k}(z))_k\in \C[[z^{\mp}]]^I$
if there exists a vector $m_0\in M$ such that
\begin{align}
  & e_{k,r}\ast m_0 = 0, \qquad \Ule^- \ast m_0 = M, \\
  & \psi^{\pm}_{k}(z)\ast m_0 = \Psi^{\pm}_{k}(z) m_0
    \quad\text{for $k\in I$}.
\end{align}
%
By using \eqref{eq:weaktri} and a standard argument, one can show that 
there is a simple {\it l\/}-highest weight module $M$ of $\Ule$ with
{\it l\/}-highest weight vector $m_0$ satisfying the above for any
$\Psi^{\pm}(z)$. Moreover, such $M$ is unique up to isomorphism.

A $\Ule$-module $M$ is said to be {\it l-integrable\/} if 
\begin{aenume}
\item $M$ has a weight space decomposition
$M = \bigoplus_{\lambda\in P} M_\lambda$ as a $\Ue$-module such that
$\dim M_\lambda < \infty$,
\item for any $m\in M$, there exists $n_0\ge 1$ such that
$e_{k,r_1}\cdots e_{k,r_n} \ast m = 
f_{k,r_1}\cdots f_{k,r_n} \ast m = 0$
for all $r_1,\dots,r_n\in\Z$, $k\in I$ and $n\ge n_0$.
\end{aenume}
For example, if $\mathfrak g$ is finite dimensional, and $M$ is a
finite dimensional module, then $M$ satisfies the above conditions
after twisting with a certain automorphism of $\Ule$
(\cite[12.2.3]{CP}).

\begin{Proposition}\label{prop:DrPol}
Assume that $\mathfrak g$ is symmetric.
The simple {\it l\/}-highest weight $\Ule$-module $M$ with 
{\it l\/}-highest weight $\Psi^{\pm}(z)$ is {\it l\/}-integrable if and
only if there exist polynomials $P_k(u)\in \C[u]$ for $k\in I$ with
$P_k(0) = 1$ such that
\begin{equation}
\label{eq:DrPol}
%
  \Psi^{\pm}_{k}(z) = \varepsilon_k^{\deg P_k}
  \left(\frac{P_k(\varepsilon_k^{-1}/z)}{P_k(\varepsilon_k/z)}\right)^{\pm},
%
\end{equation}
where $\varepsilon_k = \varepsilon^{(\alpha_k,\alpha_k)/2}$,
and
$\left(\ \right)^{\pm}
\in\C[[z^{\mp}]]$ denotes the expansion at $z = \infty$ and $0$
respectively.
\end{Proposition}

This result was announced by Drinfel'd for the Yangian \cite{Drinfeld}.
The proof of the `only if' part when $\mathfrak g$ is finite
dimensional was given by Chari-Pressley~\cite[12.2.6]{CP}. Since the
proof is based on a reduction to the case $\mathfrak g = \algsl_2$, it 
can be applied to a general Kac-Moody algebra $\mathfrak g$
(not necessarily symmetric).
The `if' part was proved by them later in \cite{CP2} when $\mathfrak
g$ is finite dimensional, again not necessarily symmetric.
As an application of the main result of this paper, we will prove the
converse for a symmetric Kac-Moody algebra $\mathfrak g$ in
\secref{sec:std}. Our proof is independent of Chari-Pressley's one.

\begin{Remark}
The polynomials $P_k$ are called {\it Drinfel'd polynomials}.
\end{Remark}

When the Drinfel'd polynomials are given by
\begin{equation*}
   P_k(u) = 
   \begin{cases}
     1 - su  & \text{if $k\neq k_0$,} \\
     1 & \text{otherwise},
   \end{cases}
\end{equation*}
for some $k_0\in I$, $s\in\C^*$, the corresponding simple {\it
l\/}-highest weight module is called a {\it l-fundamental representation}.
When $\mathfrak g$ is finite dimensional, $\Ule$ is a Hopf algebra
since Drinfel'd~\cite{Drinfeld} announced and Beck~\cite{Beck} proved
that $\Ule$ can be identified with (a quotient of) the specialized
quantized enveloping algebra associated with Cartan data of affine
type. Thus a tensor product of $\Ule$-modules is again a
$\Ule$-module. We have the following:

\begin{Proposition}[\protect{\cite[12.2.6,12.2.8]{CP}}]\label{prop:tensor}
Suppose $\mathfrak g$ is finite dimensional.

\rom{(1)} If $M$ and $N$ are simple {\it l\/}-highest weight
$\Ule$-modules with Drinfel'd polynomials $P_{k,M}$, $P_{k,N}$ such
that $M\otimes N$ is simple, then its Drinfel'd polynomial
$P_{k,M\otimes N}$ is given by
\begin{equation*}
   P_{k,M\otimes N} = P_{k,M} P_{k,N}.
\end{equation*}

\rom{(2)} Every simple {\it l\/}-highest weight $\Ule$-modules is a
subquotient of a tensor product of {\it l\/}-fundamental
representations.
\end{Proposition}

Unfortunately the coproduct is not defined for general $\mathfrak g$
as far as the author knows. Thus the above results do not make sense
for general $\mathfrak g$.

\subsection{An {\it l\/}-weight space decomposition}

Let $M$ be an {\it l\/}-integrable $\Ule$-module with the weight space 
decomposition $M = \bigoplus_{\lambda\in P} M_\lambda$. Since the
commutative subalgebra $\Ule^0$ preserves each $M_\lambda$, we can
further decompose $M$ into a sum of generalized simultaneous
eigenspaces for $\Ule^0$:
\begin{equation}\label{eq:gen_wt}
   M = \bigoplus M_{\Psi^\pm},
\end{equation}
where $\Psi^\pm(z) = (\Psi_k^\pm(z))_k\in \C[[z^\mp]]^I$ and
\begin{equation*}
   M_{\Psi^\pm} \defeq \{ m\in M
   \mid \text{$(\psi_k^\pm(z) - \Psi_k^\pm(z)\operatorname{Id})^N \ast m
   = 0$ for $k\in I$ and sufficiently large $N$}\}.
\end{equation*}
If $M_{\Psi^\pm}\neq 0$, we call $M_{\Psi^\pm}$ an {\it l-weight
space}, and the corresponding $\Psi^\pm(z)$ an {\it l-weight}.
Since the constant term of $\psi_k^\pm(z)$ is $q^{\pm h_k}$, this is a 
refinement of the weight space decomposition.
A further study of the {\it l\/}-weight space decomposition will be
given in \subsecref{subsec:FR}.

Motivated by \propref{prop:DrPol}, we introduce the following notion:
\begin{Definition}\label{def:l-dominant}
Fix $\varepsilon\in\C^*$. An {\it l\/}-weight $\Psi^\pm(z) =
(\Psi_k^\pm(z))_k$ is said to be {\it {\it l\/}-dominant\/} if there
exist a polynomial $P(u) = (P_k(u))_k\in \C[u]^I$ for with $P_k(0) =
1$ such that \eqref{eq:DrPol} holds.
\end{Definition}

Thus \propref{prop:DrPol} means that a {\it l\/}-highest weight module 
is {\it l\/}-integrable if and only if the {\it l\/}-highest weight is 
{\it l\/}-dominant.

\section{Quiver variety}\label{sec:quiver_variety}

\subsection{Notation}
Suppose that a finite graph is given and assume that there are no edge
loops (i.e., no edges joining a vertex with itself).
Let $I$ be the set of vertices and $E$ the set of edges.
Let $\bA$ be the adjacency matrix of the graph, namely
\begin{equation*}
   \bA = (\bA_{kl})_{k,l\in I}, \qquad
   \text{where $\bA_{kl}$ is the number of edges joining $k$ and $l$}.
\end{equation*}
We associate with the graph $(I, E)$ a symmetric generalized Cartan matrix 
$\bC = 2\bI - \bA$, where $\bI$ is the identity matrix.
This gives a bijection between the finite graphs without edge loops and
symmetric Cartan matrices.
We have the corresponding symmetric Kac-Moody algebra $\frak g$, the
quantized enveloping algebra $\Uq$, the quantum affine algebra $\Ua$
and the quantum loop algebra $\Ul$.
Let $H$ be the set of pairs consisting of an edge together with its
orientation. For $h\in H$, we denote by $\vin(h)$ (resp.\ $\vout(h)$)
the incoming (resp.\ outgoing) vertex of $h$.
For $h\in H$ we denote by $\overline h$ the same edge as $h$ with the
reverse orientation.
Choose and fix an orientation $\Omega$ of the graph,
i.e., a subset $\Omega\subset H$ such that
$\overline\Omega\cup\Omega = H$, $\Omega\cap\overline\Omega = \emptyset$.
The pair $(I,\Omega)$ is called a {\it quiver}.
Let us define matrices $\bA_\Omega$ and $\bA_{\overline\Omega}$ by
\begin{equation}\label{eq:Aomega}
\begin{split}
   (\bA_\Omega)_{kl} 
   & \defeq \# \{ h\in \Omega \mid \vin(h) = k, \vout(h) = l\}, \\
   (\bA_{\overline\Omega})_{kl} 
   & \defeq \# \{ h\in \overline{\Omega} \mid \vin(h) = k, \vout(h) = l\}.
\end{split}
\end{equation}
So we have $\bA = \bA_\Omega + \bA_{\overline\Omega}$,
$\lsp{t}\bA_\Omega = \bA_{\overline\Omega}$.  

Let $V = (V_k)_{k\in I}$ be a collection of finite-dimensional vector
spaces over $\C$ for each vertex $k\in I$. 
The dimension of $V$ is a vector
\[
  \dim V = (\dim V_k)_{k\in I}\in \Bbb Z_{\ge 0}^I.
\]

If $V^1$ and $V^2$ are such collections, we define vector spaces by
\begin{equation}\label{eq:LE}
  \HomL(V^1, V^2) \defeq
  \bigoplus_{k\in I} \Hom(V^1_k, V^2_k), \quad
  \HomE(V^1, V^2) \defeq
  \bigoplus_{h\in H} \Hom(V^1_{\vout(h)}, V^2_{\vin(h)})
\end{equation}

For $B = (B_h) \in \HomE(V^1, V^2)$ and 
$C = (C_h) \in \HomE(V^2, V^3)$, let us define a multiplication of $B$
and $C$ by
\[
  CB \defeq (\sum_{\vin(h) = k} C_h B_{\overline h})_k \in
  \HomL(V^1, V^3).
\]
Multiplications $ba$, $Ba$ of $a\in \HomL(V^1,V^2)$, $b\in\HomL(V^2,
V^3)$, $B\in \HomE(V^2, V^3)$ is defined in obvious manner. If
$a\in\HomL(V^1, V^1)$, its trace $\tr(a)$ is understood as $\sum_k
\tr(a_k)$.

For two collections $V$, $W$ of vector spaces
with $\bv = \dim V$, $\bw = \dim W$, we consider the vector space
given by
\begin{equation}
  \bM \equiv \bM(\bv, \bw) \defeq
  \HomE(V, V) \oplus \HomL(W, V) \oplus \HomL(V, W),
\label{def:bM}\end{equation}
where we use the notation $\bM$ unless we want to specify dimensions
of $V$, $W$.
The above three components for an element of $\bM$ will be denoted by 
$B$, $i$, $j$ respectively.
An element of $\bM$ will be called an {\it ADHM datum\/}.

Usually a point in $\bigoplus_{h\in\Omega} \Hom(V^1_{\vout(h)},
V^2_{\vin(h)})$ is called a {\it representation of the quiver\/}
$(I,\Omega)$ in literature. Thus $\HomE(V,V)$ is the product of the
space of representations of $(I,\Omega)$ and that of
$(I,\overline\Omega)$.
On the other hand, the factor $\HomL(W,V)$ or $\HomL(V,W)$ is never
appeared in literature.

\begin{Convention}\label{convention:weight}
When we relate the quiver varieties to the quantum affine algebra,
the dimension vectors will be mapped into the weight lattice in the 
following way:
\begin{equation*}
    \bv \mapsto \sum_k v_k \alpha_k, \quad
    \bw \mapsto \sum_k w_k \Lambda_k,
\end{equation*}
where $v_k$ (resp.\ $w_k$) is the $k$th component of $\bv$ (resp.\ 
$\bw$). Since $\{ \alpha_k \}$ and $\{ \Lambda_k \}$ are both
linearly independent, these maps are injective. We consider $\bv$ and
$\bw$ as elements of the weight lattice $P$ in this way hereafter.
\end{Convention}

For a collection $S = (S_k)_{k\in I}$ of subspaces of $V_k$ and
$B\in\HomE(V, V)$, we say $S$ is {\it $B$-invariant\/} if
$B_h(S_{\vout(h)}) \subset S_{\vin(h)}$.

Fix a function $\varepsilon\colon H \to \C^*$ such that
$\varepsilon(h) + \varepsilon(\overline{h}) = 0$ for all $h\in H$.
In \cite{Na-quiver,Na-alg}, it was assumed that $\varepsilon$ takes
its value $\pm 1$, but this assumption is not necessary as remarked by 
Lusztig~\cite{Lu-qv}.
For $B\in\HomE(V^1, V^2)$, let us denote
by $\varepsilon B\in \HomE(V^1, V^2)$ data given by $(\varepsilon B)_h 
= \varepsilon(h) B_h$ for $h\in H$.

Let us define a symplectic form $\omega$ on $\bM$ by
\begin{equation}
        \omega((B, i, j), (B', i', j'))
        \defeq \tr(\varepsilon B\, B') + \tr(i j' - i' j).
\label{def:symplectic}\end{equation}

Let $G$ be the algebraic group defined by
\begin{equation*}
  G \equiv G_{\bv} \defeq \prod_k \GL(V_k),
\end{equation*}
where we use the notation $G_\bv$ when we want to emphasize the
dimension.
It acts on $\bM$ by
\begin{equation}\label{eq:Gaction}
  (B, i, j) \mapsto g\cdot (B, i, j) \defeq (g B g^{-1}, g i, j g^{-1})
\end{equation}
preserving the symplectic form $\omega$. The moment map
$\mu\colon\bM\to\HomL(V, V)$ vanishing at the origin is given by
\begin{equation}\label{eq:mu}
  \mu(B, i, j) = \varepsilon B B + i j,
\end{equation}
where the dual of the Lie algebra of $G$ is identified with the Lie
algebra via the trace.
Let $\mu^{-1}(0)$ be an affine algebraic variety (not necessarily
irreducible) defined as the zero set of $\mu$.

For $(B,i,j)\in\mu^{-1}(0)$, we consider the following complex
\begin{equation}\label{eq:quiver_tangent}
        \HomL(V, V)
        \overset{\iota}{\longrightarrow}
        \HomE(V, V) \oplus \HomL(W, V) \oplus \HomL(V,W)
        \overset{d \mu}{\longrightarrow}
        \HomL(V, V),
\end{equation}
where $d\mu$ is the differential of $\mu$ at $(B,i,j)$, and $\iota$ is
given by
\begin{equation*}
        \iota(\xi) = (B \xi - \xi B) \oplus
         (-\xi i) \oplus j \xi.
\end{equation*}
If we identify $\HomE(V, V) \oplus \HomL(W, V) \oplus \HomL(V,W)$ with
its dual via the symplectic form $\omega$,
$\iota$ is the transpose of $d \mu$.

\subsection{Two quotients $\M_0$ and $\M$}\label{subsec:twoquot}

We consider two types of quotients of $\mu^{-1}(0)$ by the group $G$.
The first one is the affine algebro-geometric quotient given as
follows.  Let $A(\mu^{-1}(0))$ be the coordinate ring of the affine
algebraic variety $\mu^{-1}(0)$.  Then $\M_0$ is defined as a variety
whose coordinate ring is the invariant part of $A(\mu^{-1}(0))$:
\begin{equation}
  \M_0 \equiv \M_0(\bv,\bw) \defeq
  \mu^{-1}(0)\dslash G = \Spec A(\mu^{-1}(0))^{G}.
\end{equation}
As before, we use the notation $\M_0$ unless we need to specify the
dimension vectors $\bv$, $\bw$.
By the geometric invariant theory \cite{GIT}, this is an affine
algebraic variety.
It is also known that the geometric points of $\M_0$ are closed $G$-orbits.

For the second quotient we follow A.~King's approach \cite{King}.
Let us define a character $\chi\colon G\to \C^*$ by
$\chi(g) = \prod_k \det g_k^{-1}$ for $g = (g_k)$.
Set
\begin{equation*}
  A(\mu^{-1}(0))^{G, \chi^n} \defeq
  \{\, f\in A(\mu^{-1}(0)) \mid f(g (B, i, j)) = \chi(g)^n f(B, i, j)
  \,\} .
\end{equation*}
The direct sum with respect to $n\in\Bbb Z_{\ge 0}$ is a graded
algebra, hence we can define
\begin{equation}
  \M \equiv \M(\bv,\bw) \defeq
  \Proj \bigoplus_{n\ge 0} A(\mu^{-1}(0))^{G, \chi^n}.
\end{equation}
These are what we call {\it quiver varieties}.
\subsection{Stability Condition}\label{subsec:stability}
In this subsection, we shall give a description of the quiver variety $\M$
which is easier to deal with.
We again follow King's work \cite{King}.

\begin{Definition}\label{def:stable}
A point $(B, i, j) \in \mu^{-1}(0)$ is said to be {\it stable\/} if 
the following condition holds:
\begin{itemize}
\item[] if a collection $S = (S_k)_{k\in I}$ of subspaces of $V =
(V_k)_{k\in I}$ is $B$-invariant and contained in 
$\Ker j$, then $S = 0$.
\end{itemize}
Let us denote by $\mu^{-1}(0)^{\operatorname{s}}$ the set of stable points.
\end{Definition}
Clearly, the stability condition is invariant under the action of
$G$. Hence we may say an orbit is stable or not.

Let us lift the $G$-action on $\mu^{-1}(0)$ to the trivial line bundle
$\mu^{-1}(0) \times \C$ by $g\cdot(B,i,j,z) = (g\cdot (B,i,j),
\chi^{-1}(g)z)$.

We have the following:
\begin{Proposition}\label{prop:stable}
\rom{(1)} A point $(B,i,j)$ is stable if and only if the closure of 
$G\cdot(B,i,j,z)$ does not intersect with the zero section
of $\mu^{-1}(0) \times \C$ for $z\neq 0$.

\rom{(2)} If $(B,i,j)$ is stable, then the differential
$d\mu\colon \bM \to \HomL(V,V)$ is surjective. In particular,
$\mu^{-1}(0)^{\operatorname{s}}$ is a nonsingular variety.

\rom{(3)} If $(B,i,j)$ is stable, then $\iota$ in
\eqref{eq:quiver_tangent} is injective.

\rom{(4)} The quotient $\mu^{-1}(0)^{\operatorname{s}}/G$ has a
structure of nonsingular quasi-projective variety
of dimension $(\bv, 2\bw - \bv)$, and
$\mu^{-1}(0)^{\operatorname{s}}$ is a principal $G$-bundle
over $\mu^{-1}(0)^{\operatorname{s}}/G$.

\rom{(5)} The tangent space of $\mu^{-1}(0)^{\operatorname{s}}/G$ at
the orbit $G\cdot(B,i,j)$ is isomorphic to the middle cohomology group 
of \eqref{eq:quiver_tangent}.

\rom{(6)} The variety $\M$ is isomorphic to
$\mu^{-1}(0)^{\operatorname{s}}/G$.

\rom{(7)} $\mu^{-1}(0)^{\operatorname{s}}/G$ has a holomorphic
symplectic structure as a symplectic quotient.
\end{Proposition}

\begin{proof}
See \cite[3.ii]{Na-alg} and \cite[2.8]{Na-quiver}.
\end{proof}

\begin{Notation}
For a stable point $(B, i, j)\in\mu^{-1}(0)$, its $G$-orbit
considered as a geometric point in the quiver variety $\M$ is denoted by 
$[B, i, j]$.
If $(B, i, j)\in\mu^{-1}(0)$ has a closed $G$-orbit, then the
corresponding geometric point in $\M_0$ will be denoted also by 
$[B, i, j]$.
\end{Notation}

From the definition, we have a natural projective morphism
(see \cite[3.18]{Na-alg})
\begin{equation}\label{eq:pi_def}
      \pi\colon \M \to \M_0.
\end{equation}
If $\pi([B,i,j]) = [B^0,i^0,j^0]$, then $G\cdot(B^0, i^0, j^0)$ 
is the unique closed orbit contained in the
closure of $G\cdot(B, i, j)$.
For $x\in \M_0$, let
\begin{equation}
\label{eq:M_x}
   \M_x \defeq \pi^{-1}(x).
\end{equation}
If we want to specify the dimension, we denote the above by
$\M(\bv,\bw)_x$.
Unfortunately, this notation conflicts with the previous notation
$\M_0$ when $x = 0$.  And the central fiber $\pi^{-1}(0)$ plays an
important role later.  We shall always write $\La \equiv \La(\bv,\bw)$
for $\pi^{-1}(0)$ and not use the notation \eqref{eq:M_x} with $x =
0$.

In order to explain more precise relation between $[B,i,j]$ and
$[B^0,i^0,j^0]$, we need the following notion.

\begin{Definition}
Suppose that $(B, i, j)\in\bM$ and a $B$-invariant increasing filtration 
\[
       0 = V^{(-1)} \subset V^{(0)} \subset \cdots \subset V^{(N)} = V,
\]
with $\Ima i \subset V^{(0)}$ are given.
Then set $\gr_m V = V^{(m)}/V^{(m-1)}$ and $\gr V = \bigoplus \gr_m V$.
Let $\gr_m B$ denote the endomorphism which $B$ induces
on $\gr_m V$.
For $m = 0$, let $\gr_0 i\in \HomL(W, V^{(0)})$ be such that its
composition with the inclusion $V^{(0)} \subset V$ is $i$, and
$\gr_0 j$ be the restriction of $j$ to $V^{(0)}$.
For $m \ne 0$, set $\gr_m i = 0$, $\gr_m j = 0$.
Let $\gr (B, i, j)$ be the direct sum of $(\gr_m B, \gr_m i, \gr_m j)$
considered as data on $\gr V$.
\label{quiver:gr}\end{Definition}

\begin{Proposition}
Suppose $\pi(x) = y$.
Then there exist a representative $(B, i, j)$ of $x$ and a $B$-invariant
increasing filtration $V^{(*)}$ as in Definition~\ref{quiver:gr} such
that $\gr (B, i, j)$ is a representative of $y$ on $\gr V$.
\label{quiver:pi}\end{Proposition}
\begin{proof}
See \cite[3.20]{Na-alg}
\end{proof}

\begin{Proposition}\label{prop:homotopic}
$\La$ is a Lagrangian subvariety which is homotopic to $\M$.
\end{Proposition}
\begin{proof}
See \cite[5.5, 5.8]{Na-quiver}.
\end{proof}

\subsection{Hyper-K\"ahler structure}\label{subsec:hK}

We briefly recall hyper-K\"ahler structures on $\M$, $\M_0$. This view 
point was used for the study of $\M$, $\M_0$ in \cite{Na-quiver}.
({\bf Caution}: The following notation is different from the original
one. $K_\bv$ and $G_\bv$ were denoted by $G_\bv$ and
$G_\bv^\C$ respectively in \cite{Na-quiver}.
$\mu$ in \eqref{eq:mu} was denoted by $\mu_\C$ and the pair
$(\mu_\R,\mu)$ was denoted by $\mu$ in \cite{Na-quiver}.)

Put and fix hermitian inner products on $V$ and $W$. They together
with an orientation $\Omega$ induce a hermitian inner product and a
quaternion structure on $\bM$ (\cite[p.370]{Na-quiver}). Let $K_\bv$
be a compact Lie group defined by $K_\bv = \prod_k \U(V_k)$.
This is a maximal compact subgroup of $G_\bv$, and acts on $\bM$
preserving the hermitian and quaternion structures. The corresponding
hyper-K\"ahler moment map vanishing at the origin decomposes into
the complex part 
$\mu\colon \bM\to \bigoplus_k \gl(V_k) = \HomL(V,V)$
(defined in \eqref{eq:mu})
and
the real part
$\mu_\R\colon \bM \to \bigoplus_k \algu(V_k)$,
where
\begin{equation*}
 \mu_\R(B,i,j) = {i\over 2}\left(
        \sum_{h\in H: k = \vin(h)} B_h B_h^\dagger
         - B_{\overline h}^\dagger B_{\overline h}
         + i_k i_k^\dagger - j_k^\dagger j_k\right)_k .
\end{equation*}
\begin{Proposition}\label{prop:hK}
\rom{(1)} A $G_\bv$-orbit $[B,i,j]$ in $\mu^{-1}(0)$ intersects with
$\mu_\R^{-1}(0)$ if and only if it is closed. The map
\begin{equation*}
   \left(\mu_\R^{-1}(0)\cap \mu^{-1}(0)\right) / K_\bv
   \to \mu^{-1}(0)\dslash G_\bv = \M_0(\bv,\bw)
\end{equation*}
is a homeomorphism.

\rom{(2)} Choose a parameter $\zeta_\R = (\zeta_\R^{(k)})_k \in
\R^I$ so that $\zeta_\R^{(k)}\in \sqrt{-1}\R_{> 0}$. Then
a $G_\bv$-orbit $[B,i,j]$ in $\mu^{-1}(0)$ intersects with
with $\mu_\R^{-1}(-\zeta_\R)$ if and only if it is stable. The map
\begin{equation*}
   \left(\mu_\R^{-1}(-\zeta_\R)\cap \mu^{-1}(0)\right) / K_\bv
   \to \mu^{-1}(0)^{\operatorname{s}}/G = \M(\bv,\bw)
\end{equation*}
is a homeomorphism.
\end{Proposition}

\begin{proof}
See \cite[3.1,3.2,3.5]{Na-quiver}
\end{proof}

\subsection{}\label{subsec:M_0}

Suppose $V = (V_k)_{k\in I}$ is a collection of subspace of 
$V' = (V'_k)_{k\in I}$ and $(B,i,j)\in\mu^{-1}(0)\subset\bM(V,W)$ is
given. We can extend $(B,i,j)$ to $\bM(V',W)$ by letting it $0$ on
a complementary subspace of $V$ in $V'$. This operation induces a
natural morphism
\begin{equation}\label{eq:incl'}
   \text{$\mu^{-1}(0)$ in $\bM(\bv,\bw)$} \to
   \text{$\mu^{-1}(0)$ in $\bM(\bv',\bw)$},
\end{equation}
where $\bv' = \dim V'$. This induces a morphism
\begin{equation}\label{eq:incl}
   \M_0(\bv,\bw) \to \M_0(\bv',\bw).
\end{equation}
Moreover, we also have a map
\begin{equation*}
   \text{$\mu^{-1}(0)\cap\mu_\R^{-1}(0)$ in $\bM(\bv,\bw)$} \to
   \text{$\mu^{-1}(0)\cap\mu_\R^{-1}(0)$ in $\bM(\bv',\bw)$}.
\end{equation*}
Thus closed $G_\bv$-orbits in $\mu^{-1}(0)\subset \bM(\bv,\bw)$ is
mapped to closed $G_{\bv'}$-orbits in $\mu^{-1}(0)\subset \bM(\bv',\bw)$
by \propref{prop:hK}(1).

The following lemma was stated in \cite[p.529]{Na-alg} without proof.

\begin{Lemma}\label{lem:incl_inj}
The morphism \eqref{eq:incl} is injective.
\end{Lemma}

\begin{proof}
Suppose $x^1$, $x^2\in\M_0(\bv,\bw)$ have the same 
image under \eqref{eq:incl}. We choose representatives
$(B^1,i^1,j^1)$, $(B^2,i^2,j^2)$ which have closed $G_{\bv}$-orbit.

Let us define $S^a = (S^a_k)_{k\in I}$ ($a = 1,2$) by
\begin{equation*}
  S^a_k \defeq
  \Ima\left(\sum_{\vin(h)=k} \varepsilon(h) B^a_h + i^a_k\right).
\end{equation*}
Choose complementary subspaces $T_k^a$ of $S^a_k$ in $V_k$. We
choose a $1$-parameter subgroup $\lambda^a\colon\C^*\to G_{\bv}$ as
follows: $\lambda^a(t) = 1$ on $S^a_k$ and $\lambda^a(t) = t^{-1}$ on
$T_k^a$. Then the limit $\lambda^a(t)\cdot (B^a,i^a,j^a)$ exists and
its restriction to $T_k^a$ is $0$. Since $(B^a,i^a,j^a)$ has a closed
orbit, we may assume that the restriction of $(B^a,i^a,j^a)$ to
$T_k^a$ is $0$. Note that $S^a$ is a subspace of $V$ by the construction.

Suppose that there exists $g'\in G_{\bv'}$ such that $g'\cdot
(B^1,i^1,j^1) = (B^2,i^2,j^2)$. We want to construct $g\in G_\bv$ such
that $g\cdot (B^1,i^1,j^1) = (B^2,i^2,j^2)$.  Since we have $g'(S^1) =
S^2$, the the restriction of $g'$ to $S^1$ is invertible. Let $g$ be
an extension of the restriction $g'|_{S^1}$ to $V$ so that $T^1$ is
mapped to $T^2$. Then $g\in G_\bv$ maps $(B^1,i^1,j^1)$ to
$(B^2,i^2,j^2)$.
\end{proof}

Hereafter, we consider $\M_0(\bv,\bw)$ as a subset of
$\M_0(\bv',\bw)$. It is clearly a closed subvariety. Let
\begin{equation}
\label{eq:union}
   \M_0(\infty,\bw) \defeq \bigcup_{\bv} \M_0(\bv,\bw).
\end{equation}
If the graph is of finite type, $\M_0(\bv,\bw)$ stabilizes at some
$\bv$ (see \propref{prop:ADE} and \lemref{lem:surj}(2) below).
This is {\it not\/} true in general. However, it has no harm in
this paper. We use $\M_0(\infty,\bw)$ to simplify the notation, and do
not need any structures on it.
We can always work on individual $\M_0(\bv,\bw)$,
not on $\M_0(\infty,\bw)$.

Later, we shall also study $\M(\bv,\bw)$ for various
$\bv$ simultaneously. We introduce the following notation:
\begin{equation*}
   \Mw \defeq \bigsqcup_{\bv}\M(\bv,\bw), \qquad
   \Law \defeq \bigsqcup_{\bv}\La(\bv,\bw).
\end{equation*}
Note that there are no obvious morphisms between
$\M(\bv,\bw)$ and $\M(\bv',\bw)$ since the stability condition is {\it 
not\/} preserved under \eqref{eq:incl'}.

\subsection{Definition of $\Mreg$}

Let us introduce an open subset of $\M_0$ (possibly empty):
\begin{equation}\label{quiver:reg}
  \Mreg \equiv \Mreg(\bv,\bw) \defeq \{\, [B, i, j]\in \M_0\mid
  \text{$(B, i, j)$ has the trivial stabilizer in $G$}\,\}.
\end{equation}

\begin{Proposition}\label{prop:pi_isom}
If $[B, i, j]\in \Mreg$, then it is stable.
Moreover, $\pi$ induces an isomorphism 
$\pi^{-1}(\Mreg)\simeq \Mreg$
\end{Proposition}

\begin{proof}
See \cite[3.24]{Na-alg} or \cite[4.1(2)]{Na-quiver}.
\end{proof}

As in \subsecref{subsec:M_0}, we consider $\Mreg(\bv,\bw)$ as a subset 
of $\M_0(\bv',\bw)$ when $\bv' - \bv \in \sum_k \Z_{\ge 0}
\alpha_k$. Then we have
\begin{Proposition}\label{prop:ADE}
If the graph is of type $ADE$, then
\begin{equation*}
   \M_0(\bv',\bw) = \bigcup_{\bv} \Mreg(\bv,\bw),
\end{equation*}
where the summation runs over the set of $\bv$ such that
$\bv' - \bv \in \sum_k \Z_{\ge 0}\alpha_k$.
\end{Proposition}

\begin{proof}
See \cite[6.7]{Na-quiver}, \cite[3.28]{Na-alg}.
\end{proof}

\begin{Definition}\label{def:regular}
We say a point $x\in\M_0(\infty,\bw)$ is {\it regular\/} if it is
contained in $\Mreg(\bv,\bw)$ for some $\bv$. The above proposition
says that all points are regular if the graph is of type $ADE$. But
this is {\it not\/} true in general (see \cite[10.10]{Na-alg}).
\end{Definition}

\subsection{$G_\bw\times \C^*$-action}\label{subsec:Gw-action}
Let us define a $G_\bw\times\C^*$-action on $\M$ and $\M_0$, where
$G_\bw = \prod_{k\in I} \GL(W_k)$. ({\bf Caution}: We use the same
notation $G_\bv$ and $G_\bw$, but their roles are totally different.)

The $G_\bw$-action is simply defined by its natural action on
$\bM = \HomE(V,V)\oplus\HomL(W,V)\oplus \HomL(V,W)$. It preserves the
equation $\varepsilon B B + ij = 0$ and commutes with the $G$-action
given by \eqref{eq:Gaction}. Hence it induces an action on $\M$ and
$\M_0$.

The $\C^*$-action is slightly different from the one given in
\cite[3.iv]{Na-alg}, and we need extra data. For each pair $k,l\in I$
such that $b' = -(\alpha_k,\alpha_l) \ge 1$, we introduce and fix a
numbering $1,2,\dots,b'$ on edges joining $k$ and $l$. It induces a
numbering $h_1, \dots, h_{b'}$, $\overline{h_1}, \dots,
\overline{h_{b'}}$ on oriented edges between $k$ and $l$. Let us
define $m\colon H\to\Z$ by
\begin{equation}\label{eq:m(h)}
  m(h_p) = b'+1-2p, \quad
  m(\overline h_p) = -b'-1+2p.
\end{equation}
Then we define a $\C^*$-action on $\bM$ by
\begin{equation}\label{eq:C*action}
B_{h} \mapsto t^{m(h)+1} B_{h}, \quad
i \mapsto t i, \quad j \mapsto t j \qquad
\text{for $t\in \C^*$}.
\end{equation}
The equation $\varepsilon B B + ij = 0$ is preserved since the left
hand side is multiplied by $t^2$. It commutes with the $G$-action and
preserves the stability condition.
Hence it induces a $\C^*$-action on $\M$ and $\M_0$.
This $G_\bw\times\C^*$-action makes the projective morphism
$\pi\colon \M\to\M_0$ equivariant.

In order to distinguish this $G_\bw\times\C^*$-action from the
$G_\bv$-action \eqref{eq:Gaction}, we denote it as
\begin{equation*}
   (B,i,j) \mapsto h\star (B,i,j) \qquad (h\in G_\bw\times\C^*).
\end{equation*}

\subsection{Notation for $\C^*$-action}

For an integer $m$, we define a $\C^*$-module structure on $\C$ by 
\begin{equation}\label{eq:L(m)}
   t\cdot v \defeq t^m v \qquad\text{$t\in\C^*, v\in \C$},
\end{equation}
and denote it by $L(m)$. For a a $\C^*$-module $V$, we use the
following notational convention:
\begin{equation}\label{eq:Cst_notation}
   q^m V \defeq L(m)\otimes V.
\end{equation}
We use the same notation when $V$ is an element of $\C^*$-equivariant
$K$-theory later.

\subsection{Tautological bundles}
By the construction of $\M$, we have a natural vector bundle
\begin{equation*}
   \mu^{-1}(0)^{\operatorname{s}}\times_G V_k \to \M
\end{equation*}
associated with the principal $G$-bundle
$\mu^{-1}(0)^{\operatorname{s}}\to \M$.
For the abuse of notation, we
denote it also by $V_k$. It naturally have a structure of a
$\C^*$-equivariant vector bundle. Letting $G_\bw$ act trivially, we
make it a $G_\bw\times \C^*$-equivariant vector bundle.

The vector space $W_k$ is also considered as 
a $G_\bw\times \C^*$-equivariant vector bundle, where $G_\bw$ acts
naturally and $\C^*$ acts trivially.

We call $V_k$ and $W_k$ {\it tautological bundles}.

We consider $\HomE(V,V)$, $\HomL(W,V)$, $\HomL(V, W)$ as vector
bundles defined by the same formula as in \eqref{eq:LE}. By the
definition of tautologcial bundles, $B$, $i$, $j$ can be considered as 
sections of those bundles.
Those bundles naturally have structures of $G_\bw\times
\C^*$-equivariant vector bundles. But we modify the $\C^*$-action on
$\HomE(V,V)$ by letting $t\in \C^*$ acts by $t^{m(h)}$ on the
component \( \Hom(V_{\vout(h)}, V_{\vin(h)}) \).  This makes $B$ an
equivariant section of $\HomE(V,V)$.

We consider the following $G_\bw\times\C^*$-equivariant complex 
$C_k^\bullet \equiv C_k^\bullet(\bv,\bw)$ over 
$\M \equiv \M(\bv,\bw)$ (cf.\ \cite[4.2]{Na-alg}):
\begin{equation}\label{eq:taut_cpx}
C_k^\bullet \equiv C_k^\bullet(\bv,\bw):
\begin{CD}
  q^{-2} V_k
  @>{\sigma_k}>>
  q^{-1} \left( \displaystyle{\bigoplus_{l:k\neq l}}
     [-\langle h_k,\alpha_l\rangle]_q V_l
    \oplus W_k\right)
  @>{\tau_k}>>
  V_k,
\end{CD}
\end{equation}
where
\begin{equation*}
\sigma_k = \bigoplus_{\vin(h)=k} B_{\overline h} \oplus j_k, \qquad
\tau_k = \sum_{\vin(h)=k} \varepsilon(h) B_h + i_k.
\end{equation*}
Let us explain the factor $[-\langle h_k,\alpha_l\rangle]_q V_l$.
Set $b' = -\langle h_k,\alpha_l\rangle$.
Since the $\C^*$-action in \eqref{eq:C*action} is defined so that
\begin{equation*}
   \bigoplus_{h:\substack{\vin(h) = k\\ \vout(h)= l}} \Hom(V_k, V_l)
   = \Hom(V_k,V_l)^{\oplus b'}
\end{equation*}
has weights $b', b'-2, \dots, 2 - b'$, the above can be written as
\begin{equation*}
   \left(q^{b'} + q^{b'-2} + \dots + q^{2-b'}\right) \Hom(V_k,V_l)
   = q[b']_q \Hom(V_k,V_l)
\end{equation*}
in the notation~\eqref{eq:Cst_notation}. By the same reason
$C_k^\bullet$ is an equivariant complex.

We assign degree $0$ to the middle term. (This complex is 
the complex in \cite[4.2]{Na-alg} with a modification of the
$G_\bw\times\C^*$-action.)

\begin{Lemma}\label{lem:betasurj}
Fix a point $[B,i,j]$ and consider $C_k^\bullet$ as a complex 
of vector spaces. Then $\sigma_k$ is injective.
\end{Lemma}

\begin{proof}
See \cite[p.530]{Na-alg}. (Lemma~54 therein is a misprint of Lemma~5.2.)
\end{proof}

Note that $\tau_k$ is {\it not\/} surjective in general.
In fact, the following notion will play a crucial role later. Let $X$ be
an irreducible component of $\pi^{-1}(x)$ for $x\in\M_0$. Considering
$\tau_k$ at a generic element $[B,i,j]$ of $X$, we set
\begin{equation}\label{eq:epsilon_k}
   \varepsilon_k(X)
   \defeq \codim_{V_k} \Ima \tau_k.
\end{equation}

\begin{Lemma}\label{lem:surj}
\rom{(1)} Take and fix a point $[B,i,j]\in\M(\bv,\bw)$.
Let $\tau_k$ as in \eqref{eq:taut_cpx}.
If $[B, i, j]\in\pi^{-1}(\Mreg(\bv,\bw))$,
then we have
\begin{equation}\label{eq:decomp_surj}
     \Ima \tau_k
     = V_k \qquad
     \text{for any $k\in I$}.
\end{equation}
Moreover, the converse holds if we assume
$\pi([B,i,j])$ is regular in the sense of \defref{def:regular}. Namely
under this assumption, $[B,i,j]\in\pi^{-1}(\Mreg(\bv,\bw))$ 
if and only if \eqref{eq:decomp_surj} holds.

\rom{(2)} If $\Mreg(\bv,\bw)\neq\emptyset$, then $\bw-\bv$ is dominant.
\end{Lemma}


\begin{proof}
(1) See \cite[4.7]{Na-alg} for the first assertion.
During the proof of \cite[7.2]{Na-alg}, we have shown the second
assertion, using \cite[3.10]{Na-alg} = \propref{quiver:pi}.

(2) Consider the alternating sum of dimensions of the complex
$C_k^\bullet$. It is equal to the alternating sum of dimensions of
cohomology groups. It is nonnegative, if $\Mreg(\bv,\bw)\neq\emptyset$ 
by \lemref{lem:betasurj} and (1). On the other hand, it is equal to
\begin{equation*}
   \sum_{h:\vout(h)=k} \dim V_{\vin(h)} + W_k - 2\dim V_k
   = \langle \bw - \bv, h_k\rangle.
\end{equation*}
Thus we have the assertion.
\end{proof}


\section{Stratification of $\M_0$}\label{sec:strat}

As was shown in \cite[\S6]{Na-quiver}, \cite[3.v]{Na-alg}, there
exists a natural stratification of $\M_0$ by conjugacy classes of
stabilizers.
A local topological structure of a neiborhood of a point in a stratum
(e.g., the homology group of the fiber of $\pi$) was studied in
\cite[6.10]{Na-quiver}. We give a refinement in this section. We
define a slice to a stratum, and study a local structure as a complex
analytic space.
Our technique is based on a work of Sjamaar-Lerman~\cite{SL} in the
symplectic geometry and hence our transversal slice may not be
algebraic. It is desirable to have a purely algebraic construction of
a transversal slice, as Maffei did in a special case~\cite{Maf}.

We fix dimension vectors $\bv$, $\bw$ and denote
$\bM(\bv,\bw)$, $\M(\bv,\bw)$ by $\bM$, $\M$ in this section.

\subsection{Stratification}
\begin{Definition}\label{def:strat}(cf.\ Sjamaar-Lerman \cite{SL})
For a subgroup $\widehat G$ of $G$ denote by
$\bM_{(\widehat G)}$ the set 
of all points in $\bM$ whose stabilizer is conjugate to $\widehat G$.
A point $[(B,i,j)]\in\M_0$ is said to be of {\it $G$-orbit type\/}
$(\widehat G)$ if its representative $(B,i,j)$ is in
$\bM_{(\widehat G)}$.
The set of all points of orbit type $(\widehat G)$ is denoted by
$(\M_0)_{(\widehat G)}$.
\end{Definition}

The stratum $(\M_0)_{(1)}$ corresponding to the trivial subgroup
$1$ is $\Mreg$ by definition.
We have the following decomposition of $\M_0$:
\[
        \M_0 = \bigcup_{(\widehat G)} (\M_0)_{(\widehat G)},
\]
where the summation runs over the set of all conjugacy classes of
subgroups of $G$.

For more detailed description of $(\M_0)_{(\widehat G)}$, see
\cite[6.5]{Na-quiver}, \cite[3.27]{Na-alg}.

\subsection{Local normal form of the moment map}
Let us recall the local normal form of the moment map following
Sjamaar-Lerman~\cite{SL}.

Take $x\in\M_0$ and fix its representative $m =
(B,i,j)\in\mu^{-1}(0)$. We suppose $m$ has a closed $G$-orbit and
satisfies $\mu_\R(m) = 0$ by \propref{prop:hK}(1).
Let $\widehat G$ be the stabilizer of $m$. It is the complexification
of the stabilizer in $K = \prod \U(V_k)$ (see e.g., \cite[1.6]{Sj}).
Since $\mu(m) = 0$, the $G$-orbit $Gm = G/\widehat G$ through $m$ is
an isotropic submanifold of $\bM$.
Let $\widehat\bM$ be the quotient vector space
$(T_m Gm)^\omega/T_m Gm$, where 
$T_m Gm$ is the tangent space of the orbit $Gm$, and $(T_m Gm)^\omega$ 
is the symplectic perpendicular of $T_m Gm$ in $T_m\bM$, i.e.,
\(
   \{ v\in T_m\bM \mid
   \text{$\omega(v,w) = 0$ for all $w\in T_m Gm$}\}
\).
This is naturally a symplectic vector space.
A vector bundle $T(Gm)^\omega/T(Gm)$ over $Gm$ is called the {\it
symplectic normal bundle}.
(In general, the symplectic normal bundle of an isotropic submanifold
$S$ is defined by $TS^\omega/TS$.)
It is isomorphic to $G\times_{\widehat G}\widehat\bM$.
(In \cite[p.388]{Na-quiver}, $\widehat\bM$ was defined as the
orthogonal complement of the quaternion vector subspace spanned by
$T_m Km$ with respect to the Riemannian metric.)
The action of $\widehat G$ on $\widehat\bM$ preserves the induced
symplectic structure on $\widehat\bM$. Let $\widehat\mu\colon
\widehat\bM\to \widehat{\mathfrak g}^*$ be the corresponding moment
map vanishing at the origin.

We choose an $\operatorname{Ad}(\widehat G)$-invariant splitting
$\mathfrak g = \widehat{\mathfrak g}\oplus
\widehat{\mathfrak g}^\perp$; and its dual splitting 
$\mathfrak g^* = \widehat{\mathfrak g}^*\oplus
\widehat{\mathfrak g}^{\perp*}$.
Let us consider the natural action of $\widehat G$ on the product
$T^*G\times\widehat\bM = G\times\mathfrak g^*\times\widehat\bM$.
With the natural symplectic structure on $T^*G = G\times\mathfrak
g^*$, we have the moment map
\begin{alignat*}{2}
  \widetilde\mu\colon & G\times\mathfrak g^*\times\widehat\bM &\quad
   \to &\quad \widehat{\mathfrak g}^* \\
  & (g,\xi,\widehat m) & \mapsto &
   -\operatorname{pr}\xi + \widehat\mu(\widehat m),
\end{alignat*}
where $\operatorname{pr}\xi$ is the projection of
$\xi\in\mathfrak g^*$ to $\widehat{\mathfrak g}^*$.
Zero is a regular value of $\widetilde\mu$, hence the symplectic quotient
$\widetilde\mu^{-1}(0)/\widehat G$ is a symplectic manifold. It can be
identified with
\(
   G\times_{\widehat G} \left(
     \widehat{\mathfrak g}^{\perp*}\times \widehat\bM \right)
\)
via the map
\begin{equation*}
  G\times_{\widehat G} \left(
     \widehat{\mathfrak g}^{\perp*}\times \widehat\bM \right) \ni
   \widehat G\cdot(g,\xi, \widehat m)
     \longmapsto
   \widehat G\cdot(g,\xi+\widehat\mu(\widehat m),\widehat m)
   \in \widetilde\mu^{-1}(0)/\widehat G.
\end{equation*}
The embedding $G/\widehat G$ into
\(
   G\times_{\widehat G} \left(
     \widehat{\mathfrak g}^{\perp*}\times \widehat\bM \right)
\)
is isotropic and its symplectic normal bundle is $G\times_{\widehat
G}\widehat\bM$. Thus two embeddings of $Gm\cong G/\widehat G$, one
into $\bM$ and the other into
\(
   G\times_{\widehat G} \left(
     \widehat{\mathfrak g}^{\perp*}\times \widehat\bM \right)
\),
have the isomorphic symplectic normal bundles.

The $G$-equivariant version of Darboux-Moser-Weinstein's isotropic
embedding theorem (a special case of \cite[2.2]{SL}) says the
following:
\begin{Lemma}\label{lem:symplectic}
A neighborhood of $Gm$ \rom(in $\bM$\rom) is $G$-equivalently
symplectomorphic to a neighborhood of $G/\widehat G$ embedded as the
zero section of
\(
    G\times_{\widehat G}
      \left( \widehat{\mathfrak g}^{\perp*}\times \widehat\bM \right)
\)
with the $G$-moment map given by the formula
\[
   \mu\left(\widehat G\cdot(g,\xi,\widehat m)\right)
   = \operatorname{Ad}^*(g)\left(\xi + \widehat\mu(\widehat m)\right).
\]
\rom(Here `symplectomorphic' means that there exists a biholomorphism
intertwining symplectic structures.\rom)
\end{Lemma}

Note that Sjamaar-Lerman worked on a {\it real\/} symplectic manifold
with a {\it compact\/} Lie group action. Thus we need care to apply
their result to our situation. Darboux-Moser-Weinstein's theorem is
based on the inverse function theorem, which we have both in the
category of $C^\infty$-manifolds and in that of complex
manifolds.
A problem is that the domain of the symplectomorphism may not be
chosen so that it covers the whole $Gm$ as it is noncompact.
We can overcome this problem by taking a symplectomorphism defined in
a neighborhood of the compact orbit $Km$ first, and then extending it
to a neighborhood of $Gm$, as explained in the next three paragraphs.
This approach is based on a result in \cite{Sj}.

A subset $A$ of a $G$-space $X$ is called {\it orbitally convex\/}
with respect to the $G$-action if it is invariant under $K$ (= maximal
compact subgroup of $G$) and for all $x\in A$ and all $\xi\in\mathfrak
k$ we have that both $x$ and $\exp(\sqrt{-1}\xi)x$ are in $A$ implies
that $\exp(\sqrt{-1}t\xi)x\in A$ for all $t\in [0,1]$.
By \cite[1.4]{Sj}, if $X$ and $Y$ are complex manifolds with
$G$-actions, and if $A$ is an orbitally convex open subset of $X$ and
$f\colon A\to Y$ is a $K$-equivariant holomorphic map, then $f$ can be
uniquely extended to a $G$-equivariant holomorphic map.

Suppose that $X$ is a K\"ahler manifold with a (real) moment map
$\mu_\R\colon X\to \mathfrak k^*$ and that $x\in X$ is a point such
that $\mu_\R(x)$ is fixed under the coadjoint action of $K$. Then
\cite[Claim~1.13]{Sj} says that the compact orbit $Kx$ possesses a
basis of orbitally convex open neiborhoods.

In our situation, we have a K\"ahler metric (\subsecref{subsec:hK})
and we have assumed $\mu_\R(m) = 0$. Thus $Km$ possesses a basis of
orbitally convex open neiborhoods, and we have
\lemref{lem:symplectic}.

Now we want to study local structures of $\M_0$, $\M$ using
\lemref{lem:symplectic}. First the equation $\mu = 0$ implies
$\xi = 0$, ${\widehat\mu}(\widehat m) = 0$. Thus $\M_0$ and $\M$ are
locally isomorphic to `quotients' of 
$G\times_{\widehat G}(\{0\}\times{\widehat\mu}^{-1}(0))$ by $G$,
i.e., 'quotients' of ${\widehat\mu}^{-1}(0)$ by $\widehat G$.
Here the `quotients' are taken in the sense of the geometric invariant 
theory.
Following \propref{prop:stable}(1), we say a point $\widehat
m\in{\widehat\mu}^{-1}(0)$ is {\it stable\/} if the closure of
$\widehat G\cdot (m,z)$ does not intersect with the zero section of
${\widehat\mu}^{-1}(0)\times\C$ for $z\neq 0$.
Here we lift the $\widehat G$-action to the trivial line bundle
${\widehat\mu}^{-1}(0)\times\C$ by $\widehat g\cdot(\widehat m,z) =
(\widehat g\cdot\widehat m, \chi(g)^{-1}z)$, where $\chi$ is the
restriction of the one parameter subgroup used in
\subsecref{subsec:twoquot}.
Let ${\widehat\mu}^{-1}(0)^{\operatorname{s}}$ be the set of stable
points. As in \subsecref{subsec:stability}, we have a morphism
${\widehat\mu}^{-1}(0)^{\operatorname{s}}/\widehat G \to
{\widehat\mu}^{-1}(0)\dslash \widehat G$, which we denote by
$\widehat\pi$.
By \cite[Proposition~2.7]{Sj}, we may assume that the neighborhood of
$Gm$ in \lemref{lem:symplectic} is {\it saturated}, i.e., the closure
of the $G$-orbit of a point in the neighborhood is contained in the
neighborhood. Thus under the symplectomorphism in
\lemref{lem:symplectic}, (i)~closed $G$-orbits are mapped to closed
$\widehat G$-orbits, (ii)~the stability conditions are interchanged.

\begin{Proposition}\label{prop:localslice}
There exist a neighborhood $U$ \rom(resp.\ $U^\perp$\rom) of
$x\in\M_0$ \rom(resp.\ $0\in{\widehat\mu}^{-1}(0)\dslash\widehat
G$\rom) and biholomorphic maps $\Phi\colon U\to U^\perp$,
$\widetilde\Phi\colon \pi^{-1}(U)\to{\widehat\pi}^{-1}(U^\perp)$
such that the following diagram commutes\rom:
\begin{equation*}
\begin{CD}
   \pi^{-1}(U) @>\widetilde\Phi>\cong> {\widehat\pi}^{-1}(U^\perp)
\\
   @V{\pi}VV @VV\widehat{\pi}V \\
%
   U @>\Phi>\cong> U^\perp
\end{CD}
\end{equation*}
In particular, $\pi^{-1}(x) = \M_x$ is biholomorphic to
${\widehat\pi}^{-1}(0)$.

Furthermore, under $\Phi$, a stratum $(\M_0)_{(H)}$ of $\M_0$ is
mapped to a stratum $\left({\widehat\mu}^{-1}(0)\dslash\widehat
  G\right)_{(H)}$, which is defined as in \defref{def:strat}. \rom(If
$(\M_0)_{(H)}$ intersects with $U$, $H$ is conjugate to a subgroup of
$\widehat G$.\rom)
\end{Proposition}

The above discussion shows \propref{prop:localslice} except last
assertion.  The last assertion follows from the argument in
\cite[p.386]{SL}.

\subsection{Slice}\label{subsec:slice}
By \cite[p.391]{Na-quiver}, we have $\widehat G$-invariant splitting 
$\widehat\bM = T\times T^\perp$, where $T$ is the tangent space
$T_x(\M_0)_{(\widehat G)}$ of the stratum containing $x$, and
$\widehat G$ acts trivially on $T$. Thus we have
\begin{gather*}
   {\widehat\mu}^{-1}(0)\dslash \widehat G
    \cong
    T \times \left(T^\perp\cap \widehat\mu^{-1}(0)\right)
         \dslash \widehat G, \\
   {\widehat\mu}^{-1}(0)^{\operatorname{s}}/\widehat G
    \cong
    T \times \left(T^\perp\cap
    \widehat\mu^{-1}(0)^{\operatorname{s}}\right)/ \widehat G.
\end{gather*}
Furthermore, it was proved that
\(
    \left(T^\perp\cap \widehat\mu^{-1}(0)\right) \dslash \widehat G
\)
and
\(
    \left(T^\perp\cap
    \widehat\mu^{-1}(0)^{\operatorname{s}}\right)/ \widehat G
\)
are quiver varieties associated with a certain graph possibly
different from the original one, and possibly with edge loops.
Replacing $U^\perp$ if necessary, we may assume that $U^\perp$ is a
product of a neighborhood $U_T$ of $0$ in $T$ and that $U_{\mathfrak
S}$ of $0$ in
\(
    \left(T^\perp\cap \widehat\mu^{-1}(0)\right) \dslash \widehat G.
\)
We define a {\it transversal slice\/} to $(\M_0)_{(\widehat G)}$ at
$x$ as
\begin{equation*}
   \mathfrak S \defeq
   \Phi^{-1}\left(U^\perp\cap \left(\{ 0\} \times 
     \left(T^\perp\cap \widehat\mu^{-1}(0)\right) \dslash \widehat G
   \right)\right)
   = \Phi^{-1}\left(\{0\}\times U_{\mathfrak S}\right).
\end{equation*}
Since $\Phi$ is a local biholomorphism, this slice $\mathfrak S$
satisfies the properties in \cite[3.2.19]{Gi-book}, i.e.,
there exists a biholomorphism
\begin{equation*}
   \left(U\cap (\M_0)_{(\widehat G)}\right)\times \mathfrak S
    \xrightarrow{\cong} U
\end{equation*}
which induces biholomorphisms between factors
\begin{equation*}
   \{ x\} \times \mathfrak S \xrightarrow{\cong} \mathfrak S, \qquad
   \left(U\cap (\M_0)_{(\widehat G)}\right)\times \{x\}
    \xrightarrow{\cong} \left(U\cap (\M_0)_{(\widehat G)}\right).
\end{equation*}

\begin{Remark}
Our construction gives a slice to a stratum in $\M_{(0,\zeta_\C)}
= \mu^{-1}(-\zeta_\C)\dslash G_\bv$ for general $\zeta_\C$.
(See \cite[p.371 and Theorem~3.1]{Na-quiver} for the definition of
$\M_{(0,\zeta_\C)}$.)
In particular, the fiber $\pi^{-1}(x)$ of $\pi\colon
\M_{(\zeta_\R,\zeta_\C)}\to \M_{(0,\zeta_\C)}$ is isomorphic to the
fiber of 
\(
   \left( T^\perp\cap \widehat\mu^{-1}(0)^{\operatorname{s}}\right)/
      \widehat G 
    \to 
   \left( T^\perp\cap \widehat\mu^{-1}(0)\right)\dslash\widehat G
\)
at $0$. This is a refinement of \cite[6.10]{Na-quiver}, where an
isomorphism between homology groups were obtained.
We also remark that this gives a proof of smallness of
\[
   \pi\colon \bigsqcup_{\zeta_\C} \M_{(\zeta_\R,\zeta_\C)}
     \to \bigsqcup_{\zeta_\C}\M_{(0,\zeta_\C)}
\]
which was observed by Lusztig when $\mathfrak g$ is of type $ADE$
\cite{Lu-small}.
An essential point is, as remarked in \cite[6.11]{Na-quiver}, that 
\(
   \left( T^\perp\cap \widehat\mu^{-1}(0)^{\operatorname{s}}\right)/
      \widehat G 
\)
is diffeomorphic to an affine algebraic variety, and its homology
groups vanishes for degree greater than its complex dimension.
\end{Remark}

For our application, we only need the case when $x$ is {\it regular},
i.e., $x\in\Mreg(\bv^0,\bw)$ for some $\bv^0$.  Then, by
\cite[p.392]{Na-quiver}, \( \left(T^\perp\cap
  \widehat\mu^{-1}(0)\right) \dslash \widehat G \) and \(
\left(T^\perp\cap \widehat\mu^{-1}(0)\right)^{\operatorname{s}}/
\widehat G \) are isomorphic to the quiver varieties
$\M_0(\bv_s,\bw_s)$ and $\M(\bv_s,\bw_s)$, associated with the
original graph with dimension vector
\begin{equation*}
   \bv_s = \bv - \bv^0, \qquad
   \bw_s = \bw - \bC\bv^0,
\end{equation*}
where 
\begin{equation*}
   \bC\bv^0 
   = \sum_{k\in I} \left(2 v^0_k - a_{kl}\, v^0_l\right)\Lambda_k
   \qquad
   \text{if }\bv^0 = \sum_{k\in I} v^0_k \alpha_k
\end{equation*}
in Convention~\ref{convention:weight}.

\begin{Theorem}\label{thm:slice}
Suppose that $x\in\Mreg(\bv^0,\bw)$ as above.
Then there exist neighborhoods $U$, $U_T$, $U_{\mathfrak S}$ of
$x\in\M_0 = \M_0(\bv,\bw)$,
$0\in T$,
$0\in\M_0(\bv_s,\bw_s)$ respectively
and biholomorphic maps $U\to U_T\times U_{\mathfrak S}$,
$\pi^{-1}(U)\to U_T\times {\pi}^{-1}(U_{\mathfrak S})$
such that the following diagram commutes\rom:
\begin{equation*}
\begin{CD}
  \M \;@. \supset \; @.
   \pi^{-1}(U) @>>\cong> U_T \times \pi^{-1}(U_{\mathfrak S})
     \;@.\subset \;@. T\times \M(\bv_s,\bw_s)
\\
   @. @. @V{\pi}VV @VV{\operatorname{id}\times\pi}V @.@. \\
   \M_0\;@.\supset \; @.
   U @>>\cong> U_T\times U_{\mathfrak S}
                  \;@.\subset \;@. T\times\M_0(\bv_s,\bw_s)
\end{CD}
\end{equation*}
In particular, $\pi^{-1}(x) = \M_x$ is biholomorphic to
$\La(\bv_s,\bw_s)$.

Furthermore, a stratum of $\M_0$ is mapped to a product of $U_T$ and a 
stratum of $\M_0(\bv_s,\bw_s)$.
\end{Theorem}

\begin{Remark}\label{rem:Aslice}
Suppose that $A$ is a subgroup of $G_\bw\times\C^*$ fixing $x$. Since
the $A$-action commutes with the $G$-action, $\widehat\bM$ has an
$A$-action. The above construction can be made $A$-equivariant.
In particular, the diagram in \thmref{thm:slice} can be restricted to
a diagram for $A$-fixed points sets.
\end{Remark}

\section{Fixed point subvariety}\label{sec:Fixed}

Let $A$ be an abelian reductive subgroup of $G_\bw\times\C^*$. In this 
section, we study the $A$-fixed point subvarieties
$\M(\bv,\bw)^A$, $\M_0(\bv,\bw)^A$ of $\M(\bv,\bw)$, $\M_0(\bv,\bw)$.

\subsection{A homomorphism attached to a component of $\M(\bv,\bw)^A$}
\label{subsec:fix_hom}

Suppose that $x\in \M(\bv,\bw)^A$ is fixed by $A$. Take a
representative $(B,i,j)\in\mu^{-1}(0)^{\operatorname{s}}$ of $x$.
For every $a\in A$, there exists $\rho(a)\in G_\bv$ such that
\begin{equation}\label{eq:fixed}
   a \star (B,i,j) = \rho(a)^{-1} \cdot (B,i,j),
\end{equation}
where the left hand side is the action defined in \eqref{eq:C*action}
and the right hand side is the action defined in
\eqref{eq:Gaction}. By the freeness of $G_\bv$-action on
$\mu^{-1}(0)^{\operatorname{s}}$ (see \propref{prop:stable}), $\rho(a)$
is uniquely determined by $a$. In particular, the map $a\mapsto
\rho(a)$ is a homomorphism. 

Let $\M(\rho) \subset \M(\bv,\bw)^A$ be the set of fixed points
$x$ such that \eqref{eq:fixed} holds for some representative
$(B,i,j)$ of $x$. Note that $\M(\rho)$ depends only on the
$G_\bv$-conjugacy class of $\rho$.
Since the $G_\bv$-conjugacy class of $\rho$ is locally constant on
$\M(\bv,\bw)^A$, $\M(\rho)$ is a union of connected components of
$\M(\bv,\bw)^A$. Later we show that $\M(\rho)$ is connected under some 
assumptions (see \thmref{thm:M(rho)_conn}).
As in \propref{prop:homotopic}, we have
\begin{Proposition}\label{prop:homotopic'}
$\M(\rho)$ is homotopic to $\M(\rho)\cap\La(\bv,\bw)$.
\end{Proposition}

We regard $V$ as an $A$-module via $\rho$ and consider
the weight space corresponds to $\lambda\in \Hom(A,\C^*)$:
\begin{equation*}
V(\lambda) \defeq \{ v\in V\mid \rho(a)\cdot v = \lambda(a)v \}.
\end{equation*}
We denote by $V_k(\lambda)$ the component of $V(\lambda)$ at the
vertex $k$. We have $V = \bigoplus_\lambda V(\lambda)$.
We regard $W$ as an $A$-module via the composition
\begin{equation*}
   A \hookrightarrow G_\bw\times\C^* \xrightarrow{\text{projection}} G_\bw.
\end{equation*}
We also have the weight space decomposition
$W = \bigoplus_\lambda W(\lambda)$, $W_k =\bigoplus_\lambda W_k(\lambda)$.
We denote by $q$ the composition
\begin{equation*}
   A \hookrightarrow G_\bw\times\C^* \xrightarrow{\text{projection}} \C^*.  
\end{equation*}
Then \eqref{eq:fixed} is equivalent to
\begin{equation}\label{eq:wtmove}
B_{h}(V_{\vout(h)}(\lambda)) \subset
             V_{\vin(h)}(q^{-m(h)-1}\lambda), \quad 
i_k(W_k(\lambda)) \subset V_k(q^{-1}\lambda),\quad
   j_k(V_k(\lambda)) \subset W_k(q^{-1}\lambda),
\end{equation}
where $m(h)$ is as in \eqref{eq:m(h)}.

\begin{Lemma}\label{lem:Vwei}
If $V_k(\lambda)\neq 0$, then $W_l(q^n\lambda)\neq 0$ for some 
$n$ and $l\in I$.
\end{Lemma}

\begin{proof}
Consider $\lambda$ satisfying $W_l(q^{n}\lambda) = 0$ for any $l\in
I$, $n\in\Z$.
If we set
\begin{equation*}
   S_k \defeq \bigoplus_{\text{$\lambda$ as above}} V_k(\lambda),
\end{equation*}
then $S = (S_k)_{k\in I}$ is $B$-invariant and contained in
$\Ker j$ by \eqref{eq:wtmove}. Thus we have $S_k = 0$ by the stability 
condition.
\end{proof}

The restriction of tautological bundles $V_k$, $W_k$ to $\M(\rho)$
are bundles of $A$-modules. We have the weight decomposition
$V_k = \bigoplus V_k(\lambda)$, $W_k = \bigoplus W_k(\lambda)$. We
consider $V_k(\lambda)$, $W_k(\lambda)$ as vector bundles over
$\M(\rho)$.

Similarly, the restriction of the complex $C_k^\bullet$ 
in \eqref{eq:taut_cpx} decomposes as
$C_k^\bullet = \bigoplus_\lambda C_{k,\lambda}^\bullet$, where
\begin{equation}
\label{eq:taut_cpx_fixed}
C_{k,\lambda}^\bullet\equiv C_{k,\lambda}^\bullet(\rho):
\begin{CD}
  V_k(q^2\lambda)
  @>{\sigma_{k,\lambda}}>>
  \displaystyle{\bigoplus_{h:\vin(h)=k}}
     V_{\vout(h)}(q^{m(h)+1}\lambda)
    \oplus W_k(q\lambda)
  @>{\tau_{k,\lambda}}>>
  V_k(\lambda).
\end{CD}
\end{equation}
Here $\sigma_{k,\lambda}$, $\tau_{k,\lambda}$ are restrictions of
$\sigma_k$, $\tau_k$.
When we want to emphasize that this is a complex over $\M(\rho)$, we
denote this by $C_{k,\lambda}^\bullet(\rho)$.

The tangent space of $\M(\rho)$ at $[B,i,j]$ is the $A$-fixed part of
the tangent space of $\M$. Since the latter is the middle cohomology
group of \eqref{eq:quiver_tangent}, the former is the middle
cohomology group of the complex
\begin{equation*}
   \bigoplus_{\lambda,k} \End\left(V_k(\lambda)\right)
   \longrightarrow
   \begin{matrix}
   \bigoplus_{\lambda, h} 
     \Hom\left(V_{\vout(h)}(\lambda),
                      V_{\vin(h)}(q^{-m(h)-1}\lambda)\right) \\
     \oplus \\
   \bigoplus_{\lambda, k}
      \Hom\left(W_k(\lambda), V_k(q^{-1}\lambda)\right) \\
      \oplus \\
   \bigoplus_{\lambda, k}
      \Hom\left(V_k(\lambda), W_k(q^{-1}\lambda)\right)
   \end{matrix}
   \longrightarrow
   \bigoplus_{\lambda,k} \Hom\left(V_k(\lambda),V_k(q^{-2}\lambda)\right),
\end{equation*}
where the differentials are the restrictions of $\iota$, $d\mu$ in
\eqref{eq:quiver_tangent}. Those restrictions are injective and
surjective respectively by \propref{prop:stable}. Hence we have the
following dimension formula:
\begin{equation}\label{eq:dim_M_rho}
\begin{split}
   & \dim \M(\rho) \\
   =\; & \sum_{\lambda}\Big[
   \sum_h \dim V_{\vout(h)}(\lambda) \dim V_{\vin(h)}(q^{-m(h)-1}\lambda) \\
   & \quad
   + \sum_k \dim W_k(\lambda)
      \left( \dim V_k(q^{-1}\lambda) + \dim V_k(q\lambda) \right)
   - \dim V_k(\lambda)^2 
   - \dim V_k(\lambda)\dim V_k(q^{-2}\lambda)\Big].
\end{split}
\end{equation}

Recall that we have an isomorphism
$\pi^{-1}(\M_0^{\operatorname{reg}}) \cong \M_0^{\operatorname{reg}}$
(\propref{prop:pi_isom}).
Let 
\begin{equation}
\label{eq:defM0rhoreg}
  \Mreg(\rho) \defeq
   \pi\left( \pi^{-1}(\Mreg)\cap \M(\rho)
           \right)
  = \pi^{-1}(\Mreg)\cap \pi\left(\M(\rho)\right).
\end{equation}
By definition, $\pi^{-1}(\Mreg(\rho)) = \pi^{-1}(\Mreg)\cap \M(\rho)$
is an open subvariety of $\M(\rho)$ which is isomorphic to
$\Mreg(\rho)$ under $\pi$.



\subsection{A sufficient condition for $\M_0^A = \{ 0\}$}
Let $a = (s,\varepsilon)$ be a semisimple element in
$G_\bw\times\C^*$ and $A$ be the Zariski closure of
$\{ a^n \mid n\in\mathbb Z\}$.

\begin{Definition}\label{def:generic}
We say $a$ is {\it generic\/} if
$\M_0(\bv,\bw)^A = \{ 0\}$ for any $\bv$.
(This condition depends on $\bw$.)
\end{Definition}

\begin{Proposition}
Assume that there is at most one edge joining two vertices of $I$, and
that
\begin{equation*}
   \lambda/\lambda' \notin \{ \ve^n \mid n\in\mathbb Z\setminus\{0\}\}
\end{equation*}
for any pair of eigenvalues of $s\in G_\bw$.
\rom(The condition for the special case $\lambda = \lambda'$ implies
that $\ve$ is not a root of unity.\rom)
Then $a = (s,\ve)$ is generic.
\end{Proposition}

\begin{proof}
We prove $\M_0(\bv,\bw)^A = \{ 0\}$ by the induction on $\bv$. The
assertion is trivial when $\bv = 0$.

Take a point in $\M_0(\bv,\bw)^A$ and its representative $(B,i,j)$. As
in \eqref{eq:fixed}, there exists $g\in G_\bv$ such that
\begin{equation*}
   a \star (B,i,j) = g \cdot (B,i,j).
\end{equation*}
We decompose $V$ into eigenspaces of $g$:
\begin{equation*}
  V = \bigoplus V(\lambda), \qquad\text{where
  $V(\lambda) \defeq \{ v\in V\mid g\cdot v = \lambda(a)v \}$}.
\end{equation*}
We also decompose $W$ into eigenspaces of $s$ as $\bigoplus W(\lambda)$.
Then \eqref{eq:wtmove} holds where $q$ is replaced by $\ve$.

Choose and fix an eigenvalue $\mu$ of $s$. First suppose
$V(\ve^n\mu)\neq 0$ for some $n$. Let
\begin{equation*}
   n_0 \defeq \max\left\{ n \mid V(\ve^n\mu)\neq 0 \right\}.
\end{equation*}
Since $\varepsilon$ is not a root of unity, we have $\ve^n\mu\neq
\ve^m\mu$ for $m\neq n$. Hence the above $n_0$ is well-defined.
By \eqref{eq:wtmove} (and $m(h) = 0$ from the assumption),
we have $\Ima B_h \cap V(\ve^{n_0}\mu) = 0$.
By the assumption, we have $W(\ve^{n_0+1}\mu) = 0$, and hence
$\Ima i_k \cap V(\ve^{n_0}\mu) = 0$ again by \eqref{eq:wtmove}.
Then we may assume the restriction of $(B,i,j)$ to $V(\ve^{n_0}\mu)$
is $0$ as in the proof of \lemref{lem:incl_inj}.

Thus the data $(B,i,j)$ is defined on the smaller subspace
$V\ominus V(\ve^{n_0}\mu)$. Thus $(B,i,j) = 0$ by the induction
hypothesis.

If $V(\ve^n\mu) = 0$ for any $n$, we replace $\mu$. If we can find a
$\mu'$ so that $V(\ve^n\mu')\neq 0$ for some $n$, we are
done. Otherwise, we have $V(\ve^n\mu) = 0$ for any $n$, $\mu$, and we
have $i = j = 0$ by \eqref{eq:wtmove}.
Then we choose $\mu$, which may not be an eigenvalue of $s$, so that
$V(\mu) \neq 0$ and repeat the above argument. (This is possible since
we may assume $V\neq 0$.) We have $\Ima B_h \cap V(\ve^{n_0}\mu) = 0$
and the data $B$ is defined on the smaller subspace $V\ominus
V(\ve^{n_0}\mu)$ as above.
\end{proof}

\section{Hecke correspondence and induction of quiver varieties}

\subsection{Hecke correspondence}\label{subsec:hecke}

Take dimension vectors $\bw$, $\bv^1$, $\bv^2$ such that
$\bv^2 = \bv^1 + \alpha_k$.
Choose collections of vector spaces $W$, $V^1$, $V^2$, with
$\dim W = \bw$, $\dim V^a = \bv^a$.

Let us consider the product $\M(\bv^1,\bw)\times\M(\bv^2,\bw)$.
We denote by
$V^1_k$ (resp.\ $V^2_k$) the vector bundle
$V_k\boxtimes\shfO_{\M(\bv^2,\bw)}$
(resp.\ $\shfO_{\M(\bv^1,\bw)}\boxtimes V_k$).
A point in $\M(\bv^1,\bw)\times\M(\bv^2,\bw)$ is denoted by
$([B^1,i^1,j^1], [B^2, i^2, j^2])$.
We regard $B^a$, $i^a$, $j^a$ ($a=1,2$) as homomorphisms between
tautological bundles.

We define a three-term sequence of vector bundles
over $\M(\bv^1,\bw)\times \M(\bv^2, \bw)$ by
\begin{equation}\label{eq:hecke_complex}
  \HomL(V^1, V^2)
  \overset{\sigma}{\longrightarrow}
    q \HomE(V^1,V^2) \oplus 
    q \HomL(W, V^2) \oplus q \HomL(V^1,W)
  \overset{\tau}{\longrightarrow}
  q^2 \HomL(V^1, V^2) \oplus q^2\shfO,
\end{equation}
where
\begin{equation*}
\begin{split}
        \sigma(\xi) & = (B^2 \xi - \xi B^1) \oplus
         (-\xi i^1) \oplus j^2 \xi \\
        \tau(C\oplus a\oplus b)
         &= ( \varepsilon B^2 C + \varepsilon C B^1 + i^2 b + a j^1)
        \oplus\left( \tr(i^1 b) + \tr(a j^2)\right).
\end{split}\end{equation*}
This is a complex, that is $\tau\sigma = 0$, thanks to
the equation $\varepsilon BB+ij = 0$ and $\tr(i^1 j^2 \xi) = \tr(\xi
i^1 j^2)$.
Moreover, it is an equivariant complex with respect to the
$G_\bw\times\C^*$-action.

By \cite[5.2]{Na-alg}, $\sigma$ is injective and $\tau$ is surjective.
Hence the quotient $\Ker \tau/\Ima \sigma$ is a
$G_\bw\times\C^*$-equivariant vector bundle. Let us define an
equivariant section $s$ of $\Ker \tau/\Ima \sigma$ by
\begin{equation}
        s = \left( 0 \oplus (-i^2) \oplus j^1 \right) \bmod{\Ima \sigma},
\end{equation}
where $\tau(s) = 0$ follows from $\varepsilon BB+ij=0$ and
$\tr(B^1 B^1) = \tr(B^2 B^2) = 0$.
The point $([B^1, i^1, j^1], [B^2, i^2, j^2])$ is contained in
the zero locus $Z(s)$ of $s$ if and only if
there exists $\xi \in \HomL(V^1, V^2)$ such that
\begin{equation}
        \xi B^1 = B^2 \xi, \quad
        \xi i^1 = i^2, \quad
        j^1 = j^2 \xi.
\label{hecke:xi}\end{equation}
Moreover, $\Ker \xi$ is zero by the stability condition for $B^2$.
Hence $\Ima \xi$ is a subspace of $V^2$ with dimension $\bv^1$ which
is $B^2$-invariant and contains $\Ima i^2$.
Moreover, such $\xi$ is unique if we fix representatives
$(B^1,i^1,j^1)$ and $(B^2,i^2,j^2)$.
Hence we have an isomorphism between $Z(s)$ and the variety of all 
pairs $(B, i, j)$ and $S$ (modulo $G_{\bv^2}$-action) such that
\begin{aenume}
\item $(B, i, j)\in\mu^{-1}(0)$ is stable, and
\item $S$ is a $B$-invariant subspace containing the image of $i$ with
         $\dim S = \bv^1 = \bv^2 - \alpha_k$.
\end{aenume}

\begin{Definition}\label{def:hecke}
We call $Z(s)$ the {\it Hecke correspondence}, and denote it by
$\Pa_k(\bv^2,\bw)$. It is a $G_\bw\times\C^*$-invariant closed
subvariety.
\end{Definition}

Introducing a connection $\nabla$ on $\Ker \tau/\Ima \sigma$, we
consider the differential
\begin{equation*}
        \nabla s\colon T\M(\bv^1,\bw)\oplus T\M(\bv^2,\bw)
        \to \Ker \tau/ \Ima \sigma
\end{equation*}
of the section $s$.
Its restriction to $Z(s) = \Pa_k(\bv^2,\bw)$ is independent of the
connection.
By \cite[5.7]{Na-alg},
The differential $\nabla s$ is surjective over $\Pa_k(\bv^2,\bw)$.
Hence, $\Pa_k(\bv^2,\bw)$ is nonsingular.

By the definition, the quotient $V_k^2/V_k^1$ defines a line bundle
over $\Pa_k(\bv^2,\bw)$.

\subsection{Hecke correspondence and fixed point subvariety}
\label{subsec:HeckeFix}
Let $A$ be as in \secref{sec:Fixed} and let $\M(\rho)$ be as in
\subsecref{subsec:fix_hom} for $\rho\in\Hom(A,G_\bv)$.

Let us consider the intersection
$(\Mw^A\times\Mw^A)\cap\Pa_k(\bv^2,\bw)$.
It decomposes as
\begin{equation*}
  (\Mw^A\times\Mw^A)\cap\Pa_k(\bv^2,\bw)
  = \bigsqcup_{\rho^1,\rho^2}
     (\M(\rho^1)\times\M(\rho^2))\cap \Pa_k(\bv^2,\bw).
\end{equation*}
Take a point $([B^1,i^1,j^1], [B^2,i^2,j^2])\in
(\M(\rho^1)\times\M(\rho^2))\cap \Pa_k(\bv^2,\bw)$.
Then we have
\begin{equation*}
  a \star (B^p,i^p,j^p) = \rho^p(a)^{-1} \cdot (B^p,i^p,j^p)
   \qquad a\in A,\quad (p = 1,2),
\end{equation*}
and there exists $\xi\in\HomL(V^1,V^2)$ such that
\begin{equation*}
  \xi B^1 = B^2 \xi, \quad
  \xi i^1 = i^2, \quad
  j^1 = j^2 \xi.
\end{equation*}
By the uniqueness of $\xi$, we must have $\rho^2(a)\xi =
\xi\rho^1(a)$, that is $\xi\colon V^1\to V^2$ is $A$-equivariant.
Since $\xi$ is injective, $V^1$ can be considered as $A$-submodule of
$V^2$.

If $V^1 = \bigoplus V^1(\lambda)$, $V^2 = \bigoplus V^2(\lambda)$ are
the weight decomposition, then there exists $\lambda_0$ such that
\begin{aenume}
\item $V^1_l(\lambda) \xrightarrow{\xi} V^2_l(\lambda)$ is an isomorphism
if $\lambda\neq \lambda_0$ or $l\neq k$,
\item $V^1_k(\lambda_0) \xrightarrow{\xi} V^2_k(\lambda_0)$ is a
codimension $1$ embedding.
\end{aenume}

\subsection{}
We introduce a generalization of the Hecke correspondence.
Let us define $\Pa_k^{(n)}(\bv,\bw)$ as
\begin{equation}
\label{eq:Pa^{(n)}}
   \Pa_k^{(n)}(\bv,\bw) \defeq
   \{ (B,i,j,S)\mid
   \text{$(B,i,j)\in\bM(\bv,\bw)$, $S\subset V$ as below}\}
   /G_\bv,
\end{equation}
\begin{aenume}
\item $(B, i, j)\in\mu^{-1}(0)^{\operatorname{s}}$,
\item $S$ is a $B$-invariant subspace containing the image of $i$ with
         $\dim S = \bv - n\alpha_k$.
\end{aenume}
For $n = 1$, it is nothing but $\Pa_k(\bv,\bw)$. We consider
$\Pa_k^{(n)}(\bv,\bw)$ as a closed subvariety of
$\M(\bv-n\alpha_k,\bw)\times\M(\bv,\bw)$ by setting
\begin{align*}
  & (B^1,i^1,j^1) \defeq \text{the restriction
                               of $(B,i,j)$ to $S$},\\
  & (B^2,i^2,j^2) \defeq (B,i,j).
\end{align*}
We have a vector bundle of rank $n$ defined by $V_k^2/V_k^1$.

We shall show that $\Pa_k^{(n)}(\bv,\bw)$ is nonsingular later (see
the proof of \lemref{lem:pullbackSHecke}).

\subsection{Induction}
We recall some results in \cite[\S4]{Na-alg}. Let $Q_k(\bv,\bw)$ the
middle cohomology of the complex~\eqref{eq:taut_cpx}, i.e.,
\begin{equation*}
    Q_k(\bv,\bw) \defeq 
    \Ker \tau_k / \Ima \sigma_k
.
\end{equation*}

We introduce the following subsets of $\M(\bv,\bw)$
(cf.\ \cite[12.2]{Lu2}):
\begin{equation}\label{taut:subset}
\begin{gathered}
        \M_{k;n}(\bv,\bw) \defeq
        \left\{ [B,i,j]\in \M(\bv,\bw) \Biggm| \codim_{V_k}
        \Ima \tau_k
        = n \right\} \\
        \M_{k;\le n}(\bv,\bw) \defeq
        \bigcup_{m\le n} \M_{k;m}(\bv,\bw), \qquad
        \M_{k;\ge n}(\bv,\bw) \defeq
        \bigcup_{m\ge n} \M_{k;m}(\bv,\bw).
\end{gathered}
\end{equation}
Since $\M_{k;\le n}(\bv,\bw)$ is an open subset of $\M(\bv,\bw)$,
$\M_{k;n}(\bv,\bw)$ is a locally closed subvariety.  The restriction
of $Q_k(\bv,\bw)$ to $\M_{k;n}(\bv,\bw)$ is a
$G_\bw\times\C^*$-equivariant vector bundle of rank 
$\langle h_k, \bw - \bv\rangle + n$,
where we used Convention~\ref{convention:weight}.

Replacing $V_k$ by $\Ima\tau_k$, we have a
natural map
\begin{equation}\label{taut:induct}
  p\colon \M_{k;n}(\bv,\bw) \to \M_{k;0}(\bv - n\alpha_k,\bw).
\end{equation}
Note that the projection $\pi$ \eqref{eq:pi_def} factors through
$p$. In particular, the fiber of $\pi$ is preserved under $p$.

\begin{Proposition}\label{prop:taut_fib}
Let $G(n,Q_k(\bv - n\alpha_k,\bw)|_{\M_{k;0}(\bv-n\alpha_k,\bw)})$ be
the Grassmann bundle of $n$-planes in
the vector bundle 
obtained by restricting $Q_k(\bv - n\alpha_k,\bw)$ to
$\M_{k;0}(\bv-n\alpha_k,\bw)$.
Then
%
we have the following diagram:
\begin{equation*}
\begin{CD}
   G(n,Q_k(\bv - n\alpha_k,\bw)|_{\M_{k;0}(\bv-n\alpha_k,\bw)})
     @>\pi>> \M_{k;0}(\bv - n\alpha_k,\bw) \\
   @VV{\cong}V @| \\
   \M_{k;n}(\bv,\bw) @>p>> \M_{k;0}(\bv - n\alpha_k,\bw) \\
   @A{p_2}A{\cong}A @| \\
   \Pa_k^{(n)}(\bv,\bw)
   \cap(\M(\bv-n\alpha_k,\bw)\times\M_{k;\le n}(\bv,\bw))
   @>{p_1}>>
   \M_{k;0}(\bv-n\alpha_k,\bw),
\end{CD}
\end{equation*}
where $\pi$ is the natural projection,
$p_1$ and $p_2$ are restrictions of the projections to the first
and second factors. 
The kernel of the natural surjective homomorphism
\(
   p^* Q_k(\bv - n\alpha_k,\bw) \to Q_k(\bv,\bw)
\)
is isomorphic to the tautological vector bundle of the Grassmann
bundle of the first row, and also to the
the restriction of the vector bundle $V_k^2/V_k^1$ over
$\Pa_k^{(n)}(\bv,\bw)$ in the third row.
\end{Proposition}

\begin{proof}
The proof is essentially contained in \cite[4.5]{Na-alg}.
See also \propref{prop:taut_fib_fix} for a similar result.
\end{proof}

\subsection{Induction for fixed point subvarieties}

We consider the analogue of the results in the previous subsection for 
fixed point subvariety $\M(\rho)$. Let us use notation as in
\subsecref{subsec:fix_hom}, and suppose that $A$ is the Zariski
closure of a semisimple element $a = (s,\varepsilon)\in G_\bw\times\C^*$.

Let $Q_{k,\lambda}(\rho)$ be the middle cohomology of the
complex $C_{k,\lambda}^\bullet(\rho)$ in \eqref{eq:taut_cpx_fixed}, i.e.,
\begin{equation*}
   Q_{k,\lambda}(\rho) \defeq
   \Ker\tau_{k,\lambda}/\Ima\sigma_{k,\lambda}
\end{equation*}

Let
\begin{equation}\label{eq:fix_subset}
        \M_{k;(n_\lambda)}(\rho) \defeq
        \left\{ [B,i,j]\in \M(\rho) \Biggm| \text{$\codim_{V_k(\lambda)}
        \Ima \tau_{k,\lambda}
        = n_\lambda$ for each $\lambda$}\right\}
\end{equation}

Replacing $V_k(\lambda)$ by $\Ima\tau_{k,\lambda}$, we have a natural
map
\begin{equation*}
   p^A \colon \M_{k;(n_\lambda)}(\rho) \to \M_{k;(0)}(\rho'),
\end{equation*}
where $\rho'\colon A\to G_{\bv'}$ ($\bv' = {\bv - \sum n_\lambda
\alpha_k}$) is the homomorphism obtained from $\rho\colon A\to G_\bv$
by replacing $V_k(\lambda)$ by its codimension $n_\lambda$
subspace. Its conjugacy class is independent of the choice of the
subspace. This map is nothing but the restriction of $p$ in the
previous subsection.

For each $\lambda$,
let $G(n_\lambda,Q_{k,q^{-2}\lambda}(\rho')|_{\M_{k;(0)}(\rho')})$
denote the Grassmann bundle of $n_\lambda$-planes in the vector bundle
obtained by restricting $Q_{k,q^{-2}\lambda}(\rho')$ to
$\M_{k;(0)}(\rho')$.
Let 
\begin{equation*}
   \prod_\lambda G(n_\lambda,Q_{k,q^{-2}\lambda}(\rho')|_{\M_{k;(0)}(\rho')})
\end{equation*}
be their fiber product over $\M_{k;(0)}(\rho')$

We have the following analogue of \propref{prop:taut_fib}:
\begin{Proposition}\label{prop:taut_fib_fix}
Suppose that $\varepsilon^2 \neq 1$.
We have the following diagram:
\begin{equation*}
\begin{CD}
   {\displaystyle \prod_\lambda}
     G(n_\lambda,Q_{k,q^{-2}\lambda}(\rho')|_{\M_{k;(0)}(\rho')})
      @>\pi>> \M_{k;(0)}(\rho') \\
   @VV{\cong}V @| \\
   \M_{k;(n_\lambda)}(\rho) @>p^A>> \M_{k;(0)}(\rho'),
\end{CD}
\end{equation*}
where $\pi$ is the natural projection.
For each $\lambda$, the kernel of the natural surjective homomorphism
\(
   (p^A)^* Q_{k,q^{-2}\lambda}(\rho') \to Q_{k,q^{-2}\lambda}(\rho)
\)
is isomorphic to the tautological vector bundle of the Grassmann
bundle. Moreover, we have
\begin{gather}
   n_\lambda \ge \max\left(0, -\rank
      C_{k,\lambda}^\bullet(\rho)\right) \label{eq:n_l_ineq}
\\
   \dim \M_{k;(n_\lambda)}(\rho)
   = \dim \M(\rho) - \sum_\lambda
        n_\lambda\left(\rank C_{k,\lambda}^\bullet(\rho) + n_\lambda\right).
    \label{eq:codim_locus}
\end{gather}
\rom(Here $\rank$ of a complex means the alternating sum of dimensions
of cohomology groups.\rom)
\end{Proposition}

\begin{proof}
We have a surjective homomorphism
\(
   (p^A)^* Q_{k,q^{-2}\lambda}(\rho') \to Q_{k,q^{-2}\lambda}(\rho)
\)
of codimension $n_\lambda$ over $M_{k;(n_\lambda)}(\rho)$. This gives a 
morphism from $M_{k;(n_\lambda)}(\rho)$ to the fiber product of
Grassmann bundles. By a straightforward modification of the arguments
in \cite[4.5]{Na-alg}, one can show that it is an isomorphism.
The detail is left to the reader. The assumption $\varepsilon^2 \neq
1$ is used to distinguish $Q_{k,\lambda}$ and $Q_{k,q^{-2}\lambda}$.

Let us prove the remaining part \eqref{eq:n_l_ineq}, \eqref{eq:codim_locus}.
First note that
\begin{equation*}
\begin{split}
  \rank Q_{k,q^{-2}\lambda}(\rho')|_{\M_{k;(0)}(\rho')}
  & = \rank C_{k,q^{-2}\lambda}^\bullet(\rho') \\
  & = \rank C_{k,q^{-2}\lambda}^\bullet(\rho)
        + n_\lambda + n_{q^{-2}\lambda}.
\end{split}
\end{equation*}
Since we have an $n_\lambda$-subspace in
$Q_{k,q^{-2}\lambda}(\rho')|_{\M_{k;(0)}(\rho')}$, we must have
\begin{equation*}
   n_{q^{-2}\lambda} + \rank C_{k,q^{-2}\lambda}^\bullet \ge 0.
\end{equation*}
Replacing $q^{-2}\lambda$ by $\lambda$, we get \eqref{eq:n_l_ineq}.

Moreover, we have
\begin{equation*}
  \dim \M_{k;(n_\lambda)}(\rho)
   = \dim \M_{k;(0)}(\rho')
   + \sum_\lambda n_\lambda
    \left(\rank C_{k,q^{-2}\lambda}^\bullet(\rho) + n_{q^{-2}\lambda} \right).
\end{equation*}
On the other hand, the dimension formula~\eqref{eq:dim_M_rho} implies
\begin{equation*}
   \dim \M(\rho) - \dim \M(\rho')
   = \sum_\lambda n_\lambda \left(\rank C_{k,\lambda}^\bullet(\rho) 
           + \rank C_{k,q^{-2}\lambda}^\bullet(\rho) 
           + n_{\lambda} + n_{q^{-2}\lambda}\right).
\end{equation*}
Since $\M_{k;(0)}(\rho')$ is an open subset of $\M(\rho')$, we get
\eqref{eq:codim_locus}.
\end{proof}

Note that the inequality \eqref{eq:n_l_ineq} implies that
\begin{equation*}
   n_\lambda\left(\rank C_{k,\lambda}^\bullet(\rho)+ n_\lambda\right) \ge 0
\end{equation*}
and the equality holds if and only if
\begin{equation*}
   n_\lambda = \max\left(0, -\rank C_{k,\lambda}^\bullet(\rho)\right).
\end{equation*}

In particular, we have the following analog of \cite[4.6]{Na-alg}.
\begin{Corollary}\label{cor:codim_rank}
Suppose $\varepsilon^2\neq 1$.
On a nonempty open subset $\M(\rho)$, we have
\begin{equation*}
    \codim_{V_k(\lambda)} \Ima \tau_{k,\lambda}
   = \max\left(0, -\rank C_{k,\lambda}^\bullet(\rho)\right)
\end{equation*}
for each $\lambda$.
And the complement is a proper subvariety of $\M(\rho)$.
\end{Corollary}

As an application of this induction, we prove the following:
\begin{Theorem}\label{thm:M(rho)_conn}
Assume that $\varepsilon$ is not a root of unity and there is at most
one edge joining two vertices of $I$.
Then $\M(\rho)$ is connected if it is a nonempty set.
\end{Theorem}

\begin{proof}
We prove the assertion by induction on $\dim V$, $\dim W$. (The result is
trivial when $V = W = 0$.)

We first make a reduction to the case when
\begin{equation}\label{eq:C_k<0}
   \rank C_{k,\lambda}^\bullet(\rho) < 0 \quad\text{for some $k$, $\lambda$}.
\end{equation}

Fix a $\mu\in\Hom(A,\C^*)$ and consider
\begin{equation*}
   n_0 \defeq \max\left\{ n \mid
     \text{$V_k(q^n\mu)\neq 0$ or $W_k(q^n\mu)\neq 0$
   for some $k\in I$}\right\}.
\end{equation*}
Since $\varepsilon$ is not a root of unity, we have $q^n\mu\neq
q^m\mu$ for $m\neq n$. Hence the above $n_0$ is well-defined.

Suppose $W_k(q^{n_0}\mu)\neq 0$. By \eqref{eq:wtmove} and the choice
of $n_0$, we have $\Ima j_k\cap W_k(q^{n_0}\mu) = \{ 0\}$.
Let us replace $W_k(q^{n_0}\mu)$ by $0$. Namely we
change $(\text{restriction of $i_k$})\colon W_k(q^{n_0}\mu)\to
V_k(q^{n_0-1}\mu)$ to $0$ and all other data are unchanged.
The equation $\mu(B,i,j) = 0$ and the stability condition are
preserved by the replacement. Thus we have a morphism
\begin{equation*}
   \M(\rho) \to \M'(\rho),
\end{equation*}
where $\M'(\rho)$ is a fixed point subvariety
of $\M(\bv,\bw')$ obtained by the replacement.
(This notation will not be used elsewhere. The data $\bw$ is fixed
elsewhere.)

Conversely, we can put any homomorphism $W_k(q^{n_0}\mu)\to
V_k(q^{n_0-1}\mu)$ to get a point in $\M(\rho)$ starting from a point
in $\M'(\rho)$. This shows that $\M(\rho)$ is the total space of the
vector bundle $\Hom(W_k(q^{n_0}\mu), V_k(q^{n_0-1}\mu))$ over
$\M'(\rho)$, where $W_k(q^{n_0}\mu)$ is considered as a trivial
bundle. In particular, $\M(\rho)$ is (nonempty and) connected if and
only if $\M(\rho')$ is so. By the induction hypothesis, $\M(\rho')$ is
connected and we are done.

Thus we may assume $V_k(q^{n_0}\mu)\neq 0$.
Then $C_{k,q^{n_0}\mu}(\rho)$
consists of the last term by the choice of $n_0$. 
(Note $m(h) = 0$ under the assumption that there is at most
one edge joining two vertices of $I$.)
Hence we have
\eqref{eq:C_k<0} with $\lambda = q^{n_0}\mu$.

Now let us prove the connectedness of $\M(\rho)$ under
\eqref{eq:C_k<0}.
By \corref{cor:codim_rank}, we have
\begin{equation*}
   \dim \M_{k;(n_\lambda)}(\rho) < \dim \M(\rho)
\end{equation*}
unless $n_\lambda = \max\left(0,
-\rank C_{k,\lambda}^\bullet(\rho)\right)$ for each $\lambda$.
Hence it is enough to prove the connectedness of
$\M_{k;(n_\lambda^0)}(\rho)$ for
$n_\lambda^0 = \max\left(0,
-\rank C_{k,\lambda}^\bullet(\rho)\right)$.

Let us consider the map $p^A\colon \M_{k;(n_\lambda^0)}(\rho)\to 
\M_{k;(0)}(\rho')$.
By \eqref{eq:C_k<0}, $\dim V$ becomes smaller for
$\M_{k;(0)}(\rho')$. Hence $\M(\rho')$ is connected by the induction
hypothesis. Again by \corref{cor:codim_rank}, $\M_{k;(0)}(\rho')$ is also
connected. Since $p^A$ is a fiber product of Grassmann bundles,
$\M_{k;(n_\lambda^0)}(\rho)$ is connected.
\end{proof}

\section{Equivariant $K$-theory}

In this section, we review the equivariant $K$-theory of 
a quasi-projective variety with a group action.
See \cite[Chapter~5]{Gi-book} for further details.

\subsection{Definitions}

Let $X$ be a quasi-projective variety over $\C$.  Suppose that a
linear algebraic group $G$ acts algebraically on $X$. 
Let $K^G(X)$ be the Grothendieck group of the abelian category of
$G$-equivariant coherent sheaves on $X$.
It is a module over $R(G)$, the representation ring of $G$.

A class in $K^G(X)$ represented by a $G$-equivariant sheaf $E$ will be
denoted by $[E]$, or simply by $E$ if there is no fear of confusion.

The trivial line bundle of rank $1$, i.e., the structure sheaf, is
denoted by $\shfO_X$. If the underlying space is clear, we simply
write $\shfO$.

Let $K_G^0(X)$ be the Grothendieck group of the abelian category of
$G$-equivariant algebraic vector bundles on $X$. 
This is also an $R(G)$-module. The tensor product
$\otimes$ defines a structure of an $R(G)$-algebra on $K_G^0(X)$. Also,
$K^G(X)$ has a structure of a $K_G^0(X)$-module by the tensor
product:
\begin{equation}\label{eq:tensor}
   K_G^0(X)\times K^G(X)\ni ([E],[F]) \mapsto
   [E\otimes F] \in K^G(X).
\end{equation}

Suppose that $Y$ is a $G$-invariant closed subvariety of $X$ and let
$U = X\setminus Y$ be the complement. Two inclusions
\begin{equation*}
   Y \xrightarrow{i} X \xleftarrow{j} U.
\end{equation*}
induce an exact sequence
\begin{equation}\label{eq:exact}
\begin{CD}
   K^G(Y)  @>{i_*}>> K^G(X) @>{j^*}>> K^G(U) @>>> 0,
\end{CD}
\end{equation}
where $i_*$ is given by $[E]\mapsto [i_* E]$, and $j^*$ is given by
$[E]\mapsto [E|_U]$. (See \cite{Thom}.)

Suppose that $Y$ is a $G$-invariant closed subvariety of $X$ and that
$X$ is nonsingular. Let $K^G(X;Y)$ be the Grothendieck group of the
derived category of $G$-equivariant complexes $E^\bullet$ of
algebraic vector bundles over $X$, which are exact outside $Y$ (see
\cite[\S1]{BFM}).
We have a natural homomorphism $K^G(X;Y)\to K^G(Y)$ by setting
\begin{equation*}
  [E^\bullet] \mapsto \sum_i (-1)^i [\gr H^i(E^\bullet)].
\end{equation*}
Here $H^i(E^\bullet)$ is the $i$th cohomology sheaf of $E^\bullet$,
which is a $G$-equivariant coherent sheaf on $X$ which are supported
on $Y$.
If $\idl_Y$ is the defining ideal of $Y$, we have
$\idl_Y^N\cdot H^i(E^\bullet) = 0$ for sufficiently large $N$. Then
\begin{equation*}
   \gr H^i(E^\bullet) \defeq
   \bigoplus_j \idl_Y^j\cdot H^i(E^\bullet)/ \idl_Y^{j+1}\cdot H^i(E^\bullet)
\end{equation*}
is a sheaf on $Y$, and defines an element in $K^G(Y)$.
Conversely if a $G$-equivariant coherent sheaf
$F$ on $Y$ is given, we can take a resolution by
a finite $G$-equivariant complex of algebraic vector bundles:
\begin{equation*}
  0 \to E^{-n} \to E^{1-n} \to \dots \to E^0 \to i_* F \to 0,
\end{equation*}
where $i\colon Y\to X$ denote the inclusion. (See \cite[5.1.28]{Gi-book}.)
This shows that the homomorphism $K^G(X;Y)\to K^G(Y)$ is an
isomorphism. This relative $K$-group $K^G(X;Y)$ was not used in
\cite{Gi-book} explicitly, but many operations were defined by using it
implicitly.
When $Y = X$, $K^G(X;X)$ is isomorphic to $K_G^0(X)$. In particular,
we have an isomorphism $K_G^0(X) \cong K^G(X)$ if $X$ is
nonsingular.

We shall also use equivariant topological $K$-homology
$\Ktop^G(X)$. There are several approaches for the definition, but we
take one in \cite[5.3]{T-top}. There is a comparison map
\begin{equation*}
   K^G(X) \to \Ktop^G(X)
\end{equation*}
which satisfies obvious functorial properties.

Occasionally, we also consider the {\it higher\/} equivariant
topological $K$-homology group $\Kotop^G(X)$. (See \cite[5.3]{T-top}
again.)
In this circumstance, $\Ktop^G(X)$ may be written as $\Kztop^G(X)$.
But we do not use higher equivariant algebraic $K$-homology
$K^G_i(X)$.

Suppose that $Y$ is a $G$-invariant closed subvariety of $X$ and let
$U = X\setminus Y$ be the complement. Two inclusions
\begin{equation*}
   Y \xrightarrow{i} X \xleftarrow{j} U.
\end{equation*}
induce a natural exact hexagon
\begin{equation}\label{eq:exact'}
\begin{CD}
   \Kztop^G(Y)  @>{i_*}>> \Kztop^G(X)  @>{j^*}>> \Kztop^G(U) \\
     @AAA                 @.                       @VVV      \\
   \Kotop^G(U)  @<{j^*}<< \Kotop^G(X)  @<{i_*}<< \Kotop^G(Y),
\end{CD}
\end{equation}
for suitably defined $i_*$, $j^*$.

\subsection{Operations on $K$-theory of vector bundles}\label{subsec:opt_vec}

If $E$ is a $G$-equivariant vector bundle, its rank and dual vector
bundle will be denoted by $\rank E$ and $E^*$ respectively.

We extend $\rank$ and $*$ to operations on $K^0_G(X)$:
\begin{equation*}
   \rank\colon K^0_G(X) \to \Z^{\pi_0(X)}, \qquad
   {}^*\colon K^0_G(X) \to K^0_G(X),
\end{equation*}
where $\pi_0(X)$ is the set of the connected components of $X$.  Note
that the rank of a vector bundle may not be a constant, when $X$ has
several connected components. But we assume $X$ is connected in this
subsection for simplicity. In general, operators below can be defined
component-wisely.

If $L$ is a $G$-equivariant line bundle, we define $L^{\otimes r} =
(L^*)^{\otimes(-r)}$ for $r < 0$. Thus we have $L^{\otimes r}\otimes
L^{\otimes s} = L^{\otimes (r+s)}$ for $r,s\in\Z$.

If $E$ is a vector bundle, we define
\begin{equation*}
  \det E \defeq \Wedge^{\rank E} E, \qquad
  \Wedge_u E \defeq 
     \sum_{i=0}^{\rank E} u^i \Wedge^i E.
\end{equation*}
These operations can be extended to $K_G^0(X)$ of $G$-equivariant
algebraic vector bundles:
\begin{equation*}
  \det\colon K_G^0(X)\to K_G^0(X), \quad
  \Wedge_u\colon K_G^0(X)\to [\shfO_X] + K_G^0(X)\otimes u\Z[[u]].
\end{equation*}
This is well-defined since we have
$\det F = \det E\otimes \det G$,
$\Wedge_u F = \Wedge_u E \otimes \Wedge_u G$ for an exact sequence
$0 \to E \to F \to G \to 0$.

Note the formula
\begin{equation*}
   \Wedge_u E
   = u^{\rank E} \det E\; \Wedge_{1/u} E^*
\end{equation*}
for a vector bundle $E$.
Using this formula, we expand $\Wedge_u E$
into the Laurent expansion also at $u=\infty$:
\begin{equation*}
   u^{-\rank E}(\det E)^* \Wedge_u E
   \in [\shfO_X] + K_G^0(X)\otimes u^{-1}\Z[[u^{-1}]].
\end{equation*}

\subsection{Tor-product}\label{subsec:tor}
(cf.\ \cite[1.3]{BFM}, \cite[5.2.11]{Gi-book})
Let $X$ be a nonsingular quasi-projective variety with a
$G$-action.
Let $Y_1$, $Y_2\subset X$ are $G$-invariant closed subvariety 
of $X$.
Suppose that $E_1^\bullet$ (resp.\ $E_2^\bullet$) is a $G$-equivariant
complex of vector bundles over $X$ which is exact outside $Y_1$
(resp.\ $Y_2$).
Then we can construct a complex
\begin{equation*}
  \cdots
  \longrightarrow \bigoplus_{p+q=k} E_1^p\otimes E_2^q \longrightarrow
   \bigoplus_{p+q=k+1} E_1^p\otimes E_2^q \longrightarrow\cdots
\end{equation*}
with suitably defined differentials from the double complex
$E_1^\bullet\otimes E_2^\bullet$. It is exact outside $Y_1\cap
Y_2$. This construction defines an $R(G)$-bilinear pairing
\begin{equation*}
   K^G(X; Y_1) \times K^G(X; Y_2) \to K^G(X; Y_1\cap Y_2).
\end{equation*}
Since we assume $X$ is nonsingular, we have
$K^G(X;Y_1)\cong K^G(Y_1)$, $K^G(X;Y_2)\cong K^G(Y_2)$,
$K^G(X;Y_1\cap Y_2)\cong K^G(Y_1\cap Y_2)$. Thus we also have
an $R(G)$-bilinear pairing
\begin{equation*}
   K^G(Y_1) \times K^G(Y_2) \to K^G(Y_1\cap Y_2).
\end{equation*}
We denote these operations by $\cdot\otimes^L_X \cdot$.
(It is denoted by $\otimes$ in \cite{Gi-book}.)

\begin{Lemma}[\protect{\cite[5.4.10]{Gi-book},
\cite[Lemma~1]{Vasserot}}]\label{lem:cleanint}
Let $Y_1$, $Y_2\subset X$ be nonsingular $G$-subvarieties with conormal
bundles $T_{Y_1}^* X$, $T_{Y_2}^* X$.
Suppose that $Y \defeq Y_1\cap Y_2$ is nonsingular and
$TY_1|_Y \cap TZ_2|_Y = TY$, where $\ |_Y$ means the restriction to
$Y$.
Then for any $E_1\in K_G^0(Y_1)\cong K^G_0(Y_1)$, 
$E_2\in K_G^0(Y_2)\cong K^G_0(Y_2)$, we have
\begin{equation*}
  E_1 \otimes^L_X E_2
 = \sum_i (-1)^i \Wedge^i N \otimes E_1|_Y\otimes E_2|_Y 
 \in K_G^0(Y)\cong K^G_0(Y),
\end{equation*}
where $N \defeq T_{Y_1}^*X|_Y\cap T_{Y_2}^*X|_Y$.
\end{Lemma}

\subsection{Pull-back with support}\label{subsec:pullback}
(cf.\ \cite[1.2]{BFM}, \cite[5.2.5]{Gi-book})
Let $f\colon Y\to X$ be a $G$-equivariant morphism between {\it
nonsingular\/} $G$-varieties. Suppose that $X'$ and $Y'$ are
$G$-invariant closed subvarieties of $X$ and $Y$ respectively
satisfying $f^{-1}(X')\subset Y'$. Then the pull-back
\begin{equation*}
   E^\bullet \mapsto f^* E^\bullet
\end{equation*}
induces a homomorphism $K^G(X;X') \to K^G(Y;Y')$. Via isomorphisms
$K^G(X')\cong K^G(X;X')$, $K^G(Y')\cong K^G(Y;Y')$, we get
a homomorphism $K^G(X')\to K^G(Y')$. Note that this depends on the
ambient spaces $X$, $Y$.

Let $f\colon Y\to X$ as above.
Suppose that $X'_1, X'_2\subset X$, $Y'_1, Y'_2\subset Y$ are
$G$-invariant closed subvarieties such that $f^{-1}(X'_a) \subset Y'_a$
for $a=1,2$.
Then we have
\begin{equation}\label{eq:pull_product}
   f^*(E_1\otimes^L_X E_2) = f^*(E_1)\otimes^L_Y f^*(E_2)
\end{equation}
for $E_a\in K^G(X'_a)$ ($a=1,2$).

\subsection{Pushforward}\label{subsec:push}

Let $f\colon X\to Y$ be a proper $G$-equivariant morphism between
$G$-varieties (not necessarily nonsingular). Then we have a
pushforward homomorphism $f_*\colon K^G(X)\to K^G(Y)$ defined by
\begin{equation*}
   f_*[E] \defeq \sum (-1)^i [R^i f_* E].
\end{equation*}

Suppose further that $X$ and $Y$ are nonsingular. If $X'\subset X$,
$Y'\subset Y$ are $G$-invariant closed subvarieties, we have the
following projection formula (\cite[5.3.13]{Gi-book}):
\begin{equation}
\label{eq:projection}
   f_*(E \otimes^L_X f^* F) = f_* E \otimes^L_Y F
   \in K^G\left(f(X')\cap Y'\right)
\end{equation}
for $E\in K^G(X')$, $F\in K^G(Y')$.

\subsection{Chow group and homology group}
Let $H_*(X,\Z) = \bigoplus_k H_k(X,\Z)$ be the integral Borel-Moore
homology of $X$. Let $A_*(X) = \bigoplus_k A_k(X)$ be the Chow group
of $X$. We have a cycle map
\begin{equation*}
   A_*(X) \to H_*(X,\Z),
\end{equation*}
which has certain functorial properties (see
\cite[Chapter~19]{Fulton}).

If $Y$ is a closed subvariety of $X$ and $U = X\setminus Y$ is its
complement, then we have exact sequences which are analogue of
\eqref{eq:exact}, \eqref{eq:exact'}:
\begin{gather}
   A_k(Y) \xrightarrow{i_*} A_k(X) \xrightarrow{j^*} A_k(U) \to 0,
   \label{eq:exact2}
\\
   \cdots \to H_k(Y,\Z) \xrightarrow{i_*} H_k(X,\Z) \xrightarrow{j^*}
   H_k(U,\Z) \xrightarrow{\partial_*} H_{k-1}(Y,\Z) \to
   \cdots.\label{eq:exact3}
\end{gather}

We have operations on $A_*(X)$ and $H_*(X,\Z)$ which are analogue of
those in \subsecref{subsec:tor}, \subsecref{subsec:pullback},
\subsecref{subsec:push}. (See \cite{Fulton}.)

In the next section, we prove results for $K$-homology and Chow group
in parallel arguments. It is the reason why we avoid higher algebraic
$K$-homology. There is no analogue for the Chow group.

\section{Freeness}\label{sec:freeness}

\subsection{Properties $(S)$,$(T)$,$(T')$}

Following \cite{DLP,Lu-Base2}, we say that an algebraic variety $X$ {\it
has property\/} $(S)$ if
\begin{aenume}
\item $H_{\operatorname{odd}}(X, \mathbb Z) = 0$ and
$H_{\operatorname{even}}(X,\mathbb Z)$ is a free abelian group.
\item The cycle map $A_*(X)\to H_{\operatorname{even}}(X,\mathbb Z)$
is an isomorphism.
\end{aenume}

Similarly, we say $X$ {\it has property\/} $(T)$ if
\begin{aenume}
\item $\Kotop(X) = 0$ and $\Ktop(X) = \Kztop(X)$ is a free abelian
group.
\item The comparison map $K(X)\to \Ktop(X)$
is an isomorphism.
\end{aenume}

Suppose that $X$ is a closed subvariety of a nonsingular variety $M$.
We have a diagram (see \cite{BFM})
\begin{equation*}
\begin{CD}
   K(X)\otimes\Q @>>> A_*(X)\otimes\Q \\
     @VVV @VVV \\
   \Ktop(X)\otimes\Q @>>>
     H_{\operatorname{even}}(X,\Q),
\end{CD}
\end{equation*}
where the horizontal arrows are local Chern character homomorphisms
in algebraic and topological $K$-homologies respectively, the left
vertical arrow is a comparison map, the right vertical arrow is the
cycle map.
It is known that the upper horizontal arrow is an
isomorphism (\cite[15.2.16]{Fulton}). Thus the composite
$K(X)\otimes\Q\to H_{\operatorname{even}}(X,\Q)$ is an isomorphism if
$X$ has property $(S)$.

Assume that $X$ is nonsingular and projective. We define the bilinear
pairing $K(X)\otimes K(X)\to \mathbb Z$ by
\begin{equation}
\label{eq:pair}
  F\otimes F' \mapsto p_*(F\otimes^L_X F'),
\end{equation}
where $p$ is the canonical map from $X$ to the point.

We say that $X$ has property $(T')$ if $X$ has property $(T)$ and the
pairing~\eqref{eq:pair} is perfect.
(In \cite{Lu-Base2}, this property is called $(S')$.)

Let $G$ be a linear algebraic group. Let $X$ be an algebraic variety
with a $G$-action. We say that $X$ {\it has property\/} $(T_G)$ if
\begin{aenume}
\item $K^{G}_{1,\operatorname{top}}(X) = 0$ and
$\Ktop^{G}(X) = K^{G}_{0,\operatorname{top}}(X)$ is a
free $R_G$-module,
\item The natural map $K^G(X)\to \Ktop^{G}(X)$
is an isomorphism.
\item For a closed algebraic subgroup $H\subset G$, $H$-equivariant
$K$-theories satisfy the above properties (a), (b), and the natural
homomorphism $K^G(X)\otimes_{R(G)} R(H)\to K^H(X)$ is an isomorphism.
\end{aenume}

Suppose further that $X$ is smooth and projective. By the same formula 
as \eqref{eq:pair}, we have a bilinear pairing
$K^G(X)\otimes K^G(X)\to R(G)$. We say that $X$ {\it has property\/}
$(T_G')$ if $X$ has property $(T_G)$ and this pairing is perfect.

A finite partition of a variety $X$ into locally closed subvarieties
is said to be an {\it $\alpha$-partition\/} if the subvarieties in the
partition can be indexed $X_1, \dots, X_n$ in such a way that $X_1\cup
X_2\cup \cdots \cup X_i$ is closed in $X$ for $i=1,\dots,k$. The
following is proved in \cite[Lemma~1.8]{DLP}.

\begin{Lemma}\label{lem:alpha}
If $X$ has an $\alpha$-partition into pieces which have property
$(S)$, then $X$ has property $(S)$.
\end{Lemma}

The proof is based on exact sequences
\eqref{eq:exact2},\eqref{eq:exact3} in homology groups and Chow
groups. Since we have corresponding exact sequences
\eqref{eq:exact},\eqref{eq:exact'} in the $K$-theory, we have the
following.
\begin{Lemma}\label{lem:alpha'}
Suppose that an algebraic variety $X$ has an action of a linear
algebraic group $G$. If $X$ has an $\alpha$-partition into
$G$-invariant locally closed subvarieties which have property
$(T_G)$, then $X$ has property $(T_G)$.
\end{Lemma}

\begin{Lemma}\label{lem:vector_bundle}
Let $\pi\colon E\to X$ be a $G$-equivariant fiber bundle with affine
spaces as fibers. Suppose that $\pi$ is locally a trivial
$G$-equivariant vector bundle, i.e.\ a product of base and a vector
space with a linear $G$-action. If $X$ has property $(T_G)$
\rom(resp.\ $(S)$\rom),
then $E$ also has property $(T_G)$ \rom(resp.\ $(S)$\rom).
\end{Lemma}

\begin{proof}
We first show that $\pi^*\colon K^G(X) \to K^G(E)$ is
surjective. Choose a closed subvariety $Y$ of $X$ so that $E$ is a
trivial $G$-bundle over $U = X\setminus Y$. There is a diagram
\begin{equation*}
\begin{CD}
   K^G(Y) @>>> K^G(X) @>>> K^G(U) @>>> 0 \\
    @VV{\pi^*}V       @VV{\pi^*}V        @VV{\pi^*}V  \\
   K^G(\pi^{-1}(Y)) @>>> K^G(E) @>>> K^G(\pi^{-1}(U)) @>>> 0,
\end{CD}
\end{equation*}
with exact rows by \eqref{eq:exact}.
By a diagram chase it suffices to prove the surjectivity for the
restrictions of $E$ to $U$ and to $Y$. By repeating the process on
$Y$, it suffices to prove it for the case $E$ is a trivial
$G$-equivariant bundle.
By Thom isomorphism~\cite[4.1]{Thom} $\pi^*$ is an isomorphism if $E$
is a $G$-equivariant bundle. Thus we prove the assertion.

Let us repeat the same argument for $\pi^*\colon \Kztop^G(X)\to
\Kztop^G(E)$ and $\pi^*\colon \Kotop^G(X)\to \Kotop^G(E)$ by replacing
\eqref{eq:exact} by \eqref{eq:exact'}. By five lemma both $\pi^*$ are
isomorphisms.
In particular, we have $\Kotop^G(E) \cong \Kotop^G(X) = 0$ by
assumption.

Consider the diagram
\begin{equation*}
\begin{CD}
   K^G(X) @>\pi^*>> K^G(E) \\
   @VVV             @VVV \\
   \Ktop^G(X) @>\pi^*>> \Ktop^G(E),
\end{CD}
\end{equation*}
where the vertical arrows are comparison maps.
The left vertical arrow is an isomorphism by assumption. Thus the
right vertical arrow is also an isomorphism by the commutativity of
the diagram and what we just proved above. The condition~(c) for
$(T_G)$ can be checked in the same way, and $E$ has property $(T_G)$.

The property $(S)$ can be checked in the same way.
\end{proof}

\begin{Lemma}\label{lem:BB}
Let $X$ be a nonsingular quasi-projective variety with
$G\times\C^*$-action with a K\"ahler metric $g$ such that
\begin{aenume}
\item $g$ is complete,
\item $g$ is invariant under the maximal compact subgroup of
$G\times\C^*$,
\item there exists a moment map $f$ associated with the K\"ahler
metric $g$ and the $S^1$-action \rom(the maximal compact subgroup of
the second factor\rom), and it is proper.
\end{aenume}
Let
\begin{equation*}
   L \defeq \{ x\in X \mid
   \text{$\displaystyle\lim_{t\to\infty} t.x$ exists} \}.
\end{equation*}
If the fixed point set $X^{\C^*}$ has property
$(T_{G\times\C^*})$ \rom(resp.\ $(S)$\rom),
then both $X$ and $L$ have property
$(T_{G\times\C^*})$ \rom(resp.\ $(S)$\rom).

Furthermore, the bilinear pairing
\begin{equation}\label{eq:pairing}
   K^{G\times\C^*}(X) \times K^{G\times\C^*}(L) \ni (F, F') 
   \longmapsto p_*(F\otimes^L_X F')\in R(G\times\C^*)
\end{equation}
is nondegenerate if $X^{\C^*}$ has property $(T_{G\times\C^*}')$.
Similar intersection pairing between $A_*(X)$ and $A_*(L)$ is
nondegenerate if $X^{\C^*}$ has property $(S)$.
Here $p$ is the canonical map from $X$ to the point.
\end{Lemma}

\begin{proof}
By \cite[2.2]{Atiyah} the moment map $f$ is a Bott-Morse function, and
critical manifolds are the fixed point $X^{\C^*}$.
Let $F_1$, $F_2$, \dots be the components of $X^{\C^*}$.
By \cite[\S3]{Atiyah}, stable and unstable manifolds for the gradient
flow of $-f$ coincide with $(\pm)$-attracting sets of
Bialynicki-Birula decomposition~\cite{BiaB}:
\begin{equation*}
   S_k = \{ x\in X\mid \lim_{t\to 0} t\cdot x\in F_k \},\qquad
   U_k = \{ x\in X\mid \lim_{t\to \infty} t\cdot x \in F_k \}.
\end{equation*}
These are invariant under the $G$-action since the $G$-action commutes 
with the $\C^*$-action.

Note that results in \cite{Atiyah} are stated for compact manifolds,
but the argument can be modified to our setting. A difference is that
$\bigcup_k U_k = L$ is not $X$ unless $X$ is compact.
On the other hand,  $\bigcup_k S_k$ is $X$ since $f$ is proper.

As in \cite{Atiyah2}, we can introduce an ordering on the index set
$\{ k\}$ of components of $X^{\C^*}$ such that
$X = \bigcup S_k$ is an $\alpha$-partition and $L = \bigcup U_k$ is
an $\alpha$-partition with respect to the reversed order.

By \cite{BiaB} (see also \cite{CS} for analytic arguments), the maps
\begin{equation*}
   S_k \ni x \mapsto \lim_{t\to 0} t\cdot x\in F_k, \qquad
   U_k \ni x \mapsto \lim_{t\to \infty} t\cdot x\in F_k
\end{equation*}
are fiber bundles with affine spaces as fibers.
Furthermore, $S_k$ (resp.\ $U_k$) is locally isomorphic to a
$G\times\C^*$-equivariant vector bundle by the proof.
Thus $S_k$ and $U_k$ have properties $(S)$ and $(T_{G\times\C^*})$ by
\lemref{lem:vector_bundle}.
Hence $X$ and $L$ have properties $(S)$ and $(T_{G\times\C^*})$ by
Lemmas~\ref{lem:alpha}, \ref{lem:alpha'}.

By the argument in \cite[1.7, 2.5]{Lu-Base2}, the pairing
\eqref{eq:pairing} can be identified with a pairing
\[
   \bigoplus_k K^{G\times\C^*} (F_k) \times
     \bigoplus_k K^{G\times\C^*} (F_k) \to R(G\times\C^*)
\]
of the form
\[
   \left( \sum_k \xi_k, \sum_k \xi'_k \right)
   = \sum_{k\ge k'} (\xi_k, \xi'_{k'})_{k,k'}.
\]
for some pairing $(\ ,\ )_{k,k'}\colon
K^{G\times\C^*}(F_k)\times K^{G\times\C^*}(F_k')\to R(G\times\C^*)$
such that $(\ ,\ )_{k,k}$ is the pairing \eqref{eq:pair} for $X =
F_k$. Since $(\ ,\ )_{k,k}$ is nondegenerate for all $k$ by the assumption,
\eqref{eq:pairing} is also nondegenerate.

The proof of the statement for $A_*(X)$, $A_*(L)$ is similar. One uses 
the fact that the intersection pairing $A_*(F_k)\times A_*(F_k)\to \Z$ 
is nondegenerate under property $(S)$.
\end{proof}

\subsection{Decomposition of the diagonal}

\begin{Proposition}[\protect{cf.\ \cite{ES},\cite[5.6.1]{Gi-book}}]
\label{prop:decdiag}
Let $X$ be a nonsingular projective variety.

\rom{(1)} Let $\shfO_{\Delta X}$ be the structure sheaf of the
diagonal and $[\shfO_{\Delta X}]$ the corresponding element in
$K(X\times X)$. Assume that
\begin{equation}\label{eq:decdiag}
  [\shfO_{\Delta X}] = \sum_i \alpha_i\boxtimes \beta_i
\end{equation}
holds for some $\alpha_i$, $\beta_i\in K(X)$. Then $X$ has property
$(T')$.

\rom{(2)} Let $G$ be a linear algebraic group.
Suppose that $X$ has $G$-action and that \eqref{eq:decdiag} holds in
$K^G(X\times X)$ for some $\alpha_i$, $\beta_i\in K^G(X)$. Then $X$
has property $(T'_G)$.

\rom{(3)} Let $[\Delta X]$ be the class of the diagonal in
$A(X\times X)$. Assume that
\begin{equation}\label{eq:decdiagHom}
  [\Delta X] = \sum_i p_1^* a_i\cup p_2^*b_i
\end{equation}
holds for some $a_i$, $b_i\in A(X)$. Then $X$ has property
$(S)$.
\end{Proposition}

\begin{proof}
Let $p_a\colon X\times X\to X$ denote the projection to the $a$th factor
($a= 1,2$). Let $\Delta$ be the diagonal embedding $X \to X\times
X$. Then we have $[\shfO_{\Delta X}] = \Delta_* [\shfO_X]$. Hence
\begin{alignat*}{2}
   p_{1*}\left(p_2^* F\otimes^L_{X\times X}
    [\shfO_{\Delta X}]\right)
  & = p_{1*}\left(p_2^* F\otimes^L_{X\times X}
    \Delta_* [\shfO_{X}]\right) &&\\
  &= p_{1*}\Delta_* \left( \Delta^* p_2^* F\otimes^L_X
    [\shfO_X]\right)
  &\quad &\text{(by the projection formula)}\\
  &= F\otimes^L_X [\shfO_X] &&(p_1\circ\Delta = p_2\circ\Delta =
    \id_X)\\
  &=F&&
\end{alignat*}
If we substitute \eqref{eq:decdiag} into the above, we get
\begin{equation}\label{eq:Fspan}
  F = \sum_i p_{1*}\left(p_2^* F\otimes^L_{X\times X}
    p_1^*\alpha_i \otimes^L_{X\times X} p_2^*\beta_i\right)
    = \sum_i (F,\beta_i) \alpha_i.
\end{equation}
In particular, $K(X)$ is spanned by $\alpha_i$'s.

If $mF = 0$ for some $m\in\mathbb Z\setminus\{0\}$, then
$0 = (mF, \beta_i) = m(F,\beta_i)$. Hence we have $(F,\beta_i) =
0$. The above equality~\eqref{eq:Fspan} implies $F = 0$. This means
that $K(X)$ is torsion-free. Thus we could assume $\alpha_i$'s are
linearly independent in \eqref{eq:decdiag}. Under this assumption,
$\{\alpha_i\}$ is a basis of $K(X)$, and \eqref{eq:Fspan} implies that 
$\{ \beta_i \}$ is the dual basis.

If we perform the same computation in $\Kztop(X)\oplus\Kotop(X)$, we
get the same result.
In particular, $\{ \alpha_i \}$ is a basis of
$\Kztop(X)\oplus\Kotop(X)$. However, $\alpha_i$, $\beta_i$ are in
$\Kztop(X)$, thus we have $\Kotop(X)= 0$. We also have $K(X)\to
\Kztop(X)$ is an isomorphism. Thus $X$ has property $(T')$.

If $X$ has $G$-action and \eqref{eq:decdiag} holds in the equivariant
$K$-group, we do the same calculation in the equivariant
$K$-groups. Then the same argument shows that $X$ has property
$(T'_G)$.

The assertion for Chow groups and homology groups can be proved
in the same way.
\end{proof}

\subsection{Diagonal of the quiver variety}\label{subsec:diagquiver}

Let us recall the decomposition of the diagonal of the quiver variety
defined in \cite[Sect.~6]{Na-alg}.
In this section, we fix dimension vectors $\bv$, $\bw$ and use the
notation $\M$ instead of $\M(\bv,\bw)$.

Let us consider the product $\M\times\M$. We denote by
$V^1_k$ (resp.\ $V^2_k$) the vector bundle $V_k\boxtimes\shfO_\M$
(resp.\ $\shfO_\M\boxtimes V_k$). A point in $\M\times\M$ is denoted by
$([B^1,i^1,j^1], [B^2, i^2, j^2])$. We regard $B^a$, $i^a$, $j^a$
($a=1,2$) as homomorphisms between tautological bundles.

We consider the following 
$G_\bw\times\C^*$-equivariant complex of vector bundles over $\M\times\M$:
\begin{equation}
  \HomL(V^1, V^2)
  \overset{\sigma}{\longrightarrow}
  q\HomE(V^1,V^2) \oplus 
    q\HomL(W, V^2) \oplus q\HomL(V^1,W)
  \overset{\tau}{\longrightarrow}
    q^2\HomL(V^1, V^2),
\label{eq:diagcpx}\end{equation}
where
\begin{align*}
        \sigma(\xi) & \defeq (B^2 \xi - \xi B^1) \oplus
         (-\xi i^1) \oplus j^2 \xi \\
        \tau(C\oplus a\oplus b)
         &\defeq ( \varepsilon B^2 C + \varepsilon C B^1 + i^2 b + a j^1).
\end{align*}

It was shown that $\sigma$ is injective and $\tau$ is surjective (cf.\ 
\cite[5.2]{Na-alg}). Thus $\Ker\tau/\Ima\sigma$ is an
equivariant vector bundle. We define an equivariant section $s$ of
$\Ker\tau/\Ima\sigma$ by
\begin{equation*}
  s \defeq \left( 0 \oplus (-i^2) \oplus j^1 \right) \bmod{\Ima \sigma}.
\end{equation*}
Then $([B^1,i^1,j^1], [B^2,i^2,j^2])$ is contained in the zero locus
$Z(s)$ of $s$ if and only if there exists $\xi\in \HomL(V^1, V^2)$ such
that
\begin{equation*}
        \xi B^1 = B^2 \xi, \quad
        \xi i^1 = i^2, \quad
        j^1 = j^2 \xi.
\end{equation*}
Moreover $\xi$ is an isomorphism by the stability condition. Hence
$Z(s)$ is equal to the diagonal $\Delta\M$.
If $\nabla$ is a connection on $\Ker\tau/\Ima\sigma$, the differential 
$\nabla s\colon T(\M\times\M)\to \Ker\tau/\Ima\sigma$ is surjective on 
$Z(s)=\Delta\M$ (cf.\ \cite[5.7]{Na-alg}).
In particular, we have an exact sequence
\begin{equation*}
   0 \to \Wedge^{\max} (\Ker\tau/\Ima\sigma)^* \to \cdots
   \to (\Ker\tau/\Ima\sigma)^* \to \shfO_{\M\times\M} \to
   \shfO_{\Delta \M} \to 0.
\end{equation*}

In $K^{G_\bw\times\C^*}(\M\times\M)$, $\Ker\tau/\Ima\sigma$ is equal
to the alternating sum of terms of \eqref{eq:diagcpx} which has a form
$\sum \alpha_i \boxtimes \beta_i$ for some $\alpha_i,\beta_i\in
K^{G_\bw\times\C^*}(\M)$. Hence $\M$ satisfies the conditions of
\propref{prop:decdiag} except the projectivity.  Unfortunately, the
projectivity is essential in the proof of \propref{prop:decdiag}. (We
could not define $(p_1)_*$ otherwise.) Thus \propref{prop:decdiag} is
not directly applicable to $\M$. In order to get rid of this
difficulty, we consider the fixed point set $\M^{\C^*}$ with respect
to the $\C^*$-action.

By a technical reason, we need to use a $\C^*$-action, which is
different from \eqref{eq:C*action}. Let $\C^*$ act on $\bM$ by
\begin{equation}\label{eq:newC*action}
B_{h} \mapsto t B_{h}, \quad i \mapsto t i, \quad j \mapsto t j \qquad
\text{for $t\in \C^*$}.
\end{equation}
This induces a $\C^*$-action on $\M$ and $\M_0$ which commutes with
the previous $G_\bw\times\C^*$-action.
(If the adjacency matrix satisfies $\bA_{kl} \le 1$ for any $k, l\in I$,
then the new $\C^*$-action coincides with the old one.)
The tautological bundles $V_k$, $W_k$ become $\C^*$-equivariant vector
bundles as before.

We consider the fixed point set $\M^{\C^*}$. $[B,i,j]\in\M$ is a fixed
point if and only if there exists a homomorphism $\rho\colon
\C^*\to G_\bv$ such that
\begin{equation*}
  t\diamond (B,i,j) = \rho(t)^{-1}\cdot (B,i,j)
\end{equation*}
as in \subsecref{subsec:fix_hom}.
Here $\diamond$ denotes the new $\C^*$-action.
We decompose the fixed point set $\M^{\C^*}$ according to the
conjugacy class of $\rho$:
\begin{equation*}
\M^{\C^*} = \bigsqcup \M[\rho].  
\end{equation*}

\begin{Lemma}
$\M[\rho]$ is a nonsingular projective variety.
\end{Lemma}

\begin{proof}
Since $\M[\rho]$ is a union of connected components (possibly single
component) of the fixed point set of the $\C^*$-action on a
nonsingular variety $\M$, $\M[\rho]$ is nonsingular.

Suppose that $[B,i,j]\in \M_0$ is a fixed point of the
$\C^*$-action. It means that $(tB, ti, tj)$ lies in the closed orbit
$G\cdot(B,i,j)$. But $(tB,ti,tj)$ converges to $0$ as $t\to 0$. Hence
the closed orbit must be $\{0\}$. Since $\pi\colon\M\to\M_0$ is
equivariant, $\M^{\C^*}$ is contained in $\pi^{-1}(0)$. In particular,
$\M[\rho]$ is projective.
\end{proof}

This lemma is not true for the original $\C^*$-action.

We restrict the complex~\eqref{eq:diagcpx} to
$\M[\rho]\times\M[\rho]$. Then fibers of $V^1$ and $V^2$ become
$\C^*$-modules and hence we can take the $\C^*$-fixed part of
\eqref{eq:diagcpx}:
\begin{equation*}
  \HomL(V^1, V^2)^{\C^*}
  \overset{\sigma^{\C^*}}{\longrightarrow}
  \left(
    q \HomE(V^1,V^2) \oplus q \HomL(W, V^2) \oplus q \HomL(V^1,W)
  \right)^{\C^*}
  \overset{\tau^{\C^*}}{\longrightarrow}
    (q^2 \HomL(V^1, V^2))^{\C^*},
\end{equation*}
where $\sigma^{\C^*}$ (resp.\ $\tau^{\C^*}$) is the restriction of
$\sigma$ (resp.\ $\tau$) to the ${\C^*}$-fixed part. Then $\sigma^{\C^*}$ 
is injective and $\tau^{\C^*}$ is surjective, and
$\Ker\tau^{\C^*}/\Ima\sigma^{\C^*}$ is a vector bundle which is
the ${\C^*}$-fixed part of $\Ker\tau/\Ima\sigma$.

The section $s$ takes values in
$\Ker\tau^{\C^*}/\Ima\sigma^{\C^*} =
(\Ker\tau/\Ima\sigma)^{\C^*}$. Considering it as a section of
$\Ker\sigma^{\C^*}/\Ima\tau^{\C^*}$, we denote it by $s^{\C^*}$. The
zero locus $Z(s^{\C^*})$ is $Z(s)\cap(\M[\rho]\times\M[\rho])$ which
is the diagonal $\Delta\M[\rho]$ of
$\M[\rho]\times\M[\rho]$. Furthermore, the differential 
$\nabla s^{\C^*}\colon T(\M[\rho]\times\M[\rho])\to
(\Ker\tau/\Ima\sigma)^{\C^*}$ is surjective on $Z(s^{\C^*})=\Delta\M[\rho]$.

Our original $G_\bw\times \C^*$-action (defined in
\subsecref{subsec:Gw-action}) commutes with the new
$\C^*$-action. Thus $\M[\rho]$ has an induced $G_\bw\times
\C^*$-action. By the construction, $\Ker\tau^{\C^*}/\Ima\sigma^{\C^*}$ 
is a $G_\bw\times \C^*$-equivariant vector bundle, and $s^{\C^*}$ is
an equivariant section.

\begin{Proposition}\label{prop:fix_conn}
$\M[\rho]$ has properties $(S)$ and $(T_{G_\bw\times \C^*}')$.
Moreover, $\M[\rho]$ is connected.
\end{Proposition}

\begin{proof}
Let $\shfO_{\Delta\M[\rho]}$ be the structure sheaf of the diagonal
considered as a sheaf on $\M[\rho]\times\M[\rho]$. By the above
argument, the Koszul complex of $s^{\C^*}$ gives a resolution of
$\shfO_{\Delta\M[\rho]}$:
\begin{equation*}
  0 \to \Wedge^{\max} 
  (\Ker\tau^{\C^*}/\Ima\sigma^{\C^*})^* \to
  \cdots \to (\Ker\tau^{\C^*}/\Ima\sigma^{\C^*})^*
  \to \shfO_{\M[\rho]\times\M[\rho]}
  \to \shfO_{\Delta\M[\rho]} \to 0,
\end{equation*}
where $\max = \rank\Ker\sigma^{\C^*}/\Ima\tau^{\C^*}$.
Thus we have the following equality in the Grothendieck group
$K^{G_\bw\times\C^*}(\M[\rho]\times\M[\rho])$
\begin{equation*}
   [\shfO_{\Delta\Fi}] =
   \Wedge_{-1} [(\Ker\tau^{\C^*}/\Ima\sigma^{\C^*})^*].
\end{equation*}
Since $\sigma^{\C^*}$ is injective and $\tau^{\C^*}$ is surjective, we have
\begin{equation*}
\begin{split}
  & [\Ker\tau^{\C^*}/\Ima\sigma^{\C^*}] \\
  =\, & 
  - [ \HomL(V^1, V^2)^{\C^*} ]
  + \left[ (q\HomE(V^1,V^2) \oplus
    q \HomL(W, V^2) \oplus q\HomL(V^1,W)))^{\C^*}\right]
  - [(q^2 \HomL(V^1, V^2))^{\C^*}].
\end{split}
\end{equation*}

Each factor of the right hand side can be written in the form
$\sum_i \alpha_i\boxtimes \beta_i$ for some
$\alpha_i, \beta_i\in K^{G_\bw\times\C^*}(\M[\rho])$.
For example, the first factor is equal to
\begin{equation*}
   \HomL(V^1, V^2)^{\C^*} = \bigoplus_m \HomL(V^1(m), V^2(m)),
\end{equation*}
where $V^a(m)$ is the weight space of $V^a$, i.e.,
\begin{equation*}
   V^a(m) = \{ v\in V^a \mid t\diamond v = t^m v\}.
\end{equation*}
The remaining factors have similar description. Thus by
\propref{prop:decdiag}, $\M[\rho]$ has property
$(T_{G_\bw\times\C^*}')$.

Moreover, the above shows that $K^{G_\bw\times\C^*}(\M[\rho])$ is
generated by exterior powers of $V_k(m)$, $W_k(m)$ and its duals
(as an $R(G_\bw\times\C^*)$-algebra).
Note that these bundles have constant rank on $\M[\rho]$. If
$\M[\rho]$ have components $M_1$, $M_2$, \dots, the structure sheaf of
$M_1$ (extended to $\M[\rho]$ by setting $0$ outside) cannot be
represented by $V_k(m)$, $W_k(m)$.  This contradiction shows that
$\M[\rho]$ is connected.

The assertion for Chow groups can be proved exactly in the same
way. By the above argument, the fundamental class $[\Delta\M[\rho]]$
is the top Chern class of $\Ker\tau^{\C^*}/\Ima\sigma^{\C^*}$, which
can be represented as $\sum_i p_1^* a_i\cup p_2^* b_i$ for some $a_i$, 
$b_i\in A(X)$.
\end{proof}

\begin{Theorem}\label{thm:freeness}
$\M$ and $\La$ have properties $(S)$ and $(T_{G_\bw\times\C^*})$.
Moreover, the bilinear pairing
\begin{equation*}
   K^{G_\bw\times\C^*}(\M) \times K^{G_\bw\times\C^*}(\La) \ni (F, F') 
   \longmapsto p_*(F\otimes^L_{\M} F')\in R(G_\bw\times\C^*)
\end{equation*}
is nondegenerate. Similar pairing between $A_*(\M)$ and $A_*(\La)$ is
also nondegenerate.
Here $p$ is the canonical map from $\M$ to the point.
\end{Theorem}

\begin{proof}
We apply \lemref{lem:BB}. By \cite[2.8]{Na-quiver}, the metric on $\M$ 
defined in \subsecref{subsec:hK} is complete. By the construction, it
is invariant under $K_\bw\times S^1$, where $K_\bw = \prod \U(W_k)$ is 
the maximal compact subgroup of $G_\bw$.
(Note that the hyper-K\"ahler structure is {\it not\/} invariant under 
the $S^1$-action, but the metric is invariant.) The moment map for the
$S^1$-actions is given by
\begin{equation*}
   \frac12 \left( \sum_h \| B_h \|^2 
     + \sum_k \left( \| i_k\|^2 + \| j_k \|^2 \right) \right).
\end{equation*}
This is a proper function on $\M$. Thus \lemref{lem:BB} is
applicable. Note that we have
\(
   \La = \{ x\in \M \mid
   \text{$\displaystyle\lim_{t\to\infty} t.x$ exists} \}
\)
as in \cite[5.8]{Na-quiver}. (Though our $\C^*$-action is different
from one in \cite{Na-quiver}, the same proof works.)
\end{proof}

\subsection{Fixed point subvariety}
Let $A$ be an abelian reductive subgroup of $G_\bw\times\C^*$ as in
\secref{sec:Fixed}. Let $\M^A$ and $\La^A$ be the fixed point set in
$\M$ and $\La$ respectively. Exactly as in the previous subsection, we
have the following generalization \thmref{thm:freeness}.
\begin{Theorem}\label{thm:freeness'}
$\M^A$ and $\La^A$ have properties $(T)$ and $(S)$.
Moreover, the bilinear pairing
\begin{equation*}
   K(\M^A) \times K(\La^A) \ni (F, F') 
   \longmapsto p_*(F\otimes^L_{\M^A} F')\in \Z
\end{equation*}
is nondegenerate. Similar pairing between $A_*(\M^A)$ and $A_*(\La^A)$ is
also nondegenerate.
Here $p$ is the canonical map from $\M^A$ to the point.
\end{Theorem}

\subsection{Connectedness of $\M(\bv,\bw)$}\label{subsec:connected}

Let us consider a natural homomorphism 
\begin{equation}
\label{eq:nathom}
   R(G_\bw\times\C^*\times G_\bv) \to K^{G_\bw\times\C^*}(\M),
\end{equation}
which sends representations to bundles associated with tautological
bundles.
If we can apply \propref{prop:decdiag} to $\M$, then this homomorphism 
is surjective.
Unfortunately we can not apply \propref{prop:decdiag} since $\M$ is
not projective.
However, it seems reasonable to conjecture that the homomorphism
\eqref{eq:nathom} is surjective.
In particular, it implies that $\M$ is connected as in the proof of
\propref{prop:fix_conn}. This was stated in \cite[6.2]{Na-alg}. But the 
proof contains a gap since the function $\| s_1 \|$ may not be proper
in general.

\section{Convolution}\label{sec:convolution}


Let $X_1$, $X_2$, $X_3$ be a nonsingular quasi-projective variety, and
write $p_{ab}\colon X_1\times X_2\times X_3\to X_a\times X_b$ for the
projection ($(a,b) = (1, 2), (2,3), (1,3)$). 

Suppose $Z_{12}$ (resp.\ $Z_{23}$) is a closed subvariety of
$X_1\times X_2$ (resp.\ $X_2\times X_3$) such that the restriction of
the projection $p_{13}\colon p_{12}^{-1}(Z_{12})\cap
p_{23}^{-1}(Z_{23})\to X_1\times X_3$ is proper. Let $Z_{12}\circ
Z_{23} = p_{13}\left(p_{12}^{-1}(Z_{12})\cap
  p_{23}^{-1}(Z_{23})\right)$.  We can define the {\it convolution
product\/} \(
  \ast\colon K(Z_{12})\otimes K(Z_{23})\to K(Z_{12}\circ Z_{23})
\)
by
\begin{equation*}
  K_{12}* K_{23}
  \defeq p_{13*} \left(
  p_{12}^* K_{12}\otimes^L_{X_1\times X_2\times X_3} p_{23}^*
  K_{23}\right)
  \qquad\text{for $K_{12}\in K(Z_{12}), K_{23}\in K(Z_{23})$}.
\end{equation*}

Note that the convolution product depends on the ambient spaces $X_1$, 
$X_2$ and $X_3$. When we want to specify them, we say the convolution
product relative to $X_1$, $X_2$, $X_3$.

In this section, we study what happens when
$X_1$, $X_2$, $X_3$ are replaced by 
\begin{aenume}
\item submanifolds $S_1$, $S_2$, $S_3$ of $X_1$, $X_2$, $X_3$,
\item principal $G$-bundles $P_1$, $P_2$, $P_3$ over $X_1$, $X_2$, $X_3$.
\end{aenume}

Although we work on the non-equivariant $K$-theory, the results
extend to the case of the equivariant $K$-theory, the Borel-Moore
homology group, or any other reasonable theory in the straightforward way.

\subsection{}
Before studying the above problem, we recall the following lemma which 
will be used several times.
\begin{Lemma}\label{lem:diag}
In the above setting, we further assume that $X_1 = X_2$ and 
$Z_{12} = \operatorname{Image}\Delta_{X_1}$, where $\Delta_{X_1}$ is the
diagonal embedding $X_1\to X_1\times X_2$. Then we have
\begin{equation*}
   (\Delta_{X_1})_*[E] * K_{23} = p_2^* [E]
  \otimes K_{23}
\end{equation*}
for a vector bundle $E$ over $X_1$, where 
$p_2\colon X_2\times X_3\to X_2 = X_1$ is the projection, and
$\otimes$ in the right hand side is the tensor product
\eqref{eq:tensor} between $K^0(Z_{23})$ and $K(Z_{23})$.
\end{Lemma}

The proof is obvious from the definition, and omitted.

\subsection{Restriction of the convolution to submanifolds}
Suppose we have nonsingular closed submanifolds $S_1$, $S_2$, $S_3$ of
$X_1$, $X_2$, $X_3$ such that
\begin{equation}\label{eq:SZass}
  \left(S_1\times X_2\right) \cap Z_{12} \subset S_1\times S_2, \qquad
  \left(S_2\times X_3\right) \cap Z_{23} \subset S_2\times S_3.
\end{equation}
By this assumption, we have
\begin{equation}\label{eq:SZass'}
  (S_1\times X_3)\cap (Z_{12}\circ Z_{23}) \subset S_1\times S_3.
\end{equation}
Let $Z_{12}'$ (resp.\ $Z_{23}'$) be the intersection $(S_1\times
S_2)\cap Z_{12}$ (resp.\ $(S_2\times S_3)\cap Z_{23}$). By
\eqref{eq:SZass'}, we have $Z_{12}'\circ Z_{23}' =
(S_1\times X_3)\cap (Z_{12}\circ Z_{23})$.
We have the convolution product $\ast'\colon K(Z_{12}')\otimes
K(Z_{23}') \to K(Z_{12}'\circ Z_{23}')$ relative to $S_1, S_2, S_3$:
\begin{equation*}
  K'_{12}*' K'_{23}
  \defeq p'_{13*} \left(
  p_{12}^{\prime *} K'_{12}\otimes^L_{S_1\times S_2\times S_3}
  p_{23}^{\prime *} K'_{23}\right),
\end{equation*}
where $p_{ab}'$ is the projection
$S_1\times S_2\times S_3\to S_a\times S_b$.

We want to relate two convolution products $*$ and $*'$ via pull-back
homomorphisms. For this purpose, we consider the inclusion
$i_a\times\id_{X_b}\colon S_a\times X_b\to X_a\times X_b$,
where $i_a$ is the inclusion $S_a\hookrightarrow X_a$
($(a,b) = (1,2), (2,3), (1,3)$). By \eqref{eq:SZass}, we have
a pull-back homomorphism
\begin{equation*}
  K(Z_{12})\cong K(X_1\times X_2; Z_{12}) 
  \xrightarrow{(i_1\times\id_{X_2})^*}
  K(S_1\times X_2; Z_{12}\cap S_1\times X_2) \cong K(Z_{12}').
\end{equation*}
Similarly, we have
\begin{equation*}
  K(Z_{23})\xrightarrow{(i_2\times\id_{X_3})^*}
  K(Z_{23}'), \qquad
  K(Z_{12}\circ Z_{23})\xrightarrow{(i_1\times\id_{X_3})^*}
  K(Z_{12}'\circ Z_{23}').
\end{equation*}

\begin{Proposition}\label{prop:Sconv}
For $K_{12}\in K(Z_{12}), K_{23}\in K(Z_{23})$, we have
\begin{equation}\label{eq:Sconv}
  (i_1\times \id_{X_3})^*(K_{12} * K_{23})
  = \left((i_1\times\id_{X_2})^*K_{12}\right) *'
  \left((i_2\times\id_{X_3})^* K_{23}\right).
\end{equation}
Namely, the following diagram commutes:
\begin{equation*}
\begin{CD}
  K(Z_{12})\otimes K(Z_{23}) @>{*}>> K(Z_{12}\circ Z_{23}) \\
  @V{(i_1\times\id_{X_2})^* \otimes (i_2\times\id_{X_3})^*}VV
  @VV{(i_1\times\id_{X_3})^*}V \\
  K(Z'_{12})\otimes K(Z'_{23}) @>{*'}>> K(Z'_{12}\circ Z'_{23}).
\end{CD}
\end{equation*}
\end{Proposition}

\begin{Example}
Suppose $X_1 = X_2$, $S_1 = S_2$ and $Z_{12} =
\operatorname{Image}\Delta_{X_1}$, where $\Delta_{X_1}$ is the
diagonal embedding $X_1\to X_1\times X_2$. Then the above assumption
$S_1\times X_2 \cap Z_{12} \subset S_1\times S_2$ is satisfied, and we
have $Z_{12}' = \operatorname{Image}\Delta_{S_1}$, where
$\Delta_{S_1}$ is the diagonal embedding $S_1\to S_1\times S_2$.
If $E$ is a vector bundle over $X_1$, we have
\(
  (i_1\times\id_{X_2})^*(\Delta_{X_1})_*[E]
   = (\Delta_{S_1})_*[i_1^* E]
\)
by the base change \cite[5.3.15]{Gi-book}.
By \lemref{lem:diag}, we have 
\begin{equation*}
\begin{split}
  & (\Delta_{X_1})_*[E] * K_{23} = p_2^* [E] \otimes K_{23},\\
  & (\Delta_{S_1})_*[i^* E] *' (i_2\times\id_{X_3})^* K_{23}
  = p_2^{\prime*}i_1^* [E] \otimes
   (i_2\times\id_{X_3})^* K_{23},
\end{split}
\end{equation*}
where $p_2\colon X_2\times X_3\to X_2$,
$p_2'\colon S_2\times X_3\to S_2$
are the projections. Note
\begin{equation*}
\begin{split}
  & p_2^{\prime*}i^* [E] \otimes
  (i_2\times\id_{X_3})^* K_{23} \\
 =\; & (i_2\times\id_{X_3})^* p_2^* [E] \otimes
  (i_2\times\id_{X_3})^* K_{23}
 = (i_2\times\id_{X_3})^*\left(p_2^*[E]\otimes
  K_{23}\right)
\end{split}
\end{equation*}
by \eqref{eq:pull_product}.
Hence we have \eqref{eq:Sconv} in this case.
\end{Example}

\begin{proof}[Proof of \propref{prop:Sconv}]
In order to relate $*$ relative to $X_1$, $X_2$, $X_3$
and $*'$ relative to $S_1$, $S_2$, $S_3$,
we replace $X_a$ by $S_a$ factor by factor.

{\bf Step 1}.
First we want to replace $X_1$ by $S_1$. 
We consider the following fiber square:
\begin{equation*}
\begin{CD}
S_1\times X_2\times X_3 @>{i_1\times\id_{X_2}\times\id_{X_3}}>>
X_1\times X_2\times X_3 \\
@V{p_{13}''}VV @VV{p_{13}}V \\
S_1\times X_3 @>>{i_1\times\id_{X_3}}> X_1\times X_3,
\end{CD}
\end{equation*}
where $p_{13}''$ is the projection. We have
\begin{equation}
\label{eq:Sconv1}
\begin{split}
 & (i_1\times\id_{X_3})^* (K_{12} * K_{23}) = 
(i_1\times\id_{X_3})^* p_{13*}
\left( p_{12}^* K_{12}\otimes^L_{X_1\times X_2\times X_3} p_{23}^*
  K_{23}\right) \\ 
 =\; & p_{13*}''(i_1\times\id_{X_2}\times\id_{X_3})^*
 \left(p_{12}^* K_{12}\otimes^L_{X_1\times X_2\times X_3} p_{23}^*
  K_{23}\right) \\
 =\; & p_{13*}'' \left(
   (i_1\times\id_{X_2}\times\id_{X_3})^* p_{12}^* K_{12}
   \otimes^L_{S_1\times X_2\times X_3}
   (i_1\times\id_{X_2}\times\id_{X_3})^* p_{23}^* K_{23}\right),
\end{split}
\end{equation}
where we have used the base change (\cite[5.3.15]{Gi-book}) in the
second equality and 
\eqref{eq:pull_product} in the third equality.
If $p_{12}''\colon S_1\times X_2\times X_3\to S_1\times X_2$ denotes
the projection, we have
$p_{12}\circ (i_1\times\id_{X_2}\times\id_{X_3})
= (i_1\times\id_{X_2})\circ p_{12}''$. Hence we get
\begin{equation*}
  (i_1\times\id_{X_2}\times\id_{X_3})^* p_{12}^* K_{12}
  =  p_{12}^{\prime\prime*}(i_1\times\id_{X_2})^* K_{12}.
\end{equation*}
Similarly, we have
\begin{equation*}
  (i_1\times\id_{X_2}\times\id_{X_3})^* p_{23}^* K_{23}
  = p_{23}^{\prime\prime*} K_{23},
\end{equation*}
where $p_{23}''\colon S_1\times X_2\times X_3\to X_2\times X_3$ is the 
projection. Substituting this into \eqref{eq:Sconv1}, we obtain
\begin{equation}
\label{eq:Sconv2}
  (i_1\times\id_{X_3})^* (K_{12} * K_{23}) =
  p_{13*}'' \left( p_{12}^{\prime\prime*}(i_1\times\id_{X_2})^* K_{12}
    \otimes^L_{S_1\times X_2\times X_3}
    p_{23}^{\prime\prime*} K_{23}\right).
\end{equation}

{\bf Step 2}. Next we replace $X_2$ by $S_2$. By \eqref{eq:SZass}, we
have a homomorphism
\begin{equation*}
  (\id_{S_1}\times i_2)_*\colon K(Z_{12}')\cong
  K(S_1\times S_2; Z_{12}')\to
  K(S_1\times X_2; Z_{12}') \cong K(Z_{12}'),
\end{equation*}
which is nothing but the identity operator. We will consider
$(i_1\times\id_{X_2})^* K_{12}\in K(Z_{12}')$ as an element of
$K(S_1\times S_2; Z_{12}')$ or $K(S_1\times X_2; Z_{12}')$
interchangely. 
We consider the fiber square
\begin{equation*}
\begin{CD}
  S_1\times S_2 \times X_3 @>{p_{12}'''}>>
S_1\times S_2 \\
@V{\id_{S_1}\times i_2\times\id_{X_3}}VV @VV{\id_{S_1}\times i_2}V \\
S_1\times X_2\times X_3 @>>{p_{12}''}> S_1\times X_2,
\end{CD}
\end{equation*}
where $p_{12}'''$ is the projection.
By base change \cite[5.3.15]{Gi-book}, we get
\begin{equation*}
 p_{12}^{\prime\prime*}(i_1\times\id_{X_2})^* K_{12} = 
 p_{12}^{\prime\prime*}(\id_{S_1}\times i_2)_* (i_1\times\id_{X_2})^*
 K_{12}
 = (\id_{S_1}\times i_2\times \id_{X_3})_*
 p_{12}^{\prime\prime\prime*} (i_1\times\id_{X_2})^* K_{12}.
\end{equation*}
By the projection formula \eqref{eq:projection}, we get
\begin{equation*}
\begin{split}
  & p_{12}^{\prime\prime*}(i_1\times\id_{X_2})^* K_{12}
    \otimes^L_{S_1\times X_2\times X_3}
    p_{23}^{\prime\prime*} K_{23} \\
=\; & (\id_{S_1}\times i_2\times \id_{X_3})_* \left(
  p_{12}^{\prime\prime\prime*} (i_1\times\id_{X_2})^* K_{12}
  \otimes^L_{S_1\times S_2\times X_3} 
  (\id_{S_1}\times i_2\times \id_{X_3})^* p_{23}^{\prime\prime*}K_{23}
  \right).
\end{split}
\end{equation*}
Substituting this into \eqref{eq:Sconv2}, we have
\begin{equation}\label{eq:Sconv3}
\begin{split}
  & (i_1\times\id_{X_3})^* (K_{12} * K_{23})\\
 = \; &
  p_{13*}''(\id_{S_1}\times i_2\times \id_{X_3})_*\left(
  p_{12}^{\prime\prime\prime*} (i_1\times\id_{X_2})^* K_{12}
  \otimes^L_{S_1\times S_2\times X_3}
  (\id_{S_1}\times i_2\times \id_{X_3})^* p_{23}^{\prime\prime*}K_{23}
  \right) \\
 = \; & p_{13*}'''\left(
   p_{12}^{\prime\prime\prime*} (i_1\times\id_{X_2})^* K_{12}
   \otimes^L_{S_1\times S_2\times X_3}
  p_{23}^{\prime\prime\prime*}(i_2\times\id_{X_3})^* K_{23}
  \right)
\end{split}
\end{equation}
where
$p_{13}'''\colon S_1\times S_2\times X_3 \to S_1\times X_3$,
$p_{23}'''\colon S_1\times S_2\times X_3 \to S_2\times X_3$ are the
projections, and we have used
$p_{13}''' = p_{13}''\circ (\id_{S_1}\times i_2\times \id_{X_3})$ and
$p_{23}''\circ (\id_{S_1}\times i_2\times \id_{X_3})
= (i_2\times\id_{X_3})\circ p_{23}'''$.

{\bf Step 3}. We finally replace $X_3$ by $S_3$. By \eqref{eq:SZass},
we have a homomorphism
\begin{equation*}
  (\id_{S_2}\times i_3)_*\colon K(Z_{23}')\cong
  K(S_2\times S_3; Z_{23}')\to
  K(S_2\times X_3; Z_{23}') \cong K(Z_{23}'),
\end{equation*}
which is nothing but the identity operator. We consider the fiber
square
\begin{equation*}
\begin{CD}
  S_1\times S_2 \times S_3 @>{p_{23}'}>>
S_2\times S_3 \\
@V{\id_{S_1}\times \id_{S_2}\times i_3}VV @VV{\id_{S_2}\times i_3}V \\
S_1\times S_2\times X_3 @>>{p_{23}'''}> S_2\times X_3.
\end{CD}
\end{equation*}
By base change \cite[5.3.15]{Gi-book}, we get
\begin{equation*}
  p_{23}^{\prime\prime\prime*}(i_2\times\id_{X_3})^* K_{23}
  = p_{23}^{\prime\prime\prime*}(\id_{S_2}\times i_3)_*
  (i_2\times\id_{X_3})^* K_{23}
  = (\id_{S_1}\times \id_{S_2}\times i_3)_* p_{23}^{\prime*} 
  (i_2\times\id_{X_3})^* K_{23}.
\end{equation*}
Substituting this into \eqref{eq:Sconv3}, we obtain
\begin{equation*}
\begin{split}
  & (i_1\times\id_{X_3})^* (K_{12} * K_{23})\\
  =\; &
  p_{13*}''' \left(p_{12}^{\prime\prime\prime*}
    (i_1\times\id_{X_2})^* K_{12}
   \otimes^L_{S_1\times S_2\times X_3}
   (\id_{S_1}\times \id_{S_2}\times i_3)_* p_{23}^{\prime*} 
  (i_2\times\id_{X_3})^* K_{23}\right) \\
  =\; &
  p_{13*}''' 
  (\id_{S_1}\times \id_{S_2}\times i_3)_* \left(
    (\id_{S_1}\times \id_{S_2}\times i_3)^*
     p_{12}^{\prime\prime\prime*} (i_1\times\id_{X_2})^* K_{12}
   \otimes^L_{S_1\times S_2\times S_3} p_{23}^{\prime*} 
  (i_2\times\id_{X_3})^* K_{23}\right) \\
  =\; &
  (\id_{S_1}\times i_3)_* p_{13*}'
 \left(  p_{12}^{\prime*} (i_1\times\id_{X_2})^* K_{12}
   \otimes^L_{S_1\times S_2\times S_3} p_{23}^{\prime*} 
  (i_2\times\id_{X_3})^* K_{23}\right)
\end{split}
\end{equation*}
where we have used \eqref{eq:projection} in the second equality
and $p_{13}'''\circ(\id_{S_1}\times \id_{S_2}\times i_3) =
(\id_{S_1}\times i_3)\circ p_{13}'$,
$p_{12}'''\circ (\id_{S_1}\times \id_{S_2}\times i_3) = p_{12}'$
in the third equality. Finally, by \eqref{eq:SZass'}, the
homomorphism $(\id_{S_1}\times i_3)_*$ is nothing but the identity
operator
$K(Z_{12}'\circ Z_{23}') \to K(Z_{12}'\circ Z_{23}')$. Thus we have
the assertion.
\end{proof}

\subsection{Convolution and principal bundles}
Let $G$ be a linear algebraic group and suppose that we have
principal $G$-bundles $\pi_a\colon P_a\to X_a$ over $X_a$ for
$a=1, 2, 3$.
Consider the restriction of the principal $G$-bundle 
$\pi_a\times\id_{X_b}\colon P_a\times X_b\to X_a\times X_b$ to
$Z_{ab}$ for $(a,b) = (1,2), (2,3)$.
Then the pullback homomorphism gives a canonical isomorphism
\begin{equation}\label{eq:Zpull}
  (\pi_a\times\id_{X_b})^* \colon K(Z_{ab}) \xrightarrow{\cong}
  K^G((\pi_a\times\id_{X_b})^{-1}Z_{ab}).
\end{equation}
Similarly we have an isomorphism
\begin{equation}\label{eq:Zpull'}
  (\pi_1\times\id_{X_3})^* \colon K(Z_{12}\circ Z_{23}) \xrightarrow{\cong}
  K^G((\pi_1\times\id_{X_3})^{-1}(Z_{12}\circ Z_{23})).
\end{equation}

Let $G$ act on $P_a\times P_b$ diagonally.
We assume that there exists a closed $G$-invariant subvariety $Z_{ab}'$ of
$P_a\times P_b$ such that 
\begin{equation}
\label{eq:PZass}
\begin{split}
  & \text{the restriction of
$\id_{P_a}\times\pi_b\colon P_a\times P_b\to P_a\times X_b$ to
$Z_{ab}'$ is proper, and} \\
  & (\id_{P_a}\times\pi_b)(Z_{ab}') = (\pi_a\times\id_{X_b})^{-1}Z_{ab}
\end{split}
\end{equation}
for $(a,b) = (1,2), (2,3)$. 

Let $p_{ab}'$ be the projection $P_1\times P_2\times P_3\to P_a\times
P_b$.
Since the restriction of the projection 
$p_{13}'\colon p_{12}^{\prime-1}(Z_{12}')\cap
p_{23}^{\prime-1}(Z_{23}')\to P_1\times P_3$ is proper,
we have the convolution product
$*'\colon K^G(Z_{12}')\otimes K^G(Z_{23}') \to K^G(Z_{12}'\circ
Z_{23}')$ by
\begin{equation*}
  K'_{12}*' K'_{23}
  \defeq p'_{13*} \left(
  p_{12}^{\prime *} K'_{12}\otimes^L_{P_1\times P_2\times P_3}
  p_{23}^{\prime *} K'_{23}\right)
  \qquad\text{for $K'_{12}\in K^G(Z_{12}')$, $K'_{23}\in K^G(Z_{23}')$}.
\end{equation*}
We want to compare this convolution product with that on
$K(Z_{12})\otimes K(Z_{23})$.

By \eqref{eq:PZass}, we have $(\id_{P_a}\times\pi_b)_* K_{ab}'
\in K^G((\id_{P_a}\times\pi_b)(Z_{ab}'))
 = K^G((\pi_a\times\id_{X_b})^{-1}Z_{ab})$.
Via \eqref{eq:Zpull}, we define
\begin{equation}\label{eq:KabDef}
   K_{ab} \defeq \left((\pi_a\times\id_{X_b})^*\right)^{-1}
     (\id_{P_a}\times\pi_b)_* K_{ab}' \in K(Z_{ab}).
\end{equation}
Thus we can consider the convolution product $K_{12} * K_{23}$.

By the construction, we have 
\(
(\id_{P_1}\times\pi_3)(Z_{12}'\circ Z_{23}') = 
(\pi_1\times\id_{X_3})^{-1}(Z_{12}\circ Z_{23}).
\)
Combining with \eqref{eq:Zpull'}, we have
\begin{equation*}
   \left((\pi_1\times\id_{X_3})^*\right)^{-1}(\id_{P_1}\times\pi_3)_*
   (K_{12}' *' K_{23}') \in
   K(Z_{12}\circ Z_{23}).
\end{equation*}

\begin{Proposition}\label{prop:Pconv}
In the above setup, we have
\begin{equation}\label{eq:Pconv}
  (\pi_1\times\id_{X_3})^* (K_{12} * K_{23})
  = (\id_{P_1}\times\pi_3)_* ( K_{12}' *' K_{23}' ).
\end{equation}
Namely the following diagram commutes:
\begin{equation*}
\begin{CD}
  K^G(Z'_{12})\otimes K^G(Z'_{23}) @>{*'}>> K^G(Z'_{12}\circ Z'_{23}) \\
  @VVV
  @VVV \\
  K(Z_{12})\otimes K(Z_{23}) @>{*}>> K(Z_{12}\circ Z_{23}),
\end{CD}
\end{equation*}
where the left vertical arrow is
\begin{equation*}
\left((\pi_1\times\id_{X_2})^*\right)^{-1}
     (\id_{P_1}\times\pi_2)_* \otimes
 \left((\pi_2\times\id_{X_3})^*\right)^{-1}
     (\id_{P_2}\times\pi_3)_* ,
\end{equation*}
and the right vertical arrow is
\begin{equation*}
 \left((\pi_1\times\id_{X_3})^*\right)^{-1}
     (\id_{P_1}\times\pi_3)_*    .
\end{equation*}
\end{Proposition}

\begin{Example}
Suppose $X_1 = X_2$, $P_1 = P_2$ and $Z_{12} =
\operatorname{Image}\Delta_{X_1}$, where $\Delta_{X_1}$ is the
diagonal embedding $X_1\to X_1\times X_2$. If we take
$Z_{12}' = \operatorname{Image}\Delta_{P_1}$, where
$\Delta_{P_1}$ is the diagonal embedding
$P_1\to P_1\times P_2$, the assumption~(\ref{eq:PZass}) is
satisfied. In fact, the restriction of $\id_{P_1}\times\pi_2$
to $Z_{12}'$ is an isomorphism.
Take a vector bundle $E$ and consider $K_{12} = \Delta_{X_1*}[E]$.
By the isomorphism
$K(X_1) \xrightarrow[\cong]{\pi_1^*} K^G(P_1)$, we can define
$K_{12}' = \Delta_{P_1*}\pi_1^* [E]$. Then both
$(\id_{P_1}\times\pi_2)_* K_{12}'$ and $(\pi_1\times\id_{X_2})^*
K_{12}$
is $\Delta'_* [E]$ where $\Delta'\colon P_1\to 
(\id_{P_1}\times\pi_2)\Delta_{P_1} = (\pi_1\times\id_{X_2})^{-1}\Delta_{X_1}$
is the natural isomorphism. Hence \eqref{eq:KabDef} holds for $K_{12}$ 
and $K_{12}'$.
By \lemref{lem:diag}, we have
\begin{equation*}
\begin{split}
    & (\Delta_{X_1})_*[E] * K_{23} = p_2^* [E] \otimes K_{23},\\
  & (\Delta_{P_1})_*[\pi_1^* E] *' K_{23}'
  = p_2^{\prime*}\pi_2^* [E] \otimes K'_{23},
\end{split}
\end{equation*}
where $p_2\colon X_2\times X_3\to X_2$,
$p_2'\colon P_2\times X_3\to P_2$
are the projections. We can directly check \eqref{eq:Pconv} in this case.
\end{Example}

\begin{proof}[Proof of \propref{prop:Pconv}]
As in the proof of \propref{prop:Sconv}, we replace $X_a$ by $P_a$
factor by factor.

{\bf Step 1}. First we replace $X_1$ by $P_1$. Consider the following
fiber square:
\begin{equation*}
\begin{CD}
  P_1\times X_2\times X_3 @>{\pi_1\times\id_{X_2}\times\id_{X_3}}>>
  X_1\times X_2\times X_3 \\
  @V{p_{13}''}VV @VV{p_{13}}V \\
  P_1\times X_3 @>>{\pi_1\times\id_{X_3}}> X_1\times X_3,
\end{CD}
\end{equation*}
where $p_{13}''$ is the projection. By base change
\cite[5.3.15]{Gi-book} and \eqref{eq:pull_product}, we have
\begin{equation}\label{eq:Pconv1}
\begin{split}
  & (\pi_1\times\id_{X_3})^* (K_{12} * K_{23})
   = (\pi_1\times\id_{X_3})^* p_{13*}
\left( p_{12}^* K_{12}\otimes^L_{X_1\times X_2\times X_3} p_{23}^*
  K_{23}\right) \\
 = \; & p_{13*}'' (\pi_1\times\id_{X_2}\times\id_{X_3})^*
 \left( p_{12}^* K_{12}\otimes^L_{X_1\times X_2\times X_3} p_{23}^*
  K_{23}\right) \\
 =\; & p_{13*}'' \left( 
   (\pi_1\times\id_{X_2}\times\id_{X_3})^* p_{12}^* K_{12}
      \otimes^L_{P_1\times X_2\times X_3}
   (\pi_1\times\id_{X_2}\times\id_{X_3})^* p_{23}^* K_{23} \right) \\
 =\; & p_{13*}'' \left( 
   p_{12}^{\prime\prime*}(\pi_1\times\id_{X_2})^* K_{12}
      \otimes^L_{P_1\times X_2\times X_3}
   p_{23}^{\prime\prime*} K_{23} \right),
\end{split}
\end{equation}
where $p_{12}''\colon P_1\times X_2\times X_3\to P_1\times X_2$ and
$p_{23}''\colon P_1\times X_2\times X_3\to X_2\times X_3$ are
projections.

{\bf Step 2}.
Consider the fiber square
\begin{equation*}
\begin{CD}
  P_1\times P_2\times X_3 @>{p_{12}'''}>> P_1\times P_2 \\
  @V{\id_{P_1}\times\pi_2\times\id_{X_3}}VV 
  @VV{\id_{P_1}\times\pi_2}V \\
  P_1\times X_2\times X_3 @>>{p_{12}''}> P_1\times X_2,
\end{CD}
\end{equation*}
where $p_{12}'''$ is the projection. By base change
\cite[5.3.15]{Gi-book}, we have
\begin{equation*}
  (\id_{P_1}\times\pi_2\times\id_{X_3})_*
  p_{12}^{\prime\prime\prime*} K_{12}'
  = p_{12}^{\prime\prime*} (\id_{P_1}\times \pi_2)_* K_{12}'
  = p_{12}^{\prime\prime*}(\pi_1\times\id_{X_2})^* K_{12},
\end{equation*}
where we have used \eqref{eq:KabDef} for $(a,b) = (1,2)$.
Substituting this into \eqref{eq:Pconv1}, we get
\begin{equation}\label{eq:Pconv2}
\begin{split}
  & (\pi_1\times\id_{X_3})^* (K_{12} * K_{23}) \\
  =\; & p_{13*}''\left(
    (\id_{P_1}\times\pi_2\times\id_{X_3})_*
     p_{12}^{\prime\prime\prime*} K_{12}'
     \otimes^L_{P_1\times X_2\times X_3}
     p_{23}^{\prime\prime*} K_{23}\right) \\
  =\; & p_{13*}'' (\id_{P_1}\times\pi_2\times\id_{X_3})_*
  \left( p_{12}^{\prime\prime\prime*} K_{12}'
     \otimes^L_{P_1\times P_2\times X_3}
     (\id_{P_1}\times\pi_2\times\id_{X_3})^*
      p_{23}^{\prime\prime*}K_{23}\right)
\end{split}
\end{equation}
where we have used \eqref{eq:projection} in the second equality.
Let $p_{23}'''\colon P_1\times P_2\times X_3\to P_2\times X_3$ be the
projection.
By $p_{23}''\circ (\id_{P_1}\times\pi_2\times\id_{X_3})
= (\pi_2\times\id_{X_3})\circ p_{23}'''$, we have
\begin{equation*}
(\id_{P_1}\times\pi_2\times\id_{X_3})^*
      p_{23}^{\prime\prime*} = 
  p_{23}^{\prime\prime\prime*} (\pi_2\times\id_{X_3})^*.
\end{equation*}
We also have 
\begin{equation*}
p_{13*}'' (\id_{P_1}\times\pi_2\times\id_{X_3})_* = p_{13*}''',
\end{equation*}
where $p_{13}\colon P_1\times P_2\times X_3
\to P_1\times X_3$ is the projection. Substituting these two
equalities into \eqref{eq:Pconv2}, we get
\begin{equation}\label{eq:Pconv3}
  (\pi_1\times\id_{X_3})^* (K_{12} * K_{23})
  = p_{13*}'''\left( p_{12}^{\prime\prime\prime*} K_{12}'
            \otimes^L_{P_1\times P_2\times X_3}
            p_{23}^{\prime\prime\prime*}
            (\pi_2\times\id_{X_3})^* K_{23}\right).
\end{equation}

{\bf Step 3}. Consider the fiber square
\begin{equation*}
\begin{CD}
  P_1\times P_2\times P_3 @>{p_{23}'}>> P_2\times P_3\\
  @V{\id_{P_1}\times\id_{P_2}\times\pi_3}VV
  @VV{\id_{P_2}\times\pi_3}V \\
  P_1\times P_2\times X_3 @>{p_{23}'''}>> P_2\times X_3.
\end{CD}
\end{equation*}
By base change \cite[5.3.15]{Gi-book}, we have
\begin{equation*}
   (\id_{P_1}\times\id_{P_2}\times\pi_3)_* p_{23}^{\prime *}
   K'_{23} = p_{23}^{\prime\prime\prime*}
   (\id_{P_2}\times\pi_3)_* K'_{23}
   = p_{23}^{\prime\prime\prime*} (\pi_2\times\id_{X_3})^* K_{23},
\end{equation*}
where we have used \eqref{eq:KabDef} for $(a,b) = (2,3)$ in the
second equality.
Substituting this into \eqref{eq:Pconv3}, we have
\begin{equation*}
\begin{split}
  & (\pi_1\times\id_{X_3})^* (K_{12} * K_{23}) \\
  =\; & p_{13*}'''\left( p_{12}^{\prime\prime\prime*} K_{12}'
            \otimes^L_{P_1\times P_2\times X_3}
            (\id_{P_1}\times\id_{P_2}\times\pi_3)_* p_{23}^{\prime *}K_{23}'
            \right) \\
  =\; & p_{13*}''' (\id_{P_1}\times\id_{P_2}\times\pi_3)_*
  \left( (\id_{P_1}\times\id_{P_2}\times\pi_3)^*
            p_{12}^{\prime\prime\prime*} K_{12}'
            \otimes^L_{P_1\times P_2\times P_3}
            p_{23}^{\prime *} K_{23}'\right)
\end{split}
\end{equation*}
By $p_{13}'''\circ (\id_{P_1}\times\id_{P_2}\times\pi_3)
= (\id_{P_1}\times\pi_3)\circ p_{13}'$ and
$p_{12}'''\circ (\id_{P_1}\times\id_{P_2}\times\pi_3) = p_{12}'$, we get
\begin{equation*}
  (\pi_1\times\id_{X_3})^* (K_{12} * K_{23})
  = (\id_{P_1}\times\pi_3)_*\circ p_{13*}'\left(
    p_{12}^{\prime *} K_{12}'  \otimes^L_{P_1\times P_2\times P_3}
            p_{23}^{\prime *} K_{23}'\right).
\end{equation*}
This proves our assertion.
\end{proof}

\section{A homomorphism $\Ul\to
K^{G_\bw\times\C^*}(\Zw)\otimes_{\Z[q,q^{-1}]}\Q(q)$}\label{sec:hom}

\subsection{}\label{subsec:homdef}
Let us define an analogue of the Steinberg variety
$Z(\bv^1,\bv^2;\bw)$ by
\begin{equation}\label{eq:def_Steinberg}
   \{ (x^1, x^2)\in \M(\bv^1,\bw)\times\M(\bv^2,\bw) \mid
 \pi(x^1) = \pi(x^2)
\}.
\end{equation}
Here $\pi(x^1) = \pi(x^2)$ means that $\pi(x^1)$ is equal to
$\pi(x^2)$ if we regard both as elements of $\M_0(\infty,\bw)$ by
\eqref{eq:union}. This is a closed subvariety of
$\M(\bv^1,\bw)\times\M(\bv^2,\bw)$.

The map
$p_{13}\colon p_{12}^{-1}(Z(\bv^1,\bv^2;\bw))\cap
p_{23}^{-1}(Z(\bv^2,\bv^3;\bw)) \to \M(\bv^1,\bw)\times\M(\bv^3,\bw)$
is proper and its image is contained in $Z(\bv^1,\bv^3;\bw)$.
Hence we can define the
convolution product on the equivariant $K$-theory:
\begin{equation*}
  K^{G_\bw\times\C^*}(Z(\bv^1,\bv^2;\bw))\otimes
  K^{G_\bw\times\C^*}(Z(\bv^2,\bv^3;\bw))
  \to K^{G_\bw\times\C^*}(Z(\bv^1,\bv^3;\bw)).
\end{equation*}

Let $\prod'_{\bv^1,\bv^2} K^{G_\bw\times\C^*}(Z(\bv^1,\bv^2;\bw))$ be the
subspace of $\prod_{\bv^1,\bv^2} K^{G_\bw\times\C^*}(Z(\bv^1,\bv^2;\bw))$ 
consisting elements $(F_{\bv^1,\bv^2})$ such that
\begin{enumerate}
\item for fixed $\bv^1$, $F_{\bv^1,\bv^2} = 0$ for all but finitely
many choices of $\bv^2$,
\item for fixed $\bv^2$, $F_{\bv^1,\bv^2} = 0$  for all but finitely
many choices of $\bv^1$.
\end{enumerate}
The convolution product $*$ is well-defined on
$\prod'_{\bv^1,\bv^2} K^{G_\bw\times\C^*}(Z(\bv^1,\bv^2;\bw))$. When the
underlying graph is of type $ADE$, $\M(\bv,\bw)$ is empty for all but 
finitely many choices of $\bv$, so $\prod'$ is nothing but the direct
product $\prod$.

Let $\Zw$ denote the disjoint union
$\bigsqcup_{\bv^1,\bv^2} Z(\bv^1,\bv^2;\bw)$. When we write
$K^{G_\bw\times\C^*}(\Zw)$, we mean
$\prod'_{\bv^1,\bv^2} K^{G_\bw\times\C^*}(Z(\bv^1,\bv^2;\bw))$
as convention.

The second projection $G_\bw\times\C^*\to \C^*$ induces a homomorphism
$R(\C^*)\to R(G_\bw\times\C^*)$. Thus $R(G_\bw\times\C^*)$ is an
$R(\C^*)$-algebra. Moreover, $R(\C^*)$ is isomorphic to
$\Z[q,q^{-1}]$ where $q^m$ corresponds to $L(m)$ in \eqref{eq:L(m)}.
Thus $K^{G_\bw\times\C^*}(\Zw)$ is a
$\Z[q,q^{-1}]$-algebra.

The aim of this section and the next two sections is to define the
homomorphism from $\Ul$ into
$K^{G_\bw\times\C^*}(\Zw)\otimes_{\Z[q,q^{-1}]}\Q(q)$.
We first define the map on generators of $\Ul$, and then check the
defining relation.

\subsection{}
First we want to define the image of $q^h$, $h_{k,m}$.

Let $C_k^\bullet(\bv,\bw)$ be the $G_\bw\times\C^*$-equivariant complex 
over $\M(\bv,\bw)$ defined in \eqref{eq:taut_cpx}.
We consider $C_k^\bullet(\bv,\bw)$ as an element of
$K^{G_\bw\times\C^*}(\M(\bv,\bw))$ by identifying it with
the alternating sum
\begin{equation*}
  - [q^{-2} V_k]
  + \left[q^{-1} \left( \displaystyle{\bigoplus_{l:k\neq l}}
     [-\langle h_k,\alpha_l\rangle]_q V_l
    \oplus W_k\right) \right]
  - [V_k].
\end{equation*}

The rank of the complex~\eqref{eq:taut_cpx},
as an element of $K^{G_\bw\times\C^*}(\M(\bv,\bw))$ (see
\subsecref{subsec:opt_vec}),
is given by
\begin{equation*}
   \rank C^\bullet_k(\bv,\bw) = 
   - \sum_{l:k\neq l} (\alpha_k,\alpha_l) \dim V_l + \dim W_k - 2\dim V_k
   = \langle h_k, \bw - \bv\rangle.
\end{equation*}

Let $\Delta$ denote the diagonal embedding $\M(\bv,\bw)\to
\M(\bv,\bw)\times\M(\bv,\bw)$.

\begin{equation}\label{eq:defH}
\begin{gathered}
   q^{h} \longmapsto \sum_\bv q^{\langle h, \bw-\bv\rangle}
     \Delta_* \shfO_{\M(\bv,\bw)}, \\
   p_k^+(z)
     \longmapsto
     \sum_{\bv} \Delta_*
      \left(\Wedge_{-1/z} (C_k^\bullet(\bv,\bw))\right)^+, \\
   p_k^-(z)
     \longmapsto
     \sum_{\bv} \Delta_* 
    \left((-z)^{\rank C_k^\bullet(\bv,\bw)}\det C_k^\bullet(\bv,\bw)^*
       \Wedge_{-1/z} (C_k^\bullet(\bv,\bw))\right)^-,
\end{gathered}
\end{equation}
where $(\ )^\pm$ denotes the expansion at $z = \infty$, $0$ respectively.
Note that
\begin{equation*}
  \psi_k^\pm(z) = q^{\pm h_k} \frac{p_k^\pm(q z)}{p_k^\pm(q^{-1}z)}
  \longmapsto
  \sum_{\bv} q^{\rank C_k^\bullet(\bv,\bw)}\Delta_*
      \left(\frac{\Wedge_{-1/qz} (C_k^\bullet(\bv,\bw))}
                 {\Wedge_{-q/z} (C_k^\bullet(\bv,\bw))}\right)^\pm.
\end{equation*}

\subsection{}

Next we define the images of $e_{k,r}$ and $f_{k,r}$. They are given
by line bundles over Hecke correspondences.

Let $\bv^1$, $\bv^2$ and $\Pa_k(\bv^2,\bw)$ as in
\subsecref{subsec:hecke}.
By the definition, the quotient $V_k^2/V_k^1$ defines a line bundle
over $\Pa_k(\bv^2,\bw)$.
The generator $e_{k,r}$ is very roughly defined as the $r$th power of
$V_k^2/V_k^1$, but we need a certain modification in order to have
the correct commutation relation.

For the modification, we need to consider the
following variants of $C^\bullet_k(\bv,\bw)$:
\begin{equation}\label{eq:mod_taut_cpx}
  C^{\prime\bullet}_k(\bv,\bw) \defeq \Coker\alpha ,
  \qquad
  C^{\prime\prime\bullet}_k(\bv,\bw) \defeq V_k[-1],
\end{equation}
where $V_k[-1]$ means that we consider the complex consisting $V_k$ in 
degree $1$ and $0$ for other degrees.
Since $\alpha$ is injective, we have
\(
  C^\bullet_k(\bv,\bw) = C^{\prime\bullet}_k(\bv,\bw) +
  C^{\prime\prime\bullet}_k(\bv,\bw)
\)
in $K^{G_\bw\times\C^*}(\M(\bv,\bw))$. We have the corresponding
decomposition of the Cartan matrix $\bC$:
\begin{equation*}
  \bC = \bC' + \bC'',
  \qquad\text{where
  $\bC' \defeq \bC - \bI$, $\bC'' \defeq \bI$.}
\end{equation*}
We identify $\bC'$ (resp.\ $\bC''$) with a map given by
\begin{equation*}
   \bv = \sum v_k \alpha_k \in \bigoplus \Z\alpha_k
   \mapsto \sum v_k (\alpha_k - \Lambda_k) \quad
   \left(\text{resp.\ }\sum v_k \Lambda_k\right) \in P.
\end{equation*}

We also need matrices $\bC_\Omega$, $\bC_{\overline\Omega}$ given by
\begin{equation*}
   \bC_\Omega \defeq \bI - \bA_\Omega, \qquad
   \bC_{\overline\Omega} \defeq \bI - \bA_{\overline\Omega},
\end{equation*}
where $\bA_\Omega$, $\bA_{\overline\Omega}$ are as in
\eqref{eq:Aomega}. We also identify them with maps
$\bigoplus \Z\alpha_k \to P$ exactly as above.

Let $\omega\colon \M(\bv^1,\bw)\times\M(\bv^2,\bw)
\to \M(\bv^2,\bw)\times\M(\bv^1,\bw)$ be the exchange of factors.
Let us denote
$\omega\left(\Pa_k(\bv^2+\alpha_k,\bw)\right)
\subset \M(\bv^2+\alpha_k,\bw)\times\M(\bv^2,\bw)$ by
$\Pa_k^-(\bv^2,\bw)$.
As on $\Pa_k(\bv,\bw)$, we have a natural line bundle over
$\Pa_k^-(\bv^2,\bw)$. Let us denote it by
$V^1_k/V^2_k$.

Now we define the images of $e_{k,r}$, $f_{k,r}$ by
\begin{equation}
\label{eq:defEF}
\begin{split}
  & e_{k,r} \longmapsto 
  \sum_{\bv^2} (-1)^{-\langle h_k, \bC_{\overline\Omega}\bv^2\rangle}
  \left[ i_* \left(q^{-1}V_k^2/V_k^1\right)^{\otimes r-\langle h_k,
  \bC''\bv^2\rangle}
  \otimes
  \det C_k^{\prime\bullet}(\bv^2,\bw)^*\right], \\
  & f_{k,r} \longmapsto 
  \sum_{\bv^2} (-1)^{\langle h_k, \bw - \bC_\Omega\bv^2\rangle}
  \left[ i^-_* \left(q^{-1}V_k^1/V_k^2\right)^{\otimes r+\langle h_k,
    \bw-\bC'\bv^2\rangle}
  \otimes
  \det C_k^{\prime\prime\bullet}(\bv^2,\bw)^*\right],
\end{split}
\end{equation}
where $i\colon \Pa_k(\bv^2,\bw)\to Z(\bv^2-\alpha_k,\bv^2;\bw)$
and $i^-\colon \Pa^-_k(\bv^2,\bw)\to Z(\bv^2+\alpha_k,\bv^2;\bw)$
are the inclusions.
Hereafter, we may omit $i_*$ or $i^-_*$, hoping that it makes no confusion.

\subsection{}

\begin{Theorem}\label{thm:hommain}
The assignments \eqref{eq:defH},\eqref{eq:defEF} define a homomorphism 
\[
   \Ul \to 
   K^{G_\bw\times\C^*}(\Zw)\otimes_{\Z[q,q^{-1}]}\Q(q)
\]
of $\Q(q)$-algebras. 
\end{Theorem}

We need to check the defining relations
\eqref{eq:relCcent}$\sim$\eqref{eq:relDS}.
We do not need to consider the
relations~\eqref{eq:relCcent},\eqref{eq:relD} because we are
considering $\Ul$ instead of $\Ua$.
The relations~\eqref{eq:relHH}, \eqref{eq:relHH2} and \eqref{eq:relHH3}
follow from \lemref{lem:diag} and the fact
that $E \otimes F\cong F\otimes E$.
The remaining relations will be checked in the next two sections.

\section{Relations~(I)}\label{sec:rel1}

\subsection{Relation (\ref{eq:relHE})}

Fix a vertex $k\in I$ and take $\bv^1$, $\bv^2 = \bv^1+\alpha_k$.  Let
$i$ be the inclusion $\Pa_k(\bv^2,\bw)\to
Z(\bv^1,\bv^2;\bw)$ and let $p_1$ and $p_2$ be the projection
$\Pa_k(\bv^2,\bw)\to \M(\bv^1,\bw)$ and
$\Pa_k(\bv^2,\bw)\to \M(\bv^2,\bw)$ respectively.

By \lemref{lem:diag}, we have
\begin{equation}\label{eq:GrelHE}
\begin{split}
  & \Delta_* \Wedge_{-1/z}(C_l^\bullet(\bv^1,\bw))\ast
  i_{*} \left[\sum_{r=-\infty}^\infty
   \left[(q^{-1} V^2_k/V^1_k\right)^{\otimes r} w^{-r} \right]
  \ast \Delta_* \left(\Wedge_{-1/z} C_l^\bullet(\bv^2,\bw)\right)^{-1}\\ 
  =\; & i_* \left[\Wedge_{-1/z} p_1^*  (C_l^\bullet(\bv^1,\bw))\otimes
   \left(\Wedge_{-1/z} p_2^*  (C_l^\bullet(\bv^2,\bw))\right)^{-1}\otimes
   \sum_{r=-\infty}^\infty
           \left(q^{-1}V^2_k/V^1_k\right)^{\otimes r} w^{-r}\right].
\end{split}
\end{equation}
We have the following equality in $K^{G_\bw\times\C^*}(\Pa_k(\bv^2,\bw))$:
\begin{equation*}
   V_k^1 = V_k^2 - V_k^2/V_k^1.
\end{equation*}
Hence we have
\begin{equation*}
   \Wedge_{-1/z} p_1^* C_l^\bullet(\bv^1,\bw)\otimes
    \left(\Wedge_{-1/z} p_2^* C_l^\bullet(\bv^2,\bw)\right)^{-1}
   = 
      \Wedge_{-1/z} [\langle h_k,\alpha_l\rangle]_q
        \left(q^{-1}V^2_k/V^1_k\right).
\end{equation*}
Substituting this into \eqref{eq:GrelHE}, we get
\begin{equation*}
\begin{split}
 & \Delta_* \Wedge_{-1/z} (C_l^\bullet(\bv^1,\bw))\ast
  i_{*} \left[\sum_{r=-\infty}^\infty
   \left(q^{-1} V^2_k/V^1_k\right)^{\otimes r} w^{-r} \right]
  \ast \Delta_* \left(\Wedge_{-1/z} C_l^\bullet(\bv^2,\bw)\right)^{-1} \\
   =\; & 
   \begin{cases}
   {\displaystyle \left(1 - \frac{wq}z\right)\left(1 - \frac{w}{zq}\right)
   i_* \left[\sum_{r=-\infty}^\infty
     \left(q^{-1}V^2_k/V^1_k\right)^{\otimes r} w^{-r}\right]}
    & \text{if $k = l$,} \\
    {\displaystyle
    \prod_{p=0}^{-\langle h_k,\alpha_l\rangle-1}
    \left(1 - \frac{q^{\langle h_k,\alpha_l\rangle+1+2p}w}z\right)^{-1}
    i_* \left[\sum_{r=-\infty}^\infty
       \left(q^{-1}V^2_k/V^1_k\right)^{\otimes r} w^r\right]}
    & \text{otherwise.}
   \end{cases} \\
\end{split}
\end{equation*}
This is equivalent to \eqref{eq:relHE} for $x_l^+(w)$.
The relation \eqref{eq:relHE} for $x_l^-(w)$ can be proved in the same
way.

\subsection{Relation (\ref{eq:relEF}) for $k\neq l$}

Fix two vertices $k\neq l$. Let $\bv^1$, $\bv^2$, $\bv^3$, $\bv^4$ be
dimension vectors such that
\begin{equation*}
 \bv^2 = \bv^1 + \alpha_k = \bv^3 + \alpha_l, \quad
 \bv^4 = \bv^1 - \alpha_l = \bv^3 - \alpha_k.
\end{equation*}
We want to compute $e_{k,r} * f_{l,s}$ and $f_{l,s} * e_{k,r}$ in the
component $K^{G_\bw\times\C}(Z(\bv^1,\bv^3,\bw))$.

Let us consider the intersection
\begin{equation*}
  p_{12}^{-1}\Pa_k(\bv^2,\bw) \cap p_{23}^{-1}\Pa^-_{l}(\bv^3,\bw)
  \qquad\text{(resp.\ 
  $p_{12}^{-1}\Pa^-_{l}(\bv^2,\bw) \cap p_{23}^{-1}\Pa_{k}(\bv^3,\bw)$)}
\end{equation*}
in $\M(\bv^1,\bw)\times\M(\bv^2,\bw)\times\M(\bv^3,\bw)$
(resp.\ $\M(\bv^1,\bw)\times\M(\bv^4,\bw)\times\M(\bv^3,\bw)$).
On the intersection, we have the inclusion of restrictions of
tautological bundles
\begin{equation*}
   V^1 \subset V^2 \supset V^3
   \qquad\text{(resp.\
   $V^1 \supset V^4 \subset V^3$)}.
\end{equation*}

\begin{Lemma}\label{lem:EFint}
The above two intersections are transversal,
and there is a $G_\bw\times\C^*$-equivariant isomorphism between them
such that
\begin{aenume}
\item it is the identity operators on the factor
$\M(\bv^1,\bw)$ and $\M(\bv^3,\bw)$,
\item it induces isomorphisms $V^2_k/V^1_k\cong V^3_k/V^4_k$ and
$V^2_l/V^3_l\cong V^1_l/V^4_l$.
\end{aenume}
\end{Lemma}

\begin{proof}
See \cite[Lemmas~9.8, 9.9, 9.10 and their proofs]{Na-alg}.
\end{proof}

Since the intersection is transversal, we have
\begin{equation}\label{eq:EF}
  e_{k,r} * f_{l,s} 
  = (-1)^{-\langle h_k, \bC_{\overline\Omega}\bv^2\rangle
          +\langle h_l, \bw-\bC_\Omega\bv^3\rangle}
  p_{13*}[\mathcal L]
\end{equation}
where $\mathcal L$ is the following line bundle over
$p_{12}^{-1}\Pa_k(\bv^2,\bw) \cap p_{23}^{-1}\Pa^-_{l}(\bv^3,\bw)$:
\begin{equation*}
  \left(q^{-1}V_k^2/V^1_k\right)^{\otimes r-\langle h_k, \bC''\bv^2\rangle}
  \otimes
  \left(q^{-1}V_l^2/V^3_l\right)^{\otimes s+\langle h_l, \bw-\bC'\bv^3\rangle}
  \otimes \det C_k^{\prime\bullet}(\bv^2,\bw)^*
  \otimes \det C_l^{\prime\prime\bullet}(\bv^3,\bw)^*.
\end{equation*}
Similarly, we have
\begin{equation}\label{eq:FE}
  f_{l,s} * e_{k,r} 
  = (-1)^{\langle h_l, \bw-\bC_\Omega\bv^4\rangle
         -\langle h_k, \bC_{\overline\Omega}\bv^3\rangle}
  p_{13*}[\mathcal L']
\end{equation}
where $\mathcal L'$ is the following line bundle over
$p_{12}^{-1}\Pa^-_{l}(\bv^2,\bw) \cap p_{23}^{-1}\Pa_{k}(\bv^3,\bw)$:
\begin{equation*}
  \left(q^{-1}V_l^1/V^4_l\right)^{\otimes s+\langle h_l, \bw-\bC'\bv^4\rangle}
  \otimes
  \left(q^{-1}V_k^3/V^4_k\right)^{\otimes r-\langle h_k, \bC''\bv^3\rangle}
  \otimes \det C_l^{\prime\prime\bullet}(\bv^4,\bw)^*
  \otimes \det C_k^{\prime\bullet}(\bv^3,\bw)^*.
\end{equation*}

Let us compare \eqref{eq:EF} and \eqref{eq:FE}.
On $p_{12}^{-1}\Pa_k(\bv^2,\bw) \cap p_{23}^{-1}\Pa^-_{l}(\bv^3,\bw)$,
we have
\begin{equation*}
\begin{split}
  \det C_k^{\prime\bullet}(\bv^2,\bw)
  & = \det C_k^{\prime\bullet}(\bv^3,\bw)
  \otimes\det\left([-(\alpha_k,\alpha_l)]_q (q^{-1}V^2_l/V^3_l)\right) \\
  & = \det C_k^{\prime\bullet}(\bv^3,\bw)
  \otimes \left(q^{-1}V^2_l/V^3_l\right)^{\otimes -(\alpha_k,\alpha_l)}.
\end{split}
\end{equation*}
On the other hand, we have
\begin{equation*}
  \det C_l^{\prime\prime\bullet}(\bv^3,\bw)
  = \det C_l^{\prime\prime\bullet}(\bv^4,\bw)
\end{equation*}
on $p_{12}^{-1}\Pa^-_l(\bv^2,\bw) \cap p_{23}^{-1}\Pa_{k}(\bv^3,\bw)$.
Hence under the isomorphism in \lemref{lem:EFint}, we obtain
\begin{equation*}
  \mathcal L \cong \mathcal L'.
\end{equation*}
where we have used
\begin{equation}\label{eq:EF_FE_data}
 \langle h_l, \bC'\bv^3\rangle
= \langle h_l, \bC'\bv^4\rangle - a_{lk},
\quad
 \langle h_k, \bC''\bv^2\rangle 
= \langle h_k, \bC''\bv^3\rangle.
\end{equation}
By
\begin{equation*}
 \langle h_l, \bC_\Omega\bv^3\rangle
= \langle h_l, \bC_\Omega\bv^4\rangle - (\bA_\Omega)_{lk},
\quad
 \langle h_k, \bC_{\overline\Omega}\bv^2\rangle
= \langle h_k, \bC_{\overline\Omega}\bv^3\rangle
   - (\bA_{\overline\Omega})_{kl},
\quad
 (\bA_\Omega)_{lk} = (\bA_{\overline\Omega})_{kl},
\end{equation*}
we have
\begin{equation*}
   (-1)^{-\langle h_k, \bC_{\overline\Omega}\bv^2\rangle
          +\langle h_l, \bw-\bC_\Omega\bv^3\rangle}
   = 
   (-1)^{\langle h_l, \bw-\bC_\Omega\bv^4\rangle
         -\langle h_k, \bC_{\overline\Omega}\bv^3\rangle}
\end{equation*}
Thus we have $[e_{k,r}, f_{l,s}] = 0$.

\subsection{Relation (\ref{eq:relexE1})}
We give the proof of \eqref{eq:relexE1} for $\pm = +$, in this
subsection. The relation \eqref{eq:relexE1} for $\pm = -$ can be
proved in a similar way, and hence omitted.

Fix two vertices $k\neq l$. Let $\bv^1$, $\bv^2$, $\bv^3$, $\bv^4$ be
dimension vectors such that
\begin{equation*}
 \bv^2 = \bv^1 + \alpha_k, \quad
 \bv^4 = \bv^1 + \alpha_l, \quad
 \bv^3 = \bv^2 + \alpha_l = \bv^4 + \alpha_k
 = \bv^1 + \alpha_k + \alpha_l.
\end{equation*}
We want to compute $e_{k,r} * e_{l,s}$ and $e_{l,s} * e_{k,r}$ in the
component $K^{G_\bw\times\C}(Z(\bv^1,\bv^3,\bw))$.

Let us consider the intersection
\begin{equation*}
  p_{12}^{-1}\Pa_k(\bv^2,\bw) \cap p_{23}^{-1}\Pa_{l}(\bv^3,\bw)
  \qquad\text{(resp.\ 
  $p_{12}^{-1}\Pa_{l}(\bv^4,\bw) \cap p_{23}^{-1}\Pa_{k}(\bv^3,\bw)$)}
\end{equation*}
in $\M(\bv^1,\bw)\times\M(\bv^2,\bw)\times\M(\bv^3,\bw)$
(resp.\ $\M(\bv^1,\bw)\times\M(\bv^4,\bw)\times\M(\bv^3,\bw)$).

\begin{Lemma}\label{lem:intersect}
The above two intersections are transversal respectively.
\end{Lemma}

\begin{proof}
The proof below is modelled on \cite[9.8, 9.9]{Na-alg}.
We give the proof for $p_{12}^{-1}\Pa_k(\bv^2,\bw) \cap
p_{23}^{-1}\Pa_{l}(\bv^3,\bw)$. Then the same result for
$p_{12}^{-1}\Pa_{l}(\bv^4,\bw) \cap p_{23}^{-1}\Pa_{k}(\bv^3,\bw)$
follows by $k\leftrightarrow l$, $\bv^2 \leftrightarrow \bv^4$.

We consider the complex \eqref{eq:hecke_complex} for
$\Pa_k(\bv^2,\bw)$ and $\Pa_l(\bv^3,\bw)$:
\begin{equation*}
\begin{split}
        & \HomL(V^1, V^2)
        \overset{\sigma^{12}}{\longrightarrow}
        \HomE(V^1,V^2) \oplus \HomL(W, V^2) \oplus \HomL(V^1,W)
        \overset{\tau^{12}}{\longrightarrow}
        \HomL(V^1, V^2) \oplus \shfO, \\
        & \HomL(V^2, V^3)
        \overset{\sigma^{23}}{\longrightarrow}
        \HomE(V^2,V^3) \oplus \HomL(W, V^3) \oplus \HomL(V^2,W)
        \overset{\tau^{23}}{\longrightarrow}
        \HomL(V^2, V^3) \oplus \shfO, \\
\end{split}
\end{equation*}
where we put suffixes $12$, $23$ to distinguish endomorphisms.
We have sections $s^{12}$ and $s^{23}$ of $\Ker \tau^{12}/\Ima \sigma^{12}$
and $\Ker \tau^{23}/\Ima \sigma^{23}$ respectively.

Identifying these vector bundles and sections with those of pull-backs to
$\M(\bv^1,\bw)\times\M(\bv^2,\bw)\times\M(\bv^3,\bw)$,
we consider their zero locus 
$Z(s^{12}) = \Pa_k(\bv^2,\bw)\times\M(\bv^3,\bw)$ and
$Z(s^{23}) = \M(\bv^1,\bw)\times \Pa_{l}(\bv^3,\bw)$.

As in the proof of \cite[5.7]{Na-alg}, we consider the
transpose of $\nabla s^{12}$, $\nabla s^{23}$ via the symplectic form.
Their sum gives a vector bundle endomorphism
\begin{equation*}
\begin{split}
        {}^t (\nabla s^{12}) + {}^t (\nabla s^{22})
        \colon & \Ker{}^t \sigma^{12} / \Ima {}^t \tau^{12} \oplus
        \Ker{}^t \sigma^{23} / \Ima {}^t \tau^{23} \\
        &\qquad
        \to T\M(\bv^1,\bw)\oplus T\M(\bv^2,\bw)\oplus T\M(\bv^3,\bw).
\end{split}
\end{equation*}
It is enough to show that the kernel of
${}^t (\nabla s^{12}) + {}^t (\nabla s^{23})$ is zero at $(x^1, x^2, x^3)$.

Take representatives $(B^a , i^a , j^a )$ of $x^a$ ($a = 1,2,3$).
Then we have $\xi^{12}$, $\xi^{23}$ which satisfy \eqref{hecke:xi}
for $(B^a , i^a , j^a)$.
Suppose that
\begin{equation*}
(C^{\prime 12}, a^{\prime 12}, b^{\prime 12})\pmod{\Ima {}^t \tau^{12}}
\oplus
(C^{\prime 23}, a^{\prime 23}, b^{\prime 23})\pmod{\Ima {}^t \tau^{23}}
\end{equation*}
lies in the kernel.
Then there exist
$\gamma^a\in\HomL(V^a, V^a)$ ($a = 1,2,3$)
such that
\begin{equation}\label{eq:gamma}
\begin{gathered}
\left\{\begin{split}
        \varepsilon C^{\prime 12}\,\xi^{12}  &=
         \gamma^1\, B^1  - B^1 \, \gamma^1, \\
        b^{\prime 12} &= \gamma^1\, i^1,  \\
        -a^{\prime 12}\,\xi^{12}  &= - j^1 \,\gamma^1,
\end{split}\right.
\qquad
\left\{\begin{split}
        \varepsilon \,\xi^{23} C^{\prime 23}  &=
         \gamma^3\, B^3  - B^3 \, \gamma^3, \\
        \xi^{23}\,b^{\prime 23} &= \gamma^3\, i^3,  \\
        -a^{\prime 23} &= - j^3 \,\gamma^3,
\end{split}\right. \\
\left\{\begin{split}
                \varepsilon (\xi^{12} \, C^{\prime 12} +
                        C^{\prime 23} \, \xi^{23}) &=
                 \gamma^2\, B^2  - B^2 \, \gamma^2,\\
                \xi^{12} \, b^{\prime 12} + b^{\prime 23} &=
                 \gamma^2\,i^2,  \\
                -a^{\prime 12} - a^{\prime 23}\xi^{23} &= - j^2 \gamma^2.
\end{split}\right.
\end{gathered}
\end{equation}
Then we have
\begin{gather*}
  B^3\left(
    \xi^{23}(\gamma^2\xi^{12}-\xi^{12}\gamma^1)-\gamma^3\xi^{23}\xi^{12}
    \right)
    = \left(
    \xi^{23}(\gamma^2\xi^{12}-\xi^{12}\gamma^1)-\gamma^3\xi^{23}\xi^{12}
    \right) B^1, \\
  j^3\left(
    \xi^{23}(\gamma^2\xi^{12}-\xi^{12}\gamma^1)-\gamma^3\xi^{23}\xi^{12}
    \right) = 0.
\end{gather*}
Hence we have 
\begin{equation}
\label{eq:xigamma}
    \xi^{23}(\gamma^2\xi^{12}-\xi^{12}\gamma^1)-\gamma^3\xi^{23}\xi^{12}
     = 0
\end{equation}
by the stability condition.

Consider the equation \eqref{eq:xigamma} at the vertex $l$. Since
$k\neq l$, $\xi^{12}_l$ is an isomorphism. Hence \eqref{eq:xigamma}
implies that $\Ima\xi^{23}_l$ is invariant $\gamma^3_l$. 
Since $\Ima\xi^{23}_l$ is a codimension $1$ subspace,
the induced map $\gamma^3_l\colon
V^3_l/\Ima\xi^{23}_l \to V^3_l/\Ima\xi^{23}_l$ is a scalar which we
denote by $\lambda^{23}$. Moreover, there
exists a homomorphism
$\zeta^{23}_l\colon V^3_l \to V^2_l$ such that
\[
  \gamma^3_l - \lambda^{23}\id_{V^3_l} = \xi^{23}_l \zeta^{23}_l.
\]   
For another vertex $l'\neq l$, $\xi^{23}_{l'}$ is an isomorphism, hence we
can define $\zeta^{23}_{l'}$ so that the same equation holds also for
the vertex $l'$. Thus we have
\begin{equation}
\label{eq:zeta}
  \gamma^3 - \lambda^{23}\id_{V^3} = \xi^{23} \zeta^{23}.
\end{equation}
Substituting \eqref{eq:zeta} into \eqref{eq:gamma} and using the
injectivity of $\xi^{23}$, we get
\begin{equation*}
\left\{\begin{split}
        C^{\prime 23} &= 
        \varepsilon (\zeta^{23}B^3 - B^2 \,\zeta^{23}), \\
        a^{\prime 23} &= j^3  (\xi^{23}\zeta^{23} +
                \lambda^{23}\,\id_{V^3}),\\
        b^{\prime 23} &= (\zeta^{23}\xi^{23} 
                  + \lambda^{23}\,\id_{V^2}) i^2.
\end{split}\right.
\end{equation*}
This means that
\(
(C^{\prime 23}, a^{\prime 23}, b^{\prime 23})
= \lsp{t}{\tau^{23}}(\zeta^{23}\oplus\lambda^{23})
\).

Substituting \eqref{eq:zeta} into \eqref{eq:xigamma} and noticing
$\xi^{23}$ is injective, we obtain
\begin{equation}
\label{eq:xigamma'}
  \left(
    \gamma^2 - (\zeta^{23}\xi^{23} + \lambda^{23}\id_{V^2})
  \right)\xi^{12} 
  = \xi^{12}\gamma^1 .
\end{equation}
Thus $\Ima\xi^{12}$ is invariant under
\(
    \gamma^2 - (\zeta^{23}\xi^{23} + \lambda^{23}\id_{V^2}).
\)
Arguing as above, we can find a constant $\lambda^{12}$ and
a homomorphism $\zeta^{12}\colon V^2\to V^1$ such that
\begin{equation}
\label{eq:zeta'}
  \gamma^2 - (\zeta^{23}\xi^{23} + \lambda^{23}\id_{V^2})
  - \lambda^{12}\id_{V^2} = \xi^{12}\zeta^{12}.
\end{equation}
Substituting this equation into \eqref{eq:gamma}, we get
\begin{equation*}
\left\{\begin{split}
        C^{\prime 12} &= 
        \varepsilon (\zeta^{12}B^2 - B^1 \,\zeta^{12}), \\
        a^{\prime 12} &= j^2  (\xi^{12}\zeta^{12} +
                \lambda^{12}\,\id_{V^2}),\\
        b^{\prime 12} &= (\zeta^{12}\xi^{12} 
                  + \lambda^{12}\,\id_{V^1}) i^1.
\end{split}\right.
\end{equation*}
This means that
\(
(C^{\prime 12}, a^{\prime 12}, b^{\prime 12})
= \lsp{t}{\tau^{12}}(\zeta^{12}\oplus\lambda^{12})
\).
Hence ${}^t (\nabla s^{12}) + {}^t (\nabla s^{23})$ is injective.
\end{proof}


Let us consider the variety
$\widetilde Z_{kl}$ (resp.\ $\widetilde Z_{lk}$)
of all pairs $(B,i,j)\in \mu^{-1}(0)^{\operatorname{s}}$ and 
$S\subset V^3$
satisfying the followings:
\begin{aenume}
\item $S$ is a subspace with $\dim S = \bv^1 = \bv^3 -
\alpha_k - \alpha_l$,
\item $S$ is $B$-stable,
\item $\Ima i \subset S$,
\item the induced homomorphism
$B_h\colon V_k^3/S_k \to V_l^3/S_l$
(resp.\ $B_{\overline{h}}\colon V_l^3/S_l \to V_k^3/S_k$) is zero for
$h$ with $\vout(h) = k$, $\vin(h) = l$.
\end{aenume}
Then $\widetilde Z_{kl}/G_{\bv^3}$
(resp.\ $\widetilde Z_{lk}/G_{\bv^3}$) is isomorphic to
$
  p_{12}^{-1}\Pa_k(\bv^2,\bw) \cap p_{23}^{-1}\Pa_{l}(\bv^3,\bw)
$
(resp.\ 
$
  p_{12}^{-1}\Pa_{l}(\bv^2,\bw) \cap p_{23}^{-1}\Pa_{k}(\bv^3,\bw)
$).
The isomorphism is given by defining
\begin{align*}
  & (B^1,i^1,j^1) \defeq \text{the restriction
                               of $(B,i,j)$ to $S$},\\
  & (B^2,i^2,j^2) \defeq \text{the restriction
                               of $(B,i,j)$ to $S'$},\\
   \Big(\text{resp.\ }
  & (B^4,i^4,j^4) \defeq \text{the restriction
                               of $(B,i,j)$ to $S''$},\Big)\\
  & (B^3,i^3,j^3) \defeq (B,i,j),
\end{align*}
where $S'$ (resp.\ $S''$) is given by
\begin{equation*}
  S'_m \defeq
  \begin{cases}
     V_m &\text{if $m\neq l$,}\\
     S_l &\text{if $m = l$},
  \end{cases}
\qquad\left(\text{resp.\ }
  S''_m \defeq
  \begin{cases}
     V_m &\text{if $m\neq k$,}\\
     S_k &\text{if $m = k$},
  \end{cases}
\right).
\end{equation*}
It is also clear that the restriction of
$p_{13}$ to
$
  \widetilde Z_{kl}/G_{\bv^3}
$
(resp.\ 
$
  \widetilde Z_{lk}/G_{\bv^3}
$)
is an isomorphism onto its image.
Hereafter, we identify
$
  \widetilde Z_{kl}/G_{\bv^3}
$
(resp.\ 
$
  \widetilde Z_{lk}/G_{\bv^3}
$)
with the image. Then
$
  \widetilde Z_{kl}/G_{\bv^3}
$
and
$
  \widetilde Z_{lk}/G_{\bv^3}
$
are closed subvarieties of $Z(\bv^1,\bv^3;\bw)$.
Let
$i_{kl}\colon \widetilde Z_{kl}/G_{\bv^3}\to Z(\bv^1,\bv^3;\bw)$
(resp.\ 
$i_{lk}\colon \widetilde Z_{lk}/G_{\bv^3}\to Z(\bv^1,\bv^3;\bw)$)
denote the inclusion.


The quotient $V_k^3/S_k$ (resp.\ $V_l^3/S_l$) forms a line bundle
over
\(
  \widetilde Z_{kl}/G_{\bv^3} 
  =  p_{12}^{-1}\Pa_k(\bv^2,\bw) \cap p_{23}^{-1}\Pa_{l}(\bv^3,\bw)
\)
(resp.\ 
\(
  \widetilde Z_{lk}/G_{\bv^3}
   = p_{12}^{-1}\Pa_{l}(\bv^2,\bw) \cap p_{23}^{-1}\Pa_{k}(\bv^3,\bw)
\)
).
By the above consideration, $e_{k,r}\ast e_{l,s}$
(resp.\ $e_{l,s}\ast e_{k,r}$) is represented by
{\allowdisplaybreaks[4]
\begin{equation}\label{eq:E_krE_ls}
\begin{split}
   & 
  \begin{aligned}[m]
   (-1)^{\langle h_k, \bC_{\overline{\Omega}}\bv^2\rangle
        +\langle h_l, \bC_{\overline{\Omega}}\bv^3\rangle}
   i_{kl*}\Bigl[
   & (q^{-1}V^3_k/S_k)^{\otimes r-\langle h_k, \bC''\bv^2\rangle}
     \otimes (q^{-1}V^3_l/S_l)^{\otimes s-\langle h_l, \bC''\bv^3\rangle}\\
   &\qquad\qquad\qquad
   \otimes \det C_k^{\prime\bullet}(\bv^2,\bw)^*
   \otimes \det C_l^{\prime\bullet}(\bv^3,\bw)^*\Bigr]
  \end{aligned}
 \\
   \Biggl(
   \text{resp.\ }&
   \begin{aligned}[m]
    (-1)^{\langle h_l, \bC_{\overline{\Omega}}\bv^4\rangle
         +\langle h_k, \bC_{\overline{\Omega}}\bv^3\rangle}
    i_{lk*}\Bigl[
    & (q^{-1}V^3_l/S_l)^{\otimes s-\langle h_l, \bC''\bv^4\rangle}
     \otimes (q^{-1}V^3_k/S_k)^{\otimes r-\langle h_k, \bC''\bv^3\rangle}\\
   &\qquad\qquad\qquad
   \otimes \det C_l^{\prime\bullet}(\bv^4,\bw)^*
   \otimes \det C_k^{\prime\bullet}(\bv^3,\bw)^*\Bigr]
   \end{aligned}
\Biggr).
\end{split}
\end{equation}
Note} that we have
{\allowdisplaybreaks[4]
\begin{equation}
\label{eq:notice}
\begin{gathered}
   \langle h_k, \bC_{\overline{\Omega}}\bv^2\rangle
  +\langle h_l, \bC_{\overline{\Omega}}\bv^3\rangle
   = \langle h_l, \bC_{\overline{\Omega}}\bv^4\rangle
    +\langle h_k, \bC_{\overline{\Omega}}\bv^3\rangle
   \pm (\alpha_k,\alpha_l),\\
   \langle h_k, \bC''\bv^2\rangle = \langle h_k, \bC''\bv^3\rangle, \qquad 
   \langle h_l, \bC''\bv^3\rangle = \langle h_l, \bC''\bv^4\rangle, \\
   \det C_k^{\prime\bullet}(\bv^2,\bw)
   = \det C_k^{\prime\bullet}(\bv^3,\bw)\otimes
   (q^{-1}V^3_l/S_l)^{\otimes (\alpha_k,\alpha_l)},\\
   \det C_l^{\prime\bullet}(\bv^4,\bw)
   = \det C_l^{\prime\bullet}(\bv^3,\bw)\otimes
   (q^{-1}V^3_k/S_k)^{\otimes (\alpha_k,\alpha_l)}.
\end{gathered}
\end{equation}
}

Set $b' = -(\alpha_k,\alpha_l)$.
We consider
\begin{equation*}
  \bigoplus_{\vout(h)=k, \vin(h)=l} B_{\overline h}
\end{equation*}
as a section of the vector bundle
$
q [b']_q \Hom(V^3_l/S_l, V^3_k/S_k)
$
over
$
  \widetilde Z_{kl}/G_{\bv^3}
$.
Let us denote it by $s_{kl}$.
Similarly
$
\bigoplus_{\vout(h)=k, \vin(h)=l} B_{h}
$
is a section (denoted by $s_{lk}$) of the vector bundle
$
q [b']_q \Hom(V^3_k/S_k, V^3_l/S_l)
$
over
$
  \widetilde Z_{lk}/G_{\bv^3}
$.
\begin{Lemma}\label{lem:transzero}
The section $s_{kl}$ \rom(resp.\ $s_{lk}$\rom) is transversal to the
zero section
\rom(if it vanishes somewhere\rom).
\end{Lemma}

\begin{proof}
Fix a subspace $S\subset V^3$ with $\dim S = \bv^1$. 
Let $P$ be the parabolic subgroup of $G_{\bv^3}$ consisting elements
which preserve $S$.
We also fix a complementary subspace $T$. Thus we have
$V^3 = S\oplus T$.
We will check the assertion for $s_{kl}$.
The assertion for $s_{lk}$ follows if we exchange $k$ and $l$.

We consider
\begin{equation*}
  \widetilde\bM
  \defeq \left\{ (B,i,j)\in\bM(\bv^3,\bw)\left|
      \begin{aligned}[m]
        & B(S) \subset S,\; \Ima i\subset S, \\
        & \text{$B_h\colon V^3_k/S_k \to V^3_l/S_l$
        is $0$ for $h$ with $\vout(h) = k$, $\vin(h) = l$}
      \end{aligned}
    \right\}\right..
\end{equation*}
It is a linear subspace of $\bM(\bv^3,\bw)$. Let 
$\widetilde\mu\colon \widetilde\bM\to \HomL(V^3,S)$ be the composition
of the restriction of the moment map
$\mu\colon \bM(\bv^3,\bw)\to \HomL(V^3,V^3)$ to $\widetilde\bM$
and the projection
$\HomL(V^3,V^3)\to \HomL(V^3,S)$.
Let $\widetilde\mu^{-1}(0)^{\operatorname{s}}$ denote the set of
$(B,i,j)\in\widetilde\mu^{-1}(0)$ which is stable. It is preserved
under the action of $P$ and we have a $G_{\bv^3}$-equivariant isomorphism
\begin{equation*}
   G_{\bv^3}\times_{P} \widetilde\mu^{-1}(0)^{\operatorname{s}}
   \cong \widetilde Z_{kl}.
\end{equation*}
Note that the $\HomL(V^3,T)$-part of the moment map
$\mu$ vanishes on $\widetilde\bM$ thanks to the definition of
$\widetilde\bM$.

The assertion follows if we check
\begin{equation*}
   d\widetilde\mu\oplus \Pi \colon \widetilde\bM \to
    \HomL(V^3,S) \oplus q [b']_q \Hom(V^3_l/S_l, V^3_k/S_k)
\end{equation*}
is surjective at $(B,i,j)\in\widetilde\mu^{-1}(0)^{\operatorname{s}}$.
Here $\Pi\colon\widetilde\bM \to
q [b']_q \Hom(V^3_l/S_l, V^3_k/S_k)$ is the natural projection.
Thus it is enough to show that
\(
   d\widetilde\mu\colon \widetilde\bM\cap \Ker\Pi\to \HomL(V^3,S)
\)
is surjective.

Suppose that $\zeta\in\HomL(S,V^3)$ is orthogonal to
$d\widetilde\mu(\widetilde\bM\cap \Ker\Pi)$, namely
\begin{equation*}
   \tr\left(\varepsilon(\delta B\,B + B\,\delta B)\zeta
     + \delta i\,j\,\zeta + i\,\delta j\,\zeta\right) = 0
\end{equation*}
for any $(\delta B, \delta i, \delta j)\in \widetilde\bM\cap \Ker\Pi$.
Hence we have
\begin{equation*}
   0 = j\zeta \in \HomL(S,W), \qquad
   0 = B\zeta - \zeta B|_S \in \HomL(S,V^3),
\end{equation*}
where $B|_S$ is the restriction of $B$ to $S$. Therefore the image of
$\zeta$ is invariant under $B$ and contained in $\Ker j$. By the
stability condition, we have $\zeta = 0$. Thus we have proved the
assertion.
\end{proof}

Let $Z(s_{kl})$ (resp.\ $Z(s_{lk})$) be the zero locus of
$s_{kl}$ (resp.\ $s_{lk}$). By \lemref{lem:transzero}, we have the
following exact sequence (Koszul complex) on
$\widetilde Z_{kl}/G_{\bv^3}$ (resp.\ $\widetilde Z_{lk}/G_{\bv^3}$):
\begin{align*}
&
\begin{aligned}[m]
   & 0 \to \Wedge^{\max} (q [b']_q \Hom(V^3_l/S_l, V^3_k/S_k))^*
   \to \cdots
   \to (q [b']_q \Hom(V^3_l/S_l, V^3_k/S_k))^* \\
   & \qquad\qquad\qquad\qquad\qquad\qquad\qquad\qquad
   \to \shfO_{\widetilde Z_{kl}/G_{\bv^3}} \to
   \shfO_{Z(s_{kl})} \to 0
\end{aligned}
\\
\Biggl(\text{resp.\ }
&
\begin{aligned}[m]
   & 0 \to \Wedge^{\max} (q [b']_q \Hom(V^3_k/S_k, V^3_l/S_l))^*
   \to \cdots
   \to (q [b']_q \Hom(V^3_k/S_k, V^3_l/S_l))^* \\
   & \qquad\qquad\qquad\qquad\qquad\qquad\qquad\qquad
   \to \shfO_{\widetilde Z_{lk}/G_{\bv^3}} \to
   \shfO_{Z(s_{lk})} \to 0
\end{aligned}
\Biggr).
\end{align*}
Hence we have the following equality in
$K^{G_{\bw}\times\C^*}(\widetilde Z_{kl}/G_{\bv^3})$
(resp.\ $K^{G_{\bw}\times\C^*}(\widetilde Z_{lk}/G_{\bv^3})$):
\begin{equation*}
\begin{split}
   & \shfO_{Z(s_{kl})}
   = \Wedge_{-1} (q [b']_q \Hom(V^3_l/S_l, V^3_k/S_k))^*
\\
\Big(
\text{resp.\ }
   & \shfO_{Z(s_{lk})}
   = \Wedge_{-1} (q [b']_q \Hom(V^3_k/S_k, V^3_l/S_l))^* 
\Big).
\end{split}
\end{equation*}
Both $Z(s_{kl})$ and $Z(s_{lk})$ consist of all pairs
$(B,i,j)\in\mu^{-1}(0)^{\operatorname{s}}$ and
$S\subset V^3$ satisfying the followings
\begin{aenume}
\item $S$ is a subspace with $\dim S = \bv^1 = \bv^3 -
\alpha_k - \alpha_l$,
\item $S$ is $B$-stable,
\item $\Ima i \subset S$,
\item the induced homomorphism
$B\colon V/S \to V/S$ is zero,
\end{aenume}
modulo the action of $G_{\bv^3}$. In particular, we have
$Z(s_{kl}) = Z(s_{lk})$.
Hence we have
\begin{equation*}
   i_{kl*}\left[
     \Wedge_{-1} (q [b']_q \Hom(V^3_l/S_l, V^3_k/S_k))^*\right]
  = i_{lk*}\left[
     \Wedge_{-1} (q [b']_q \Hom(V^3_k/S_k, V^3_l/S_l))^*\right].
\end{equation*}
This implies
\begin{equation}
\begin{split}
    & i_{kl*}\left[(q^{-1}V^3_k/S_k)^{\otimes r} \otimes
      (q^{-1} V^3_l/S_l)^{\otimes s} \otimes
    \Wedge_{-1} (q [b']_q \Hom(V^3_l/S_l, V^3_k/S_k))^*\right] \\
  =\; & i_{lk*}\left[(q^{-1}V^3_k/S_k)^{\otimes r} \otimes
        (q^{-1}V^3_l/S_l)^{\otimes s} \otimes
    \Wedge_{-1} (q [b']_q \Hom(V^3_k/S_k, V^3_l/S_l))^*\right]
\end{split}
\end{equation}
by the projection formula~\eqref{eq:projection}.
Multiplying this equality by $z^{-r} w^{-s}$ and taking sum with
respect to $r$ and $s$, we get
\begin{equation*}
\begin{split}
   & \prod_{p=1}^{b'} (1 - q^{b'-2p}\frac{w}{z})
   \sum_{r,s=-\infty}^\infty
    i_{kl*}\left[(q^{-1}V^3_k/S_k)^{\otimes r} \otimes
      (q^{-1} V^3_l/S_l)^{\otimes s}\right] z^{-r} w^{-s}\\
  =\; & \prod_{p=1}^{b'} (1 - q^{b'-2p}\frac{z}{w})
   \sum_{r,s=-\infty}^\infty
    i_{lk*}\left[(q^{-1}V^3_k/S_k)^{\otimes r}\otimes
      (q^{-1} V^3_l/S_l)^{\otimes s}\right] z^{-r} w^{-s}.
\end{split}
\end{equation*}
Comparing this with \eqref{eq:E_krE_ls} and using \eqref{eq:notice},
we get \eqref{eq:relexE1}.

\subsection{Relation (\ref{eq:relDS})}
We give the proof of \eqref{eq:relDS} for $\pm = +$, in this
subsection, assuming other relations.
(The relations~\eqref{eq:relEF} with $k = l$ and \eqref{eq:relexE2}
will be checked in the next section, but its proof is independent of
results in this subsection.)
The relation \eqref{eq:relDS}
for $\pm = -$ can be proved in a similar way, and hence omitted.

By the proof of \cite[9.3]{Na-alg}, operators $e_{k,0}$, $f_{k,0}$
acting on $K^{G_\bw\times\C^*}(\Zw)$ are locally nilpotent. (See also
\lemref{lem:lint} below.) It is known that the constant term of
\eqref{eq:relDS}, i.e.,
\begin{equation}\label{eq:DS0}
\sum_{p=0}^{b}(-1)^p \begin{bmatrix} b \\ p\end{bmatrix}_{q_k}
                     e_{k,0}^p e_{l,0} e_{k,0}^{b-p} = 0
\end{equation}
can be deduced from the other relations and the local nilpotency of
$e_{k,0}$, $f_{k,0}$ (see e.g., \cite[4.3.2]{Gi-book} for the proof
for $q = 1$). Thus our task is to reduce 
\begin{equation}\label{eq:genDS}
  \sum_{\sigma\in S_b}
   \sum_{p=0}^{b}(-1)^p 
   \begin{bmatrix} b \\ p\end{bmatrix}_{q_k}
   e_{k,r_{\sigma(1)}}\cdots e_{k, r_{\sigma(p)}} e_{l,s}
   e_{k, r_{\sigma(p+1)}}\cdots e_{k,r_{\sigma(b)}} = 0
\end{equation}
to \eqref{eq:DS0}.
This reduction was done by Grojnowski \cite{Gr_aff}, but we reproduce
it here for the sake of completeness.

For $p\in \{0,1,\dots, b\}$, let $\bv^0,\dots,\bv^{b+1}$ be dimension
vectors with
\begin{equation*}
   \bv^i = 
   \begin{cases}
     \bv^{i-1} + \alpha_k & \text{if $i\neq p+1$}, \\
     \bv^{i-1} + \alpha_l & \text{if $i = p+1$}.
   \end{cases}
\end{equation*}
Let
\begin{equation*}
   \pi_{ij}\colon
   \M(\bv^0,\bw)\times\cdots\times \M(\bv^{b+1},\bw) \to
   \M(\bv^{i},\bw)\times\M(\bv^j,\bw)
\end{equation*}
be the projection. Let
\begin{multline}\label{eq:inter}
   \widehat\Pa_p \defeq
   \pi_{12}^{-1}(\Pa_k(\bv^1,\bw))\cap \cdots\cap
   \pi_{p-1,p}^{-1}(\Pa_k(\bv^p,\bw))\cap \\
   \pi_{p,p+1}^{-1}(\Pa_l(\bv^{p+1},\bw))\cap
   \pi_{p+1,p+2}^{-1}(\Pa_k(\bv^{p+2},\bw))\cap\cdots\cap
   \pi_{b,b+1}^{-1}(\Pa_k(\bv^{b+1},\bw)).
\end{multline}
This is equal to
\begin{equation*}
   \{ (B,i,j,V^0\subset\cdots\subset V^{b+1}) \mid
   \text{as below}\} /G_{\bv^{b+1}}
\end{equation*}
\begin{aenume}
\item $(B,i,j)\in \mu^{-1}(0)^{\operatorname{s}}$
   (in $\bM(\bv^{b+1},\bw)$)
\item $V^i$ is a $B$-invariant subspace containing the image of $i$
with $\dim V^i = \bv^i$.
\end{aenume}
In particular, we have line bundles $V^i/V^{i-1}$ on $\widehat\Pa_p$
($i = 1,\dots,b+1$).
By the definition, there exists a line bundle $\mathfrak
L_p(r_1,\dots,r_b;s)\in K^{G_\bw\times\C^*}(\widehat\Pa_p)$ such that
\begin{equation*}
   e_{k,r_1}\cdots e_{k,r_p} e_{l,s} e_{k,r_{p+1}}\cdots e_{k,r_b}
    = \pi_{0,b+1 *} \mathfrak L_p(r_1,\dots,r_b;s).
\end{equation*}
Moreover, we have
\begin{multline}\label{eq:LL}
   \mathfrak L_p(r_1,\dots,r_b;s) =
   \left(q^{-1} V^1_k/V^0_k \right)^{\otimes r_1}
    \otimes \cdots \otimes
   \left(q^{-1} V^p_k/V^{p-1}_k \right)^{\otimes r_p} \otimes
   \left(q^{-1} V^{p+1}_l/V^{p}_l \right)^{\otimes s} \otimes \\
   \left(q^{-1} V^{p+2}_k/V^{p+1}_k \right)^{\otimes r_{p+1}}
    \otimes\cdots\otimes
   \left(q^{-1} V^{b+1}_k/V^{b}_k \right)^{\otimes r_{b}} \otimes
   \mathfrak L_p(0,\dots,0;0).
\end{multline}

Now consider the symmetrization.
By \eqref{eq:LL}, we have
\begin{equation*}
  \sum_{\sigma \in S_b}
   \mathfrak L_p(r_{\sigma(1)},\dots,r_{\sigma(b)};s)
   = T(V_k^{b+1}/V_k^0)\otimes
     \left(q^{-1} V^{p+1}_l/V^{p}_l \right)^{\otimes s}
     \otimes \mathfrak L_p(0,\dots,0;0)
\end{equation*}
for some tensor product $T(V_k^{b+1}/V_k^0)$ of
exterior products of the bundle $V_k^{b+1}/V_k^0$ and its dual.
(In the notation in \subsecref{subsec:rank1} below,
$T(V_k^{b+1}/V_k^0)$ corresponds to the symmetric function
$\sum_{\sigma\in S_b} x_1^{r_{\sigma(1)}}\cdots x_b^{r_{\sigma(b)}}$.)
Note that we have $q^{-1} V^{p+1}_l/V^{p}_l = q^{-1}
V^{b+1}_l/V^{0}_l$.
Thus
\(
   T(V_k^{b+1}/V_k^0)\otimes
     \left(q^{-1} V^{b+1}_l/V^{0}_l \right)^{\otimes s}
\)
can be considered as a vector bundle over $\pi_{0,b+1}(\widehat\Pa_p)$.
Then the projection formula implies that
\begin{equation*}
  \sum_{\sigma \in S_b} \pi_{0,b+1*}
   \mathfrak L_p(r_{\sigma(1)},\dots,r_{\sigma(b)};s)
   = T(V_k^{b+1}/V_k^0)\otimes
     \left(q^{-1} V^{b+1}_l/V^{0}_l \right)^{\otimes s}
     \otimes \pi_{0,b+1*}\mathfrak L_p(0,\dots,0;0).
\end{equation*}
Noticing that
\(
   T(V_k^{b+1}/V_k^0)\otimes
     \left(q^{-1} V^{b+1}_l/V^{0}_l \right)^{\otimes s}
\)
is independent of $p$, we can derive \eqref{eq:genDS} from \eqref{eq:DS0}.

\section{Relations~(II)}\label{sec:rel2}
The purpose of this section is to check the relations~\eqref{eq:relEF}
with $k = l$ and \eqref{eq:relexE2}.
Our strategy is the following.
We first reduce the computation of the convolution product to the case
of the graph of type $A_1$ using results in \secref{sec:convolution}
and introducing modifications of quiver varieties and Hecke
correspondences.
Then we perform the computation using the
explicit description of the equivariant $K$-theory for quiver
varieties for the graph of type $A_1$.

In this section we fix a vertex $k$.

\subsection{Modifications of quiver varieties}\label{subsec:modQui}
We take a collection of vector space 
$V = (V_l)_{l\in I}$ with $\dim V = \bv$.
Let $\mu_k$ be the $\Hom(V_k,V_k)$-component of
$\mu\colon \bM(\bv,\bw)\to \HomL(V,V)$, i.e.,
\begin{equation*}
  \mu_k(B,i,j) \defeq
  \sum_{\vin(h)=k} \varepsilon(h) B_h B_{\overline h} + i_k j_k.
\end{equation*}
Let
\begin{equation*}
\begin{split}
  & \Mp(\bv,\bw) \\
 \defeq \; &\left\{ (B,i,j)\in \mu_k^{-1}(0) \left|
   \text{$\bigoplus_{\vout(h)=k} B_h \oplus j_k\colon
   V_k \to \bigoplus_{\vout(h)=k} V_{\vin(h)}\oplus W_k$ is injective}
   \right\}\right/\! \GL(V_k).
\end{split}
\end{equation*}
This is a product of the quiver variety for the graph of type $A_1$
and the affine space:
\begin{equation*}
  \Mp(\bv,\bw) = \M(v_k,N) \times \bM'(\bv,\bw),
\end{equation*}
where
\begin{gather*}
  v_k = \dim V_k, \qquad 
  N = - \sum_{l:l\neq k} (\alpha_k,\alpha_l)\dim V_l+\dim W_k, \\
  \bM'(\bv,\bw) \defeq
  \bigoplus_{h\colon \substack{\vin(h)\neq k\\ \vout(h)\neq k}}
  \Hom(V_{\vout(h)}, V_{\vin(h)}) \oplus
  \bigoplus_{l\colon l\neq k} \Hom(W_l, V_l)\oplus \Hom(V_l,W_l).
\end{gather*}
Moreover, the variety $\M(v_k,N)$ is isomorphic to the cotangent
bundle of the Grassmann manifold $G(v_k,N)$ of $v_k$-dimensional
subspaces in the $N$-dimensional space. (See \cite[Chap.~7]{Na-quiver}
for detail.)
The isomorphism is given as follows: $\M(v_k,N)$ is the set of
$\GL(v_k,\C)$-orbit of $i\colon\C^N \to \C^{v_k}$, $j\colon
\C^{v_k}\to \C^N$ such that
\begin{aenume}
\item $ij = 0$,
\item $j$ is injective.
\end{aenume}
The action is given by $(i,j)\mapsto (gi, jg^{-1})$. Then
\begin{equation*}
  \M(v_k,N)\ni G\cdot (i,j) \mapsto
  (\operatorname{Image}j\subset\C^N)
\end{equation*}
defines a map $\M(v_k,N)\to G(v_k,N)$. And the linear map
\begin{equation*}
   ji\colon \C^N/\operatorname{Image}j \longrightarrow
   \operatorname{Image}j 
\end{equation*}
defines a cotangent vector at $\operatorname{Image}j$.

Let $\mu^{-1}(0)^{\operatorname{s}}$ be as in \defref{def:stable} and
$\mu_k^{-1}(0)^{\operatorname{s}}$ be the set of stable points in
$\mu_k^{-1}(0)$. Although the stability condition~\eqref{def:stable}
was defined only for $(B,i,j)\in\mu^{-1}(0)$, it can be defined
for any $(B,i,j)\in\bM(\bv,\bw)$.
Let
\begin{equation*}
  \Mp^{\circ}(\bv,\bw) \defeq \mu_k^{-1}(0)^{\operatorname{s}}/\GL(V_k),
   \qquad
  \Mpp(\bv,\bw) \defeq \mu^{-1}(0)^{\operatorname{s}}/\GL(V_k).
\end{equation*}
We have a natural action of 
\[
G_\bv' \defeq \prod_{l:l \neq k} \GL(V_l)
\]
on $\Mp(\bv,\bw)$, $\Mp^\circ(\bv,\bw)$ and $\Mpp(\bv,\bw)$.
We have the following relations between these varieties:
\begin{aenume}
\item $\Mp^{\circ}(\bv,\bw)$ is an open subvariety
of $\Mp(\bv,\bw)$,
\item $\Mpp(\bv,\bw)$ is a nonsingular closed subvariety of
$\Mp^\circ(\bv,\bw)$ (defined by the equation
$\mu_l = 0$ for $l\neq k$),
\item  $\Mpp(\bv,\bw)$ is a principal 
$G_\bv'$-bundle over $\M(\bv,\bw)$.
\end{aenume}

The vector space $V_k$ defines 
on $\Mp(\bv,\bw)$, $\Mp^\circ(\bv,\bw)$, $\Mpp(\bv,\bw)$ and
$\M(\bv,\bw)$. We denote all of them by $V_k$ for brevity, hoping that
it makes no confusion.

\subsection{Modifications of Hecke correspondences}\label{subsec:modHecke}
Fix $n \in \Z_{> 0}$.
Take collections of vector spaces $V^1 = (V^1_l)_{l\in I}$, $V^2 =
(V^2_l)_{l\in I}$ whose dimension vectors $\bv^1$, $\bv^2$ satisfy
$\bv^2 = \bv^1 + n\alpha_k$. 
(For the proof of \thmref{thm:hommain}, it is enough to consider the case
$n = 1$. But we study general $n$ for a later purpose.)
These data will be fixed throughout this
subsection, and we use the following notation:
\begin{gather*}
\Mp_1 = \Mp(\bv^1,\bw),\; \Mp_2 = \Mp(\bv^2,\bw),\;
  \Mp_1^\circ = \Mp^\circ(\bv^1,\bw),\;
  \Mp_2^\circ = \Mp^\circ(\bv^2,\bw),\; \\
  \Mpp_1 = \Mpp(\bv^1,\bw),\; \Mpp_2 = \Mpp(\bv^2,\bw),\;
  \M_1 = \M(\bv^1,\bw), \; \M_2 = \M(\bv^2,\bw).
\end{gather*}

Consider $\Mp_1$ and $\Mp_2$. These varieties are products of quiver
varieties for the graph of type $A_1$ and the affine space. We fix an
isomorphism $V^1_l \cong V^2_l$ for $l\neq k$. Then we have
identifications $G'_{\bv^1}\cong G'_{\bv^2}$ and
$\bM'(\bv^1,\bw)\cong \bM'(\bv^2,\bw)$. We write them
$G'$ and $\bM'$ respectively for brevity.
We write $V_l$ for $V^1_l$ and $V^2_l$ for $l\neq k$.
Let us define a subvariety $\widetilde\Pa_k^{(n)}\subset
\Mp_1\times\Mp_2$ as the product of the Hecke
correspondence for the graph of type $A_1$ and the diagonal for the
affine space. Namely
\begin{equation*}
\begin{matrix}
   & \widetilde\Pa_{k}^{(n)} & \defeq &
   \Pa_{k}^{(n)}(v^2_k,N)\times\Delta\bM' \\
   & \cap                       &        &  \cap                    \\
   & \Mp_1\times\Mp_2 & = &
   \M(v^2_k-n,N)\times\M(v^2_k,N)\times
   \bM'\times\bM',
\end{matrix}
\end{equation*}
where $v_k^2 = \dim V_k^2$,
$N = - \sum_{l:l\neq k} (\alpha_k,\alpha_l)\dim V^1_l+\dim W_k
= - \sum_{l:l\neq k} (\alpha_k,\alpha_l)\dim V^2_l+\dim W_k$.
Here $\Pa_{k}^{(n)}(v^2_k,N)\subset\M(v^2_k-n,N)\times\M(v^2_k,N)$
is the generalization of the Hecke correspondence introduced in
\eqref{eq:Pa^{(n)}}. Since the graph is of type $A_1$,
it is isomorphic to the conormal bundle of
\begin{equation*}
   O^{(n)}(v_k^2,N) \defeq
   \{ (V^1_k, V^2_k) \in G(v^2_k-n, N)\times G(v^2_k, N) \mid
   V^1_k \subset V^2_k \}.
\end{equation*}
The quotient $V^2_k/V^1_k$ defines a vector bundle over
$\Pa_{k}^{(n)}(v^2_k,N)$ of rank $n$.

We have
\begin{equation}
\label{eq:stab}
  \widetilde\Pa_{k}^{(n)}\cap
  \left(\Mp^\circ_1\times\Mp_2\right)
  \subset \Mp^\circ_1\times\Mp^\circ_2, \qquad
  \widetilde\Pa_{k}^{(n)}\cap
  \left(\Mp_1\times\Mp^\circ_2\right)
  \subset \Mp^\circ_1\times\Mp^\circ_2.
\end{equation}
The latter inclusion is obvious from the definition of stability, and
the former one follows from the argument in \cite[Proof of 4.5]{Na-alg}.
Let
\begin{equation*}
  \widetilde\Pa_{k}^{(n)\circ} \defeq
  \widetilde\Pa_{k}^{(n)}\cap
  \left(\Mp^\circ_1\times\Mp^\circ_2\right).
\end{equation*}

For $([B^1,i^1,j^1], [B^2,i^2,j^2])\in \widetilde\Pa_{k}^{(n)}$,
the first factor $[B^1,i^1,j^1]$ satisfies $\mu(B^1,i^1,j^1) = 0$ if and
only if the other factor $[B^2,i^2,j^2]$ also satisfies
$\mu(B^2,i^2,j^2) = 0$. This implies that
\begin{equation}
\label{eq:SHecke}
  \widetilde\Pa_{k}^{(n)\circ}\cap
  \left(\Mpp_1\times\Mp^\circ_2\right)
  \subset \Mpp_1\times\Mpp_2, \qquad
  \widetilde\Pa_{k}^{(n)\circ}\cap
  \left(\Mp^\circ_1\times\Mpp_2\right)
  \subset \Mpp_1\times\Mpp_2.
\end{equation}
Let
\begin{equation*}
  \widehat\Pa_{k}^{(n)}
  \defeq \widetilde\Pa_{k}^{(n)\circ}\cap
  \left(\Mpp_1\times\Mpp_2\right).
\end{equation*}
The quotient $V^2_k/V^1_k$ defines vector bundles over
$\widetilde\Pa_{k}^{(n)}$, $\widetilde\Pa_{k}^{(n)\circ}$, and
$\widehat\Pa_{k}^{(n)}$. For brevity, all are simply denoted by
$V^2_k/V^1_k$.

Let us denote by $i_a$ the inclusion $\Mpp_a\subset \Mp^\circ_a$ for
$a=1,2$.  By \eqref{eq:SHecke}, the inclusion map $i_1\times
\id_{\Mp^\circ(\bv^2,\bw)}\colon
\Mpp(\bv^1,\bw)\times\Mp^\circ(\bv^2,\bw)\to
\Mp^\circ(\bv^1,\bw)\times\Mp^\circ(\bv^2,\bw)$ induces the
pull-back homomorphism with support
\begin{equation*}
  (i_1\times \id_{\Mp^\circ_2})^*
  \colon K^{G'\times G_\bw\times\C^*}
  (\widetilde\Pa_{k}^{(n)\circ})\to
  K^{G'\times G_\bw\times\C^*}(\widehat\Pa_{k}^{(n)}).
\end{equation*}
Similarly, we have
\begin{equation*}
  (\id_{\Mp^\circ_1}\times i_2)^*
  \colon K^{G'\times G_\bw\times\C^*}
  (\widetilde\Pa_{k}^{(n)\circ})\to
  K^{G'\times G_\bw\times\C^*}(\widehat\Pa_{k}^{(n)}).
\end{equation*}

\begin{Lemma}\label{lem:pullbackSHecke}
We have
\begin{equation*}
  (i_1\times \id_{\Mp^\circ_2})^*
  [\shfO_{\widetilde\Pa_{k}^{(n)\circ}}]
  = [\shfO_{\widehat\Pa_{k}^{(n)}}], \qquad
  (\id_{\Mp^\circ_2}\times i_2)^*
  [\shfO_{\widetilde\Pa_{k}^{(n)\circ}}]
  = [\shfO_{\widehat\Pa_{k}^{(n)}}].
\end{equation*}
More generally, if $T(V^2_k/V^1_k)$ denotes a tensor product of
exterior products of the bundle $V^2_k/V^1_k$ and its dual, we have
\begin{equation*}
\begin{split}
  & (i_1\times \id_{\Mp^\circ_2})^*
  [T\!\left(V^2_k/V^1_k\right)]
  = [T\!\left(V^2_k/V^1_k\right)\otimes \shfO_{\widehat\Pa_{k}^{(n)}}],\\
  & (\id_{\Mp^\circ_1}\times i_2)^*
  [T\!\left(V^2_k/V^1_k\right)]
  = [T\!\left(V^2_k/V^1_k\right)\otimes \shfO_{\widehat\Pa_{k}^{(n)}}].
\end{split}
\end{equation*}
\end{Lemma}

\begin{proof}
The latter statement follows from the first statement and the
formula~\eqref{eq:pull_product}. Thus it is enough to check the first
statement. And the first statement follows from the transversality of
intersections \eqref{eq:SHecke} in $\Mp^\circ_1\times\Mp^\circ_2$.

Let $\mu'\colon \bM(\bv,\bw)\to \bigoplus_{l\neq k}\gl(V_l)$ be the
$\bigoplus_{l\neq k}\gl(V_l)$-part of the moment map $\mu$.
It induces a map
$\Mp^\circ_a\to \bigoplus_{l\neq k}\gl(V_l)$ for $a=1, 2$. Let us
denote it by $\mu_a'$. Thus we have $\Mpp_a = \mu_a^{\prime-1}(0)$.
Composing $\mu_a'$ with the projection 
$p_a\colon\Mp^\circ_1\times\Mp^\circ_2\to \Mp^\circ_a$, we have
a map $\mu_a'\circ p_a\colon \Mp^\circ_1\times\Mp^\circ_2\to
\bigoplus_{l\neq k}\gl(V_l)$ for $a=1,2$. We denote it also by
$\mu_a'$ for brevity.
It is enough to show that
the restriction of the differential $d\mu_a'$ to
$T\widetilde\Pa_k^{(n)\circ}$ is surjective on
$\widehat\Pa_k^{(n)} =
\widetilde\Pa_k^{(n)\circ}\cap(\Mp^\circ_1\times\Mpp_2) =
\widetilde\Pa_k^{(n)\circ}\cap(\Mpp_1\times\Mp^\circ_2)$. 

We consider the homomorphisms $\sigma_k$, $\tau_k$ defined in
\eqref{eq:taut_cpx} where $(B,i,j)$ is replaced by $(B^a,i^a,j^a)$
($a=1,2$). We denote them by $\sigma^a_k$ and $\tau^a_k$ respectively.

Take a point
$([B^1,i^1,j^1],[B^2,i^2,j^2])\in\widehat\Pa_k^{(n)}$. Then
\begin{equation*}
   B^1_h = B^2_h \quad (\vin(h)\neq k, \vout(h)\neq k), \quad
   i^1_l = i^2_l, \quad j^1_l = j^2_l \quad(l\neq k),
\end{equation*}
and there exists $\xi_k\colon V_k^1\to V_k^2$ such that
\begin{equation*}
   \xi_k B^1_h = B^2_h, \quad B^1_{\overline h} \xi_k = B^2_{\overline h}
    \quad(\vin(h) = k),\quad
   \xi_k i^1_k = i^2_k, \quad
   j^1_k = j^2_k \xi_k.
\end{equation*}
The tangent space $T\widetilde\Pa_k^{(n)\circ}$ at is isomorphic to
the space of $(\delta B^1, \delta i^1, \delta j^1, \delta B^2, \delta
i^2, \delta j^2) \in \bM(\bv^1,\bw)\times\bM(\bv^2,\bw)$ such that
\begin{equation}
\label{eq:TP}
{\allowdisplaybreaks[4]
\begin{gathered}
   \delta B^1_h = \delta B^2_h
    \quad (\text{if $\vin(h)\neq k$, $\vout(h)\neq k$}), \\
    \delta i_l^1 = \delta i_l^2,\quad \delta j_l^1 = \delta j_l^2
    \quad (\text{for $l\neq k$}),\\
   \tau_k^a\, \delta\sigma_k^a + \delta\tau_k^a\, \sigma_k^a = 0
   \qquad(a=1,2),\\
    \xi_k\circ \delta\tau_k^1 = \delta\tau_k^2
   , \qquad
     \delta\sigma_k^1 = \delta\sigma_k^2\circ\xi_k
\end{gathered}
}
\end{equation}
modulo the image of
\begin{equation}\label{eq:Parabolic}
\begin{split}
   & \{ (\delta\xi_k^1, \delta \xi_k^2) \in \gl(V_k^1)\times \gl(V^2_k)\mid
    \delta\xi_k^2\, \xi_k = \xi_k \delta\xi_k^1 \} \\
   & \qquad \longmapsto 
  \left\{
   \begin{aligned}[m]
    &\delta\xi_k^1 B^1_h, \; - B^1_{\overline h}\delta\xi_k^1, \;
    \delta\xi_k^2 B^2_h, \; - B^2_{\overline h}\delta\xi_k^2
     \quad(\vin(h) = k), \\
    & \delta\xi_k^1 i_k^1, \; -j_k^1 \delta\xi_k^1,\;
    \delta\xi_k^2 i_k^2, \; -j_k^2 \delta\xi_k^2\\
    & \text{other components are $0$},
   \end{aligned}
   \right.
\end{split}
\end{equation}
where
\begin{equation*}
  \delta\sigma_k^a
    = \bigoplus_{\vin(h)=k}\delta B_{\overline{h}}^a\oplus\delta j_k^a,
  \quad
  \delta\tau_k^a
    = \bigoplus_{\vin(h)=k}\varepsilon(h)\delta B_h^a + \delta i_k^a,
  \quad (a=1,2),
\end{equation*}
and we have used the identification $V^1_l \cong V^2_l$ for $l\neq
k$.

Now suppose that $(\zeta_l)_{l\neq k} \in \bigoplus_{l\neq k}
\gl(V_l)$ is orthogonal to the image of
$d\mu'_a|_{T\widetilde\Pa_k^{(n)\circ}}$.
Putting $\zeta_k = 0$, we consider $\zeta = (\zeta_l)$ as an
element of $\HomL(V^a,V^a)$. 
Then
\begin{equation*}
   \tr\left(\varepsilon \delta B^a(B^a\zeta - \zeta B^a) 
     + \delta i^a\, j^a\zeta + \zeta i^a\, \delta j^a\right) = 0
\end{equation*}
for any $(\delta B^1, \delta i^1, \delta j^1, \delta B^2, \delta
i^2, \delta j^2)\in T\widetilde\Pa_k^{(n)\circ}$.
Since the image of \eqref{eq:Parabolic} lies in the kernel
of $d\mu'_a$, the above equality holds for any
$(\delta B^1, \delta i^1, \delta j^1, \delta B^2, \delta
i^2, \delta j^2)$ satisfying \eqref{eq:TP}.

Taking $(\delta B^1, \delta i^1, \delta j^1, \delta B^2, \delta
i^2, \delta j^2)$ from $\Delta\bM'(\bv,\bw)$, we find
\begin{equation}
\label{eq:notk}
\begin{gathered}
   B_h^a\zeta_{\vout(h)} = \zeta_{\vin(h)} B_h^a\qquad
    \text{if $\vin(h)\neq k$, $\vout(h)\neq k$},\\
   \zeta_l i^2_l = 0, \quad j_l^2 \zeta_l = 0\qquad 
   \text{if $l\neq k$}.
\end{gathered}
\end{equation}

Next taking 
$(\delta B^1, \delta i^1, \delta j^1, \delta B^2,
\delta i^2, \delta j^2)$ from the other component (data related to the
vertex $k$), we get
\begin{equation*}
   \sigma_k^a \circ
   \left(\sum_{\vin(h)=k} \varepsilon(h) B^a_h \zeta_{\vout(h)}, 0\right)
   = \left(
   \bigoplus_{\vin(h)=k} \zeta_{\vout(h)} B^a_{\overline{h}} \oplus 0
    \right)\circ \tau_k^a
   \in
   \End\left(\bigoplus_{\vin(h)=k} V_{\vout(h)}\oplus W_k\right).
\end{equation*}
Comparing $\Hom(V_{\vout(h)}, W_k)$-components, we find
\begin{equation}\label{eq:kk}
   \varepsilon(h) j_k B^a_h \zeta_{\vout(h)} 
   = 0.
\end{equation}
Comparing $\Hom(V_{\vout(h)}, V_{\vout(h')})$-components, we have
\begin{equation}\label{eq:kkk}
   B^a_{\overline{h'}}\varepsilon(h) B^a_h \zeta_{\vout(h)}
   = \zeta_{\vout(h')} B^a_{\overline{h'}}\varepsilon(h) B^a_h.
\end{equation}
If we define
\begin{equation*}
   S_l \defeq
   \begin{cases}
      \Ima \zeta_l & \text{if $l\neq k$}, \\
      \sum_{\vin(h)=k} \Ima\left(B^a_h \zeta_{\vout(h)}\right)
      &\text{if $l = k$},
   \end{cases}
\end{equation*}
\eqref{eq:notk}, \eqref{eq:kk} and \eqref{eq:kkk} implies that $S =
(S_l)$ is $B^a$-invariant and contained in $\Ker j$. Thus $S = 0$ by
the stability condition. In particular, we have $\zeta = 0$.  This
means that $d\mu'_a|_{T\widetilde\Pa_k^{(n)\circ}}$ is surjective.
\end{proof}

Let $\pi_a\colon \Mpp_a\to \M_a$ be the projection
($a = 1,2$). Then we have
\begin{equation}\label{eq:PHecke}
\begin{gathered}
 \text{the restriction of
$\id_{\Mpp_1}\times\pi_2\colon \Mpp_1\times\Mpp_2
\to \Mpp_1\times \M_2$ to
$\widehat\Pa_k^{(n)}$ is proper,} \\
 (\id_{\Mpp_1}\times\pi_2)(\widehat\Pa_k^{(n)})
 = (\pi_1\times\id_{\M_2})^{-1}(\Pa_k^{(n)}), \\
 \text{the restriction of
$\pi_1\times\id_{\Mpp_2}\colon \Mpp_1\times\Mpp_2
\to \M_1\times \Mpp_2$ to
$\widehat\Pa_k^{(n)}$ is proper,} \\
 (\pi_1\times\id_{\Mpp_2})(\widehat\Pa_k^{(n)})
 = (\id_{\M_1}\times\pi_2)^{-1}(\Pa_k^{(n)}).
\end{gathered}
\end{equation}
By these properties, we have homomorphisms
\begin{equation*}
\begin{split}
  & (\pi_1\times\id_{\M_2})^{*-1}
  \left(\id_{\Mpp_1}\times\pi_2\right)_*\colon
  K^{G'\times G_\bw\times\C^*}(\widehat\Pa_k^{(n)})\to
   K^{G_\bw\times\C^*}(\Pa_k^{(n)}), \\
  & (\id_{\M_1}\times\pi_1)^{*-1}
  \left(\pi_1\times\id_{\Mpp_2}\right)_*\colon
  K^{G'\times G_\bw\times\C^*}(\widehat\Pa_k^{(n)})\to
   K^{G_\bw\times\C^*}(\Pa_k^{(n)}).
\end{split}
\end{equation*}

\begin{Lemma}\label{lem:pushforwardPHecke}
We have
\begin{equation*}
\begin{split}
 & (\pi_1\times\id_{\M_2})^{*-1}
  \left(\id_{\Mpp_1}\times\pi_2\right)_*
  [\shfO_{\widehat\Pa_k^{(n)}}]
  = [\shfO_{\Pa_k^{(n)}}], \\
 & (\id_{\M_1}\times\pi_2)^{*-1}
  \left(\pi_1\times\id_{\Mpp_2}\right)_*
  [\shfO_{\widehat\Pa_k^{(n)}}]
  = [\shfO_{\Pa_k^{(n)}}].
\end{split}
\end{equation*}
More generally, if $T(V^2_k/V^1_k)$ denotes a tensor product of
exterior products of the bundle $V^2_k/V^1_k$ and its dual, we have
\begin{equation*}
\begin{split}
 & (\pi_1\times\id_{\M_2})^{*-1}
  \left(\id_{\Mpp_1}\times\pi_2\right)_*
  [(\pi_1\times\pi_2)^* T\left( V^2_k/V^1_k\right)]
  = [T\left( V^2_k/V^1_k\right)], \\
 & (\id_{\M_1}\times\pi_2)^{*-1}
  \left(\pi_1\times\id_{\Mpp_2}\right)_*
  [(\pi_1\times\pi_2)^* T\left(V^2_k/V^1_k\right)]
  = [T\left(V^2_k/V^1_k\right)].
\end{split}
\end{equation*}
\end{Lemma}

\begin{proof}
The latter statement follows from the former one together with the
projection formula~\eqref{eq:projection}. Thus it is enough to prove
the former statement.

By definition, $(\pi_1\times\id_{\M_2})^{-1}(\Pa_k^{(n)})$ consists of
$(\GL(V^1_k)\cdot(B^1,i^1,j^1), G_{\bv^2}\cdot (B^2,i^2,j^2)) \in
\Mpp_1\times\M_2$ such that there exists $\xi\in\HomL(V^1, V^2)$
satisfying \eqref{hecke:xi}. We fix representatives
$(B^1,i^1,j^1)$, $(B^2,i^2,j^2)$. Then the above $\xi$ is uniquely
determined. Recall that we have chosen the identification
$V^1_l\cong V^2_l$ for $l\neq k$ over $\Mpp_1\times\Mpp_2$.
Let us define $\xi'\in \HomL(V^2,V^2)$ by
\begin{equation*}
  \xi'_l \defeq 
  \begin{cases}
    \id & \text{if $l = k$} \\
    \xi_l & \text{otherwise}.
  \end{cases}
\end{equation*}
We define a new datum
\begin{equation*}
(B^3, i^3, j^3) \defeq
  (\xi^{\prime-1}B^2 \xi', \xi^{\prime-1}i^2, j^2\xi').
\end{equation*}
By definition, we have
\begin{gather*}
  \xi_k B^1_h = B^3_h, \quad
  B^1_{\overline{h}} = B^3_{\overline{h}}\xi_k
  \quad(\vin(h) = k), \quad
  \xi_k i^1_k = i^3_k , \quad j^1_k = j^3_k\xi_k, \\
   B^1_h = B^3_h \quad (\vin(h)\neq k, \vout(h)\neq k), \quad
   i^1_l = i^3_l, \quad j^1_l = j^3_l \quad(l\neq k).
\end{gather*}
Hence $(\GL(V^1_k)\cdot (B^1,i^1,j^1), \GL(V^2_k)\cdot (B^3,i^3,j^3))$ is
contained in $\widehat\Pa_k$. Moreover, $\GL(V^2_k)\cdot (B^3,i^3,j^3)$ is
independent of the choice of the representative $(B^2,i^2,j^2)$. Thus
we have defined a map $(\pi_1\times\id_{\M_2})^{-1}(\Pa_k^{(n)})\to
\widehat\Pa_k^{(n)}$ by
\begin{equation*}
(\GL(V^1_k)\cdot(B^1,i^1,j^1), G_{\bv^2}\cdot (B^2,i^2,j^2))
\mapsto 
(\GL(V^1_k)\cdot(B^1,i^1,j^1), \GL(V^2_k)\cdot (B^3,i^3,j^3)),
\end{equation*}
which is the inverse of the restriction of
$\id_{\Mp_1}\times \pi_2$. In particular, this implies
\begin{equation*}
  \left(\id_{\Mpp_1}\times\pi_2\right)_*
  [\shfO_{\widehat\Pa_k^{(n)}}] =
  [\shfO_{(\pi_1\times\id_{\M_2})^{-1}(\Pa_k^{(n)})}].
\end{equation*}

Since $\pi_1\times\id_{\M_2}\colon 
(\pi_1\times\id_{\M_2})^{-1}(\Pa_k^{(n)})\to \Pa_k^{(n)}$ is a principal
$G'$-bundle, we have
\begin{equation*}
  (\pi_1\times\id_{\M_2})^*[\shfO_{\Pa_k^{(n)}}] = 
[\shfO_{(\pi_1\times\id_{\M_2})^{-1}(\Pa_k^{(n)})}].
\end{equation*}
Thus we have proved the first equation. The second equation can be
proved in a similar way.
\end{proof}

\subsection{Reduction to rank $1$ case}
First consider the relation~\eqref{eq:relEF} for $k = l$.
Let $\bv^1$, $\bv^2$, $\bv^3$, $\bv^4$ be dimension vectors such that
\begin{equation*}
 \bv^1 = \bv^3, \quad
 \bv^2 = \bv^1 + \alpha_k, \quad
 \bv^4 = \bv^1 - \alpha_k.
\end{equation*}
We want to compute $x_{k}^+(z) * x_{k}^-(w)$ and $x_{k}^-(w) *
x_{k}^+(z)$ in the component $K^{G_\bw\times\C}(Z(\bv^1,\bv^3,\bw))$,
and then compare it with the right hand side of \eqref{eq:relEF} with
$k = l$ in the same component.

Let $\Mpp(\bv^a,\bw)$, $\Mp(\bv^a,\bw)$, $\Mp^\circ(\bv^a,\bw)$,
$G'_{\bv^a}$ and $\bM'(\bv^a,\bw)$ be as in \subsecref{subsec:modQui}.
Let $\Pa_k(\bv^a,\bw)$, $\widehat\Pa_k(\bv^a,\bw)$,
$\widetilde\Pa_k(\bv^a,\bw)$, $\widetilde\Pa_k^\circ(\bv^a,\bw)$ be
the Hecke correspondence and its modifications introduced in
\subsecref{subsec:modHecke}. (We drop the superscript $(n)$ and write
the dimension vector $\bv^a$, $\bw$.) Let $\omega$ be the exchange of
factors as before.
Let $\widehat Z(\bv^a,\bv^b;\bw)$, $\widetilde Z(\bv^a,\bv^b;\bw)$,
$\widetilde Z^\circ(\bv^a,\bv^b;\bw)$ be subvarieties in
$\Mpp(\bv^a,\bw)\times\Mpp(\bv^b,\bw)$,
$\Mp(\bv^a,\bw)\times\Mp(\bv^b,\bw)$,
$\Mp^\circ(\bv^a,\bw)\times\Mp^\circ(\bv^b,\bw)$ defined in the same
way as $Z(\bv^a,\bv^b;\bw)$.

We have the following commutative diagram:
\begin{equation*}
\begin{CD}
{
   K^{G_\bw\times\C^*}(\Pa_k(\bv^2,\bw))\times
   K^{G_\bw\times\C^*}(\omega \Pa_k(\bv^2,\bw))} @>>>
{
   K^{G_\bw\times\C^*}(Z(\bv^1,\bv^3;\bw))} \\
      @AAA @AAA \\
{
   K^{G'\times G_\bw\times\C^*}(\widehat\Pa_k(\bv^2,\bw))\times
   K^{G'\times G_\bw\times\C^*}(\omega\widehat\Pa_k(\bv^2,\bw))} @>>>
{
   K^{G'\times G_\bw\times\C^*}(\widehat Z(\bv^1,\bv^3;\bw))} \\
      @AAA @AAA \\
{
   K^{G'\times G_\bw\times\C^*}(\widetilde\Pa_k^\circ(\bv^2,\bw))\times
   K^{G'\times G_\bw\times\C^*}(\omega\widetilde\Pa_k^\circ(\bv^2,\bw))} @>>>
{
   K^{G'\times G_\bw\times\C^*}(\widetilde Z^\circ(\bv^1,\bv^3;\bw))} \\
      @AAA @AAA \\
{
   K^{G'\times G_\bw\times\C^*}(\widetilde\Pa_k(\bv^2,\bw))\times
   K^{G'\times G_\bw\times\C^*}(\omega\widetilde\Pa_k(\bv^2,\bw))} @>>>
{
   K^{G'\times G_\bw\times\C^*}
   (\widetilde Z(\bv^1,\bv^3;\bw))}.
\end{CD}
\end{equation*}
The horizontal arrows are convolution products relative to
\begin{enumerate}
\item $\M(\bv^1,\bw)$, $\M(\bv^2,\bw)$,
      $\M(\bv^3,\bw)$,
\item $\Mpp(\bv^1,\bw)$, $\Mpp(\bv^2,\bw)$,
      $\Mpp(\bv^3,\bw)$,
\item $\Mp^\circ(\bv^1,\bw)$, $\Mp^\circ(\bv^2,\bw)$,
  $\Mp^\circ(\bv^3,\bw)$,
\item $\Mp(\bv^1,\bw)$, $\Mp(\bv^2,\bw)$,
      $\Mp(\bv^3,\bw)$.
\end{enumerate}
The vertical arrows between the first and the second rows are
homomorphisms given in \propref{prop:Pconv}. The arrows between the
second and the third are homomorphisms given in \propref{prop:Sconv}.
By the property (\ref{eq:SHecke}, \ref{eq:PHecke}) and
\begin{aenume}
\item $\widetilde Z^\circ(\bv^1,\bv^3;\bw)\cap
  \left(\Mpp(\bv^1,\bw)\times\Mp^\circ(\bv^2,\bw)\right)
  \subset \widehat Z(\bv^1,\bv^3;\bw)$,
\item the restriction of $\pi_1\times\id_{\Mpp(\bv^3,\bw)}\colon
  \Mpp(\bv^1,\bw)\times\Mpp(\bv^3,\bw)\to \M(\bv^1,\bw)\times
  \Mpp(\bv^3,\bw)$ to $\widehat Z(\bv^1,\bv^3;\bw)$ is proper,
\item $(\pi_1\times\id_{\Mpp(\bv^3,\bw)})(\widehat Z(\bv^1,\bv^3;\bw))
\subset (\id_{\M(\bv^1,\bw)}\times\pi_3)^{-1}(Z(\bv^1,\bv^3;\bw))$,
\end{aenume}
those homomorphisms can be defined.
Finally the arrows between the third and the fourth are restriction to
open subvarieties.

The commutativity for the first and the second squares follow from
\propref{prop:Pconv},~\ref{prop:Sconv} respectively. The last square
is also commutative since $\Mp^\circ(\bv^a,\bw)$ is an open subvariety
of $\Mp(\bv^a,\bw)$ and since we have \eqref{eq:stab}.

Recall that the modified Hecke correspondences in the last row is the
product of the Hecke correspondence for type $A_1$ and the diagonal
$\Delta\bM'(\bv^1,\bw)$.
Under the composite of vertical homomorphisms, $e_{k,r}$, $f_{k,s}$ at
the upper left are the images of the exterior products of the
corresponding elements for type $A_1$ and
$\shfO_{\Delta\bM'(\bv^1,\bw)}$ at the lower left, except the
following two differences:
\begin{aenume}
\item the groups acting varieties are different,
\item the sign factors in \eqref{eq:defEF}, which involve the
orientation $\Omega$, are different.
\end{aenume}
For the quiver varieties of type $A_1$, the group is
\begin{equation*}
   \GL\left(\bigoplus_{h:\vout(h)=k} V_{\vin(h)}\oplus W_k\right)
   \times \C^*
   \cong \GL_N(\C)\times\C^*.
\end{equation*}
But, if we define a homomorphism
$G'\times G_\bw\times \C^* \to \GL_N(\C)\times\C^*$ by
\begin{equation*}
  \left((g_l)_{l\in I:l\neq k}, (h_{l'})_{l'\in I}, q\right)
  \longmapsto
  \left(\bigoplus_{h:\vout(h)=k} q^{m(h)} g_{\vin(h)} \oplus
    h_k, q\right),
\end{equation*}
we have an induced homomorphism in equivariant $K$-groups:
$K^{\GL_N(\C)\times\C^*}(\ )\to
 K^{G'\times G_\bw\times \C^*}(\ )$.
(Here $m(h)$ is as in \eqref{eq:m(h)}.)
It is compatible with the convolution produce, hence it is enough to
check the relation in $K^{\GL_N(\C)\times\C^*}(\ )$.

Furthermore, the sign factor cansels out in $e_{k,r} * f_{k,s}$.
Thus the above differences make no effect when we check the
relation~\eqref{eq:relEF}.

By the commutativity of the diagram, $e_{k,r} * f_{k,s}$ is the image of
the corresponding element in the lower right.

We have a similar commutative diagram to compute $f_{k,s} * e_{k,r}$.
Hence the commutator $[e_{k,r}, f_{k,s}]$ is the image of the
corresponding commutator in the lower right.
In the next section, we will check the relation~\eqref{eq:relEF} for
type $A_1$. In particular, the commutator in the lower right is
represented by tautological bundles, considered as an element of the
$K$-theory of the diagonal $\Delta\Mp(\bv^1,\bw)$. Note that
$\shfO_{\Delta\Mp(\bv^1,\bw)}$ is mapped to
$\shfO_{\Delta\M(\bv^1,\bw)}$ by Examples in \secref{sec:convolution},
and that the tautological bundles on $\Mp(\bv^1,\bw)$ are restricted
to tautological bundles on $\M(\bv^1,\bw)$. Hence we have exactly the
same relation \eqref{eq:relEF} for general case.

Similarly, we can reduce the check of the relation~\eqref{eq:relexE2}
to the case of type $A_1$.

\subsection{Rank $1$ case}\label{subsec:rank1}
In this subsection, we check the relation when the graph is of type
$A_1$. This calculation is essentially the same as one by
Vasserot~\cite{Vasserot}, but we reproduce it here for the convenience
of the reader. (Remark that our $\C^*$-action is different from one
in \cite{Vasserot}. The definition of $e_{k,r}$, etc.\ is also
different.) We drop the subscript $k\in I$ as usual.

We prepare several notations.
For $a,b\in\Z$, let
\begin{equation*}
   [a,b] \defeq 
\begin{cases}
  \{ a, a+1, \dots, b \} & \text{if $b \ge a$}, \\
  \emptyset & \text{otherwise}.
\end{cases}
\end{equation*}
Let
\begin{equation*}
  \mathbf R \defeq 
\Z[q,q^{-1}][x_1^{\pm},\dots, x_{N}^\pm].
\end{equation*}
For a partition $I = (I_1, I_2)$ of the set $\{1, \dots, N\}$ into $2$
subsets, let $S_I = S_{I_1} \times S_{I_2}$ be the subgroups of $S_N$
consisting of permutations which preserve each subset.
For a subgroup $G\subset S_N$, let $\mathbf R^{G}$ be the subring of
$\mathbf R$ consisting of elements which are fixed by the action of $G$.
If $J$ is another partition of $\{1,\dots,N\}$, we define the 
symmetrizer ${\mathfrak S}_I^J\colon
{\mathcal R}^{S_I\cap S_J}\to {\mathcal R}^{S_J}$ by
\begin{equation*}
   f \mapsto \sum_{\sigma\in S_J/S_I\cap S_J} \sigma(f),
\end{equation*}
where $\mathcal R$ is the quotient field of $\mathbf R$.
%
For each $v\in [0, N]$, let $[v]$ be the partition
\(
   ([1, v], [v+1,N])
\).
If $I = (I_1, I_2)$ is a partition of $\{1,\dots, N\}$ into $2$
subsets and $k\in I_1$ (resp.\ $k\in I_2$), we define a new partion
$\tau_k^+(I)$ (resp.\ $\tau_k^-(I)$) by
\begin{equation*}
   \tau_k^+(I) \defeq (I_1\setminus \{k\}, I_2\cup \{ k\}), \quad
   \left(\text{resp.\
   $\tau_k^-(I)\defeq (I_1\cup \{ k\}, I_2\setminus \{k\})$}\right).
\end{equation*}
If $I = (I_1, I_2)$ is a partition and $f\in \mathbf R^{S_{[v]}}$, we
define
\begin{equation*}
   f(x_I) \defeq f(x_{i_1}, \dots, x_{i_v},
   x_{j_1}, \dots, x_{j_{N-v}}),
\end{equation*}
where $I_1 = \{ i_1, \dots, i_v \}$, $I_2 = \{ j_1, \dots, j_{N-v} \}$.

Let $\M(v,N)$ be the quiver variety for the graph of type $A_1$
with dimension vectors $v$, $N$. It is isomorphic to the cotangent
bundle of the Grassmann variety of $v$-dimensional subspaces of an
$N$-dimensional space. Let $G(v,N)$ denote the Grassmann variety
contained in $\M(v,N)$ as $0$-section. Let $Z(v^1,v^2;N)$ be the
analogue of Steinberg's variety as before. The following lemma is
crucial.
\begin{Lemma}[\protect{\cite[Lemma~13]{Vasserot},
\cite[Claim~7.6.7]{Gi-book}}]\label{lem:Vass1}
The representation of
$\bigoplus K^{\GL_N(\C)\times\C^*}(\linebreak[4]Z(v^1,v^2;N))$ on
$\bigoplus K^{\GL_N(\C)\times\C^*}(T^*G(v,N))$ by convolution is faithful.
\end{Lemma}

Thus it is enough to check the relation in $\bigoplus
K^{\GL_N(\C)\times\C^*}(T^*G(v,N))$.

Let $\Pa^{(n)}(v,N)\subset\M(v-n,N)\times\M(v,N)$ be as in
\eqref{eq:Pa^{(n)}}. It is the conormal bundle of
\begin{equation*}
   O^{(n)}(v,N) \defeq
   \{ (V^1, V^2) \in G(v-n, N)\times G(v, N) \mid
   V^1 \subset V^2 \}.
\end{equation*}
We denote the projections for $\Pa(v,N)$ by $p_1$, $p_2$, and the
projections for $O(v,N)$ by $P_1$, $P_2$. Note that both $P_1$, $P_2$
are smooth and proper. Let $\pi_1$, $\pi_2$ denote the projections
$T^*G(v-n,N)\to G(v-1,N)$, $T^*G(v,N)\to G(v,N)$.

\begin{Lemma}[\protect{\cite[Corollary~4]{Vasserot}}]\label{lem:Vass2}
For $E\in K^{\GL_N(\C)\times\C^*}(G(v,N))$
\rom(resp.\ $E\in K^{\GL_N(\C)\times\C^*}(G(v-n,N))$\rom), we have
\begin{equation*}
\begin{split}
   & [\shfO_{\Pa^{(n)}(v,N)}] * \pi_2^* E
   = \sum_i (-q^{-2})^i \pi_1^* P_{1*} \left( [\Wedge^i TP_1]
     \otimes P_2^* E\right), \\
   \Bigg(\text{resp.\ }
   & [\shfO_{\Pa^{(n)}(v,N)}] * \pi_1^* E
   = \sum_i (-q^{-2})^i \pi_2^* P_{2*} \left( [\Wedge^i TP_2]
     \otimes P_1^* E\right),\Bigg)
\end{split}
\end{equation*}
where $TP_1$ \rom(resp.\ $TP_2$\rom) is the relative tangent bundle
along the fibers of $P_1$ \rom(resp.\ $P_2$\rom). 
\end{Lemma}

\begin{proof}
As is explained in \cite[Corollary~4]{Vasserot}, the result follows
from \lemref{lem:cleanint}.
The factor $q^{-2}$ is introduced to make the differential in the
Koszul complex {\it equivariant}.
\end{proof}

By the Thom isomorphism~\cite[5.4.17]{Gi-book}, $\pi^*\colon
K^{\GL_N(\C)\times\C^*}(G(v,N))\to K^{\GL_N(\C)\times\C^*}(T^*G(v,N))$
is an isomorphism. Moreover, we have the following explicit
description of the $K$-group of the Grassmann variety (cf.\ 
\cite[6.1.6]{Gi-book}):
\begin{equation*}
   K^{\GL_N(\C)\times\C^*}(G(v,N)) \cong
   R(\C^*\times \GL_v(\C)\times \GL_{N-v}(\C)) \cong
   \mathbf R^{S_{[v]}},
\end{equation*}
where $S_{[v]} = S_v\times S_{N-v}$ acts as permutations of
$x_1,\dots, x_v$ and $x_{v+1},\dots, x_N$. If $E$ denotes the
tautological rank $v$ vector bundle over $G(v,N)$ and $Q$ denotes
the quotient bundle $\shfO^{\oplus N}/E$, the isomorphism is given by
\begin{gather*}
  e_i(x_1,\dots, x_v) \mapsto \Wedge^i E, \qquad
  e_i(x_1^{-1},\dots, x_v^{-1}) \mapsto \Wedge^i E^*, \\
  e_i(x_{v+1},\dots, x_N) \mapsto \Wedge^i Q, \quad
  e_i(x_{v+1}^{-1},\dots, x_N^{-1}) \mapsto \Wedge^i Q^*
\end{gather*}
where $e_i$ denotes the $i$th elementary symmetric polynomial.

The tautological vector bundle $V$ is isomorphic to $q E$, and
$W$ is isomorphic to the trivial bundle $\shfO^{\oplus N}$.
Let $C^\bullet(v,N)$ 
(resp.\ $C^{\prime\bullet}(v,N)$, $C^{\prime\prime\bullet}(v,N)$)
be the complex~\eqref{eq:taut_cpx} 
(resp.\ \eqref{eq:mod_taut_cpx}) over $\M(v,N)$.
In the description above, we have
\begin{equation}\label{eq:genP}
\begin{gathered}
  \Wedge_{-1/z} C^\bullet(v,N) =
  \left( \prod_{u\in [1,v]}(1 - z^{-1} q x_u) \right)^{-1}
  \prod_{t\in [v+1,N]} (1 - z^{-1} q^{-1} x_t), \\
  \det C^{\prime\bullet}(v,N) = \prod_{t\in [v+1,N]} q^{-1} x_t,
  \quad
  \det C^{\prime\prime\bullet}(v,N) = \prod_{u\in [1,v]} q^{-1} x_u^{-1}.
\end{gathered}
\end{equation}

We also have
\begin{equation*}
   K^{\GL_N(\C)\times\C^*}(O^{(n)}(v,N)) \cong
   \mathbf R^{S_{[v-n]}\cap S_{[v]}},
\end{equation*}
where $S_{[v-n]}\cap S_{[v]} \cong S_{v-n}\times S_n\times S_{N-v}$
acts as permutations of $x_1,\dots, x_{v-n}$, $x_{v-n+1},\dots, x_v$,
and $x_{v+1},\dots, x_N$.
The natural vector bundle $V^2/V^1$ is $q (x_{v-n+1} + \cdots +
x_{v})$.
The relative tangent bundles $TP_1$, $TP_2$ are
\begin{equation*}
    [TP_1] = \sum_{t=v+1}^N \sum_{k=v-n+1}^v \frac{x_t}{x_k}, \qquad
    [TP_2] = \sum_{u=1}^{v-n} \sum_{k=v-n+1}^v \frac{x_k}{x_u}.
\end{equation*}

\begin{Lemma}[\protect{\cite[Proposition~6]{Vasserot}}]\label{lem:Vass3}
\rom{(1)} The pullback homomorphisms
$P_1^*\colon K^{\GL_N(\C)\times\C^*}(G(v-n,N))\to
K^{\GL_N(\C)\times\C^*}(O(v,N))$, $P_2^*\colon
K^{\GL_N(\C)\times\C^*}(G(v,N))\to K^{\GL_N(\C)\times\C^*}(O(v,N))$
are identified with the natural homomorphisms
\begin{equation*}
  \mathbf R^{S_{[v-n]}} \to 
  \mathbf R^{S_{[v-n]}\cap S_{[v]}}, \quad
  \mathbf R^{S_{[v]}} \to 
  \mathbf R^{S_{[v-n]}\cap S_{[v]}}
\end{equation*}
respectively.

\rom{(2)} The pushforward homomorphisms
$P_{1*}\colon K^{\GL_N(\C)\times\C^*}(O(v,N))\to
K^{\GL_N(\C)\times\C^*}(G(v-n,N))$, $P_{2*}\colon
K^{\GL_N(\C)\times\C^*}(O(v,N))\to K^{\GL_N(\C)\times\C^*}(G(v,N))$
are identified with the natural homomorphisms
\begin{equation*}
\begin{split}
  \mathbf R^{S_{[v-n]}\cap S_{[v]}} \ni f & \mapsto
  {\mathfrak S}_{[v]}^{[v-n]}
  \left(f \prod_{t=v+1}^{N}\prod_{k=v-n+1}^v
    ( 1 - \frac{x_k}{x_t})^{-1}\right), \\
  \mathbf R^{S_{[v-n]}\cap S_{[v]}} \ni f & \mapsto
  {\mathfrak S}_{[v-n]}^{[v]}
  \left(f \prod_{u=1}^{v-1}\prod_{k=v-n+1}^v
    ( 1 - \frac{x_u}{x_k})^{-1}\right)
\end{split}
\end{equation*}
respectively.
\rom(The right hand sides are a priori in ${\mathcal R}$, but they are
in fact in $\mathbf R$.\rom)
\end{Lemma}

Using above lemmas, we can write down the operators $x^+(z)$ explicitly:
\begin{equation*}
\begin{split}
   x^+(z) f 
   &=  (-1)^{N-v}
   {\mathfrak S}_{[v]}^{[v-1]}
    \left(f \sum_{r=-\infty}^\infty\left(\frac{x_v}{z}\right)^r
    x_v^{-v} 
    \prod_{t\in [v+1,N]} q x_t^{-1}
    \left( 1 - \frac{x_v}{x_t}\right)^{-1} 
    \left( 1 - q^{-2} \frac{x_t}{x_v}\right) \right) \\
&=
    \sum_{k\in [v,N]}
          f(x_{\tau_k^-[v-1]})
    \sum_{r=-\infty}^\infty\left(\frac{x_k}{z}\right)^r x_k^{-N}
    \prod_{t\in [v, N]\setminus \{k\}}
    \frac{q x_k - q^{-1} x_t}{x_k - x_t},
\end{split}
\end{equation*}
for $f\in \mathbf R^{S_{[v]}}$.
Similarly,
\begin{equation*}
\begin{split}
   x^-(w) g 
 &= (-1)^{1-v}
 {\mathfrak S}_{[v-1]}^{[v]}
  \left(g \sum_{s=-\infty}^\infty\left(\frac{x_v}{w}\right)^s
    {x_v^{N-v+1}} 
    \prod_{u\in [1,v-1]}  q x_u
     \left( 1 - \frac{x_u}{x_v}\right)^{-1} 
     \left( 1 - q^{-2} \frac{x_v}{x_u}\right) \right) \\
&= \sum_{l\in [1,v]}
     g(x_{\tau_l^+[v]})
  \sum_{s=-\infty}^\infty\left(\frac{x_l}{w}\right)^s x_l^{N}
    \prod_{u\in [1,v]\setminus\{l\}}
    \frac{q^{-1} x_l - q x_u}{x_l - x_u}
\end{split}
\end{equation*}
for $g\in\mathbf R^{S_{[v-1]}}$.

Let us compare $x^+(z)x^+(w)$ with $x^+(w)x^+(z)$ in the component
$K^{\GL_N(\C)\times\C^*}(T^*G(v,N)) \to
K^{\GL_N(\C)\times\C^*}(T^*G(v-2,N))$:
\begin{equation}\label{eq:x^{(2)}}
\begin{split}
   x^+(z)x^+(w) f & =
   - \sum_{l\in [v-1,N]}\sum_{k\in [v-1,N]\setminus\{l\}}
          f(x_{\tau_k^-\tau_l^-[v-2]})
    \sum_{r=-\infty}^\infty\left(\frac{x_k}{w}\right)^r 
    \sum_{s=-\infty}^\infty\left(\frac{x_l}{z}\right)^s 
    x_k^{-N}x_l^{-N}
\\
    &\qquad\qquad \times
    \prod_{t\in [v-1, N]\setminus \{k,l\}}
    \frac{q x_k - q^{-1} x_t}{x_k - x_t}   
    \prod_{u\in [v-1, N]\setminus \{l\}}
    \frac{q x_l - q^{-1} x_u}{x_l - x_u},   \\
   x^+(w)x^+(z) f & =
   - \sum_{k\in [v-1,N]}\sum_{l\in [v-1,N]\setminus\{k\}}
          f(x_{\tau_l^-\tau_k^-[v-2]})
    \sum_{r=-\infty}^\infty\left(\frac{x_k}{w}\right)^r 
    \sum_{s=-\infty}^\infty\left(\frac{x_l}{z}\right)^s
    x_k^{-N}x_l^{-N}
\\
    &\qquad\qquad\times
    \prod_{t\in [v-1, N]\setminus \{k\}}
    \frac{q x_k - q^{-1} x_t}{x_k - x_t}   
    \prod_{u\in [v-1, N]\setminus \{k,l\}}
    \frac{q x_l - q^{-1} x_u}{x_l - x_u}.
\end{split}
\end{equation}
Hence we have
\begin{equation*}
  x^+(w) x^+(z) = \frac{q w - q^{-1}z}{q^{-1}w - qz} x^+(z) x^+(w).
\end{equation*}
The relation~\eqref{eq:relexE2} for $x^-(z)$, $x^-(w)$ can be proved
in the same way.

Let us compare $x^-(w)x^+(z)$ and $x^+(z)x^-(w)$ in the component
$K^{\GL_N(\C)\times\C^*}(T^*G(v,N)) \to
K^{\GL_N(\C)\times\C^*}(T^*G(v,N))$
\begin{equation*}
\begin{split}
  x^-(w)x^+(z) f
  =& \sum_{l\in [1,v]} \sum_{k\in [v+1,N]\cup \{l\}}
  f(x_{\tau_k^-\tau_l^+[v]})
  \sum_{r,s=-\infty}^\infty \left(\frac{x_k}{z}\right)^r
                            \left(\frac{x_l}{w}\right)^s
  \left(\frac{x_l}{x_k}\right)^N
\\
  &\qquad\qquad\times
  \prod_{t\in ([v+1,N]\cup \{l\}) \setminus \{k\}}
     \frac{q x_k - q^{-1} x_t}{x_k - x_t}
  \prod_{u\in [1,v]\setminus\{l\}}
     \frac{q^{-1} x_l - q x_u}{x_l - x_u}, \\
  x^+(z)x^-(w) f
  =& \sum_{k\in [v+1,N]}\sum_{l\in [1,v]\cup\{k\}} 
  f(x_{\tau_l^+\tau_k^-[v]})
  \sum_{r,s=-\infty}^\infty \left(\frac{x_k}{z}\right)^r
                            \left(\frac{x_l}{w}\right)^s
  \left(\frac{x_l}{x_k}\right)^N
\\
  &\qquad\qquad\times
  \prod_{t\in [v+1,N] \setminus \{k\}}
     \frac{q x_k - q^{-1} x_t}{x_k - x_t}
  \prod_{u\in ([1,v]\cup\{k\}) \setminus\{l\}}
    \frac{q^{-1} x_l - q x_u}{x_l - x_u}.
\end{split}
\end{equation*}
Terms with $k\neq l$ cancel out for $x^+(z)x^-(w)f$ and $x^-(w)x^+(z)f$. Thus
\begin{equation*}
\begin{split}
  \left[ x^+(z), x^-(w) \right]
  = \;
   & \sum_{k\in [v+1,N]} 
     \sum_{r,s=-\infty}^\infty \left(\frac{x_k}{z}\right)^r
                               \left(\frac{x_k}{w}\right)^s
\prod_{t\in [v+1,N] \setminus \{k\}}
     \frac{q x_k - q^{-1} x_t}{x_k - x_t}
  \prod_{u\in [1,v]}
    \frac{q^{-1} x_k - q x_u}{x_k - x_u}\\
& - \sum_{l\in [1,v]} 
     \sum_{r,s=-\infty}^\infty \left(\frac{x_l}{z}\right)^r
                               \left(\frac{x_l}{w}\right)^s
\prod_{t\in [v+1,N]}
     \frac{q x_l - q^{-1} x_t}{x_l - x_t}
  \prod_{u\in [1,v]\setminus \{l\}}
    \frac{q^{-1} x_l - q x_u}{x_l - x_u}.
\end{split}
\end{equation*}

Let
\begin{equation*}
\begin{split}
  A(x) & \defeq \prod_{u\in [1,v]} (x - x_u)\;
                \prod_{t\in [v+1,N]} (x - x_t), \\
  B(x) & \defeq \prod_{u\in [1,v]} (q^{-1} x - q x_u)\;
                \prod_{t\in [v+1,N]} (q x - q^{-1} x_t).
\end{split}
\end{equation*}
Then we have
\begin{equation*}
  \left[ x^+(z), x^-(w)\right]
  = \frac{1}{q - q^{-1}}\sum_{m\in [1,N]}
    \sum_{r,s=-\infty}^\infty \left(\frac{x_m}{z}\right)^r
                               \left(\frac{x_m}{w}\right)^s
    \frac{x_m^{-1}B(x_m)}{A'(x_m)}.
\end{equation*}
Applying the residue theorem to
\(
   \frac{1}{q - q^{-1}}  
   \sum_{r,s=-\infty}^\infty \left(\frac{x}{z}\right)^r
                             \left(\frac{x}{w}\right)^s
   \frac{B(x)}{A(x)} \frac{dx}{x},
\)
we get
\begin{equation*}
    \left[ x^+(z), x^-(w)\right]
  = \frac{1}{q - q^{-1}}\left(
    \sum_{r=-\infty}^\infty \left(\frac{z}{w}\right)^r
    \left(\frac{B(z)}{A(z)}\right)^+
    - \sum_{r=-\infty}^\infty \left(\frac{z}{w}\right)^r
    \left(\frac{B(z)}{A(z)}\right)^-\right),
\end{equation*}
where
\(
\left(\frac{B(z)}{A(z)}\right)^\pm \in \C[[z^{\mp}]]
\)
denotes the Laurent expansion of $\frac{B(z)}{A(z)}$ at $z = \infty$ and
$0$ respectively. Since
\begin{equation*}
  \frac{B(z)}{A(z)} = q^{N-2v}
  \frac{\Wedge_{-1/(qz)}C^\bullet(v,N)}{\Wedge_{-q/z}C^\bullet(v,N)}
\end{equation*}
by \eqref{eq:genP}, we have completed the proof of \thmref{thm:hommain}.

\section{Integral structure}\label{sec:int}
In this section, we compare $\Uli$ with $K^{G_\bw\times\C^*}(\Zw)$. In 
the case of the affine Hecke algebra, the equivariant $K$-group of the 
Steinberg variety is isomorphic to the integral form of the affine
Hecke algebra (see \cite[7.2.5]{Gi-book}).
We shall prove a weaker form of the corresponding result for quiver
varieties in this section.

\subsection{rank $1$ case}
We first consider the case when the graph is of type $A_1$. We drop
the subscript $k$. We use the notation in \subsecref{subsec:rank1}.
We also consider $\omega(\Pa^{(n)}(v,N))$ where
$\Pa^{(n)}(v,N)$ is as in \eqref{eq:Pa^{(n)}} and
$\omega\colon \M(v-n,N)\times \M(v,N)\to \M(v,N)\times\M(v-n,N)$ is
the exchange of factors. We identify its equivariant $K$-group with
$\mathbf R^{S_{[v-n]}\cap S_{[v]}}$ as in \subsecref{subsec:rank1}.
In particular, the vector bundle $V^1/V^2$ is identified with
$q(x_{v-n+1}+\cdots+x_v)$.

\begin{Lemma}\label{lem:rank1_Grass}
\rom{(1)} Let $p_1 < \dots < p_s$ be an increasing sequence of integers
and let $n_1, \dots, n_s$ be a sequence of positive integers such that
$\sum n_i = n$.
Let $\lambda$ be the partition 
\[
   \left((p_2-p_1)^{n_2} \cdots (p_s-p_1)^{n_s}\right).
\]
Then for $g \in K^{\GL_N(\C)\times\C^*}(T^*G(v-n,N))$, we have
\begin{equation*}
\begin{split}
  & f_{p_1}^{(n_1)} f_{p_2}^{(n_2)} \cdots f_{p_s}^{(n_s)} g \\
  =\; & \pm q^L 
   \sum_{\{l_1,\dots,l_n\}}
    g(x_{\tau_{l_1}^+\cdots\tau_{l_n}^+[v]})
      {(x_{l_1}\cdots x_{l_n})}^{N+p_1}
      P_\lambda(x_{l_1},\dots, x_{l_n};q^2)
      \prod_{\substack{i=1,\dots, n\\
                      u\in [1,v]\setminus \{l_1,\dots, l_n\}}}
       \frac{q^{-1} x_{l_i} - q x_u}{x_{l_i} - x_u}
\end{split}
\end{equation*}
for some $L\in\Z$. 
Here $P_\lambda$ is the Hall-Littlewood polynomial \rom(see
\cite[III(2.1)]{Mac}\rom), and the summation runs over the set of
unordered $n$-tuples $\{ l_1,\dots,l_n\} \subset [1,v]$ such that
$l_i\neq l_j$ for $i\neq j$.

\rom{(2)} Let us consider a tensor product $T(V^1/V^2)$ of exterior
products of the bundle $V^1/V^2$ and its dual over
$\omega(\Pa^{(n)}(v,N))$, and denote by 
\(
   T(x_{v-n+1}, \dots, x_v) \in 
   \Z[x_{v-n+1}^{\pm},\dots, x_v^\pm]^{S_n}\subset 
   \mathbf R^{S_{[v-n]}\cap S_{[v]}}
\)
the corresponding element in the equivariant $K$-group.
Then for $g\in \mathbf R^{S_{[v-n]}}$, we have the following formula
\begin{equation*}
\begin{split}
    & \left[T(V^1/V^2)\otimes
        \det C^{\prime\prime\bullet}(v-n,N)^{\otimes-n}\right] \ast g \\
    =\; & 
    \pm \sum_{\{l_1,\dots,l_n\}}
    g(x_{\tau_{l_1}^+\cdots\tau_{l_n}^+[v]})\,
    (x_{l_1}\dots x_{l_n})^{v-n}\,
     T(x_{l_1}, \dots, x_{l_n})
      \prod_{\substack{i=1,\dots, n\\
                      u\in [1,v]\setminus \{l_1,\dots, l_n\}}}
        \frac{q^{-1} x_{l_i} - q x_u}{x_{l_i} - x_u},
\end{split}
\end{equation*}
where the summation runs over the set of unordered $n$-tuples $\{
l_1,\dots,l_n\} \subset [1,v]$ such that $l_i\neq l_j$ for $i\neq j$.
\end{Lemma}

\begin{proof}
(1) Generalizing \eqref{eq:x^{(2)}}, we have the following formula
for $f_{r_1}f_{r_2}\dots f_{r_n}\colon 
K^{\GL_N(\C)\times\C^*}(G(v-n,N)) \to
K^{\GL_N(\C)\times\C^*}(G(v,N))$:
\begin{equation*}
\begin{split}
   f_{r_1}f_{r_2}\dots f_{r_n} g =&
   \pm \sum_{(l_1, \dots, l_n)}
         g(x_{\tau_{l_1}^+\cdots\tau_{l_n}^+[v]})
         x_{l_1}^{r_1}\cdots x_{l_n}^{r_n}
    {(x_{l_1}\cdots x_{l_n})}^{N}
\\
   &\qquad\qquad\times  
  \prod_{\substack{i=1,\dots, n\\
                   t\in [1,v]\setminus \{l_1,\dots, l_n\}}}
    \frac{q^{-1} x_{l_i} - q x_u}{x_{l_i} - x_u}
    \prod_{i > j}
    \frac{q^{-1} x_{l_i} - q x_{l_j}}{x_{l_i} - x_{l_j}},
\end{split}
\end{equation*}
where the summation runs over the set of {\it ordered\/} $n$-tuples
$(l_1, \dots, l_n)$ such that $l_i\in [1,v]$, $l_i\neq l_j$ for $i\neq
j$.


Choose $r_1 \le r_2 \le \dots \le r_n$ so that
\begin{equation*}
   (r_1,r_2,\dots,r_n) = (
     \underbrace{p_1,\dots,p_1}_{\text{$n_1$ times}}, 
     \underbrace{p_2,\dots,p_2}_{\text{$n_2$ times}}, \dots).
\end{equation*}
Consider the following term appeared in the above formula:
\begin{equation*}
\begin{split}
   & \sum_{\sigma\in S_n} 
    x_{l_{\sigma(1)}}^{r_1}\cdots x_{l_{\sigma(n)}}^{r_n}
    \prod_{i > j} 
    \frac{q^{-1} x_{l_{\sigma(i)}} - q x_{l_{\sigma(j)}}}
    {x_{l_{\sigma(i)}} - x_{l_{\sigma(j)}}} \\
  =\; & (x_{l_1}\cdots x_{l_n})^{r_1}
   \sum_{\sigma\in S_n} 
    x_{l_{\sigma(1)}}^0x_{l_{\sigma(2)}}^{r_2-r_1}\cdots 
     x_{l_{\sigma(n)}}^{r_n-r_1}
    \prod_{i > j} 
    \frac{q^{-1} x_{l_{\sigma(i)}} - q x_{l_{\sigma(j)}}}
    {x_{l_{\sigma(i)}} - x_{l_{\sigma(j)}}}.
\end{split}
\end{equation*}
By \cite[III(2.1)]{Mac} it is equal to
\begin{equation*}
   (x_{l_1}\cdots x_{l_n})^{r_1} q^{-n(n-1)/2} v_\lambda(q^2)
   P_\lambda(x_{l_1},\dots, x_{l_n};q^2),
\end{equation*}
where $P_\lambda$ is the Hall-Littlewood polynomial and
\begin{equation*}
   v_\lambda(q^2) = q^{n_1(n_1-1)/2} [n_1]_q ! 
   q^{n_2(n_2-1)/2} [n_2]_q ! \cdots q^{n_s(n_s-1)/2} [n_s]_q !\ .
\end{equation*}
Thus we have the assertion.

(2) By Lemmas~\ref{lem:Vass2}, \ref{lem:Vass3} we have
\makeatletter\tagsleft@false
\begin{align}
   & \left[T(V^1/V^2)\otimes
           \det C^{\prime\prime\bullet}(v-n,N)^{\otimes-n}\right] \ast g
      \notag\\
  =\; & {\mathfrak S}_{[v-n]}^{[v]}
    \left(g\, T(x_{v-n+1},\dots, x_v)
    \prod_{u\in [1,v-n]} (q x_u)^n \prod_{l\in [v-n+1,v]}
     \left( 1 - \frac{x_u}{x_l}\right)^{-1} 
     \left( 1 - q^{-2} \frac{x_l}{x_u}\right) \right) \notag\\
  =\; & \pm \sum_{\{l_1,\dots, l_n\}}
     g(x_{\tau_{l_1}^+\cdots\tau_{l_n}^+[v]})\,
     (x_{l_1}\dots x_{l_n})^{v-n}\,
     T(x_{l_1}, \dots, x_{l_n})
     \prod_{\substack{i=1,\dots, n\\
                      u\in [1,v]\setminus \{l_1,\dots, l_n\}}}
       \frac{q^{-1} x_{l_i} - q x_u}{x_{l_i} - x_u}. \tag*{\qed}
\end{align}
\makeatother
\renewcommand{\qed}{}
\end{proof}

\subsection{}
Let $K^{G_\bw\times\C^*}(\Zw)/\!\operatorname{torsion}$
be
\begin{equation*}
   \operatorname{Image}\left(K^{G_\bw\times\C^*}(\Zw)
     \to
   K^{G_\bw\times\C^*}(\Zw)\otimes_{\Z[q,q^{-1}]}\Q(q) 
   \right).
\end{equation*}
(It seems reasonable to conjecture that
$K^{G_\bw\times\C^*}(Z(\bv^1,\bv^2;\bw))$ is free over
$R(G_\bw\times\C^*)$ since it is true for type $A_n$.
But I do not know how to prove it in general.)

\begin{Theorem}\label{thm:homint}
The homomorphism in \thmref{thm:hommain} induces
a homomorphism $\Uli\to
K^{G_\bw\times\C^*}(\Zw)/\!\operatorname{torsion}$.
\end{Theorem}

\begin{Remark}
The homomorphism is neither injective nor surjective. It is likely
that there exists a surjective homomorphism from a modification of
$\Uli$ to
$K^{G_\bw\times\C^*}(Z^{\operatorname{reg}}(\bw))/\!\operatorname{torsion}$
for a suitable subset $Z^{\operatorname{reg}}(\bw)$ of $\Zw$,
as in \cite[9.5, 10.15]{Na-alg}.
\end{Remark}

\begin{proof}[Proof of \thmref{thm:homint}]
It is enough to check that
$e_{k,r}^{(n)}$, $f_{k,r}^{(n)}$, $q^h$
and the coefficients of $p_k^\pm(z)$ are mapped to
$K^{G_\bw\times\C^*}(\Zw)$.
For $q^h$ and the coefficients of $p_k^\pm(z)$, the assertion is clear 
from the definition.

For $e_{k,r}^{(n)}$ and $f_{k,r}^{(n)}$, we can use
a reduction to rank $1$ case as in \secref{sec:rel2}.
Namely, it is enough to show the assertion when the graph is of type
$A_1$.

Now if the graph is of type $A_1$, \lemref{lem:rank1_Grass}
together with \lemref{lem:Vass1}
shows that 
$f_{r}^{(n)}$ is represented by a certain line bundle over
$\omega(\Pa^{(n)}(v,N))$ extended to $Z(v,v-n;N)$ by $0$.
We leave the proof for $e_r^{(n)}$ as an exercise. The only thing we
need is to write down an analogue of \lemref{lem:rank1_Grass} for
$e_r^{(n)}$. It is straightforward.
\end{proof}

\subsection{The module $K^{G_\bw\times\C^*}(\Law)$}

By \thmref{thm:homint} and \thmref{thm:freeness},
$K^{G_\bw\times\C^*}(\Law)$ is a $\Uli$-module. We show that it is an
{\it l\/}-highest weight module in this subsection.

\begin{Lemma}\label{lem:e's}
Let $\Pa_k^{(n)}(\bv,\bw)$ be as in \eqref{eq:Pa^{(n)}} and 
$\omega\colon \M(\bv-n\alpha_k,\bw)\times \M(\bv,\bw)
\to \M(\bv,\bw)\times \M(\bv-n\alpha_k,\bw)$ denote the exchange of
factors.
Let $T(V^1_k/V^2_k)$ be a tensor product of exterior products of the
vector bundle $V^1_k/V^2_k$ and its dual over
$\omega(\Pa_k^{(n)}(\bv,\bw))$.
Let us consider it as an element of
$K^{G_\bw\times\C^*}(Z(\bv,\bv-n\alpha_k;\bw))$.
Then
\(\left[
  T(V_k^1/V_k^2)\otimes
  \det C_k^{\prime\prime\bullet}(\bv-n\alpha_k,\bw)^{\otimes -n}
\right]\)
can be written as a linear combination
\rom(over $\Z[q,q^{-1}]$\rom) of elements of the form
\begin{equation*}
   f_{k,p_1}^{(n_1)} f_{k,p_2}^{(n_2)} \cdots
   f_{k,p_s}^{(n_s)} \ast [\shfO_{\Delta\M(\bv-n\alpha_k,\bw)}]
   \qquad(n_1+n_2+\cdots+n_s = n,\quad \text{$p_i$ distinct}),
\end{equation*}
where $\Delta\M(\bv-n\alpha_k,\bw)$ is the diagonal in
$\M(\bv-n\alpha_k,\bw)\times\M(\bv-n\alpha_k,\bw)$.
\end{Lemma}

\begin{proof}
As in \secref{sec:rel2}, we may assume that the graph is of type
$A_1$. Now \lemref{lem:rank1_Grass} together with the fact that
Hall-Littlewood polynomials form a basis of symmetric polynomials
implies the assertion.
\end{proof}

\begin{Proposition}\label{prop:Lw}
Let $[0]\in K^{G_\bw\times\C^*}(\La(0,\bw))$ be the class represented
by the structure sheaf of $\M(0,\bw)\cong \La(0,\bw) =
\operatorname{point}$. Then
\begin{equation*}
 K^{G_\bw\times\C^*}(\Law)
 = \Uli^-\ast \left(R(G_\bw\times\C^*)[0]\right).
\end{equation*}
\end{Proposition}

\begin{proof}
The following proof is an adaptation of proof of \cite[10.2]{Na-alg},
which was inspired by \cite[3.6]{Lu-affine} in turn.

We need the following notation:
\begin{equation*}
\begin{gathered}
  \La_{k;n}(\bv,\bw) \defeq \La(\bv,\bw)\cap \M_{k;n}(\bv,\bw), \qquad
  \La_{k;\le n}(\bv,\bw) \defeq \La(\bv,\bw)\cap \M_{k;\le n}(\bv,\bw), \\
  \La_{k;\ge n}(\bv,\bw) \defeq \La(\bv,\bw)\cap \M_{k;\ge n}(\bv,\bw).
\end{gathered}
\end{equation*}

We prove $K^{G_\bw\times\C^*}(\La(\bv,\bw))\subset
\Uli^-\ast\left(R(G_\bw\times\C^*)[0]\right)$ by induction on the
dimension vector $\bv$.
When $\bv = 0$, the result is trivial since
$K^{G_\bw\times\C^*}(\operatorname{point}) =
R(G_\bw\times\C^*)$.
Consider $\La(\bv,\bw)$ and suppose that
\begin{equation}
\label{eq:indass}
\text{if $\bv - \bv'\in \bigoplus \Z_{\ge 0}\alpha_k \setminus \{0\}$,
then $K^{G_\bw\times\C^*}(\La(\bv',\bw))\subset
\Uli^-\ast\left(R(G_\bw\times\C^*)[0]\right)$.}
\end{equation}

Take $[E]\in K^{G_\bw\times\C^*}(\La(\bv,\bw))$. We want to show
$[E]\in \Uli^-\ast\left(R(G_\bw\times\C^*)[0]\right)$.
We may assume that the
support of $E$ is contained in an irreducible component of
$\La(\bv,\bw)$ without loss of generality. In fact, suppose that
$\operatorname{Supp}E\subset X\cup Y$ such that $X$ is an irreducible
component.
Since $Y$ is a closed subvariety of $X\cup Y$ and since $X\cap Y$ is a
closed subvariety of $X$, we have the diagram
\begin{equation*}
\begin{CD}
   K^{G_\bw\times\C^*}(Y) @>{i_*}>>
   K^{G_\bw\times\C^*}(X\cup Y) @>{j^*}>>
   K^{G_\bw\times\C^*}(X\setminus Y) @>>> 0 \\
   @AAA @A{i''_*}AA @| \\
   K^{G_\bw\times\C^*}(X\cap Y) @>{i'_*}>>
   K^{G_\bw\times\C^*}(X) @>{j^{\prime*}}>>
   K^{G_\bw\times\C^*}(X\setminus Y) @>>> 0,
\end{CD}
\end{equation*}
where the first and the second row are exact by \eqref{eq:exact}.
Thus there exists $[E']\in K^{G_\bw\times\C^*}(X)$ such that
$j^{\prime*}[E'] = j^*[E]$. Then $j^*([E] - i''_*[E']) = 0$, therefore
there exists $E''\in K^{G_\bw\times\C^*}(Y)$ such that
$[E] = i_*[E''] + i''_*[E']$. By the induction on the number of
irreducible components in the support, we may assume that the support
of $E$ is contained in an irreducible component, which is denoted by
$X_E$.

Let us consider $\varepsilon_k$ defined in \eqref{eq:epsilon_k}. If
$\varepsilon_k(X_E) = 0$ for all $k\in I$, $X_E$ must be $\La(0,\bw)$
by \lemref{lem:surj}. We have nothing to prove in this case.  Thus
there exists $k$ such that $\varepsilon_k(X_E) > 0$. Set $n =
\varepsilon_k(X_E)$.
By the descending induction on $\varepsilon_k$, we may assume that
\begin{equation}
\label{eq:indass2}
\text{if
$\operatorname{Supp}(E')\subset \La_{k;\ge n+1}(\bv,\bw)$, then
$[E']\in \Uli^-\ast\left(R(G_\bw\times\C^*)[0]\right)$.}
\end{equation}

Since $\La_{k;\le n}(\bv,\bw)$ is an open
subvariety of $\La(\bv,\bw)$, we have an exact sequence
\begin{equation*}
   K^{G_\bw\times\C^*}(\La_{k;\ge n+1}(\bv,\bw))
    \xrightarrow{a_*}
   K^{G_\bw\times\C^*}(\La(\bv,\bw))
    \xrightarrow{b^*}
   K^{G_\bw\times\C^*}(\La_{k;\le n}(\bv,\bw))
    \to 0
\end{equation*}
by \eqref{eq:exact}.
Consider $b^*[E]$.
By \eqref{eq:indass2}, it is enough to show that
\begin{equation}
\label{eq:to_show}
\text{there exists $[\widetilde E]\in\Uli^-\ast[0]$ such that
$b^*[E] = b^*[\widetilde E]$.}
\end{equation}

Since $X_E\cap \La_{k;\le n}(\bv,\bw)\subset \La_{k;n}(\bv,\bw)$, the
support of $b^*(E)$ is contained in $\La_{k;n}(\bv,\bw)$.
We have a map
\begin{equation*}
   P\colon \La_{k;n}(\bv,\bw) \to
     \La_{k;0}(\bv-n\alpha_k,\bw)
\end{equation*}
which is the restriction of the map \eqref{taut:induct}.
Recall that this map is a Grassmann bundle (see
\propref{prop:taut_fib}).
Let us denote its tautological bundle by $S$.
Then $b^*[E]$ can be written as a linear combination of elements of
the form
\begin{equation*}
   [T(S)]\otimes P^* [E_0],
\end{equation*}
where $T(S)$ is a tensor product of exterior powers of the tautological
bundle $S$ and 
$E_0\in K^{G_\bw\times\C^*}(\La_{k;0}(\bv-n\alpha_k,\bw))$.
Since the homomorphism
\[
b^{\prime*}\colon K^{G_\bw\times\C^*}(\La(\bv-n\alpha_k,\bw))\to
K^{G_\bw\times\C^*}(\La_{k;0}(\bv-n\alpha_k,\bw))
\]
is surjective by \eqref{eq:exact}, there exists
$[E_1]\in K^{G_\bw\times\C^*}(\La(\bv-n\alpha_k,\bw))$ such
that $b^{\prime*}[E_1] = [E_0]$.

Consider $\Pa_k^{(n)\prime}(\bv-n\alpha_k,\bw)\cap(
\La_{k;\le n}(\bv,\bw)\times\La(\bv-n\alpha_k,\bw))$.
By \propref{prop:taut_fib}, it is isomorphic to $\La_{k;n}(\bv,\bw)$
and the map $P$ can be identified with the projection to the second
factor. Moreover, the tautological bundle $S$ is identified with the
restriction of the natural vector bundle $V_k^1/V_k^2$.
Hence we have
\begin{equation*}
     [T(S)]\otimes P^* [E_0]
   = b^* \left( T(V_k^1/V_k^2) \ast [E_1]\right),
\end{equation*}
where $T(V_k^1/V_k^2)$ is considered as an element of
$K^{G_\bw\times\C^*}(Z(\bv,\bv-n\alpha_k;\bw))$.
By \lemref{lem:e's}, $T(V_k^1/V_k^2)$ can be written as a linear
combination of elements
\begin{equation*}
   f_{k,p_1}^{(n_1)}f_{k,p_2}^{(n_2)}\dots f_{k,p_s}^{(n_s)}\ast
   [p_2^*\det C_k^{\prime\prime\bullet}(\bv-n\alpha_k,\bw)^{\otimes n}],
\end{equation*}
where $p_2\colon \Delta\M(\bv-n\alpha_k,\bw)\to \M(\bv-n\alpha_k,\bw)$ 
is the projection.
By \eqref{eq:indass},
\[
   [\det C_k^{\prime\prime\bullet}(\bv-n\alpha_k,\bw)^{\otimes n}
\otimes E_1]\in \Uli^-\ast\left(R(G_\bw\times\C^*)[0]\right).
\]
Hence $T(V_k^1/V_k^2) \ast [E_1]\in
\Uli^-\ast\left(R(G_\bw\times\C^*)[0]\right)$.
Thus we have shown \eqref{eq:to_show}.
\end{proof}

\section{Standard modules}\label{sec:std}

In this section, we start the study of the representation theory of
$\Uli$ using $K^{G_\bw\times\C^*}(\Zw)$ and
the homomorphism in (\ref{thm:homint}). We shall define certain
modules called {\it standard modules}, and study their properties.
Results in this section holds even if $\varepsilon$ is a root of unity.

Note that $R(G_\bw\times\C^*)$ is contained in the center of
$K^{G_\bw\times\C^*}(\Zw)$ by $R(G_\bw\times\C^*)\ni \rho \mapsto
\rho\otimes\sum_\bv [\shfO_{\Delta\M(\bv,\bw)}]$. Hence a
$K^{G_\bw\times\C^*}(\Zw)/\!\operatorname{torsion}$-module $M$ (over
$\C$), which is {\it l\/}-integrable as a $\Uli$-module, decomposes as
$M = \bigoplus M_\chi$ where $\chi$ is a homomorphism from
$R(G_\bw\times\C^*)$ to $\C$ and $M_\chi$ is the corresponding
simultaneous generalized eigenspace, i.e., some powers of the kernel
of $\chi$ acts as $0$ on $M_\chi$.
Such a homomorphism $\chi$ is given by the evaluation of the character
at a semisimple element $a = (s,\varepsilon)$ in
$G_\bw\times\C^*$. (This gives us a bijection between homomorphisms and
semisimple elements.)

What is the meaning of the choice of $a = (s,\varepsilon)$ when we
consider $M$ as a $\Uli$-module~?
The role of $\varepsilon$ is clear. It is a specialization $q \to
\varepsilon$, and we get $\Ule$-modules. It will becomes clear later
that $s$ corresponds to the Drinfel'd polynomials by
\begin{equation*}
   P_k(u) = \text{(a normalization of) the characteristic polynomial
   of $s_k$}.
\end{equation*}

\subsection{Fixed data}

Let $a = (s,\varepsilon)$ be a semisimple element in $G_\bw\times
\C^*$ and $A$ be the Zariski closure of $\{ a^n \mid n\in\Z\}$.
Let $\chi_a\colon R(A) \to \C$ be the homomorphism given by the
evaluation at $a$. Considering $\C$ as an $R(A)$-module by this
evaluation homomorphism, we denote it by $\C_a$.
Via the homomorphism $R(G_\bw\times\C^*)\to R(A)$, we consider $\C_a$
also as an $R(G_\bw\times\C^*)$-module.
We consider $R(A)$ as a $\Z[q,q^{-1}]$-algebra, where
$R(G_\bw\times\C^*)$ is a $\Z[q,q^{-1}]$-algebra as in
\subsecref{subsec:homdef}.

Let $\M(\bw)^A$, $\M_0(\infty,\bw)^A$ be the fixed point subvarieties
of $\M(\bw)$, $\M_0(\infty,\bw)$ respectively.
Let us take a point $x\in \M_0(\infty,\bw)^A$ which is regular, i.e.,
$x\in \M_0^{\operatorname{reg}}(\bv^0,\bw)$ for some $\bv^0$.

These data $x$, $a$ will be fixed through this section.

\subsection{Definition}
As in \eqref{eq:M_x}, let $\M(\bv,\bw)_x$ denote the inverse image of
$x\in \M_0^{\operatorname{reg}}(\bv^0,\bw)\subset
\M_0(\infty,\bw)$ under the map
$\pi\colon \M(\bv,\bw)\to \M_0(\bv,\bw)\hookrightarrow
\M_0(\infty,\bw)$.
It is invariant under the $A$-action.
Let $\Mw_x$ be $\bigsqcup_\bv \M(\bv,\bw)_x$. We set
\begin{equation*}
    K^A(\Mw_x) \defeq \bigoplus_\bv K^A(\M(\bv,\bw)_x)
\end{equation*}
as convention.

Let $K^{A}(\Zw)/\!\operatorname{torsion}$ be
\[
   \operatorname{Image}\left(K^{A}(\Zw)
     \to
   K^{A}(\Zw)\otimes_{\Z[q,q^{-1}]}\Q(q) 
   \right).
\]
Let
\begin{equation}\label{eq:standmod}
  M_{x,a} \defeq K^A(\Mw_x)\otimes_{R(A)}\C_a.
\end{equation}
By \thmref{thm:freeness} together with \thmref{thm:slice},
$K^A(\M(\bv,\bw)_x)$ is a free $R(A)$-module.
Thus the $K^A(\Zw)$-module structure on $K^A(\Mw_x)$ descends to
a $K^A(\Zw)/\!\operatorname{torsion}$-module structure.
Hence $M_{x,a}$ is a $\Ule$-module via the composition of
\begin{equation}\label{eq:homA}
\begin{split}
   \Ule & \to
   K^{G_\bw\times\C^*}(\Zw)
    /\!\operatorname{torsion}\otimes_{R(G_\bw\times\C^*)}\C_a \\
   & \to K^{A}(\Zw)
       /\!\operatorname{torsion}\otimes_{R(A)}\C_a.
\end{split}
\end{equation}
We call $M_{x,a}$ the {\it standard module}.

It has a decomposition
$M_{x,a} = \bigoplus K^A(\M(\bv,\bw)_x)\otimes_{R(A)}\C_a$, and
each summand is a weight space:
\begin{equation}\label{eq:weightdec}
  q^h\ast v = \varepsilon^{\langle h, \bw - \bv\rangle} v
  \qquad\text{for $v\in K^A(\M(\bv,\bw)_x)\otimes_{R(A)}\C_a$}.
\end{equation}
Thus $M_{x,a}$ has the weight decomposition as a $\Ue$-module.

In the remainder of this section, we study properties of
$M_{x,a}$. The first one is the following.
\begin{Lemma}\label{lem:lint}
As a $\Ule$-module, $M_{x,a}$ is {\it l\/}-integrable.
\end{Lemma}
\begin{proof}
The assertion is proved exactly as \cite[9.3]{Na-alg}.
Note that the regularity assumption of $x$ is {\it not\/} used here.
\end{proof}

\subsection{Highest weight vector}\label{subsec:hwvector}

Recall that $\pi\colon \M(\bv^0,\bw) \to \M_0(\bv^0,\bw)$ is an
isomorphism on $\pi^{-1}(\Mreg(\bv^0,\bw))$ (\propref{prop:pi_isom}).
Under this isomorphism, we can consider $x$ as a point in
$\M(\bv^0,\bw)$.  Then $\M(\bv^0,\bw)_x$ consists of the single point
$x$, thus we have a canonical generator of $K^A(\M(\bv^0,\bw)_x)$.  We
denote it by $[x]$.

Since $x$ is fixed by $A$, the fibers $(V_k)_x$, $(W_k)_x$ of tautological 
bundles at $x$ are $A$-modules. Then the restriction of the complex
$C_k^\bullet(\bv^0,\bw)$ to $x$ can be considered as a complex of
$A$-modules. In particular, it defines an element in $R(A)$.
Let us denote it by $C_{k,x}^\bullet$.

Let us spell out $C_{k,x}^\bullet$ more explicitly.
Since $x$ is fixed by $A$, we have a homomorphism $\rho\colon A\to
G_{\bv^0}$ by \subsecref{subsec:fix_hom}. It is uniquely determined by
$x$ up to the conjugacy. Then a virtual $G_{\bv^0}\times G_\bw\times
\C^*$-module
\begin{equation*}
  q^{-1} \left( \displaystyle{\bigoplus_{l}}
     [-\langle h_k,\alpha_l\rangle]_q V_l^0 \oplus W_k\right)
\end{equation*}
can considered as a virtual $A$-module via
$\rho\times(\operatorname{inclusion})\colon A\to G_{\bv^0}\times
G_\bw\times \C^*$. Its isomorphism class is independent of $\rho$ and
coincides with $C_{k,x}^\bullet$. Note that the first and third terms
in \eqref{eq:taut_cpx} are absorbed in the term $l = k$.

\begin{Proposition}\label{prop:stdhw}
The standard module $M_{x,a}$ is an {\it l\/}-highest weight module with
{\it l\/}-highest weight 
$P_k(u) \defeq \chi_a\left(\Wedge_{-u} C_{k,x}^\bullet\right)$. Namely,
the followings hold\rom:

\rom{(1)} $P_k(u)$ is a polynomial in $u$ of degree
$\langle h_k, \bw - \bv^0\rangle$.

\rom{(2)}
\begin{gather*}
  x_k^+(z)\ast [x] = 0, \quad 
  q^h \ast [x] = \varepsilon^{\langle h, \bw - \bv^0\rangle}[x], \\
  p_k^+(z)\ast [x] = 
    P_k(1/z)[x], \quad
  p_k^-(z)\ast [x] = (-z)^{\rank C_{k,x}^\bullet}P_k(1/z)\;
  \chi_{a} \left((\det C_{k,x}^{\bullet})^*\right) [x].
\end{gather*}

\rom{(3)}
\(
 M_{x,a}
 = \Ule^-\ast [x].
\)
\end{Proposition}

\begin{proof}
(1) If we restrict the complex $C_k^\bullet(\bv^0,\bw)$ to $x$,
$\tau_k$ is surjective and $\sigma_k$ is injective by
\lemref{lem:betasurj}.
Thus $C_{k,x}^\bullet$ is represented by a genuine $A$-module, and
$\chi_a\left(\Wedge_{-u} C_{k,x}^\bullet\right)$ is a
polynomial in $u$. The degree is equal to
$\langle h_k, \bw - \bv^0\rangle$ by the definition of $C_{k,x}$.

(2) The first equation is the consequence of
$\M(\bv-\alpha^k,\bw)_x = \emptyset$, which follows from \lemref{lem:surj}.
The remaining equations follows from the definition and
\lemref{lem:diag}.

(3) The assertion is proved exactly as in \propref{prop:Lw}. Note that 
the assumption $x\in\Mreg(\bv^0,\bw)$ is used here in order to apply
\lemref{lem:surj}.
\end{proof}

\begin{Remark}
$P_k(u)$ is the Drinfel'd polynomial attached to the
simple quotient of $M_{x,a}$, which we will study later.
\end{Remark}

We give a proof of \propref{prop:DrPol} as promised:
\begin{proof}[Proof of \propref{prop:DrPol}]
It is enough to show that there exists a simple {\it l\/}-integrable {\it
l\/}-highest module with given Drinfel'd polynomials $P_k(u)$. We can
construct it as the quotient of the standard module $M_{0,a}$ by the
unique maximal proper submodule. (The uniqueness can be proved as in the
case of Verma modules.) Here the parameter $a = (s,\varepsilon)$ is
chosen so that $P_k(u) = \chi_a\left(\Wedge_{-u} q^{-1} W_k\right)$,
i.e., $P_k(u)$ is a normalization of the characteristic polynomial of
$s_k$.
\end{proof}

\subsection{Localization}\label{subsec:local}

Let $R(A)_a$ denote the localization of $R(A)$ with respect to
$\Ker\chi_a$.

Let $\Zw^A$ denote the fixed point set of $A$ on
$\Zw$, and let
\(
   i\colon\Mw^A\times\Mw^A
   \to \Mw\times\Mw
\)
be the inclusion. Note that it induces an inclusion
\(
   \Zw^A \to \Zw
\)
which we also denote by $i$.
By the concentration theorem~\cite{Thom}
\begin{equation*}
  i_* \colon K^A(\Zw^A)\otimes_{R(A)} R(A)_a
  \to K^A(\Zw)\otimes_{R(A)} R(A)_a
\end{equation*}
is an isomorphism. Let
\begin{multline*}
  i^*\colon K^A(\Zw)
  \cong K^A(\Mw\times\Mw;\Zw) \\
  \longrightarrow
  K^A(\Mw^A\times\Mw^A;\Zw^A)
  \cong 
  K^A(\Zw^A)
\end{multline*}
be the pullback with support map. Then $i^* i_*$ is given by
multiplication by $\Wedge_{-1} N^*\boxtimes\Wedge_{-1} N^*$, where $N$
is the normal bundle of $\Mw^A$ in $\Mw$. By \cite[5.11.3]{Gi-book},
$\Wedge_{-1} N^*$ becomes invertible in the
localized $K$-group. Thus $i^*$ is an isomorphism on 
the localized $K$-group.
As in \cite[5.11.10]{Gi-book}, we introduce a correction factor to $i^*$
\begin{equation*}
   r_a \defeq (1\boxtimes (\Wedge_{-1} N^*)^{-1})\circ i^*
   \colon K^A(\Zw^A)\otimes_{R(A)} R(A)_a
   \to K^A(\Zw^A)\otimes_{R(A)} R(A)_a.
\end{equation*}
Then $r_a$ is an algebra isomorphism with respect to the
convolution.

Since $A$ acts trivially on $\Zw^A$, we have
\begin{equation}\label{eq:triv}
   K^A(\Zw^A) \cong
   K(\Zw^A)\otimes R(A).
\end{equation}
Thus we have the evaluation map
\begin{equation*}
  \operatorname{ev}_a\colon 
  K^A(\Zw^A)\otimes_{R(A)} R(A)_a
   \cong 
  K(\Zw^A)\otimes R(A)_a
   \to
   K(\Zw^A)\otimes \C,
\end{equation*}
by sending $F\otimes (f/g)$ to $F\otimes (\chi_a(f)/\chi_a(g))$.

By the bivariant Riemann-Roch theorem \cite[5.11.11]{Gi-book},
\begin{equation*}
   \operatorname{RR} \defeq
   (1\boxtimes \operatorname{Td}_{\Mw^A})\cup
   \operatorname{ch}\colon 
   K(\Zw^A)\to H_*(\Zw^A,\Q)
\end{equation*}
is an algebra homomorphism with respect to the convolution. Here
$\operatorname{Td}_{\Mw^A}$ is the Todd genus
of $\Mw^A$.

Composing \eqref{eq:homA} with all these homomorphisms, we have a
homomorphism
\begin{equation}\label{eq:homloc}
   \Ule \to
   H_*(\Zw^A, \C).
\end{equation}
Note that the torsion part in \eqref{eq:homA}
disappears in the right hand side of \eqref{eq:triv} after tensoring
with $R(A)_a$.

We have similar $\C$-linear maps for $\Mw_x$:
\begin{equation}\label{eq:StdIsom}
\begin{split}
   M_{x,a} = & \,K^A(\Mw_x)\otimes_{R(A)}\C_a
   \xrightarrow[\cong]{i^*}
     K^A(\Mw_x^A)\otimes_{R(A)}\C_a \\
   \xrightarrow[\cong]{\operatorname{ev}_a}
   & K(\Mw_x^A)\otimes \C
   \xrightarrow[\cong]{\operatorname{ch}}
     H_*(\Mw_x^A,\C),
\end{split}
\end{equation}
where $i^*$ is an isomorphism by the concentration theorem~\cite{Thom}
and the invertibility of $\Wedge_{-1} N^*$ in the localized
$K$-homology group, $\operatorname{ev}_a$ is an isomorphism since $A$
acts trivially on $\Mw_x^A$, and $\operatorname{ch}$ is an isomorphism
by \thmref{thm:freeness'} and \thmref{thm:slice}.  The composition is
compatible with the $\Ule$-module structure, where $H_*(\Mw_x^A,\C)$
is a $\Ule$-module via the convolution together with
\eqref{eq:homloc}.

Recall that we have decomposition
\begin{equation*}
   \Mw^A = \bigsqcup_\rho \M(\rho),
\end{equation*}
where $\rho$ runs the set of homomorphisms $A\to G_\bv$ (with various
$\bv$) (\subsecref{subsec:fix_hom}). Let
\begin{equation*}
   \M(\rho)_x \defeq \M(\rho)\cap \Mw^A_x.
\end{equation*}
Thus we have the canonical decomposition
\begin{equation}\label{eq:genwei}
   M_{x,a} = H_*(\Mw_x^A,\C)
           = \bigoplus_\rho H_*(\M(\rho)_x,\C).
\end{equation}
Each summand $H_*(\M(\rho)_x,\C)$ in \eqref{eq:genwei} is an
{\it l\/}-weight space with respect to the $\Ule$-action in the
sense that operators $\psi_k^\pm(z)$ acts on $H_*(\M(\rho)_x,\C)$ as
scalars plus nilpotent transformations. More precisely, we have
\begin{Proposition}\label{prop:genweight}
\rom{(1)}
Let $V_k$ be the tautological vector bundle over $\M(\bv,\bw)$.
Viewing $\Wedge_{u} V_k$ as an element of $K^A(\Delta\M(\bv,\bw))[u]$,
we consider it as an operator on $M_{x,a}$. Then we have
\begin{equation*}
\begin{split}
   & H_*(\M(\rho)_x,\C) \\ =\, &
   \{ m\in M_{x,a}\mid\left(
    \Wedge_{u} V_k - \chi_a(\Wedge_{u} V_k)\operatorname{Id}\right)^N
   \ast m = 0\;  \text{for $k\in I$ and sufficiently large $N$}\},
\end{split}
\end{equation*}
where $\chi_a(\Wedge_{u} V_k)$ is the evaluation at $a$ of
$\Wedge_{u} V_k$, considered as an $A$-module via $\rho\colon A\to G_\bv$.

\rom{(2)}
Let us consider
\begin{equation*}
  C_k^\bullet(\bv,\bw) = q^{-1} \left( \displaystyle{\bigoplus_{l}}
     [-\langle h_k,\alpha_l\rangle]_q V_l \oplus W_k\right)
\end{equation*}
as a virtual $A$-module via
$\rho\times(\operatorname{inclusion})\colon A\to G_\bv\times
G_\bw\times\C^*$.
Then operators $\psi_k^\pm(z)$ act on $H_*(\M(\rho)_x,\C)$ by
\begin{equation}\label{eq:genweight}
   \varepsilon^{\rank C_k^\bullet(\bv,\bw)}
   \chi_a\left(\frac{\Wedge_{-1/qz}C_k^\bullet(\bv,\bw)}
        {\Wedge_{-q/z}C_k^\bullet(\bv,\bw)}\right)^\pm
\end{equation}
plus nilpotent transformations.
\end{Proposition}

\begin{proof}
(2) follows from (1). We show (1).

Note that $\shfO_{\Delta\M(\bv,\bw)}$ is mapped to
$\shfO_{\Delta\M(\bv,\bw)^A}$ under $r_a$. And
$\shfO_{\Delta\M(\bv,\bw)^A}$ is mapped to
the fundamental class $[\Delta\M(\bv,\bw)^A]$ under
$\operatorname{RR}$.
Combining with the projection formula~\eqref{eq:projection}, we find
that the operator $\Wedge_u V_k$ is mapped to
\begin{equation*}
   \left(\operatorname{ch}\,\circ\, \operatorname{ev}_a\,\circ\, i^*
   \Wedge_u V_k\right) \cap [\Delta \M(\bv,\bw)^A]
\end{equation*}
under the homomorphism~\eqref{eq:homloc}.
Thus as an operator on $H_*(\M(\rho)_x, \C)$, it is equal to
\begin{equation}\label{eq:genweight'}
  m \mapsto 
    (j^*\circ \operatorname{ch}\,\circ\, \operatorname{ev}_a\,\circ\, i^*
    \Wedge_u V_k)\cap m,
\end{equation}
where $j\colon \M(\rho)_x\to \M(\bv,\bw)^A$ is the inclusion.

Now, on a connected space $X$, any $\alpha\in H^*(X,\C)$ acts on
$H_*(X,\C)$ as a scalar plus nilpotent operator, where the scalar is
$H^0(X,\C) (\cong \C)$-part of $\alpha$.
In our situation, the $H^0$-part of \eqref{eq:genweight'} is given by
$\chi_a(\Wedge_u V_k)$.
(Although we do not prove $\M(\rho)_x$ is connected, the $H^0$-part is 
the same on any component.)

Furthermore, $\chi_a(\Wedge_u V_k)$ determines all eigenvalues of the
operator $a$ acting on $V_k$. Hence it determines the conjugacy class
of the homomorphism $\rho\colon A\to G_\bv$. Thus the generalized
eigenspace of $\Wedge_u V_k$ with the eigenvalue $\chi_a(\Wedge_u
V_k)$ coincides with $H_*(\M(\rho)_x,\C)$.
\end{proof}

\subsection{Frenkel-Reshetikhin's $q$-character}\label{subsec:FR}
In this subsection, we study Frenkel-Reshetikhin's $q$-character for
the standard module $M_{x,a}$. The result is a simple application of
\propref{prop:genweight}. Results in this subsection will not be use
in the rest of the paper.

We assume $\mathfrak g$ is of type $ADE$ in this subsection.

Let us recall the definition of $q$-character.
It is a map from the Grothendieck group of finite dimensional
$\Ule$-modules $M$. As we shall see later in
\subsecref{subsec:std_decomp}, standard modules $M_{0,a}$ ($x = 0$ is
fixed, $\bw$ and $a = (s,\varepsilon)\in G_\bw\times\C^*$ are moving)
give a basis of the Grothendieck group, thus it is enough to define
the $q$-character for standard modules $M_{0,a}$. We decompose $M =
M_{0,a}$ as
\begin{equation*}
   M = \bigoplus M_{\Psi^\pm},
\end{equation*}
as in \eqref{eq:gen_wt}.
Moreover, by \propref{prop:genweight}, $\Psi_k^\pm(z)$ have the form
\begin{equation}\label{eq:Q_k}
  \Psi_k^\pm(z) = \varepsilon^{\deg Q_k - \deg R_k}
  \frac{Q_k(1/\varepsilon z) R_k(\varepsilon/z)}
  {R_k(1/\varepsilon z)Q_k(\varepsilon/z)}
\end{equation}
where $Q_k(u)$, $R_k(u)$ are polynomials in $u$ with constant term
$1$.
(Compare with \cite[Proposition~1]{FR}. Note $u = 1/z$.)
Suppose
\begin{equation*}
   \frac{Q_k(u)}{R_k(u)}
   = \frac{\prod_{r} (1 - u a_{kr})}{\prod_s (1 - u b_{ks})}.
\end{equation*}
Then the $q$-character is defined by
\begin{equation*}
   \chi_q(M_{0,a}) \defeq
   \sum_{\Psi^\pm(z)}
   \dim V_{\Psi^\pm_k(z)}\,
   \prod_{k\in I}\prod_{r} Y_{k,a_{kr}} \prod_{s} Y_{k,b_{ks}}^{-1},
\end{equation*}
where $Y_{k,a_{kr}}$, $Y_{k,b_{ks}}$ are formal variables and $\chi_q$ 
takes value in $\Z[Y_{k,c}^{\pm}]_{k\in I, c\in \C^*}$.
($\chi_q$ should not be confused with $\chi_a$.)

Let
\begin{equation*}
  A_{k,a} \defeq
  Y_{k,a\varepsilon} Y_{k,a\varepsilon^{-1}}
  \prod_{h: \vin(h) = k} Y_{\vout(h), aq^{m(h)}}^{-1}.
\end{equation*}

\begin{Proposition}[\protect{cf.\
Conjecture 1 in \cite{FR}}]
Let $M_{0,a}$ be a standard module with $x = 0$.
Suppose that $P_k(u)$ in \propref{prop:stdhw} equals to
\begin{equation*}
   P_k(u) = \prod_{i=1}^{n_k} ( 1 - u a_i^{(k)})
\end{equation*}
for $k\in I$. Then the $q$-character of $M_{0,a}$ has the following
form\rom:
\begin{equation*}
   \prod_{k\in I} \prod_{i=1}^{n_k}
     Y_{k,a_i^{(k)}}\left(1+\sum M_p'\right),
\end{equation*}
where each $M_p'$ is a product of
$A_{l,c}^{-1}$ with $c\in \bigcup a_i^{(k)} \varepsilon^{\Z}$.
\end{Proposition}

\begin{proof}
By \propref{prop:genweight}, $H_*(\M(\rho)_0,\C)$ is a geneneralized
eigenspace for $\psi_k^\pm(z)$ for a homomorphism $\rho\colon A\to
G_\bv$. Thus it is enough to study the eigenvalue.
We consider $V_k$, $W_k$ as $A$-modules via
$\rho\times(\operatorname{inclusion})\colon A\to G_\bv\times
G_\bw\times\C^*$ as before. Let $V_k(\lambda)$, $W_k(\lambda)$ be
weight space as in \subsecref{subsec:fix_hom}.

By the definition of $C_k^\bullet(\bv,\bw)$ and $P_k(u)$, we have
\begin{equation*}
\begin{gathered}
  \chi_a\left(\frac{\Wedge_{-1/qz}C_k^\bullet(\bv,\bw)}
        {\Wedge_{-q/z}C_k^\bullet(\bv,\bw)}\right)
  = \chi_a\left(
     \frac{\Wedge_{-1/qz}q^{-1}W_k}{\Wedge_{-q/z}q^{-1}W_k}
     \prod_l \frac{\Wedge_{-1/qz}q^{-1}
     [-\langle h_k, \alpha_l\rangle]_q V_l}
        {\Wedge_{-q/z}q^{-1}
     [-\langle h_k, \alpha_l\rangle]_q V_l}\right), \\
   P_k(u) = \chi_a\left(\Wedge_{-u} q^{-1} W_k\right).
\end{gathered}
\end{equation*}
By \propref{prop:genweight} we have
\begin{equation*}
   \frac{Q_k(u)}{R_k(u)} = P_k(u)\; 
    \chi_a\left(
     \prod_l \Wedge_{-u}q^{-1}
     [-\langle h_k, \alpha_l\rangle]_q V_l \right),
\end{equation*}
where $Q_k(u)$, $R_k(u)$ are defined by \eqref{eq:Q_k}.

Let $\{ c^{(l)}_t \}$ be the set of eigenvalues of $a\in A$ on $V_l$
counted with multiplicities. Then we have
\begin{equation*}
   \chi_a(\Wedge_{-u} q^{-1} [-\langle h_k, \alpha_l\rangle]_q V_l)
   = 
   \begin{cases}{\displaystyle
     \left(\prod_t (1 - u c^{(k)}_t)
             (1 - u \varepsilon^{2}c^{(k)}_t)\right)^{-1}}
        & \text{if $k=l$}, \\
     {\displaystyle \prod_{h: \substack{\vin(h) = k\\ \vout(h) = l}}
     \prod_t (1 - u \varepsilon^{m(h)+1} c^{(l)}_t)}
        & \text{otherwise}.
   \end{cases}
\end{equation*}
Thus we have
\begin{equation*}
   \chi_q(M_{0,a})
  = \sum_{\rho} \dim H_*(\M(\rho),\C)
 \prod_{k\in I}\prod_{i=1}^{n_k} Y_{k,a_i^{(k)}}
 \prod_t A_{k,\varepsilon^{-1}c^{(k)}_t}^{-1}.
\end{equation*}

Note that the term for $\rho$ with $\bv = 0$ has the contribution
\begin{equation*}
  \prod_{k\in I}\prod_{i=1}^{n_k} Y_{k,a_i^{(k)}},
\end{equation*}
and any other terms are monomials of
$A_{k,\varepsilon^{-1}c^{(k)}_t}^{-1}$ which are not constant.

Moreover, we have $c^{(k)}_t\in \bigcup a_i^{(k)} \varepsilon^{\Z}$ by
\lemref{lem:Vwei}. This completes the proof.
\end{proof}

\section{Simple modules}\label{sec:simple}

The purpose of this section is to study simple modules of $\Ule$. Our
discussion relies on Ginzburg's classification of simple modules of
the convolution algebra \cite[Chapter~9]{Gi-book}. (See also
\cite{Lu-cusp}.)
He applied his classification to the affine Hecke algebra.
However, unlike the case of the affine Hecke algebra, his
classification does not directly imply a classification of simple
modules of $\Ule$, and we need an extra argument. A difficulty lies in
the fact that the homomorphism $\Ule\to H_*(\Zw^A,\C)$ in
\eqref{eq:homloc} is not necessarily isomorphism. Our additional input
is \propref{prop:stdhw}(3). In order to illustrate its usage, we first
consider the special case when $a = (s,\varepsilon)$ is generic in the
first subsection. In this case, Ginzburg's classification becomes
trivial. Then we shall review Ginzburg's classification in
\subsecref{subsec:Ginzburg}, and finally we shall study general case
in the last subsection.

We preserve the setup in \secref{sec:std}.

\subsection{}\label{subsec:generic}
Let us identify $e_{k,r}$, $f_{k,r}$ with their image under
\eqref{eq:homloc}.
Let $1_\rho\in H_*(\Zw^A)$ denote the fundamental
class $[\Delta\M(\rho)]$ of the diagonal of $\M(\rho)\times\M(\rho)$.

\begin{Lemma}\label{lem:ETaut}
Let us consider
$(\M(\rho^1)\times\M(\rho^2))\cap \Pa_k(\bv^2,\bw)\subset
(\M(\bv^1,\bw)^A\times\M(\bv^2,\bw)^A)\cap\Pa_k(\bv^2,\bw)$,
and let $\lambda_0$ be the weight of $A$ determined by $\rho^1$ and
$\rho^2$ as in \subsecref{subsec:HeckeFix}.
Then we have the following equality in
$H_*((\M(\rho^1)\times\M(\rho^2))\cap \Pa_k(\bv^2,\bw),\C)$\rom:
\begin{equation*}
\begin{split}
   & p_2^* \operatorname{ch}\left(\Wedge_{-u} V^2_l(\lambda)\right)
     \cap (1_{\rho^1}x_k^+(z)1_{\rho^2})  \\
   = 
   \; &\begin{cases}
      p_1^*\operatorname{ch}\left(\Wedge_{-u} V^1_l(\lambda)\right)\cap
          (1_{\rho^1}x_k^+(z)1_{\rho^2}) & \text{if $l\neq k$ or 
                                 $\lambda \neq \lambda_0$}, \\
      (1-\frac{uq}z)p_1^*\operatorname{ch}\left(\Wedge_{-u}
          V^1_l(\lambda)\right)\cap (1_{\rho^1} x_k^+(z) 1_{\rho^2})
        & \text{if $l = k$, $\lambda = \lambda_0$}.
   \end{cases}
\end{split}
\end{equation*}
\end{Lemma}

\begin{proof}
We have the following equality in
$K^0((\M(\rho^1)\times\M(\rho^2))\cap \Pa_k(\bv^2,\bw))$:
\begin{equation*}
   p_2^* V^2_l(\lambda) = 
   \begin{cases}
      p_1^* V^1_l(\lambda)
        & \text{if $l\neq k$ or $\lambda \neq \lambda_0$}, \\
      p_1^* V^1_k(\lambda_0) + (V^2/V^1)
        & \text{if $l = k$, $\lambda = \lambda_0$},
   \end{cases}
\end{equation*}
where $V^2/V^1$ is (the restriction of) the natural line bundle over
$\Pa_k(\bv^2,\bw)$. The assertion follows immediately.
\end{proof}

\begin{Theorem}\label{thm:generic}
Suppose that $a = (s,\varepsilon)$ is generic in the sense of
\defref{def:generic}. \rom(Hence, 
$\Mw_0^A = \Mw^A$.\rom)
Then the standard module
\begin{equation*}
   M_{0,a} = K^A(\Mw^A)\otimes_{R(A)}\otimes\C_a
   \cong H_*(\Mw^A,\C)
\end{equation*}
is a simple $\Ule$-module. Its Drinfel'd polynomial is given by
\begin{equation*}
   P_k(u) = \det( 1 - u\varepsilon^{-1}s_k),
\end{equation*}
where $s_k$ is the $\GL(W_k)$-component of $s\in G_\bw$.
Moreover, $M_{0,a}$ is isomorphic to a tensor product of
{\it l\/}-fundamental representations when $\mathfrak g$ is finite
dimensional.
\end{Theorem}

\begin{proof}
Recall that we have a distinguished vector (we denote it by $[0]$) in
the standard module $M_{0,a}$ (\subsecref{subsec:hwvector}).
It has the properties listed in \propref{prop:stdhw}.
In particular, it is the eigenvector for $p_k^\pm(z)$, and the
eigenvalues are given in terms of $P_k(u)$ therein. In the present
setting, $P_k(u)$ is equal to $\det( 1 - u\varepsilon^{-1}s_k)$.

Let
\begin{equation*}
   M_{0,a}^\circ \defeq \{ m\in M_{0,a} \mid
    \text{$e_{k,r}\ast m = 0$ for any $k\in I$, $r\in\Z$} \}.
\end{equation*}
We have $[0]\in M_{0,a}^\circ$.
We want to show that any nonzero submodule $M'$ of $M_{0,a}$ is
$M_{0,a}$ itself.
The weight space decomposition (as a
$\Ue$-module)~\eqref{eq:weightdec} of $M_{0,a}$ induces that of
$M'$.
Since the set of weights of $M'$ is bounded from $\bw$ with respect to
the dominance order, there exists a maximal weight of $M'$. Then a
vector in the corresponding weight space is killed by all
$e_{k,r}$ by the maximality. Thus
%
%
$M'$ contains a nonzero vector
$m\in M_{0,a}^\circ$. Hence it is enough to show that $M_{0,a}^\circ =
\C[0]$ since we have already shown that $M_{0,a} = \Ule^-\ast[0]$ in
\propref{prop:stdhw}(3).

Let us consider the operator
\begin{equation*}
   \left[ \Delta_* \Wedge_u V_l \right]
   \in K^A(Z(\bv,\bv;\bw))
\end{equation*}
where $\Delta\colon \M(\bv,\bw)\to Z(\bv,\bv;\bw)$. If we consider such
operators for various $\bv$, $l\in I$, they form a commuting family.
Moreover, $M_{0,a}^\circ$ is invariant under them since we have the
relation
\begin{equation*}
   e_{k,r} \ast \left[ \Delta_* \Wedge_u V^2_l \right]
   = 
   \begin{cases}
   \left[ \Delta_* \Wedge_u V^1_l \right] \ast e_{k,r}
      & \text{if $k\neq l$}, \\
   \left[ \Delta_* \Wedge_u V^1_k \right] \ast (e_{k,r}+uq e_{k,r+1})
      & \text{if $k = l$},
   \end{cases}
\end{equation*}
where $V^1_k$, $V^2_k$ are tautological bundles over
$\M(\bv^1,\bw)$, $\M(\bv^2,\bw)$ respectively.
(Here $\bv^2 = \bv^1 + \alpha_k$.)
Thus $M_{0,a}^\circ$ is a direct sum of generalized eigenspaces for
$\Delta_* \Wedge_u V_l$. Let us take a direct summand
$M_{0,a}^{\circ\circ}$. By \propref{prop:genweight}(1),
$M_{0,a}^{\circ\circ}$ is contained in $H_*(\M(\rho),\C)$ for some
$\rho\colon A\to G_\bv$. If we can show $M_{0,a}^{\circ\circ} =
\C[0]$, then we get $M_{0,a}^\circ = \C[0]$ since
$M_{0,a}^{\circ\circ}$ is an {\it arbitrary\/} direct summand.

Since $a$ is generic, we have $\M_0(\infty,\bw)^A = \{0\}$. Hence
\begin{equation*}
   \Zw^A = \Mw^A\times \Mw^A,
\end{equation*}
and $\Mw^A$ is a nonsingular projective variety
(having possibly infinitely many components). By the Poincar\'e
duality, the intersection pairing
\begin{equation*}
   (\ ,\ )\colon 
   H_*(\Mw^A,\C)\otimes H_*(\Mw^A,\C) \to \C
\end{equation*}
is nondegenerate.

Let $\lsp{t}f_{k,r}$ denote the transpose of
$f_{k,r}$ with respect to the pairing $(\ ,\ )$, namely
\begin{equation*}
   (f_{k,r}\ast m, m') = (m, \lsp{t}f_{k,r}\ast m')
   \qquad\text{for $m$, $m'\in H_*(\Mw^A,\C)$}.
\end{equation*}
By the definition of the convolution, $\lsp{t}f_{k,r}$ is equal to
$\omega_* f_{k,r}$ where
$\omega\colon \Mw^A\times\Mw^A
\to \Mw^A\times\Mw^A$ is the 
map exchanging the first and second factors
and
$\omega_*$ is the induced homomorphism on
$H_*(\Mw^A\times \Mw^A,\C)$.

Let us consider $1_\rho f_{k,r} 1_{\rho'}$, where $\rho$ is as above
and $\rho'$ is any other homomorphism. It is just the projection of
$f_{k,r}$ to the component $H_*(\M(\rho)\times \M(\rho'),\C)$. We have
\begin{equation*}
   1_{\rho'}\,\lsp{t}f_{k,r}\,1_{\rho}
   = 1_{\rho'}\, \omega_* f_{k,r}\, 1_{\rho}
   = (p_1^*\alpha \cup p_2^* \beta) \cap 1_{\rho'}\, e_{k,r'}\, 1_{\rho}
\end{equation*}
for some $r'\in\Z$, and $\alpha\in H_*(\M(\rho),\C)$,
$\beta\in H_*(\M(\rho'),\C)$.
Those $\alpha$ and $\beta$ come from asymmetry in the defintion of
$f_{k,r}$, $e_{k,r}$ and in the homomorphism \eqref{eq:homloc}.
We do not give their explicit forms, though it is possible.
What we need is that $\beta$ is written
%
%
by tensor powers of exterior products of $V_k^2(\lambda)$ for various
$k$, $\lambda$.
Thus we can write
\begin{equation*}
   (p_1^*\alpha \cup p_2^* \beta) \cap
   1_{\rho'}\, e_{k,r'}\, 1_{\rho}
   = \sum_{r''} p_1^* \alpha_{r''} \cap 1_{\rho'}\, e_{k,r''}\, 1_{\rho}
\end{equation*}
for some $\alpha_{r''}\in H_*(\M(\rho),\C)$ by \lemref{lem:ETaut}.
Therefore, for $m\in M_{0,a}^{\circ\circ}$, we have
\begin{equation*}
   (f_{k,r}\ast m', m) 
   = (1_\rho\, f_{k,r}\, 1_{\rho'}\ast m', m) 
   = \left(m', \sum_{r''} p_1^* \alpha_{r''} \cap
            1_{\rho'}\,e_{k,r''}\, 1_{\rho} \ast m\right)
   = 0
\end{equation*}
for any $k\in I$, $r\in\Z$, $\rho'$, $m'\in H_*(\M(\rho'),\C)$.
Here we have used $1_\rho\ast m= m$, $1_{\rho'}\ast m' = m'$,
$e_{k,r''}\ast m = 0$.
Since $H_*(\Mw^A,\C) = \Ule^-\ast[0]$, we have
one of the followings:
\begin{aenume}
\item $(m',m) = 0$ for any $m'\in H_*(\Mw^A,\C)$,
\item $m\in \C[0]$.
\end{aenume}
The first case is excluded by the nondegeneracy of $(\ ,\ )$.
Thus we have $m \in \C[0]$.

Let us prove the last assertion. First consider the case $\bw =
\Lambda_k$ for some $k$. If $\varepsilon$ is not a root of
unity, $a = (s,\varepsilon)\in\C^*\times\C^*$ is generic for the
quiver variety $\M(\Lambda_k)$. Hence the above shows that the
standard module for $\M(\Lambda_k)$ is simple, and hence gives an
{\it l\/}-fundamental representation.

Let us return to the case for general $\bw = \sum_k w_k \Lambda_k$.
Let $a_k^{1},\dots, a_k^{w_k}$ be eigenvalues of $s_k$ counted with
multiplicities.
By \propref{prop:tensor}, it is enough to show that
\begin{equation}\label{eq:wantprod}
  \dim M_{0,a} = \prod_k \prod_{i=1}^{w_k} \dim M_{0,a_k^i}(\Lambda_k), 
\end{equation}
where $M_{0,a_k}(\Lambda_k)$ is the standard module for
$\M(\Lambda_k)$ with $a_k = (s_k,\varepsilon)$.
Since $\Mw^A$ has no odd homology
groups (\thmref{thm:freeness'}), we have
\begin{equation*}
   \dim M_{0,a} = \operatorname{Euler}(\Mw^A),
\end{equation*}
where $\operatorname{Euler}(\ )$ denotes the topological Euler
number. By a property of the Euler number, we have
\begin{equation*}
    \operatorname{Euler}(\Mw^A) =  \operatorname{Euler}(\Mw).
\end{equation*}
If we take a maximal torus $T$ of $G_\bw$, the fixed point set
$\Mw^T$ is isomorphic to
\(
   \prod_k \M(\Lambda_k)^{w_k}.   
\)
Hence we have
\begin{equation*}
    \operatorname{Euler}(\Mw) = \prod_k
    \operatorname{Euler}(\M(\Lambda_k))^{w_k}.
\end{equation*}
Since we have
\begin{equation*}
   \dim M_{0,a_k^i}(\Lambda_k) = \operatorname{Euler}(\M(\Lambda_k)),
\end{equation*}
we get \eqref{eq:wantprod}.
\end{proof}

\subsection{Simple modules of the convolution algebra}\label{subsec:Ginzburg}

We briefly recall Ginzburg's classification of simple modules of
the convolution algebra \cite[\S8.6]{Gi-book}. (See also
\cite{Lu-cusp}.)

Let $X$ be a complex algebraic variety. We consider the derived
category of complexes with constructible cohomology sheaves, and
denote it by $D^b(X)$.  We use the notation in \cite{Gi-book}. For
example, we put
\(
   \Ext^k_{D^b(X)} (A, B) \defeq \Hom_{D^b(X)}(A, B[k]),
\)
\(
    \Ext^*_{D^b(X)}(A,B) \defeq \bigoplus_k \Ext^k_{D^b(X)}(A,B).
\)
$\Ext^*_{D^b(X)}(A,A)$ is an algebra by the Yoneda product.
The Verdier duality operator is denoted by ${}^\vee$.
Given graded vector spaces $V$, $W$, we write $V\ngpiso W$ if there
exists a linear isomorphism which does {\it not\/} necessarily
preserve the gradings. We will also use the same notation to denote
two objects are quasi-isomorphic up to a shift in the derived category.

Let $f\colon M\to X$ is a projective morphism between algebraic
varieties $M$, $X$, and assume that $M$ is nonsingular. Then we are in
the setting for the convolution in \secref{sec:convolution} with
$X_1 = X_2 = X_3$ and $Z_{12} = Z_{23} = Z$, where
\begin{equation*}
   Z \defeq M\times_X M 
   = \{ (m^1, m^2)\in M\times M\mid f(m^1) = f(m^2) \}.
\end{equation*}
Since $Z\circ Z = Z$, we have the convolution product
\begin{equation*}
   H_*(Z,\C)\otimes H_*(Z,\C) \to H_*(Z,\C).  
\end{equation*}
Let $\mathcal A$ be the algebra $H_*(Z,\C)$.
Set $M_x = f^{-1}(x)$. Then the convolution defines an
$\mathcal A$-module structure on $H_*(M_x,\C)$.
More generally, if $Y$ is a locally closed subset of $X$, then 
$H_*(f^{-1}(Y),\C)$ has an $\mathcal A$-module structure via
convolution.

By \cite[8.6.7]{Gi-book}, we have an algebra isomorphism, which does not
necessarily preserve gradings,
\begin{equation*}
   \mathcal A = H_*(Z,\C)
   \ngpiso \Ext^*_{D^b(X)}(f_*\C_M, f_*\C_M),
\end{equation*}
where $\C_M$ is the constant sheaf on $M$.

We apply the decomposition theorem~\cite{BBD} to $f_*\C_M$.
There exists an isomorphism in $D^b(X)$:
\begin{equation}\label{eq:decomp}
   f_*\C_M \cong \bigoplus_{\phi,k} L_{\phi,k}\otimes P_\phi[k],
\end{equation}
where $\{ P_\phi\}$ is the set of isomorphism classes of simple
perverse sheaves on $X$ such that some shift is a direct summand of
$f_*\C_M$. We thus have an isomorphism
\begin{equation*}
   \mathcal A \ngpiso
   \bigoplus_{i,j,k,\phi,\psi}
   \Hom_\C(L_{\phi,i}, L_{\psi,j})\otimes \Ext^k_{D^b(X)}(P_\phi, P_\psi).
\end{equation*}
Set $L_\phi\defeq \bigoplus_k L_{\phi,k}$ and
\begin{equation*}
   \mathcal A_k \defeq
   \bigoplus_{\phi,\psi}
   \Hom_\C(L_{\phi}, L_{\psi})\otimes \Ext^k_{D^b(X)}(P_\phi, P_\psi)
\end{equation*}
so that $\mathcal A = \bigoplus \mathcal A_k$. By definition,
$\mathcal A_k \cdot \mathcal A_l \subset \mathcal A_{k+l}$ under the
multiplication of $\mathcal A$.
By a property of perverse sheaves, we have
$\mathcal A_k = 0$ for $k < 0$
and $\Ext^0_{D^b(X)}(P_\phi, P_\psi) =
\C\delta_{\phi\psi}\operatorname{id}$.
Hence,
\begin{equation}\label{eq:convalg}
   \mathcal A = \mathcal A_0 \oplus \bigoplus_{k > 0} \mathcal A_k;
   \qquad
   \mathcal A_0 \cong \bigoplus_\phi \End(L_{\phi}).
\end{equation}
In particular, the projection $\mathcal A\to \mathcal A_0$ is an
algebra homomorphism.
Furthermore, $\mathcal A_0$ is a semisimple algebra.
And the kernel of the projection, i.e., $\bigoplus_{k > 0}\mathcal A_k$,
consists of nilpotent elements, thus it
is precisely the radical of $\mathcal A$.
In particular,
\begin{equation*}
   \{ L_\phi \}_{\phi}
\end{equation*}
is a complete set of mutually non-isomorphic simple
$\mathcal A$-modules.

For $x\in X$, let $i_x\colon \{x\} \to X$ denote the inclusion. Then
$H^*(i_x^! f_*\C_M)$ is an $\Ext^*_{D^b(X)}(f_*\C_M, \linebreak[3]
f_*\C_M)$-module.
More generally, if $i_Y\colon Y\hookrightarrow X$ is a locally closed
embedding, then the hyper-cohomology groups $H^*(Y,i_Y^! f_*\C_M)$ and
$H^*(Y,i_Y^*, f_*\C_M)$ are $\Ext^*_{D^b(X)}(f_*\C_M, f_*\C_M)$-modules.
It is known \cite[8.6.16, 8.6.35]{Gi-book} that $H^*(Y, i_Y^! f_*\C_M)$
is isomorphic to $H_*(f^{-1}(Y),\C)$ as an $\mathcal A \cong
\Ext^*_{D^b(X)}(f_*\C_M, f_*\C_M)$-module.

For $C\in D^b(\{x\})$, we write $H^k(C)$ instead of $H^k(\{x\}, C)$,
and $H^*(C)$ instead of $H^*(\{x\}, C)$.
By applying $H^*(i_x^!\bullet)$ to \eqref{eq:decomp}, we get an
isomorphism
\begin{equation*}
   H_*(M_x,\C) \ngpiso
   \bigoplus_{\phi,k} L_{\phi}\otimes H^k(i_x^!P_\phi).
\end{equation*}
Let 
\[
 M_{\ge k}(i_x^!f_*\C_M) \defeq
  \bigoplus_{k'\ge k} L_{\phi}\otimes H^{k'}(i_x^!P_\phi).
\]
By definition, we have $\mathcal A_k\cdot M_{\ge l}(i_x^!f_*\C_M)
\subset M_{\ge k+l}(i_x^!f_*\C_M)$ under the $\mathcal A$-module on
$H^*(i_x^!f_*\C_M)$. In particular, $M_{\ge k}(i_x^!f_*\C_M)$ is
an $\mathcal A$-submodule for each $k$. Hence
\begin{equation*}
   \operatorname{gr} M(i_x^!f_*\C_M)
   \defeq \bigoplus_k  M_{\ge k}(i_x^!f_*\C_M)/ M_{\ge k+1}(i_x^!f_*\C_M)
\end{equation*}
is an $\mathcal A$-module, on which $\bigoplus_{k > 0} \mathcal A_k$
acts as $0$. By definition, 
\begin{equation*}
  \operatorname{gr} M(i_x^!f_*\C_M)
  \cong \bigoplus_\phi L_\phi\otimes H^*(i_x^!P_\phi),
\end{equation*}
where the $\mathcal A$-module structure on the right hand side is
given by
$a\colon \xi\otimes\xi'\mapsto a\xi\otimes\xi'$.
Thus we have
\begin{Theorem}
In the Grothendieck group of $\mathcal A$-modules of finite dimension
over $\C$, we have
\begin{equation*}
   H_*(M_x,\C) = \bigoplus_\phi L_\phi\otimes H^*(i_x^!P_\phi),
\end{equation*}
where the $\mathcal A$-module structure on the right hand side is
given by
$a\colon \xi\otimes\xi'\mapsto a\xi\otimes\xi'$.
\end{Theorem}

\begin{proof}
Since $\operatorname{gr} M(i_x^!f_*\C_M)$ is equal to
$M(i_x^!f_*\C_M)$ in the Grothendieck group, the assertion follows
from the discussion above.
\end{proof}

\subsection{}\label{subsec:std_decomp}
In this subsection, we assume that the graph is of type $ADE$, and
$\varepsilon$ is not a root of unity.
We apply the results in the previous subsection to our quiver
varieties.

Recall
\begin{equation*}
   \Zw^A 
   = \{ (x^1, x^2)\in \Mw^A\times\Mw^A
   \mid \pi^A(x^1) = \pi^A(x^2) \},
\end{equation*}
where $\pi^A\colon\Mw^A\to \M_0(\infty,\bw)^A$ is the
restriction of $\pi\colon\Mw\to \M_0(\infty,\bw)$.
Thus the results in the previous subsection are applicable to this setting.
We have an algebra isomorphism
\begin{equation*}
   H_*(\Zw^A,\C) \ngpiso
   \Ext^*_{D^b(\M_0(\infty,\bw)^A)}
   (\pi^A_*\C_{\M(\bw)^A}, \pi^A_*\C_{\M(\bw)^A}).
\end{equation*}
Let us denote this algebra by $\mathcal A$ as in the previous
subsection.

Since the graph is of type $ADE$, we have
$\M_0(\infty,\bw) = \bigsqcup_\bv \Mreg(\bv,\bw)$ by \propref{prop:ADE}.
Thus we have the stratification
$\M_0(\infty,\bw)^A = \bigsqcup \Mreg(\rho)$.
Since the restriction
\[
    \pi^A|_{(\pi^A)^{-1}(\Mreg(\rho))}
\colon (\pi^A)^{-1}(\Mreg(\rho))\to
\Mreg(\rho)
\]
is a locally trivial topological fibration by \thmref{thm:slice},
all the complexes in the right hand side of \eqref{eq:decomp}
(applied to $f = \pi^A$, $M = \Mw^A$)
have locally constant cohomology sheaves along each stratum
$\Mreg(\rho)$. Since $\Mreg(\rho)$ is irreducible by
\thmref{thm:M(rho)_conn}, it implies that $P_\phi$ is the
intersection cohomology complex
\(
     IC(\Mreg(\rho),\phi)
\)
associated with an irreducible local system $\phi$ on
$\Mreg(\rho)$. Thus we have
\begin{equation}
\label{eq:DecompSecond}
   \pi^A_*\C_{\M(\bw)^A} \cong
   \bigoplus_{(\rho,\phi,k)} L_{(\rho,\phi,k)}\otimes 
   IC(\Mreg(\rho),\phi)[k]
\end{equation}
for some finite dimensional vector space $L_{(\rho,\phi,k)}$.
Let $L_{(\rho,\phi)} \defeq \bigoplus_k L_{(\rho,\phi,k)}$.
By a discussion in the previous subsection, $\{ L_{(\rho,\phi)} \}$ is 
a complete set of mutually non-isomorphic simple
\(
    \mathcal A
\)-modules. Via the homomorphism~\eqref{eq:homloc}, $L_{(\rho,\phi)}$
is considered also as a $\Ule$-module.

\begin{Theorem}\label{thm:std_decomp}
Assume $\varepsilon$ is not a root of unity.

\rom{(1)} Simple perverse sheaves $P_{\phi}$ whose shift appear in a
direct summand of $\pi^A_*\C_{\M(\bw)^A}$ are the intersection
cohomology complexes $IC(\Mreg(\rho))$
associated with the constant local system $\C_{\Mreg(\rho)}$ on various
$\Mreg(\rho)$.

\rom{(2)} Let us denote the constant local system $\C_{\Mreg(\rho)}$
by $\C_\rho$ for simplicity.
Then $L_{(\rho,\C_\rho)}$ is nonzero if and only if
$\Mreg(\rho)\neq\emptyset$. Moreover, there is a bijection between
the set $\{\rho\mid L_{(\rho,\C_{\rho})}\neq 0\}$
and
the set of {\it l\/}-weights of $M_{0,a}$ which are {\it l\/}-dominant.

\rom{(3)} The simple $\mathcal A = H_*(Z(\bw)^A,\C)$-module 
$L_{(\rho,\C_\rho)} = \bigoplus_k L_{(\rho,\C_\rho,k)}$
is also simple as a $\Ule$-module, and its Drinfel'd
polynomial is
$P_k(u) = \chi_a\left(\Wedge_{-u} C_{k,x}^\bullet\right)$
in \propref{prop:stdhw} for $x\in \Mreg(\rho)$.

\rom{(4)} $L_{(\rho,\C_\rho)}$ is the simple quotient of $M_{x,a}$,
where $x$ is a point in a stratum $\Mreg(\rho)$.


\rom{(5)} Standard modules $M_{x,a}$ and $M_{y,a}$ are isomorphic as
$\Ule$-modules if and only if $x$ and $y$ are contained in the same
stratum.
\end{Theorem}

\begin{proof}
We use the transversal slice in \subsecref{subsec:slice}.
The idea to use transversal slices is taken from \cite[\S8.5]{Gi-book}.

Choose and fix a point $x\in \M_0(\infty,\bw)^A$. Suppose that $x$ is
contained in a stratum $\Mreg(\rho_x)$ for some $\rho_x$. We first
show
\begin{Claim}
If $\C_{\rho_x}$ denote the constant local system on
$\Mreg(\rho_x)$, the corresponding vector space
$L_{(\rho_x,\C_{\rho_x})}$ is nonzero.
\end{Claim}

If we restrict $\pi^A$ to the component $\M(\rho_x)$, then we have
\begin{equation}\label{eq:DecompThird}
   \pi^A_*\C_{\M(\rho_x)} \ngpiso
   \bigoplus_{(\rho,\phi)} L_{(\rho,\phi)}'\otimes 
   IC(\Mreg(\rho),\phi),
\end{equation}
where $L_{(\rho,\phi)}'$ is a direct summand of $L_{(\rho,\phi)}$.
The summation runs over the set of pairs $(\rho,\phi)$ such that
$\Mreg(\rho)$ is contained in $\pi^A(\M(\rho_x))$.
(In fact, \eqref{eq:DecompSecond} was obtained by applying the
decomposition theorem to each component $\M(\rho_x)$ and taking direct
sum.) 
If we restrict \eqref{eq:DecompThird} to the open stratum $\Mreg(\rho_x)$
of $\pi^A(\M(\rho_x))$, the right hand side of
\eqref{eq:DecompThird} becomes
\begin{equation*}
   \bigoplus_\phi L_{(\rho_x,\phi)}'\otimes \phi,
\end{equation*}
where the summation runs over the set of isomorphism classes of
irreducible local systems $\phi$ on
$\Mreg(\rho_x)$.
On the other hand, $\pi^A$ induces an isomorphism between
$(\pi^A)^{-1}(\Mreg(\rho_x))$ and $\Mreg(\rho_x)$ by
\propref{prop:pi_isom}. This means that the restriction of the left
hand side of \eqref{eq:DecompThird} is the constant local system
$\C_{\rho_x}$. Hence we have $L_{(\rho_x,\C_{\rho_x})}' \cong \C$, and
$L_{(\rho_x,\C_{\rho_x})}$ is nonzero. This is the end of the proof of
the claim.

The claim implies the first assertion of (2).
Let us prove the latter assertion of (2).
Suppose $\Mreg(\rho)\neq\emptyset$. Then we have
$\M(\rho)\neq\emptyset$ and $H_*(\M(\rho)_0,\C)\neq 0$ by
\propref{prop:homotopic'}.
By \propref{prop:genweight}, the corresponding {\it l\/}-weight space
is nonzero, where the {\it l\/}-weight $\Psi^\pm(z) =
(\Psi_k^\pm(z))_k$ is given by \eqref{eq:genweight}. Furthermore,
since $C_k^\bullet(\bv,\bw)$ can be represented by a genuine
$A$-module over a point in $\Mreg(\rho)$ by \lemref{lem:betasurj},
$\chi_a\left(\Wedge_{-u} C_k^\bullet(\bv,\bw)\right)$ is a polynomial
in $u$. Thus $\Psi^\pm(z)$ is {\it l\/}-dominant.

Conversely suppose that we have the {\it l\/}-weight space with the
{\it l\/}-weight \eqref{eq:genweight} is nonzero. Since $\varepsilon$
is {\it not\/} a root of unity, the $\varepsilon$-analogue of the
Cartan matrix $[-\langle h_k,\alpha_l\rangle]_\varepsilon$ is
invertible. Hence \eqref{eq:genweight} determines $\chi_a(\Wedge_{u}
V_k)$. Thus the {\it l\/}-weight space is precisely
$H_*(\M(\rho)_0,\C)$ by \propref{prop:genweight}(1). In particular, we
have $H_*(\M(\rho)_0,\C)\neq 0$, and hence $\M(\rho)\neq\emptyset$.
Furthermore, if we decompose $C_k^\bullet(\bv,\bw)$ into
$\bigoplus_\lambda C_{k,\lambda}^\bullet(\rho)$ as in
\subsecref{subsec:fix_hom}, we have \( \rank
C_{k,\lambda}^\bullet(\rho) \ge 0 \) since the {\it l\/}-weight
\eqref{eq:genweight} is {\it l\/}-dominant.  By
\corref{cor:codim_rank}, $\tau_{k,\lambda}$ is surjective for any $k$,
$\lambda$ on a nonempty open subset of $\M(\rho)$. By
\lemref{lem:surj} (plus the subsequent remark), $\Mreg(\rho)$ is
nonempty. This shows the latter half of (2).

Let $\bv_x$ denote the dimension vector corresponding to $\rho_x$,
i.e., $\Mreg(\rho_x)\subset \Mreg(\bv_x,\bw)$. Take a transversal
slice to $\Mreg(\bv_x,\bw)$ at $x$ as in \subsecref{subsec:slice}.
Let $S$ be its intersection with $\M_0(\infty,\bw)^A$. Since the
transversal slice in \subsecref{subsec:slice} can be made
$A$-equivariant (\remref{rem:Aslice}), it is a transversal slice to
$\Mreg(\rho_x)$ (at $x$) in $\M_0(\infty,\bw)^A$. Let $\widetilde S
\defeq (\pi^A)^{-1}(S)$.
Let $\varepsilon\colon S\to
\M_0(\infty,\bw)^A$, $\widetilde\varepsilon \colon \widetilde S\to
\Mw^A$ be the inclusions.

The stratification
$\M_0(\bw)^A = \bigsqcup \Mreg(\rho)$
induces by restriction a stratification
$S = \bigsqcup S_\rho$ where $S_\rho =
\Mreg(\rho)\cap S$.
Any intersection complex $IC(\Mreg(\rho),\phi)$
restricts (up to shift) to the intersection complex
$IC(S_\rho,\phi|_{S_\rho})$ by transversality. Here $\phi|_{S_\rho}$
is the restriction of $\phi$ to $S_{\rho}$. Taking $\varepsilon^!$ of
\eqref{eq:DecompSecond}, we get
\begin{equation}\label{eq:DecompFourth}
   \varepsilon^! \left(\pi^A_*\C_{\M(\bw)^A}\right)
   \ngpiso \bigoplus_{(\rho,\phi)} L_{(\rho,\phi)}\otimes 
   IC(S_\rho,\phi|_{S_\rho}).
\end{equation}

Let $i^S_x\colon \{x\}\to S$ be the inclusion. It induces two
pull-back homomorphisms $i^{S!}_x$, $i^{S*}_x$, and there is a
natural morphism $i^{S!}_x E \to i^{S*}_x E$ for any $E\in
D^b(S)$. We apply these functors to both hand sides of
\eqref{eq:DecompFourth} and take cohomology groups.
By a property of intersection cohomology
sheaves (see \cite[8.5.3]{Gi-book}), the homomorphism
\begin{equation}\label{eq:cpt_usl}
   H^*(i^{S!}_x IC(S_\rho,\phi|_{S_\rho})) \to
   H^*(i^{S*}_x IC(S_\rho,\phi|_{S_\rho}))
\end{equation}
is zero unless $S_\rho = \{x\}$ (or equivalently
$\rho = \rho_x$),
in which case it is a quasi-isomorphism.
Thus
\begin{equation}\label{eq:Ima}
   \Ima\left[
      H^*(i^{S!}_x \varepsilon^! \pi^A_* \C_{\M(\bw)^A})
    \to 
      H^*(i^{S*}_x \varepsilon^! \pi^A_* \C_{\M(\bw)^A})
    \right]
   \ngpiso \bigoplus_{\phi} L_{(\rho_x,\phi)}\otimes 
   \phi_{x},
\end{equation}
where the summation runs over isomorphism classes of irreducible local 
systems on $\Mreg(\rho_x)$,
and $\phi_{x}$ is the fiber of the local system $\phi$ at $x$.
Moreover, \eqref{eq:cpt_usl} is a homomorphism of
$\mathcal A$-modules, and \eqref{eq:Ima} is an isomorphism of
$\mathcal A$-modules, where the module structure on the right
hand side is given by $a\colon \xi\otimes\xi'\mapsto a\xi\otimes \xi'$.

On the other hand, we have 
\[
   H^*(i^{S!}_x \varepsilon^! \pi^A_* \C_{\M(\bw)^A})
   = H^*(i_x^! \pi^A_* \C_{\widetilde S})
   \ngpiso H_*(\M(\bw)^A_x,\C).
\]
As shown in \eqref{eq:StdIsom}, the right hand side is isomorphic to
the standard module $M_{x,a}$. Thus the left hand side of
\eqref{eq:Ima} is a quotient of $M_{x,a}$, and it is indecomposable by
\propref{prop:stdhw}(3). Thus the right hand side of \eqref{eq:Ima}
consists of at most single direct summand. Since we have already shown
that $L_{(\rho_x,\C_{\rho_x})}\neq 0$ in the claim, we get
\(
   L_{(\rho_x,\phi)} = 0
\)
if $\phi$ is a {\it nonconstant\/} irreducible local system. Since
$x$ was an arbitrary point, we have the statement~(1).

Let us prove (3).
For the proof, we need a further study of \eqref{eq:Ima}. By the above
discussion, we have
\begin{equation}\label{eq:Ima'}
      \Ima\left[
      H^*(i^{S!}_x \varepsilon^! \pi^A_* \C_{\M(\bw)^A})
    \to 
      H^*(i^{S*}_x \varepsilon^! \pi^A_* \C_{\M(\bw)^A})
    \right]
   \ngpiso L_{(\rho_x,\C_{\rho_x})}.
\end{equation}
By the base change theorem, we have 
\(
  \varepsilon^! \left(\pi^A_*\C_{\M(\bw)^A}\right)
   = \pi^S_* \widetilde\varepsilon^!\C_{\M(\bw)^A}
\)
 where $\pi^S$ is the restriction of
$\pi^A$ to $\widetilde S$. Further, we have
\(
   \widetilde\varepsilon^!\C_{\M(\bw)^A} \ngpiso \C_{\widetilde S}
\)
since $\widetilde S$ is a nonsingular submanifold of
$\M(\bw)^A$. Applying the Verdier duality, we have
\begin{equation*}
   \Hom\left(H^*(i^{S!}_x \varepsilon^! \pi^A_*\C_{\M(\bw)^A}),
     \C \right) 
  \ngpiso
   H^*((i^{S!}_x \pi^S_* \C_{\widetilde S})^\vee)
  \ngpiso
   H^*(i^{S*}_x \pi^S_* \C_{\widetilde S}).
\end{equation*}
Hence \eqref{eq:Ima'} becomes
\begin{equation}\label{eq:Ima''}
    \Ima\left[
      M_{x,a}
    \to 
      M_{x,a}^*
    \right]
   \ngpiso L_{(\rho_x,\C_{\rho_x})},
\end{equation}
where $M_{x,a}^*$ is the dual space of $M_{x,a}$ as a complex vector space.
Let us introduce an $\mathcal A$-module on $M_{x,a}^*$ by
\begin{equation*}
   \langle a \ast h, \xi\rangle 
    = \langle h, (\omega_* a)\ast \xi\rangle, \quad 
     a\in \mathcal A, \ h\in M_{x,a}^*, \ \xi\in M_{x,a},
\end{equation*}
where $\langle\ ,\ \rangle$ denote the dual pairing,
$\omega\colon \Zw^A\to \Zw^A$ is the exchange of two factors of
$\Zw^A = \Mw^A\times_{\M_0(\infty,\bw)^A}\Mw^A$,
and $\omega_*$ is the induced homomorphism on
$\mathcal A = H_*(\Zw^A,\C)$.
Then \eqref{eq:Ima''} is compatible with $\mathcal A$-module
structures (cf.\ \cite[paragraphs preceding 8.6.25]{Gi-book}).

The decomposition \eqref{eq:genwei} induces a similar one for
$M_{x,a}^*$:
\begin{equation*}
   M_{x,a}^* = \bigoplus_\rho \Hom(H_*(\M(\rho)_x,\C),\C).
\end{equation*}
The homomorphism $M_{x,a}\to M_{x,a}^*$ respects the decomposition,
and induces a decomposition on \eqref{eq:Ima''}.

Recall that we have the distinguished vector $[x]$ in $M_{x,a}$. The
component $H_*(\M(\rho_x),\C)$ of $M_{x,a}$ is $1$-dimensional space
$\C[x]$. (See \subsecref{subsec:hwvector}.) By the above discussion,
$[x]$ is {\it not\/} annihilated by the above homomorphism
$M_{x,a} \to M_{x,a}^*$. Thus we may consider
$[x]$ also as an element of $M_{x,a}^*$.

We want to show that any nonzero $\Ule$-submodule $L'$ of
$L_{(\rho_x,\C_{\rho_x})}$ is $L_{(\rho_x,\C_{\rho_x})}$ itself.
Our strategy is the same as the proof of \thmref{thm:generic}.
Since we already show that $L_{(\rho_x,\C_{\rho_x})}$ is a quotient
of $M_{x,a}$, \propref{prop:stdhw}(3) implies
$L_{(\rho_x,\C_{\rho_x})} = \Ule^-\ast [x]$. Thus it is enough to show 
that $L'$ contains $[x]$. To show this, consider
\begin{equation*}
   M_{x,a}^{*\circ} \defeq \{ m^*\in M^*_{x,a} \mid
    \text{$e_{k,r}\ast m^* = 0$ for any $k\in I$, $r\in\Z$} \}.
\end{equation*}
By the argument as in the proof of \thmref{thm:generic}, $L'$ contains
a nonzero vector in $M_{x,a}^{*\circ}$. Hence it is enough to show
that $M_{x,a}^{*\circ} = \C[x]$.

As in the proof of \thmref{thm:generic}, $M_{x,a}^{*\circ}$ is a
direct sum of generalized eigenspaces for $\Delta_*\Wedge_u V_l$. Let
us choose and fix a direct summand $M_{x,a}^{*\circ\circ}$ contained
in $H_*(\M(\rho)_x,\C)$.
Then $m^*\in M_{x,a}^{*\circ\circ}$ satisfies
\begin{equation*}
  \langle f_{k,r}\ast m, m^*\rangle = 0
\end{equation*}
for any $k$, $r$.
Since $M_{x,a} = \Ule^-\ast [x]$ by \propref{prop:stdhw}(3), the above 
equation implies that $m^*\in \Hom(\C[x],\C)$. Thus we get
$M_{x,a}^{*\circ} = \C[x]$ as desired.

We have shown the statement~(4) during the above discussion.

Let us prove (5). Since $\pi^A$ is a locally trivial topological
fibration on each stratum $\Mreg(\rho)$, $M_{x,a}$ and $M_{y,a}$ are
isomorphic if both $x$ and $y$ is contained in
$\Mreg(\rho)$. Conversely, if $M_{x,a}$ and $M_{y,a}$ are
isomorphic as $\Ule$-modules, the corresponding {\it l\/}-highest
weights $\chi_a\left(\Wedge_{-u} C_{k,x}^\bullet\right)$
and $\chi_a\left(\Wedge_{-u} C_{k,y}^\bullet\right)$ are equal.
Since $\chi_a\left(\Wedge_{-u} C_{k,x}^\bullet\right)$ determines the
homomorphism $\rho$ as in the proof of (2), $x$ and $y$ are in the
same stratum.
\end{proof}

\begin{Remark}\label{rem:root of unity}
The assumption that $\varepsilon$ is {\it not\/} a root of unity is
used to apply \thmref{thm:M(rho)_conn} and to have the invertibility
of the $\varepsilon$-analogue of the Cartan matrix.
It seems likely that \thmref{thm:M(rho)_conn} holds even if
$\varepsilon$ is a root of unity.
The latter condition was used to parametrize the index set of $\rho$
(i.e., \thmref{thm:std_decomp}(2)). But one should have a similar
parametrization if one replace a notion of {\it l\/}-weights in a
suitable way.
Thus \thmref{thm:std_decomp} should hold even if $\varepsilon$ is a
root of unity, if one replace the statement~(2).
\end{Remark}

Let $\mathcal P = \{ P(u) = (P_k(u))_k \}$ be the set of {\it
l\/}-weights of $M_{0,a}$, which are {\it l\/}-dominant.  Since the
index set $\{ \rho \}$ of the stratum coincides with $\mathcal P$, we
may write $L_{(\rho,\C_\rho)}$ as $L(P)$, when $\Mreg(\rho)$
corresponds to $P\in\mathcal P$. The standard module $M_{x,a}$ depends
only on the stratum containing $x$, so we may also write $M_{x,a}$ as
$M(P)$. We have an analogue of the Kazhdan-Lusztig multiplicity
formula:
\begin{Theorem}
Assume $\varepsilon$ is not a root of unity.

For $x\in \M_0(\infty,\bw)^A$, let $i_x\colon \{x\} \to
\M_0(\infty,\bw)^A$ denote the inclusion.
Let $P\in\mathcal P$ be the {\it l\/}-weight corresponding to the
stratum $\Mreg(\rho_x)$ containing $x$.
In the Grothendieck group of finite dimensional
$\Ule$-modules, we have
\begin{equation*}
   M(P) = \bigoplus_{Q\in\mathcal P} L(Q) \otimes
     H^*(i_x^! IC(\Mreg(\rho_Q)),
\end{equation*}
where $\Mreg(\rho_Q)$ is the stratum corresponding to $Q\in\mathcal
P$, and $IC(\Mreg(\rho_Q))$ is the intersection cohomology complex
attached to $\Mreg(\rho_Q)$ and the constant local system
$\C_{\Mreg(\rho_Q)}$.
Here the $\Ule$-module structure on the right hand side is given by
$a\colon \xi\otimes\xi'\mapsto a\xi\otimes \xi'$.
\end{Theorem}

This follows from \thmref{thm:std_decomp} and a result in the
previous subsection.

\begin{Remark}
By \cite[8.7.8]{Gi-book} and \thmref{thm:freeness'} with \thmref{thm:slice}, 
$H^{d_Q + n}(i_x^! IC(\Mreg(\rho_Q))$ 
vanish for all odd $n$, where $d_Q$ is the dimension of $\Mreg(\rho_Q)$.
\end{Remark}

\section{The $\Ue$-module structure}\label{sec:Ue}

In this section, we assume the graph is of type $ADE$.
The result of this section holds even if $\varepsilon$ is a root of
unity, if we replace the simple module $L(\Lambda)$ by the
corresponding Weyl module (see \cite[11.2]{CP} for the definition).

\subsection{}
For a given $\bw\in\bigoplus\Z_{\ge 0}\Lambda_k$,
let $\mathcal V_{\bv^0}(\bw)$ be the finite set consisting of all
$\bv\in\bigoplus \Z_{\ge 0}\alpha_k$ such that $\bw - \bv$ is
dominant and the weight space with weight $\bw - \bv^0$ is
nonzero in the simple highest weight $\Uq$-module $L(\bw-\bv)$.

Let $\mathcal V(\bw)$ be the union of all $\mathcal V_{\bv^0}(\bw)$ for
various $\bv^0$.
It is the set consisting of all $\bv$ such that $\bw - \bv$ is 
dominant.

Since the graph is of type $ADE$, we have
$\M_0(\infty,\bw) = \bigsqcup_{\bv} \Mreg(\bv,\bw)$.
Since $\Mreg(\bv,\bw)$ is isomorphic to an open subvariety of
$\M(\bv,\bw)$, $\Mreg(\bv,\bw)$ is irreducible if $\M(\bv,\bw)$ is
connected.
Although we do not know whether $\M(\bv,\bw)$ is connected or not (see
\subsecref{subsec:connected}), we consider the intersection cohomology
complex $IC(\Mreg(\bv,\bw))$ attached to $\Mreg(\bv,\bw)$ and the
constant local system $\C_{\Mreg(\bv,\bw)}$. It may not be a simple
perverse sheaf if $\M(\bv,\bw)$ is not connected.

We prove the following in this section:
\begin{Theorem}\label{thm:UqmoduleStr}
As a $\Ue$-module, we have the following decomposition:
\begin{equation*}
  \operatorname{Res}M_{x,a} = \bigoplus_{\bv\in\mathcal V(\bw)}
      H^*(i_x^! IC(\Mreg(\bv,\bw)))\otimes L(\bw - \bv),
\end{equation*}
where $i_x\colon\{x\}\to \M_0(\infty,\bw)$ is the inclusion,
and $\Ue$ acts trivially on $H^*(i_x^! IC(\Mreg(\bv^0,\bw)))$.
\end{Theorem}

\begin{Remark}
By \cite[8.7.8]{Gi-book} and \thmref{thm:freeness} with
\thmref{thm:slice}, $H^{\operatorname{odd}}(i_x^! IC(\Mreg(\bv,\bw)))$
vanishes.
\end{Remark}

\subsection{Reduction to $\varepsilon=1$}
Suppose that $x$ is contained in a stratum $\Mreg(\bv,\bw)$. Take a
representative $(B,i,j)$ of $x$ and define $\rho(a)$ as in
\eqref{eq:fixed}. Here $a$ is fixed and we do not consider $A$. We
choose $S\in\mathfrak g_\bw = \operatorname{Lie} G_\bw$,
$R\in\mathfrak g_\bv = \operatorname{Lie} G_\bv$, $E\in\C$ so that
$\exp S = s$, $\exp R = \rho(a)$, $\exp E = \varepsilon$, where $a =
(s,\varepsilon)$.
Let $a_t = (\exp tS, \exp tE)$ for $t\in\C$. Then we have
\begin{equation*}
   a_t \ast (B,i,j) = \exp(tR)^{-1}\cdot (B,i,j)
\end{equation*}
from \eqref{eq:fixed}. If $A_x$ denotes the stabilizer of $x$ in
$G_\bw\times\C^*$, the above equation means that $a_t\in A_x$.

Let us consider a ${\mathbf U}_{\exp tE}({\mathbf L}\mathfrak
g)$-module 
\[
   M_t \defeq K^{A_x}(\Mw_x)\otimes_{R(A_x)} \C_{a_t}
\]
parametrized by $t\in\C$, where $\C_{a_t}$ is an $R(A_x)$-module
given by the evaluation at $a_t$ as in \secref{sec:std}.
When $t = 1$, we can replace $A_x$ by $A$ by \thmref{thm:freeness} and
\thmref{thm:slice}, hence the module $M_{t=1}$ coincides with
$M_{x,a}$.
Moreover, it depends continuously on $t$ also by
\thmref{thm:freeness}.

Let us consider $M_t$ as a ${\mathbf U}_{\exp tE}(\mathfrak g)$-module
by the restriction. Since finite dimensional ${\mathbf U}_{\exp
tE}(\mathfrak g)$-modules are classified by discrete data (highest
weights), it is independent of $t$.
(Simple modules $L(\Lambda)$ of ${\mathbf U}_{\exp tE}(\mathfrak g)$
depends continuously on $t$.)
Thus it is enough to decompose $M_t$ when $t = 0$, i.e., $s = 1$,
$\varepsilon = 1$.
By \thmref{thm:freeness} and \thmref{thm:slice}, $K^{A_x}(\Mw_x)$ is
specialized to $H_*(\Mw_x,\C)$ at $s=1$, $\varepsilon=1$. Thus our
task now becomes the decomposition of $H_*(\Mw_x,\C)$ into simple
$\mathfrak g$-modules.

\subsection{}

When $Y$ is pure dimensional, we denote by $H_{\topdeg}(Y,\C)$ the top
degree part of $H_*(Y,\C)$, that is the subspace spanned by the
fundamental classes of irreducible components of $Y$. 
Suppose that $Y$ has several connected components $Y_1, Y_2, \dots$
such that each $Y_i$ is pure dimensional, but $\dim Y_i$ may change
for different $i$. Then we define $H_{\topdeg}(Y,\C)$ as $\bigoplus
H_{\topdeg}(Y_i,\C)$. Note that the degree $\topdeg$ may differ for
different $i$ since the dimensions are changing.

By \cite[9.4]{Na-alg}, there is a homomorphism
\begin{equation*}
   \Ueone \to
     H_{\topdeg}(\Zw,\C).
\end{equation*}
In fact, it is the restriction of the homomorphism in
\eqref{eq:homloc} for $A = \{1\}$, $\varepsilon = 1$,
composed with the projection
\begin{equation*}
  H_*(\Zw,\C)
  \to H_{\topdeg}(\Zw,\C).
\end{equation*}
For each $\bv$, we take a point $x_\bv\in \Mreg(\bv,\bw)$.
(By \lemref{lem:surj}(2), $\bw - \bv$ is dominant if
$\Mreg(\bv,\bw)$ is nonempty.)
By \cite[10.2]{Na-alg}, $H_{\topdeg}(\Mw_{x_{\bv}},\C)$ is the
simple highest weight module $L(\bw-\bv)$ via this homomorphism.
(In fact, we have already proved a similar result, i.e.,
\propref{prop:stdhw}.)

\begin{Proposition}
Consider the map $\pi\colon \M(\bv^0,\bw)\to \M_0(\bv^0,\bw)$.
Then $\pi$, as a map into $\pi(\M(\bv,\bw))$, is semi-small
and all strata are relevant, namely
\begin{equation*}
 2 \dim \M(\bv^0,\bw)_{x_{\bv}} = \codim \M_0(\bv,\bw)
 \quad \text{for $x_{\bv}\in \Mreg(\bv,\bw)$},
\end{equation*}
where $\codim$ is the codimension in $\pi(\M(\bv,\bw))$.
\end{Proposition}

\begin{proof}
See \cite[6.11]{Na-quiver} and \cite[10.11]{Na-alg}.
\end{proof}

\begin{Proposition}
We have
\begin{equation}\label{eq:decompQ}
  \pi_*\left(\C_{\M(\bv^0,\bw)}[\dim\M(\bv^0,\bw)]\right)
  = \bigoplus_{\bv\in\mathcal V_{\bv^0}(\bw)}
          H_{\topdeg}(\M(\bv^0,\bw)_{x_{\bv}},\C)\otimes IC(\Mreg(\bv,\bw)),
\end{equation}
where $x_{\bv}$ is taken from $\Mreg(\bv,\bw)$.
\rom(By \thmref{thm:slice} $H_{\topdeg}(\M(\bv^0,\bw)_{x_{\bv}},\C)$
is independent of choice of $x_{\bv}$.\rom)
\end{Proposition}

\begin{proof}
By the decomposition theorem for a semi-small map
\cite[8.9.3]{Gi-book}, the left hand side of \eqref{eq:decompQ}
decomposes as
\begin{equation*}
   \pi_*\left(\C_{\M(\bv^0,\bw)}[\dim\M(\bv^0,\bw)]\right)
  = \bigoplus_{\bv,\alpha,\phi}
  L_{(\bv,\alpha,\phi)} \otimes IC(\Mreg(\bv,\bw)^{\alpha},\phi),
\end{equation*}
where $\Mreg(\bv,\bw)^{\alpha}$ is a component of
$\Mreg(\bv,\bw)$ and $IC(\Mreg(\bv,\bw)^{\alpha},\phi)$ is the
intersection complex associated with an irreducible local system
$\phi$ on $\Mreg(\bv,\bw)^{\alpha}$.
Moreover, by \cite[8.9.9]{Gi-book}, we have
\begin{equation}\label{eq:Htop}
   H_{\topdeg}(\M(\bv^0,\bw)_{x_{\bv}},\C) = \bigoplus_\phi
     L_{(\bv,\alpha,\phi)},
\end{equation}
where $\phi$ runs over the set of irreducible local systems on the
component of $\Mreg(\bv,\bw)^{\alpha}$ containing $x_{\bv}$.
But as argued in the proof of \thmref{thm:std_decomp},
the indecomposability of $H_{\topdeg}(\Mw_{x_{\bv}},\C)$
implies that no intersection complex associated with a {\it
nontrivial\/} local system appears in the summand.
Moreover the left hand side of \eqref{eq:Htop} is independent of the
choice of the component by \thmref{thm:slice}.
Thus we can combine the summation over $\alpha$ together as
\[
   \bigoplus_{\alpha,\phi} L_{\bv,\alpha,\phi} \otimes
       IC(\Mreg(\bv,\bw)^{\alpha},\phi)
   = H_{\topdeg}(\M(\bv^0,\bw)_{x_{\bv}},\C)\otimes IC(\Mreg(\bv,\bw)).
\]

Our remaining task is to identify the index set of $\bv$.
The fundamental class $[\M(\bv^0,\bw)_{x_{\bv}}]$ is nonzero 
if $\M(\bv^0,\bw)_{x_{\bv}}$ is nonempty. Thus
$\M(\bv^0,\bw)_{x_{\bv}}$ is nonempty if and only if
\[
    H_{\topdeg}(\M(\bv^0,\bw)_{x_{\bv}},\C)\neq 0.
\]
By \cite[10.2]{Na-alg} and the construction,
$H_{\topdeg}(\M(\bv^0,\bw)_{x_{\bv}},\C)$ is isomorphic to the weight
space of weight $\bw - \bv^0$ in $L(\bw-\bv)$. Thus it is nonzero if
and only if $\bv\in \mathcal V_{\bv^0}(\bw)$.
\end{proof}

Take $x\in \M(\infty,\bw)$ and consider the inclusion
$i_x\colon \{x\} \to \M_0(\infty,\bw)$.
Applying $H^*(i_x^!\bullet)$ to \eqref{eq:decompQ} and then
summing up with respect to $\bv^0$, we get
\begin{equation}\label{eq:decomp'}
   H_*(\Mw_x,\C)
   = \bigoplus_{\bv\in\mathcal V(\bw)}
        H_{\topdeg}(\Mw_{x_{\bv}},\C)\otimes
        H^*(i_x^!IC(\Mreg(\bv,\bw))).
\end{equation}
By the convolution product, $H_{\topdeg}(\Mw_{x_{\bv}},\C)$ is a
module of $H_{\topdeg}(\Zw,\C)$.
By \cite[\S8.9]{Gi-book}, the decomposition \eqref{eq:decomp'} is
compatible with the module structure, where
$H_{\topdeg}(\Zw,\C)$ acts on
$H_{\topdeg}(\Mw_{x_{\bv}},\C)\otimes H^*(i_x^!IC(\Mreg(\bv,\bw)))$ by
$z\colon \xi\otimes\xi' \mapsto z\xi\otimes\xi'$.
This completes the proof of \thmref{thm:UqmoduleStr}.

\end{document}